\documentclass[11pt]{amsart}

\usepackage{
    amsmath,
    amsfonts,
    amssymb,
    amsthm,
    amscd,
    gensymb,      
    graphicx,
    comment,      
    etoolbox,     
    mathtools,    
    enumitem,
    mathdots,     
    stmaryrd,     
    txfonts,      
    booktabs,     
    stackrel,     
}
\usepackage[usenames,dvipsnames]{xcolor}

\makeatletter
\@namedef{subjclassname@2020}{%
  \textup{2020} Mathematics Subject Classification}
\makeatother

\usepackage[T1]{fontenc}
\usepackage{bbm}                     
\usepackage[colorlinks=true, linkcolor=blue, citecolor=blue, urlcolor=blue, breaklinks=true]{hyperref}

\DeclareFontFamily{OT1}{pzc}{}
\DeclareFontShape{OT1}{pzc}{m}{it}{<-> s * [1.10] pzcmi7t}{}
\DeclareMathAlphabet{\mathpzc}{OT1}{pzc}{m}{it}

\usepackage[capitalize]{cleveref}   

\crefname{defin}{Definition}{Definitions}
\crefname{eg}{Example}{Examples}
\crefname{lem}{Lemma}{Lemmas}
\crefname{theo}{Theorem}{Theorems}
\crefname{equation}{}{}
\crefname{enumi}{}{}
\newcommand{\creflastconjunction}{, and\nobreakspace}  

\usepackage{tikz}
\usetikzlibrary{arrows.meta}
\usetikzlibrary{decorations.markings}
\usetikzlibrary{calc}

\newcommand\braidup{to[out=up,in=down]}
\newcommand\braiddown{to[out=down,in=up]}

\newcommand\dotlabel[1]{$\scriptstyle{#1}$}
\newcommand\token[4][black]{
    \filldraw[#1] (#2) circle (1.5pt) node[anchor=#3] {\dotlabel{#4}}
}
\newcommand\singdot[2][black]{
    \filldraw[fill=white, draw=#1] (#2) circle (1.5pt)
}
\newcommand\multdot[4][black]{
    \filldraw[fill=white, draw=#1] (#2) circle (1.5pt) node[anchor=#3] {{\color{#1} \dotlabel{#4}}}
}
\newcommand\heartdec[4]{ 
    \draw (#1)++(0,-0.03) node[color=white] {$\scriptstyle{\varheartsuit}$} node[color=violet] {$\scriptstyle{\heartsuit}$};
    \draw (#1)++(0,-0.03) node[anchor=#2] {\dotlabel{(#4,#3)}}
}
\newcommand\diamdec[4]{ 
    \draw (#1) node[color=white] {$\scriptstyle{\Diamondblack}$} node[color=violet] {$\scriptstyle{\diamondsuit}$};
    \draw (#1) node[anchor=#2] {\dotlabel{(#4,#3)}}
}

\newcommand\teleport[2]{
    \draw (#1) -- (#2);
    \filldraw[black] (#1) circle (1.5pt);
    \filldraw[black] (#2) circle (1.5pt);
}
\newcommand\telecolor[4]{
    \draw (#3) to (#4);
    \filldraw[#1] (#3) circle (1.5pt);
    \filldraw[#2] (#4) circle (1.5pt)
}
\newcommand\spear[5][0.65]{
    \draw[thick,densely dotted] (#4) to (#5);
    \draw[->] ($ (#4) !#1-0.01! (#5) $) to ($ (#4) !#1! (#5) $);
    \singdot[#2]{#4};
    \singdot[#3]{#5}
}
\newcommand\dart[5][0.65]{
    \draw[thick,-] (#4) to (#5);
    \draw[->] (#4) to ($ (#4) !#1! (#5) $);
    \singdot[#2]{#4};
    \singdot[#3]{#5}
}
\newcommand\tripdart[6]{
    \draw[thick,-] (#4) -- (#5) -- (#6);
    \draw[->] (#4) to ($ (#4) !0.65! (#5) $);
    \draw[->] (#5) to ($ (#5) !0.65! (#6) $);
    \singdot[#1]{#4};
    \singdot[#2]{#5};
    \singdot[#3]{#6}
}
\newcommand\boxbell[5][0.65]{
    \draw[thick,-] (#4) to (#5);
    \filldraw[fill=white,draw=#2] (#4)++(-0.06,-0.06) rectangle ++(0.12,0.12);
    \filldraw[fill=white,draw=#3] (#5)++(-0.06,-0.06) rectangle ++(0.12,0.12)
}
\newcommand\boxdart[5][0.65]{
    \draw[thick,-] (#4) to (#5);
    \draw[->] (#4) to ($ (#4) !#1! (#5) $);
    \filldraw[fill=white,draw=#2] (#4)++(-0.06,-0.06) rectangle ++(0.12,0.12);
    \filldraw[fill=white,draw=#3] (#5)++(-0.06,-0.06) rectangle ++(0.12,0.12)
}
\newcommand\tripboxdart[6]{
    \draw[thick,-] (#4) -- (#5) -- (#6);
    \draw[->] (#4) to ($ (#4) !0.65! (#5) $);
    \draw[->] (#5) to ($ (#5) !0.65! (#6) $);
    \filldraw[fill=white,draw=#1] (#4)++(-0.06,-0.06) rectangle ++(0.12,0.12);
    \filldraw[fill=white,draw=#2] (#5)++(-0.06,-0.06) rectangle ++(0.12,0.12);
    \filldraw[fill=white,draw=#3] (#6)++(-0.06,-0.06) rectangle ++(0.12,0.12)
}

\newcommand\bubrightblank[2][black]{
    \draw[->,#1] (#2)++(0,0.2) arc(90:-270:0.2)
}
\newcommand\bubright[4][black]{
    \draw[->,#1] (#2)++(0,0.2) arc(90:-270:0.2);
    \filldraw[#1] (#2)++(0.2,0) circle (1.5pt) node[anchor=west] {\dotlabel{#3}};
    \filldraw[fill=white, draw=#1] (#2)++(-0.2,0) circle (1.5pt) node[anchor=east] {{\color{#1} \dotlabel{#4}}}
}
\newcommand\bubleftblank[2][black]{
    \draw[->,#1] (#2)++(0,0.2) arc(90:450:0.2)
}
\newcommand\bubleft[4][black]{
    \draw[->,#1] (#2)++(0,0.2) arc(90:450:0.2);
    \filldraw[#1] (#2)++(-0.2,0) circle (1.5pt) node[anchor=east] {\dotlabel{#3}};
    \filldraw[fill=white, draw=#1] (#2)++(0.2,0) circle (1.5pt) node[anchor=west] {{\color{#1} \dotlabel{#4}}}
}

\newcommand\plusrightblank[2][black]{
    \draw[->,#1] (#2)++(0,0.2) arc(90:-270:0.2);
    \node at (#2) {{\color{#1} \dotlabel{+}}}
}
\newcommand\plusright[4][black]{
    \draw[->,#1] (#2)++(0,0.2) arc(90:-270:0.2);
    \node at (#2) {{\color{#1} \dotlabel{+}}};
    \filldraw[#1] (#2)++(0.2,0) circle (1.5pt) node[anchor=west] {\dotlabel{#3}};
    \filldraw[fill=white, draw=#1] (#2)++(-0.2,0) circle (1.5pt) node[anchor=east] {{\color{#1} \dotlabel{#4}}}
}
\newcommand\plusleftblank[2][black]{
    \draw[->,#1] (#2)++(0,0.2) arc(90:450:0.2);
    \node at (#2) {{\color{#1} \dotlabel{+}}}
}
\newcommand\plusleft[4][black]{
    \draw[->,#1] (#2)++(0,0.2) arc(90:450:0.2);
    \node at (#2) {{\color{#1} \dotlabel{+}}};
    \filldraw[#1] (#2)++(-0.2,0) circle (1.5pt) node[anchor=east] {\dotlabel{#3}};
    \filldraw[fill=white, draw=#1] (#2)++(0.2,0) circle (1.5pt) node[anchor=west] {{\color{#1} \dotlabel{#4}}}
}
\newcommand\plusgenright[2][black]{
    \draw[->,#1] (#2)++(0,0.2) arc(90:-270:0.2);
    \node at (#2) {\dotlabel{w}}
}
\newcommand\plusgenleft[2][black]{
    \draw[->,#1] (#2)++(0,0.2) arc(90:450:0.2);
    \node at (#2) {\dotlabel{w}}
}

\newcommand\minusrightblank[2][black]{
    \draw[->,#1] (#2)++(0,0.2) arc(90:-270:0.2);
    \node at (#2) {{\color{#1} \dotlabel{-}}}
}
\newcommand\minusright[4][black]{
    \draw[->,#1] (#2)++(0,0.2) arc(90:-270:0.2);
    \node at (#2) {{\color{#1} \dotlabel{-}}};
    \filldraw[#1] (#2)++(0.2,0) circle (1.5pt) node[anchor=west] {\dotlabel{#3}};
    \filldraw[fill=white, draw=#1] (#2)++(-0.2,0) circle (1.5pt) node[anchor=east] {{\color{#1} \dotlabel{#4}}}
}
\newcommand\minusleftblank[2][black]{
    \draw[->,#1] (#2)++(0,0.2) arc(90:450:0.2);
    \node at (#2) {{\color{#1} \dotlabel{-}}}
}
\newcommand\minusleft[4][black]{
    \draw[->,#1] (#2)++(0,0.2) arc(90:450:0.2);
    \node at (#2) {{\color{#1} \dotlabel{-}}};
    \filldraw[#1] (#2)++(-0.2,0) circle (1.5pt) node[anchor=east] {\dotlabel{#3}};
    \filldraw[fill=white, draw=#1] (#2)++(0.2,0) circle (1.5pt) node[anchor=west] {{\color{#1} \dotlabel{#4}}}
}
\newcommand\minusgenright[2][black]{
    \draw[->,#1] (#2)++(0,0.2) arc(100:90:0.2);
    \draw[densely dotted,thick,#1] (#2)++(0,0.2) arc(90:-270:0.2);
    \node at (#2) {\dotlabel{w}}
}
\newcommand\minusgenleft[2][black]{
    \draw[->,#1] (#2)++(0,0.2) arc(80:90:0.2);
    \draw[densely dotted,thick,#1] (#2)++(0,0.2) arc(90:450:0.2);
    \node at (#2) {\dotlabel{w}}
}

\newcommand\intleft[2]{ 
    \filldraw[->,draw=#1,fill=white] (#2)++(-0.18,0) arc (-180:180:0.18)
}
\newcommand\intright[2]{
    \filldraw[->,draw=#1,fill=white] (#2)++(0.18,0) arc (0:-360:0.18)
}
\newcommand\intleftsm[2]{
    \filldraw[white] (#2) circle (0.1);
    \draw[#1,->] (#2)++(0,-0.1) to[out=0,in=-90] ++(0.1,0.1) to[out=90,in=0] ++(-0.1,0.1) to[out=180,in=90] ++(-0.09,-0.14);
    \draw[#1] (#2)++(-0.09,-0.02) to[out=-84,in=0] ++(0.09,-0.08)
}
\newcommand\intrightsm[2]{
    \filldraw[white] (#2) circle (0.1);
    \draw[#1,->] (#2)++(0,-0.1) to[out=180,in=-90] ++(-0.1,0.1) to[out=90,in=180] ++(0.1,0.1) to[out=0,in=90] ++(0.09,-0.14);
    \draw[#1] (#2)++(0.09,-0.02) to[out=-96,in=0] ++(-0.09,-0.08)
}

\newcommand\cbubble[3][black]{
    \begin{tikzpicture}[centerzero]
        \bubright[#1]{0,0}{#2}{#3};
    \end{tikzpicture}
}
\newcommand\ccbubble[3][black]{
    \begin{tikzpicture}[centerzero]
        \bubleft[#1]{0,0}{#2}{#3};
    \end{tikzpicture}
}
\newcommand\pluscbubble[3][black]{
    \begin{tikzpicture}[centerzero]
        \plusright[#1]{0,0}{#2}{#3};
    \end{tikzpicture}
}
\newcommand\plusccbubble[3][black]{
    \begin{tikzpicture}[centerzero]
        \plusleft[#1]{0,0}{#2}{#3};
    \end{tikzpicture}
}
\newcommand\minuscbubble[3][black]{
    \begin{tikzpicture}[centerzero]
        \minusright[#1]{0,0}{#2}{#3};
    \end{tikzpicture}
}
\newcommand\minusccbubble[3][black]{
    \begin{tikzpicture}[centerzero]
        \minusleft[#1]{0,0}{#2}{#3};
    \end{tikzpicture}
}
\newcommand\tokup[1][a]{
    \begin{tikzpicture}[centerzero]
        \draw[->] (0,-0.2) -- (0,0.2);
        \token{0,0}{west}{#1};
    \end{tikzpicture}
}
\newcommand\tokdown[1][a]{
    \begin{tikzpicture}[centerzero]
        \draw[<-] (0,-0.2) -- (0,0.2);
        \token{0,0}{west}{#1};
    \end{tikzpicture}
}
\newcommand\dotup{
    \begin{tikzpicture}[centerzero]
        \draw[->] (0,-0.2) -- (0,0.2);
        \singdot{0,0};
    \end{tikzpicture}
}
\newcommand\poscross{
    \begin{tikzpicture}[centerzero]
        \draw[->] (0.2,-0.2) -- (-0.2,0.2);
        \draw[wipe] (-0.2,-0.2) -- (0.2,0.2);
        \draw[->] (-0.2,-0.2) -- (0.2,0.2);
    \end{tikzpicture}
}
\newcommand\negcross{
    \begin{tikzpicture}[centerzero]
        \draw[->] (-0.2,-0.2) -- (0.2,0.2);
        \draw[wipe] (0.2,-0.2) -- (-0.2,0.2);
        \draw[->] (0.2,-0.2) -- (-0.2,0.2);
    \end{tikzpicture}
}
\newcommand\rightcup{
    \begin{tikzpicture}[centerzero]
        \draw[->] (-0.2,0.2) -- (-0.2,0) arc(180:360:0.2) -- (0.2,0.2);
    \end{tikzpicture}
}
\newcommand\rightcap{
    \begin{tikzpicture}[centerzero]
        \draw[->] (-0.2,-0.2) -- (-0.2,0) arc(180:0:0.2) -- (0.2,-0.2);
    \end{tikzpicture}
}
\newcommand\leftcup{
    \begin{tikzpicture}[centerzero]
        \draw[<-] (-0.2,0.2) -- (-0.2,0) arc(180:360:0.2) -- (0.2,0.2);
    \end{tikzpicture}
}
\newcommand\leftcap{
    \begin{tikzpicture}[centerzero]
        \draw[<-] (-0.2,-0.2) -- (-0.2,0) arc(180:0:0.2) -- (0.2,-0.2);
    \end{tikzpicture}
}

\tikzset{anchorbase/.style={>=To,baseline={([yshift=-0.5ex]current bounding box.center)}}}
\tikzset{ 
    centerzero/.style={>=To,baseline={([yshift=-0.5ex](#1))}},
    centerzero/.default={0,0}
}
\tikzset{wipe/.style={white,line width=3pt}}

\tikzset{->-/.style={decoration={
  markings,
  mark=at position #1 with {\arrow{>}}},postaction={decorate}}}
\tikzset{-<-/.style={decoration={
  markings,
  mark=at position #1 with {\arrow{<}}},postaction={decorate}}}

\newcommand\C{\mathbb{C}}
\newcommand\kk{\Bbbk}
\newcommand\KK{\mathbb{K}}
\newcommand\N{\mathbb{N}}
\newcommand\OO{\mathbb{O}}
\newcommand\tOO{\tilde{\OO}}
\newcommand\Q{\mathbb{Q}}
\newcommand\R{\mathbb{R}}
\newcommand\Z{\mathbb{Z}}

\newcommand\cI{\mathcal{I}}

\newcommand\ba{\mathbf{a}}
\newcommand\bb{\mathbf{b}}
\newcommand\bc{\mathbf{c}}
\newcommand\bC{\mathbf{C}}
\newcommand\bsf{\mathbf{f}}
\newcommand\bsg{\mathbf{g}}
\newcommand\bF{\mathbf{F}}
\newcommand{\B}{\mathbf{B}}             
\newcommand{\BA}{{\B_A}}                
\newcommand\one{\mathbbm{1}}

\newcommand\gr{\mathrm{gr}}
\newcommand\new{\mathrm{new}}
\newcommand\old{\mathrm{old}}
\newcommand\op{\mathrm{op}}              
\newcommand\rev{\mathrm{rev}}
\newcommand\tr{\mathrm{tr}}

\newcommand\psmod{\textup{psmod-}}
\newcommand\smod{\textup{smod-}}

\newcommand\even{{\bar{0}}}
\newcommand\odd{{\bar{1}}}

\newcommand\cA{\mathpzc{A}}
\newcommand\cB{\mathpzc{B}}
\newcommand\cC{\mathpzc{C}}
\newcommand\cV{\mathpzc{V}}

\newcommand\cH{\mathpzc{H}}             
\newcommand\Heis{\mathpzc{Heis}}        
\newcommand\rHeis{\mathrm{Heis}}        
\newcommand\QAW{\mathpzc{AW}}          
\newcommand\invQAW{\overline{\mathpzc{AW}}}     
\newcommand\QAWA{\mathrm{AW}}          
\newcommand\QW{\mathpzc{W}}            
\newcommand\QWA{\mathrm{W}}            
\newcommand\SEnd{\mathpzc{SEnd}}        
\newcommand\SVec{\mathpzc{SVec}}        

\newcommand\fS{\mathfrak{S}}

\newcommand\sQ{\mathsf{Q}}

\newcommand\Laurent[1]{
    (\!(#1)\!)
}
\newcommand\red[1]{{\color{red} #1 }}
\newcommand\blue[1]{{\color{blue} #1}}
\newcommand\upblue{\blue{\uparrow}}
\newcommand\upred{\red{\uparrow}}
\newcommand\downblue{\blue{\downarrow}}
\newcommand\downred{\red{\downarrow}}
\newcommand\barodot{\, \overline{\odot}\, }
\newcommand\cDelt[1][l]{\Delta_{\blue{#1}|\red{m}}}
\newcommand\cocenter[1]{
    \dot{#1}
}

\DeclareMathOperator{\Add}{Add}
\DeclareMathOperator{\Aut}{Aut}
\DeclareMathOperator{\coind}{coind}
\DeclareMathOperator{\End}{End}
\DeclareMathOperator{\END}{END}
\DeclareMathOperator{\Ev}{Ev}
\DeclareMathOperator{\flip}{flip}
\DeclareMathOperator{\Hom}{Hom}
\DeclareMathOperator{\HOM}{HOM}
\DeclareMathOperator{\id}{id}
\DeclareMathOperator{\im}{im}       
\DeclareMathOperator{\ind}{ind}
\DeclareMathOperator{\Kar}{Kar}     
\DeclareMathOperator{\res}{res}
\DeclareMathOperator{\Span}{Span}
\DeclareMathOperator{\Sym}{Sym}

\newtheorem{theo}{Theorem}[section]

\newtheorem{lem}[theo]{Lemma}
\newtheorem*{lem*}{Lemma}
\newtheorem{cor}[theo]{Corollary}

\theoremstyle{definition}
\newtheorem{defin}[theo]{Definition}
\newtheorem{rem}[theo]{Remark}

\numberwithin{equation}{section}
\allowdisplaybreaks

\setenumerate[1]{label=(\alph*)}  
\setenumerate[2]{label=(\roman*)}

\newtoggle{comments}
\newtoggle{details}
\newtoggle{detailsnote}

\toggletrue{detailsnote}   

\iftoggle{comments}{%
  \newcommand{\acomments}[1]{
    \ \\
    {\color{red}
      \textbf{AS:} #1
    }
    \\
    }
}{%
  \newcommand{\acomments}[1]{}
}

\iftoggle{details}{%
  \newcommand{\details}[1]{
      \ \\
      {\color{OliveGreen}
        \textbf{Details:} #1
      }
      \\
  }
}{%
  \newcommand{\details}[1]{}
}

\begin{document}

\title{Quantum Frobenius Heisenberg categorification}

\author{Jonathan Brundan}
\address[J.B.]{
    Department of Mathematics \\
    University of Oregon \\
    Eugene, OR, USA
}
\email{brundan@uoregon.edu}

\author{Alistair Savage}
\address[A.S.]{
    Department of Mathematics and Statistics \\
    University of Ottawa \\
    Ottawa, ON, Canada
}
\urladdr{\href{https://alistairsavage.ca}{alistairsavage.ca}, \textrm{\textit{ORCiD}:} \href{https://orcid.org/0000-0002-2859-0239}{orcid.org/0000-0002-2859-0239}}
\email{alistair.savage@uottawa.ca}

\author{Ben Webster}
\address[B.W.]{
    Department of Pure Mathematics, University of Waterloo \&
    Perimeter Institute for Theoretical Physics\\
    Waterloo, ON, Canada
}
\email{ben.webster@uwaterloo.ca}

\thanks{J.B.\ supported in part by NSF grant DMS-1700905.}
\thanks{A.S.\ supported by Discovery Grant RGPIN-2017-03854 from the Natural Sciences and Engineering Research Council of Canada.}
\thanks{B.W.\ supported by Discovery Grant RGPIN-2018-03974 from the Natural Sciences and Engineering Research Council of Canada. This research was supported in part by Perimeter Institute for Theoretical Physics. Research at Perimeter Institute is supported by the Government of Canada through the Department of Innovation, Science and Economic Development Canada and by the Province of Ontario through the Ministry of Research, Innovation and Science.}

\begin{abstract}
    We associate a diagrammatic monoidal category $\Heis_k(A;z,t)$,
    which we call the \emph{quantum Frobenius Heisenberg category}, to
    a symmetric Frobenius superalgebra $A$, a central charge $k \in
    \Z$, and invertible parameters $z,t$ in some ground ring. When $A$ is trivial, i.e.\ it equals the ground ring,
    these categories recover the quantum Heisenberg categories
    introduced in our previous work, and when the central charge $k$ is
    zero they yield generalizations of the affine HOMFLY-PT skein
    category.  By exploiting some natural categorical actions of
    $\Heis_k(A;z,t)$ on generalized cyclotomic quotients, we prove a basis theorem for morphism spaces.
\end{abstract}


\ifboolexpr{togl{comments} or togl{details}}{%
  {\color{magenta}DETAILS OR COMMENTS ON}
}{%
}

\maketitle
\thispagestyle{empty}


\section{Introduction}

Fix a commutative ground ring $\kk$ and parameters $z,t \in \kk^\times$.  In this paper, to any central charge $k \in \Z$ and symmetric Frobenius superalgebra $A$, we associate  a  monoidal supercategory $\Heis_k(A;z,t)$, which we call the \emph{quantum Frobenius Heisenberg category}.  The case $A = \kk$ recovers the \emph{quantum Heisenberg categories} $\Heis_k(z,t)$ of \cite{BSW-qheis}, which extend the $q$-deformed Heisenberg category introduced in \cite{LS13}.  In all other cases, the categories are new.  In this way, the current paper can be viewed as a generalization of all of these earlier treatments of quantum Heisenberg categories.  Indeed, in \cite{BSW-qheis} the proofs of many diagrammatic lemmas were omitted, deferring them to the current paper, which gives these proofs in the more general setting.
Quantum Frobenius Heisenberg categories can also be viewed as quantum analogues of the \emph{Frobenius Heisenberg categories} introduced in \cite{RS17,Sav19} and further developed in \cite{BSW-foundations}.  In turn, these were generalizations of the \emph{Heisenberg categories} of \cite{MS18,Bru18}, which were based on the original construction of \cite{Kho14}.

Following the approach of \cite{BSW-qheis}, we give three equivalent definitions of $\Heis_k(A;z,t)$.  All three take as their starting point the \emph{quantum affine wreath product category} $\QAW(A;z)$, whose morphism spaces are the \emph{quantum affine wreath product algebras} introduced in \cite{RS20}; see \cref{sec:QAWC}.  The next step is to adjoin a right dual $\downarrow$ to the generating object $\uparrow$ of $\QAW(A;z)$.  This involves adding right cups and caps satisfying zigzag relations.  In the first two approaches, described in \cref{sec:first,sec:second}, the final step is to require that a certain morphism in the additive envelope is invertible, and then impose one more relation on certain bubbles.  The only difference between these two definitions is that the morphism to be inverted involves different crossings (positive or negative).  The third definition of $\Heis_k(A;z,t)$, given in \cref{sec:third}, is the most explicit.  Here we adjoin additional generating morphisms, the left cups and caps, subject to further relations which are often useful in practice.  The proof that all three approaches yield isomorphic categories involves a systematic ``unpacking'' of the inversion relations.

The categories $\Heis_k(A;z,t)$ have many desirable properties, analogous to those of the other Heisenberg categories mentioned above.  For example, they are strictly pivotal, and we have curl relations, bubble slide relations, and infinite Grassmannian relations.  When $k \ne 0$, they act naturally on categories of modules for \emph{quantum cyclotomic wreath product algebras}, a fact which is easiest to see in the inversion relation approach of the first two definitions of $\Heis_k(A;z,t)$; see \cref{sec:action}.  We also prove a basis theorem (\cref{basis-thm}) for the morphisms spaces in $\Heis_k(A;z,t)$.  Following the methods of \cite{BSW-K0,BSW-qheis,BSW-foundations}, our proof of this basis theorem uses a \emph{categorical comultiplication} (see \cref{sec:comult}) to construct a sufficiently large module category starting from the actions on modules over the cyclotomic wreath product algebras for generic parameters. We also introduce certain {\em generalized cyclotomic quotients}; see \cref{sec:GCQ}.

The central charge zero quantum Heisenberg category $\Heis_0(z,t)$ from \cite{BSW-qheis} is the \emph{affine HOMFLY-PT skein category} from \cite[\S4]{Bru17}.  Omitting the dot generator, we obtain the \emph{HOMFLY-PT skein category}, which is the ribbon category underpinning the definition of the HOMFLY-PT invariant of an oriented link.  The central charge zero quantum Frobenius Heisenberg category $\Heis_0(A;z,t)$ yields a Frobenius algebra analog of the affine HOMFLY-PT skein category. By analogy with the $A = \kk$ case, we expect that $\Heis_0(A;z,t)$ should act on categories of modules for some natural Frobenius algebra analogue of the quantized enveloping algebra $U_q(\mathfrak{gl}_n)$. On omitting the dot generator from $\Heis_0(A;z,t)$, we obtain the {\em Frobenius HOMFLY-PT skein category}, from which one can define a Frobenius algebra analogue of the HOMFLY-PT link invariant.

The name \emph{Heisenberg category} comes from Khovanov's original conjecture that his Heisenberg category categorifies a central reduction of the infinite rank Heisenberg algebra.  This conjecture was proved in \cite[Th.~1.1]{BSW-K0} and then extended in \cite[Th.~10.5]{BSW-foundations} to show that the Grothendieck ring of the Frobenius Heisenberg category is isomorphic to a \emph{lattice Heisenberg algebra} associated to the Frobenius algebra $A$.  We expect that the Grothendieck ring of the quantum Frobenius Heisenberg category is also isomorphic to this lattice Heisenberg algebra.  The obstruction to proving this conjecture, present even in the $A=\kk$ case, is that we do not know how to compute the split Grothendieck group of the quantum affine wreath product algebra.

There is an alternative framework for decategorification using \emph{trace} (or zeroth Hochschild homology) in place of the Grothendieck ring.  In \cite{CLLSS18}, the trace of the $q$-deformed Heisenberg category of \cite{LS13} was related to the elliptic Hall algebra.  It is thus natural to expect that the trace of $\Heis_k(A;z,t)$ should be related to a Frobenius algebra generalization of the elliptic Hall algebra.  This would be the quantum analogue of the situation for Frobenius Heisenberg categories, where a computation of the trace in \cite{ReeksS20} suggested a Frobenius algebra generalization of the W-algebra $W_{1+\infty}$, extending results from \cite{CLLS18}.

\section{The quantum affine wreath product category\label{sec:QAWC}}

A fundamental ingredient in the definition of the quantum Frobenius Heisenberg category is the quantum affine wreath product algebra introduced in \cite{RS20}.  In this section we recall the definition of this algebra, together with a few key results needed in the current paper.  We then switch our point of view to monoidal categories, introducing categories whose endomorphism spaces are quantum affine wreath product algebras.

\subsection*{Monoidal supercategories}

We fix a commutative ground ring $\kk$.  All vector spaces, algebras, categories, and functors will be assumed to be linear over $\kk$ unless otherwise specified.  Unadorned tensor products denote tensor products over $\kk$.  Almost everything in the paper will be enriched over the category $\SVec$ of vector superspaces with parity-preserving morphisms.  We write $\bar{v}$ for the parity of a homogeneous vector $v$ in a vector superspace.  When we write formulae involving parities, we assume the elements in question are homogeneous; we then extend by linearity.

For superalgebras $A = A_\even \oplus A_\odd$ and $B = B_\even \oplus B_\odd$, multiplication in the superalgebra $A \otimes B$ is defined by
\begin{equation}
    (a' \otimes b) (a \otimes b') = (-1)^{\bar a \bar b} a'a \otimes bb'
\end{equation}
for homogeneous $a,a' \in A$, $b,b' \in B$.  The \emph{opposite} superalgebra $A^\op$ is a copy $\{a^\op : a \in A\}$ of the vector superspace $A$ with multiplication defined from
\begin{equation}\label{mups}
    a^\op  b^\op := (-1)^{\bar a \bar b} (ba)^{\op}.
\end{equation}
The \emph{center} of $A$ is the subalgebra
\begin{equation}\label{center}
    Z(A) := \Span_\kk \{ a \in A \text{ homogeneous} : ab = (-1)^{\bar a \bar b}ba \text{ for all } b \in A\}
\end{equation}
The \emph{cocenter} $C(A)$, which, in general, is merely a vector superspace not an superalgebra, is
\begin{equation}\label{cocenter}
    C(A) := A / \Span_\kk \{ab - (-1)^{\bar{a} \bar{b}} ba : a,b \in A \text{ homogeneous} \}
\end{equation}

Throughout this paper we will work with \emph{strict monoidal supercategories}, in the sense of \cite{BE17}.  We refer the reader to \cite[\S 2]{BSW-foundations} for a summary of this topic well adapted to the current work, or to \cite{BE17} for a thorough treatment.  We summarize here a few crucial properties that play an important role in the present paper.

A \emph{supercategory} means a category enriched in $\SVec$.  Thus, its morphism spaces are vector superspaces and composition is parity-preserving.  A \emph{superfunctor} between supercategories induces a parity-preserving linear map between morphism superspaces.  For superfunctors $F,G \colon \cA \to \cB$, a \emph{supernatural transformation} $\alpha \colon F \Rightarrow G$ of \emph{parity $r\in\Z/2$} is the data of morphisms $\alpha_X\in \Hom_{\cB}(FX, GX)_r$ for each $X \in \cA$ such that $Gf \circ \alpha_X = (-1)^{r \bar f}\alpha_Y\circ Ff$ for each homogeneous $f \in \Hom_{\cA}(X, Y)$.  Note when $r$ is odd that $\alpha$ is \emph{not} a natural transformation in the usual sense due to the sign. A \emph{supernatural transformation} $\alpha \colon F \Rightarrow G$ is of the form $\alpha = \alpha_\even + \alpha_\odd$, with each $\alpha_r$ being a supernatural transformation of parity $r$.

In a \emph{strict monoidal supercategory}, morphisms satisfy the \emph{super interchange law}:
\begin{equation}\label{interchange}
    (f' \otimes g) \circ (f \otimes g')
    = (-1)^{\bar f \bar g} (f' \circ f) \otimes (g \circ g').
\end{equation}
We denote the unit object by $\one$ and the identity morphism of an object $X$ by $1_X$.  We will use the usual calculus of string diagrams, representing the horizontal composition $f \otimes g$ (resp.\ vertical composition $f \circ g$) of morphisms $f$ and $g$ diagrammatically by drawing $f$ to the left of $g$ (resp.\ drawing $f$ above $g$).  Care is needed with horizontal levels in such diagrams due to the signs arising from the super interchange law:
\begin{equation}\label{intlaw}
    \begin{tikzpicture}[anchorbase]
        \draw (-0.5,-0.5) -- (-0.5,0.5);
        \draw (0.5,-0.5) -- (0.5,0.5);
        \filldraw[fill=white,draw=black] (-0.5,0.15) circle (5pt);
        \filldraw[fill=white,draw=black] (0.5,-0.15) circle (5pt);
        \node at (-0.5,0.15) {$\scriptstyle{f}$};
        \node at (0.5,-0.15) {$\scriptstyle{g}$};
    \end{tikzpicture}
    \quad=\quad
    \begin{tikzpicture}[anchorbase]
        \draw (-0.5,-0.5) -- (-0.5,0.5);
        \draw (0.5,-0.5) -- (0.5,0.5);
        \filldraw[fill=white,draw=black] (-0.5,0) circle (5pt);
        \filldraw[fill=white,draw=black] (0.5,0) circle (5pt);
        \node at (-0.5,0) {$\scriptstyle{f}$};
        \node at (0.5,0) {$\scriptstyle{g}$};
    \end{tikzpicture}
    \quad=\quad
    (-1)^{\bar f\bar g}\
    \begin{tikzpicture}[anchorbase]
        \draw (-0.5,-0.5) -- (-0.5,0.5);
        \draw (0.5,-0.5) -- (0.5,0.5);
        \filldraw[fill=white,draw=black] (-0.5,-0.15) circle (5pt);
        \filldraw[fill=white,draw=black] (0.5,0.15) circle (5pt);
        \node at (-0.5,-0.15) {$\scriptstyle{f}$};
        \node at (0.5,0.15) {$\scriptstyle{g}$};
    \end{tikzpicture}
    \ .
\end{equation}

If $\cA$ is a supercategory, the category $\SEnd_\kk(A)$ of superfunctors $\cA \to \cA$ and supernatural transformations is a strict monoidal supercategory.  A \emph{module supercategory} over a strict monoidal supercategory $\cC$ is a supercategory $\cA$ together with a monoidal superfunctor $\cC \rightarrow \SEnd_\kk(\cA)$.

If $R$ is a superalgebra, we let $\smod R$ denote the supercategory of right $R$-supermodules and let $\psmod R$ denote the supercategory of finitely-generated projective right $R$-supermodules.

\subsection*{Quantum affine wreath product algebras}

Fix $z \in \kk^\times$.  (In fact, we do not need the assumption that $z$ is invertible until we introduce quantum Frobenius Heisenberg categories in \cref{sec:first}, but we make this assumption from the beginning to be uniform.)  Let $A$ be a symmetric Frobenius superalgebra with even trace map $\tr \colon A \to \kk$.  Thus
\[
    \tr(ab) = (-1)^{\bar{a}\bar{b}} \tr(ba),\quad a,b \in A.
\]
The definition of Frobenius superalgebra gives that $A$ possesses a homogeneous basis $\BA$ and a left dual basis $\{b^\vee : b \in \BA\}$ such that
\begin{equation} \label{chkdef}
    \tr(b^\vee c) = \delta_{b,c},\quad b,c \in \BA.
\end{equation}
It follows that, for all $a \in A$, we have
\begin{gather} \label{adecomp}
    a
    = \sum_{b \in \BA} \tr(b^\vee a)b
    = \sum_{b \in \BA} \tr(ab) b^\vee,
    \\ \label{beam}
    \sum_{b \in \BA} ab \otimes b^\vee
    = b \otimes b^\vee a
    ,\quad
    \sum_{b \in \BA} (-1)^{\bar{a} \bar{b}} b a \otimes b^\vee
    = \sum_{b \in \BA} (-1)^{\bar{a} \bar{b}} b \otimes a b^\vee.
\end{gather}
\details{
    For the second equality in \cref{beam}, we have
    \[
        \sum_{b \in \BA} (-1)^{\bar{a}\bar{b}} b a \otimes b^\vee
        \overset{\cref{adecomp}}{=} \sum_{b,c \in \BA} (-1)^{\bar{a} \bar{b}} \tr(c^\vee ba) c \otimes b^\vee
        = \sum_{b,c \in \BA} (-1)^{\bar{a} \bar{c}} c \otimes \tr(a c^\vee b) b^\vee
        \overset{\cref{adecomp}}{=} \sum_{c \in \BA} (-1)^{\bar{a}\bar{c}} c \otimes a c^\vee.
    \]
}
Note that $\bar{b} = \overline{b^\vee}$, and that the dual basis to $\{b^\vee : b \in \BA\}$ is given by
\begin{equation} \label{doubledual}
    (b^\vee)^\vee = (-1)^{\bar{b}} b.
\end{equation}

For the remainder of the paper we adopt the following summation convention: \emph{any expression involving both the symbols $b$ and $b^\vee$ includes an implicit sum over $b \in \BA$}.  We adopt an analogous convention when $b$ is replaced by $a$, $c$, etc.  Thus, for instance,
\[
    a b \otimes b^\vee = \sum_{b \in \BA} a b \otimes b^\vee = ac \otimes c^\vee.
\]
For any homogeneous $a \in A$, we define
\begin{equation}\label{adag}
    a^\dagger
    := (-1)^{\bar a \bar b} z b a  b^\vee
    = z \sum_{b \in \BA} (-1)^{\bar{a}\bar{b}} b a b^\vee,
\end{equation}
which is well-defined independent of the choice of the basis $\BA$.

\begin{defin}[{\cite[Def.~2.1]{RS20}}]\label{tiger}
    For $n \in \N$, $n \ge 2$, the \emph{quantum wreath product algebra} (or \emph{Frobenius Hecke algebra}) $\QWA_n(A;z)$ is the free product
    $A^{\otimes n} * \langle \sigma_i : 1 \le i \le n-1 \rangle$
    of the superalgebra $A^{\otimes n}$ and the
    free associative superalgebra with even generators $\sigma_1,\dots,\sigma_{n-1}$, modulo the relations
    \begin{align}
        \sigma_i \sigma_j &= \sigma_j \sigma_i, &1 \le i,j \le n-1,\ |i-j| > 1, \\
        \sigma_i \sigma_{i+1} \sigma_i &= \sigma_{i+1} \sigma_i \sigma_{i+1}, &1 \le i \le n-2, \\
        \sigma_i^2 &= z \tau_i \sigma_i + 1, &1 \le i \le n-1, \\
        \sigma_i \ba &= s_i(\ba) \sigma_i, &\ba \in A^{\otimes n},\ 1 \le i \le n-1,
    \end{align}
    where
    \begin{equation} \label{tau}
        \tau_i := 1^{\otimes (n-i-1)} b \otimes b^\vee \otimes 1^{\otimes (i-1)},\quad
        1 \le i,j \le n-1,
    \end{equation}
    and $s_i(\ba)$ denotes the action of the simple transposition $s_i$ on $\ba$ by superpermutation of the factors:
    \[
        s_i(a_n \otimes \dotsb \otimes a_1)
        =(-1)^{\overline{a_i}\, \overline{a_{i+1}}} a_n \otimes \dotsb \otimes a_{i+2} \otimes a_i \otimes a_{i+1} \otimes a_{i-1} \otimes \dotsb \otimes a_1.
    \]
    We define
    \begin{equation}
         a^{(i)} := 1^{\otimes (n-i)} \otimes a \otimes 1^{\otimes (i-1)},\quad
         a \in A,\ 1 \le i \le n.
    \end{equation}
    Note that here, and throughout the paper, we number factors in the tensor product \emph{from right to left}.  It is straightforward to verify that $\tau_i$ does not depend on the choice of basis $\BA$.  We adopt the conventions that $\QWA_1(A;z) := A$ and $\QWA_0(A;z) := \kk$.
\end{defin}

\begin{defin}[{\cite[Def.~2.5]{RS20}}]
    For $n \in \N$, $n \ge 1$, we define the \emph{quantum affine wreath product algebra} (or \emph{affine Frobenius Hecke algebra}) $\QAWA_n(A;z)$ to be the free product $\kk[x_1^{\pm 1},\dotsc,x_n^{\pm 1}] * \QWA_n(A;z)$
    of the Laurent polynomial algebra viewed as a superalgebra with all $x_i$ being even and the superalgebra $\QWA_n(a;z)$ from \cref{tiger},
    modulo the relations
    \begin{align}
        \sigma_i x_j &= x_j \sigma_i, & 1 \le i \le n-1,\ 1 \le j \le n,\ j \ne i,i+1, \\
        \sigma_i x_i \sigma_i &= x_{i+1}, &1 \le i \le n-1, \\
        x_i \ba &= \ba x_i, &1 \le i \le n,\ \ba \in A^{\otimes n}.
    \end{align}
    We adopt the convention that $\QAWA_0(A;z) := \kk$.
\end{defin}

\begin{rem}
When $A = \kk = \C(q)$ and $z = q-q^{-1}$, then $\QWA_n(A;z)$ (resp.\ $\QAWA_n(A;z)$) is the Iwahori--Hecke algebra (resp.\ affine Hecke algebra) of type $A_{n-1}$.  More generally, let $C_d$ be a cyclic group of order $d$.  If $\kk = \C(q)$, $z=(q-q^{-1})/d$, and $A = \kk C_d$, with trace map given by projection onto the identity element of the group, then $\QWA_n(\kk C_d;z)$ (resp.\ $\QAWA_n(\kk C_d;z)$) is the Yokonuma--Hecke algebra (resp.\ affine Yokonuma--Hecke algebra) from \cite{CPd14}.
\end{rem}

\begin{rem}\label{Frobop}
    The opposite superalgebra $A^\op$ is again a symmetric Frobenius superalgebra, with the trace map being the same underlying linear map as for $A$.  It follows that the elements $\tau_i$ in $\QWA_n(A^\op)$ or $\QAWA_n(A^\op)$ are given by the exactly the same formula \cref{tau} as for $\QWA_n(A)$ or $\QAWA_n(A)$.
\end{rem}

\subsection*{The quantum wreath product category}

Define the \emph{quantum wreath product category} $\QW(A;z)$ to be the strict $\kk$-linear monoidal supercategory generated by one object $\uparrow$ and morphisms
\begin{equation}
    \begin{tikzpicture}[centerzero]
        \draw[->] (0.2,-0.2) -- (-0.2,0.2);
        \draw[wipe] (-0.2,-0.2) -- (0.2,0.2);
        \draw[->] (-0.2,-0.2) -- (0.2,0.2);
    \end{tikzpicture}
    \ ,\
    \begin{tikzpicture}[centerzero]
        \draw[->] (-0.2,-0.2) -- (0.2,0.2);
        \draw[wipe] (0.2,-0.2) -- (-0.2,0.2);
        \draw[->] (0.2,-0.2) -- (-0.2,0.2);
    \end{tikzpicture}
    \ \colon
    \uparrow \otimes \uparrow \to \uparrow \otimes \uparrow,
    \qquad
    \begin{tikzpicture}[anchorbase]
        \draw[->] (0,-0.2) -- (0,0.2);
        \token{0,0}{west}{a};
    \end{tikzpicture}
    \colon \uparrow \to \uparrow,\quad a \in A,
\end{equation}
where the crossings are even and the parity of the morphism $\tokup$ is the same as the parity of $a$.  We refer to the generators above as a \emph{positive crossing}, a \emph{negative crossing}, and a \emph{token}, respectively.  The relations are as follows:
\begin{gather} \label{tokrel}
    \begin{tikzpicture}[anchorbase]
        \draw[->] (0,0) -- (0,0.7);
        \token{0,0.35}{west}{1};
    \end{tikzpicture}
    =
    \begin{tikzpicture}[anchorbase]
        \draw[->] (0,0) -- (0,0.7);
    \end{tikzpicture}
    ,\quad
    \lambda\
    \begin{tikzpicture}[anchorbase]
        \draw[->] (0,0) -- (0,0.7);
        \token{0,0.35}{west}{a};
    \end{tikzpicture}
    + \mu\
    \begin{tikzpicture}[anchorbase]
        \draw[->] (0,0) -- (0,0.7);
        \token{0,0.35}{west}{b};
    \end{tikzpicture}
    =
    \begin{tikzpicture}[anchorbase]
        \draw[->] (0,0) -- (0,0.7);
        \token{0,0.35}{west}{\lambda a + \mu b};
    \end{tikzpicture}
    ,\quad
    \begin{tikzpicture}[anchorbase]
        \draw[->] (0,0) -- (0,0.7);
        \token{0,0.2}{east}{b};
        \token{0,0.45}{east}{a};
    \end{tikzpicture}
    =
    \begin{tikzpicture}[anchorbase]
        \draw[->] (0,0) -- (0,0.7);
        \token{0,0.35}{west}{ab};
    \end{tikzpicture}
    \ ,\quad a,b \in A,\ \lambda,\mu \in \kk,
    \\ \label{braid}
    \begin{tikzpicture}[anchorbase]
        \draw[->] (0.2,-0.5) to[out=up,in=down] (-0.2,0) to[out=up,in=down] (0.2,0.5);
        \draw[wipe] (-0.2,-0.5) to[out=up,in=down] (0.2,0) to[out=up,in=down] (-0.2,0.5);
        \draw[->] (-0.2,-0.5) to[out=up,in=down] (0.2,0) to[out=up,in=down] (-0.2,0.5);
    \end{tikzpicture}
    \ =\
    \begin{tikzpicture}[anchorbase]
        \draw[->] (-0.2,-0.5) -- (-0.2,0.5);
        \draw[->] (0.2,-0.5) -- (0.2,0.5);
    \end{tikzpicture}
    \ =\
    \begin{tikzpicture}[anchorbase]
        \draw[->] (-0.2,-0.5) to[out=up,in=down] (0.2,0) to[out=up,in=down] (-0.2,0.5);
        \draw[wipe] (0.2,-0.5) to[out=up,in=down] (-0.2,0) to[out=up,in=down] (0.2,0.5);
        \draw[->] (0.2,-0.5) to[out=up,in=down] (-0.2,0) to[out=up,in=down] (0.2,0.5);
    \end{tikzpicture}
    \ ,\quad
    \begin{tikzpicture}[anchorbase]
        \draw[->] (0.4,-0.5) -- (-0.4,0.5);
        \draw[wipe] (0,-0.5) to[out=up, in=down] (-0.4,0) to[out=up,in=down] (0,0.5);
        \draw[->] (0,-0.5) to[out=up, in=down] (-0.4,0) to[out=up,in=down] (0,0.5);
        \draw[wipe] (-0.4,-0.5) -- (0.4,0.5);
        \draw[->] (-0.4,-0.5) -- (0.4,0.5);
    \end{tikzpicture}
    \ =\
    \begin{tikzpicture}[anchorbase]
        \draw[->] (0.4,-0.5) -- (-0.4,0.5);
        \draw[wipe] (0,-0.5) to[out=up, in=down] (0.4,0) to[out=up,in=down] (0,0.5);
        \draw[->] (0,-0.5) to[out=up, in=down] (0.4,0) to[out=up,in=down] (0,0.5);
        \draw[wipe] (-0.4,-0.5) -- (0.4,0.5);
        \draw[->] (-0.4,-0.5) -- (0.4,0.5);
    \end{tikzpicture}
    \ ,\quad
    \begin{tikzpicture}[centerzero]
        \draw[->] (0.3,-0.4) -- (-0.3,0.4);
        \draw[wipe] (-0.3,-0.4) -- (0.3,0.4);
        \draw[->] (-0.3,-0.4) -- (0.3,0.4);
        \token{-0.15,-0.2}{east}{a};
    \end{tikzpicture}
    \ =\
    \begin{tikzpicture}[centerzero]
        \draw[->] (0.3,-0.4) -- (-0.3,0.4);
        \draw[wipe] (-0.3,-0.4) -- (0.3,0.4);
        \draw[->] (-0.3,-0.4) -- (0.3,0.4);
        \token{0.15,0.2}{west}{a};
    \end{tikzpicture}
    \ ,\quad
    \begin{tikzpicture}[centerzero]
        \draw[->] (-0.3,-0.4) -- (0.3,0.4);
        \draw[wipe] (0.3,-0.4) -- (-0.3,0.4);
        \draw[->] (0.3,-0.4) -- (-0.3,0.4);
        \token{-0.15,-0.2}{east}{a};
    \end{tikzpicture}
    \ =\
    \begin{tikzpicture}[centerzero]
        \draw[->] (-0.3,-0.4) -- (0.3,0.4);
        \draw[wipe] (0.3,-0.4) -- (-0.3,0.4);
        \draw[->] (0.3,-0.4) -- (-0.3,0.4);
        \token{0.15,0.2}{west}{a};
    \end{tikzpicture}
    \ ,
    \\ \label{preskein}
    \begin{tikzpicture}[centerzero]
        \draw[->] (0.3,-0.3) -- (-0.3,0.3);
        \draw[wipe] (-0.3,-0.3) -- (0.3,0.3);
        \draw[->] (-0.3,-0.3) -- (0.3,0.3);
    \end{tikzpicture}
    \ -\
    \begin{tikzpicture}[centerzero]
        \draw[->] (-0.3,-0.3) -- (0.3,0.3);
        \draw[wipe] (0.3,-0.3) -- (-0.3,0.3);
        \draw[->] (0.3,-0.3) -- (-0.3,0.3);
    \end{tikzpicture}
    = z\
    \begin{tikzpicture}[centerzero]
        \draw[->] (-0.2,-0.3) -- (-0.2,0.3);
        \draw[->] (0.2,-0.3) -- (0.2,0.3);
        \token{-0.2,0}{east}{b};
        \token{0.2,0}{west}{b^\vee};
    \end{tikzpicture}.
\end{gather}
Of course \cref{preskein} should be interpreted using the usual summation convention; we call it the \emph{Frobenius skein relation}.  The relations \cref{tokrel} imply that the map
\[
    A \to \End_{\QW(A;z)}(\uparrow),\quad
    a \mapsto \tokup,
\]
is a superalgebra homomorphism.  Using \cref{braid,intlaw}, it also follows that
\begin{equation} \label{couch}
    \begin{tikzpicture}[centerzero]
        \draw[->] (0.3,-0.4) -- (-0.3,0.4);
        \draw[wipe] (-0.3,-0.4) -- (0.3,0.4);
        \draw[->] (-0.3,-0.4) -- (0.3,0.4);
        \token{0.15,-0.2}{west}{a};
    \end{tikzpicture}
    \ =\
    \begin{tikzpicture}[centerzero]
        \draw[->] (0.3,-0.4) -- (-0.3,0.4);
        \draw[wipe] (-0.3,-0.4) -- (0.3,0.4);
        \draw[->] (-0.3,-0.4) -- (0.3,0.4);
        \token{-0.147,0.194}{east}{a};
    \end{tikzpicture}
    \ ,\quad
    \begin{tikzpicture}[centerzero]
        \draw[->] (-0.3,-0.4) -- (0.3,0.4);
        \draw[wipe] (0.3,-0.4) -- (-0.3,0.4);
        \draw[->] (0.3,-0.4) -- (-0.3,0.4);
        \token{0.15,-0.2}{west}{a};
    \end{tikzpicture}
    \ =\
    \begin{tikzpicture}[centerzero]
        \draw[->] (-0.3,-0.4) -- (0.3,0.4);
        \draw[wipe] (0.3,-0.4) -- (-0.3,0.4);
        \draw[->] (0.3,-0.4) -- (-0.3,0.4);
        \token{-0.147,0.194}{east}{a};
    \end{tikzpicture}
    \ ,\quad
    \begin{tikzpicture}[centerzero]
        \draw[->] (-0.2,-0.4) -- (-0.2,0.4);
        \draw[->] (0.2,-0.4) -- (0.2,0.4);
        \token{-0.2,0.15}{east}{a};
        \token{0.2,-0.15}{west}{b};
    \end{tikzpicture}
    \ = (-1)^{\bar{a} \bar{b}}
    \begin{tikzpicture}[centerzero]
        \draw[->] (-0.2,-0.4) -- (-0.2,0.4);
        \draw[->] (0.2,-0.4) -- (0.2,0.4);
        \token{-0.2,-0.15}{east}{a};
        \token{0.2,0.15}{west}{b};
    \end{tikzpicture}
    \ .
\end{equation}

We introduce the \emph{teleporters}
\begin{equation}\label{brexit}
    \begin{tikzpicture}[anchorbase]
        \draw[->] (0,-0.4) --(0,0.4);
        \draw[->] (0.5,-0.4) -- (0.5,0.4);
        \teleport{0,0}{0.5,0};
    \end{tikzpicture}
    =
    \begin{tikzpicture}[anchorbase]
        \draw[->] (0,-0.4) --(0,0.4);
        \draw[->] (0.5,-0.4) -- (0.5,0.4);
        \teleport{0,0.2}{0.5,-0.2};
    \end{tikzpicture}
    =
    \begin{tikzpicture}[anchorbase]
        \draw[->] (0,-0.4) --(0,0.4);
        \draw[->] (0.5,-0.4) -- (0.5,0.4);
        \teleport{0,-0.2}{0.5,0.2};
    \end{tikzpicture}
    := z\
    \begin{tikzpicture}[anchorbase]
        \draw[->] (0,-0.4) --(0,0.4);
        \draw[->] (0.5,-0.4) -- (0.5,0.4);
        \token{0,0.15}{east}{b};
        \token{0.5,-0.15}{west}{b^\vee};
    \end{tikzpicture}
    \overset{\cref{doubledual}}{=} z\
    \begin{tikzpicture}[anchorbase]
        \draw[->] (0,-0.4) --(0,0.4);
        \draw[->] (0.5,-0.4) -- (0.5,0.4);
        \token{0,-0.15}{east}{b^\vee};
        \token{0.5,0.15}{west}{b};
    \end{tikzpicture}\ .
\end{equation}
We do not insist that the tokens in a teleporter are drawn at the same horizontal level, the convention when this is not the case being that $b$ is on the higher of the tokens and $b^\vee$ is on the lower one.  We will also draw teleporters in larger diagrams.  When doing so, we add a sign of $(-1)^{y \bar b}$ in front of the $b$ summand in \cref{brexit}, where $y$ is the sum of the parities of all morphisms in the diagram vertically between the tokens labeled $b$ and $b^\vee$.  For example,
\[
    \begin{tikzpicture}[anchorbase]
        \draw[->] (-1,-0.4) -- (-1,0.4);
        \draw[->] (-0.5,-0.4) -- (-0.5,0.4);
        \draw[->] (0,-0.4) -- (0,0.4);
        \draw[->] (0.5,-0.4) -- (0.5,0.4);
        \token{-1,0}{east}{a};
        \token{0,0.1}{west}{c};
        \teleport{-0.5,0.2}{0.5,-0.2};
    \end{tikzpicture}
    =
    (-1)^{(\bar a + \bar c) \bar b}z\
    \begin{tikzpicture}[anchorbase]
        \draw[->] (-1,-0.4) -- (-1,0.4);
        \draw[->] (-0.5,-0.4) -- (-0.5,0.4);
        \draw[->] (0,-0.4) -- (0,0.4);
        \draw[->] (0.5,-0.4) -- (0.5,0.4);
        \token{-1,0}{east}{a};
        \token{0,0.1}{west}{c};
        \token{-0.5,0.2}{east}{b};
        \token{0.5,-0.2}{west}{b^\vee};
    \end{tikzpicture}
    \ .
\]
This convention ensures that one can slide the endpoints of teleporters along strands:
\[
    \begin{tikzpicture}[anchorbase]
        \draw[->] (-1,-0.4) -- (-1,0.4);
        \draw[->] (-0.5,-0.4) -- (-0.5,0.4);
        \draw[->] (0,-0.4) -- (0,0.4);
        \draw[->] (0.5,-0.4) -- (0.5,0.4);
        \token{-1,0}{east}{a};
        \token{0,0.1}{west}{c};
        \teleport{-0.5,0.2}{0.5,-0.2};
    \end{tikzpicture}
    =
    \begin{tikzpicture}[anchorbase]
        \draw[->] (-1,-0.4) -- (-1,0.4);
        \draw[->] (-0.5,-0.4) -- (-0.5,0.4);
        \draw[->] (0,-0.4) -- (0,0.4);
        \draw[->] (0.5,-0.4) -- (0.5,0.4);
        \token{-1,0}{east}{a};
        \token{0,0.1}{east}{c};
        \teleport{-0.5,-0.2}{0.5,0.2};
    \end{tikzpicture}
    =
    \begin{tikzpicture}[anchorbase]
        \draw[->] (-1,-0.4) -- (-1,0.4);
        \draw[->] (-0.5,-0.4) -- (-0.5,0.4);
        \draw[->] (0,-0.4) -- (0,0.4);
        \draw[->] (0.5,-0.4) -- (0.5,0.4);
        \token{-1,0}{east}{a};
        \token{0,0.1}{west}{c};
        \teleport{-0.5,0.2}{0.5,0.2};
    \end{tikzpicture}
    =
    \begin{tikzpicture}[anchorbase]
        \draw[->] (-1,-0.4) -- (-1,0.4);
        \draw[->] (-0.5,-0.4) -- (-0.5,0.4);
        \draw[->] (0,-0.4) -- (0,0.4);
        \draw[->] (0.5,-0.4) -- (0.5,0.4);
        \token{-1,0}{east}{a};
        \token{0,0.1}{west}{c};
        \teleport{-0.5,-0.2}{0.5,-0.2};
    \end{tikzpicture}
    \ .
\]
Using teleporters, the Frobenius skein relation \cref{preskein} can be written as
\begin{equation} \label{skein}
    \begin{tikzpicture}[centerzero]
        \draw[->] (0.3,-0.3) -- (-0.3,0.3);
        \draw[wipe] (-0.3,-0.3) -- (0.3,0.3);
        \draw[->] (-0.3,-0.3) -- (0.3,0.3);
    \end{tikzpicture}
    \ -\
    \begin{tikzpicture}[centerzero]
        \draw[->] (-0.3,-0.3) -- (0.3,0.3);
        \draw[wipe] (0.3,-0.3) -- (-0.3,0.3);
        \draw[->] (0.3,-0.3) -- (-0.3,0.3);
    \end{tikzpicture}
    \ = \
    \begin{tikzpicture}[anchorbase]
        \draw[->] (-0.2,-0.3) -- (-0.2,0.3);
        \draw[->] (0.2,-0.3) -- (0.2,0.3);
        \teleport{-0.2,0}{0.2,0};
    \end{tikzpicture}
\end{equation}
It follows from \cref{beam} that tokens can ``teleport'' across teleporters (justifying the terminology) in the sense that, for $a \in A$, we have
\begin{equation} \label{teleport}
    \begin{tikzpicture}[anchorbase]
        \draw[->] (0,-0.5) --(0,0.5);
        \draw[->] (0.5,-0.5) -- (0.5,0.5);
        \token{0.5,-0.25}{west}{a};
        \teleport{0,0}{0.5,0};
    \end{tikzpicture}
    =
    \begin{tikzpicture}[anchorbase]
        \draw[->] (0,-0.5) --(0,0.5);
        \draw[->] (0.5,-0.5) -- (0.5,0.5);
        \token{0,0.25}{east}{a};
        \teleport{0,0}{0.5,0};
    \end{tikzpicture}
    \ ,\qquad
    \begin{tikzpicture}[anchorbase]
        \draw[->] (0,-0.5) --(0,0.5);
        \draw[->] (0.5,-0.5) -- (0.5,0.5);
        \token{0,-0.25}{east}{a};
        \teleport{0,0}{0.5,0};
    \end{tikzpicture}
    =
    \begin{tikzpicture}[anchorbase]
        \draw[->] (0,-0.5) --(0,0.5);
        \draw[->] (0.5,-0.5) -- (0.5,0.5);
        \token{0.5,0.25}{west}{a};
        \teleport{0,0}{0.5,0};
    \end{tikzpicture}
    \ ,
\end{equation}
where the strings can occur anywhere in a diagram (i.e.\ they do not need to be adjacent).  The endpoints of teleporters slide through crossings and they can teleport too.  For example we have
\begin{equation} \label{laser}
  \begin{tikzpicture}[anchorbase]
    \draw[->] (0.4,-0.4) to (-0.4,0.4);
    \draw[wipe] (-0.4,-0.4) to (0.4,0.4);
    \draw[->] (-0.4,-0.4) to (0.4,0.4);
    \draw[->] (0.8,-0.4) to (0.8,0.4);
    \teleport{-0.17,0.17}{0.8,0.17};
  \end{tikzpicture}
  \ =\
  \begin{tikzpicture}[anchorbase]
    \draw[->] (0.4,-0.4) to (-0.4,0.4);
    \draw[wipe] (-0.4,-0.4) to (0.4,0.4);
    \draw[->] (-0.4,-0.4) to (0.4,0.4);
    \draw[->] (0.8,-0.4) to (0.8,0.4);
    \teleport{0.2,-0.2}{0.8,-0.2};
  \end{tikzpicture}
  \ ,\qquad
  \begin{tikzpicture}[anchorbase]
    \draw[->] (-0.4,-0.5) to (-0.4,0.5);
    \draw[->] (0,-0.5) to (0,0.5);
    \draw[->] (0.4,-0.5) to (0.4,0.5);
    \teleport{-0.4,-0.25}{0,-0.25};
    \teleport{0,0}{0.4,0};
  \end{tikzpicture}
  \ =\
  \begin{tikzpicture}[anchorbase]
    \draw[->] (-0.4,-0.5) to (-0.4,0.5);
    \draw[->] (0,-0.5) to (0,0.5);
    \draw[->] (0.4,-0.5) to (0.4,0.5);
    \teleport{-0.4,-0.2}{0.4,0.3};
    \teleport{0,-.1}{0.4,-.1};
  \end{tikzpicture}
  \ =\
  \begin{tikzpicture}[anchorbase]
    \draw[->] (-0.4,-0.5) to (-0.4,0.5);
    \draw[->] (0,-0.5) to (0,0.5);
    \draw[->] (0.4,-0.5) to (0.4,0.5);
    \teleport{-0.4,0.15}{0.4,0.15};
    \teleport{0,-.1}{0.4,-.1};
  \end{tikzpicture}
    \ .
\end{equation}
For more discussion of teleporters, see \cite[$\S$4]{BSW-foundations}, which already adopted all of these conventions.

The objects of $\QW(A;z)$ are $\{\uparrow^{\otimes n} : n \in \N\}$.  There are no nonzero morphisms $\uparrow^{\otimes m} \to \uparrow^{\otimes n}$ for $m \ne n$.  Furthermore, we have an isomorphism of superalgebras
\[
    \QWA_n(A;z) \xrightarrow{\cong} \End_{\QW(A;z)} (\uparrow^{\otimes n})
\]
sending $a^{(i)}$, $a \in A$, to a token labelled $a$ on the $i$-th strand and $\sigma_i$ to a positive crossing of the $i$-th and $(i+1)$-st strands.  Here and throughout the paper, we number strings \emph{from right to left}.

\subsection*{The quantum affine wreath product category}

We define the \emph{quantum affine wreath product category} $\QAW(A;z)$ to be the strict $\kk$-linear monoidal supercategory obtained by adjoining to $\QW(A;z)$ an additional invertible even morphism
\[
    \begin{tikzpicture}[anchorbase]
        \draw[->] (0,0) -- (0,0.6);
        \singdot{0,0.3};
    \end{tikzpicture}
    \colon \uparrow \to \uparrow
\]
and imposing the additional relations
\begin{equation} \label{QAWC}
    \begin{tikzpicture}[centerzero]
        \draw[->] (0.3,-0.4) -- (-0.3,0.4);
        \draw[wipe] (-0.3,-0.4) -- (0.3,0.4);
        \draw[->] (-0.3,-0.4) -- (0.3,0.4);
        \singdot{0.171,-0.228};
    \end{tikzpicture}
    \ =\
    \begin{tikzpicture}[centerzero]
        \draw[->] (-0.3,-0.4) -- (0.3,0.4);
        \draw[wipe] (0.3,-0.4) -- (-0.3,0.4);
        \draw[->] (0.3,-0.4) -- (-0.3,0.4);
        \singdot{-0.15,0.2};
    \end{tikzpicture}
    \ ,\quad
    \begin{tikzpicture}[centerzero]
        \draw[->] (-0.3,-0.4) -- (0.3,0.4);
        \draw[wipe] (0.3,-0.4) -- (-0.3,0.4);
        \draw[->] (0.3,-0.4) -- (-0.3,0.4);
        \singdot{-0.171,-0.228};
    \end{tikzpicture}
    \ =\
    \begin{tikzpicture}[centerzero]
        \draw[->] (0.3,-0.4) -- (-0.3,0.4);
        \draw[wipe] (-0.3,-0.4) -- (0.3,0.4);
        \draw[->] (-0.3,-0.4) -- (0.3,0.4);
        \singdot{0.15,0.2};
    \end{tikzpicture}
    \ ,\quad
    \begin{tikzpicture}[centerzero]
        \draw[->] (0,-0.4) -- (0,0.4);
        \token{0,0.15}{west}{a};
        \singdot{0,-0.15};
    \end{tikzpicture}
    \ =\
    \begin{tikzpicture}[anchorbase]
        \draw[->] (0,-0.4) -- (0,0.4);
        \token{0,-0.15}{west}{a};
        \singdot{0,0.15};
    \end{tikzpicture}
    \ ,\quad a \in A.
\end{equation}
Note that, in fact, one only needs to impose one of the first two relations in \cref{QAWC}; the other then follows by composing on the top and bottom by the appropriate crossings.  It also follows that endpoints of teleporters pass through dots:
\[
  \begin{tikzpicture}[anchorbase]
    \draw[->] (-0.2,-0.5) to (-0.2,0.5);
    \draw[->] (0.2,-0.5) to (0.2,0.5);
    \singdot{-0.2,0};
    \teleport{-0.2,-0.25}{0.2,0};
  \end{tikzpicture}
  \ =\
  \begin{tikzpicture}[anchorbase]
    \draw[->] (-0.2,-0.5) to (-0.2,0.5);
    \draw[->] (0.2,-0.5) to (0.2,0.5);
    \singdot{-0.2,0};
    \teleport{-0.2,0.25}{0.2,0};
  \end{tikzpicture}
  \ .
\]

The objects of $\QAW(A;z)$ are $\{\uparrow^{\otimes n} : n \in \N\}$.  There are no nonzero morphisms $\uparrow^{\otimes m} \to \uparrow^{\otimes n}$ for $m \ne n$.  Furthermore, we have an isomorphism of superalgebras
\[
    \QAWA_n(A;z) \xrightarrow{\cong} \End_{\QAW(A;z)} (\uparrow^{\otimes n}),
\]
sending $a^{(i)}$, $a \in A$, to a token labeled $a$ on the $i$-th strand, $x_i$ to a dot on the $i$-th strand, and $\sigma_i$ to a positive crossing of the $i$-th and $(i+1)$-st strands.

\begin{rem}
    Throughout this paper, we assume that the Frobenius superalgebra $A$ is symmetric for simplicity of exposition.  We expect that all the results can be extended to the case where $A$ is a general Frobenius superalgebra.  In this setup, the Nakayama automorphism $\psi$ will appear in some of the relations.  For example, the last relation in \cref{QAWC} becomes
    \begin{equation}
        \begin{tikzpicture}[centerzero]
            \draw[->] (0,-0.4) -- (0,0.4);
            \token{0,0.15}{west}{a};
            \singdot{0,-0.15};
        \end{tikzpicture}
        \ =\
        \begin{tikzpicture}[anchorbase]
            \draw[->] (0,-0.4) -- (0,0.4);
            \token{0,-0.15}{west}{\psi(a)};
            \singdot{0,0.15};
        \end{tikzpicture}
    \end{equation}
    and one must also apply a power of $\psi$ to tokens as they travel over the left caps and cups to be introduced below.  We refer the reader to \cite{Sav19} for a treatment of Frobenius Heisenberg categories in this level of generality.  The main case of interest that our assumption on $A$ excludes is the Clifford superalgebra, where the quantum affine wreath product algebra would correspond to the affine Hecke--Clifford superalgebra.  However, this case is more naturally treated by considering an odd affinization of the quantum wreath product algebra and the resulting Heisenberg category.  We hope to explore this in future work.
\end{rem}

\section{First approach\label{sec:first}}

In this section we give our first definition of the \emph{quantum Frobenius Heisenberg category} $\Heis_k(A;z,t)$.  We will give two more, equivalent definitions in \cref{sec:second,sec:third}.


We begin by adjoining a right dual $\downarrow$ to the object $\uparrow$.  Thus we have additional generating morphisms
\[

    \ \colon \uparrow \otimes \downarrow \to \one
\]
subject to the right adjunction relations
\begin{equation} \label{rightadj}
    %
\ .
\end{equation}

We define the downward dots and tokens
\begin{equation}
    %
    \ ,\quad a \in A.
\end{equation}
Using \cref{rightadj}, we immediately have that
\begin{equation}
    A^\op \to \End(\downarrow),\quad a^\op \mapsto \tokdown,
\end{equation}
is a homomorphism of superalgebras, i.e.
\begin{equation} \label{antimatter}
    %
    ,\quad a,b \in A,
\end{equation}
We introduce teleporters on downward strings as in \cref{brexit} (with the orientation of strings reversed); see \cref{Frobop}.  Tokens teleport across these new downward teleporters in just the same way as in \cref{teleport}.  We may also draw teleporters with endpoints on oppositely oriented strings, interpreting them as usual by putting the label $b$ on the higher of the tokens and $b^\vee$ on the lower one, summing over all $b \in \BA$. The way that dots teleport across these is slightly different; for example, we have that
\[
    %
    \ .
\]

We define the positive and negative right crossings by
\begin{equation} \label{rcross}
    %
    \ .
\end{equation}
We then define positive and negative downwards crossings by
\begin{equation}
    %
    \ .
\end{equation}
It follows that the rightwards and downwards Frobenius skein relations hold:

\noindent\begin{minipage}{0.5\linewidth}
    \begin{equation} \label{rskein}
        %
        \ .
    \end{equation}
\end{minipage}\par\vspace{\belowdisplayskip}

\noindent In addition, for all $a \in A$, we have
\begin{align} \label{rcslide}
    %
    \ .
\end{align}
Furthermore, attaching right caps to the top and right cups to the bottom of the relations \cref{braid,couch,QAWC} gives the following relations for all $a \in A$:
\begin{gather} \label{spit1}
    %
\]
denote the composition of $n$ dots if $n \ge 0$ and the composition of $|n|$ inverse dots if $n < 0$.

\begin{lem} \label{dotslide}
    The following relations hold for $n \in \Z$:
    \begin{align} \label{teapos}
        %
            & \text{if } n < 0.
        \end{cases}
    \end{align}
\end{lem}

\begin{proof}
    For $n \ge 0$, this follows from repeated application of \cref{QAWC,skein}.  The cases for $n < 0$ then follow from the $n > 0$ cases by composing on the top and bottom with an appropriate number of negative dots.
\end{proof}

Using \cref{rdotcross}, we have rightwards and downwards analogues of \cref{dotslide}.  Similarly, after proving the leftwards dot slide relation \cref{ldotcross} below, we also have a leftwards analogue.  In what follows, we will simply refer to these rotated relations using the equation numbers of \cref{dotslide}.


Now we fix the index set
\[
    J := \{ \star \} \cup \{(r,b) : 0 \le r \le |-k| - 1,\ b \in \BA \}
\]
and we fix a total order on this set such that $\star < (r,b)$ for all $r,b$.  We can then naturally speak of $1 \times J$, $J \times 1$, and $J \times J$ matrices.

\begin{defin} \label{def1}
    The \emph{quantum Frobenius Heisenberg category} $\Heis_k(A;z,t)$ is the strict $\kk$-linear monoidal supercategory obtained from $\QAW(A;z,t)$ by adjoining a right dual $\downarrow$ to $\uparrow$, together with the relation that the following matrix of morphisms is invertible in the additive envelope of $\Heis_k(A;z,t)$:
    \begin{equation} \label{invrel}
        \begin{aligned}
            \left[
                \begin{tikzpicture}[anchorbase]
                    \draw[->] (0,0) -- (0.6,0.6);
                    \draw[wipe] (0.6,0) -- (0,0.6);
                    \draw[<-] (0.6,0) -- (0,0.6);
                \end{tikzpicture}
                \
                \begin{tikzpicture}[anchorbase]
                    \draw[->] (-0.2,-0.4) -- (-0.2,0) arc (180:0:0.2) -- (0.2,-0.4);
                    \multdot{-0.2,0}{east}{r};
                    \token{-0.2,-0.2}{east}{b^\vee};
                \end{tikzpicture}
                \ ,\ 0 \le r \le k-1,\ b \in \BA
            \right]^T
            \ &\colon \uparrow \downarrow \to \downarrow \uparrow \oplus \one^{\oplus k \dim A} \quad \text{if } k \ge 0,
            \\
            \left[
                \begin{tikzpicture}[anchorbase]
                    \draw[->] (0,0) -- (0.6,0.6);
                    \draw[wipe] (0.6,0) -- (0,0.6);
                    \draw[<-] (0.6,0) -- (0,0.6);
                \end{tikzpicture}
                \quad
                \begin{tikzpicture}[anchorbase]
                    \draw[->] (-0.2,0.4) -- (-0.2,0) arc (180:360:0.2) -- (0.2,0.4);
                    \multdot{0.2,0.2}{west}{r};
                    \token{0.2,0}{west}{b^\vee};
                \end{tikzpicture}
                \ ,\ 0 \le r \le -k-1,\ b \in \BA
            \right]
            \ &\colon \uparrow \downarrow \oplus \one^{\oplus (-k \dim A)} \to \downarrow \uparrow \quad \text{if } k < 0.
        \end{aligned}
    \end{equation}
    (The above matrices are of size $J \times 1$ and $1 \times J$ in the cases $k \ge 0$ and $k < 0$, respectively.)  We denote the matrix entries of the two-sided inverse of \cref{invrel} as follows:
    \begin{equation} \label{invrel2}
        \begin{aligned}
            \left[
                \begin{tikzpicture}[anchorbase]
                    \draw[->] (0.6,0) -- (0,0.6);
                    \draw[wipe] (0,0) -- (0.6,0.6);
                    \draw[<-] (0,0) -- (0.6,0.6);
                \end{tikzpicture}
                \
                \begin{tikzpicture}[centerzero={(0,-0.1)}]
                    \draw[<-] (-0.2,0.2) -- (-0.2,0) arc (180:360:0.2) -- (0.2,0.2);
                    \diamdec{0,-0.2}{north}{b}{r};
                \end{tikzpicture},\
                0 \le r \le k-1,\ b \in \BA
            \right]
            &=
            \left(
                \left[
                    \begin{tikzpicture}[anchorbase]
                        \draw[->] (0,0) -- (0.6,0.6);
                        \draw[wipe] (0.6,0) -- (0,0.6);
                        \draw[<-] (0.6,0) -- (0,0.6);
                    \end{tikzpicture}
                    \
                    \begin{tikzpicture}[anchorbase]
                        \draw[->] (-0.2,-0.4) -- (-0.2,0) arc (180:0:0.2) -- (0.2,-0.4);
                        \multdot{-0.2,0}{east}{r};
                        \token{-0.2,-0.25}{east}{b^\vee};
                    \end{tikzpicture}
                    \ ,\ 0 \le r \le k-1,\ b \in \BA
                \right]^T
            \right)^{-1}
            \text{ if } k \ge 0,
            \\
            \left[
                \begin{tikzpicture}[anchorbase]
                    \draw[<-] (0,0) -- (0.6,0.6);
                    \draw[wipe] (0.6,0) -- (0,0.6);
                    \draw[->] (0.6,0) -- (0,0.6);
                \end{tikzpicture}
                \quad
                \begin{tikzpicture}[centerzero={(0,0.1)}]
                    \draw[<-] (-0.2,-0.2) -- (-0.2,0) arc(180:0:0.2) -- (0.2,-0.2);
                    \heartdec{0,0.2}{south}{b}{r};
                \end{tikzpicture}
                ,\ 0 \le r \le -k-1,\ b \in \BA
            \right]^T
            &=
            \left[
                \begin{tikzpicture}[anchorbase]
                    \draw[->] (0,0) -- (0.6,0.6);
                    \draw[wipe] (0.6,0) -- (0,0.6);
                    \draw[<-] (0.6,0) -- (0,0.6);
                \end{tikzpicture}
                \quad
                \begin{tikzpicture}[anchorbase]
                    \draw[->] (-0.2,0.4) -- (-0.2,0) arc (180:360:0.2) -- (0.2,0.4);
                    \multdot{0.2,0.2}{west}{r};
                    \token{0.2,0}{west}{b^\vee};
                \end{tikzpicture}
                \ ,\ 0 \le r \le -k-1,\ b \in \BA
            \right]^{-1}
            \text{ if } k < 0.
        \end{aligned}
    \end{equation}
    We extend the definition of the decorated left cups and caps by linearity in the second argument of the label.  In other words, for $a \in A$, we define
    \[
        \begin{tikzpicture}[centerzero]
            \draw[<-] (-0.2,0.2) -- (-0.2,0) arc (180:360:0.2) -- (0.2,0.2);
            \diamdec{0,-0.2}{north}{a}{r};
        \end{tikzpicture}
        \, = \tr(b^\vee a)\
        \begin{tikzpicture}[centerzero]
            \draw[<-] (-0.2,0.2) -- (-0.2,0) arc (180:360:0.2) -- (0.2,0.2);
            \diamdec{0,-0.2}{north}{b}{r};
        \end{tikzpicture}
        \ ,\quad \text{if } k > 0,
        \qquad
        \begin{tikzpicture}[centerzero={(0,0.2)}]
            \draw[<-] (-0.2,-0.2) -- (-0.2,0) arc(180:0:0.2) -- (0.2,-0.2);
            \heartdec{0,0.2}{south}{a}{r};
        \end{tikzpicture}
        \ = \tr(b^\vee a)\
        \begin{tikzpicture}[centerzero={(0,0.2)}]
            \draw[<-] (-0.2,-0.2) -- (-0.2,0) arc(180:0:0.2) -- (0.2,-0.2);
            \heartdec{0,0.2}{south}{b}{r};
        \end{tikzpicture}
        \ ,\quad \text{if } k < 0.
    \]
    We define the left cup and cap by
    \begin{equation} \label{leftwards}
        \begin{tikzpicture}[anchorbase]
            \draw[<-] (-0.2,0.2) -- (-0.2,0) arc (180:360:0.2) -- (0.2,0.2);
        \end{tikzpicture}
        :=
        \begin{cases}
            -\displaystyle{\frac{t^{-1}}{z}}\
            \begin{tikzpicture}[centerzero={(0,-0.2)}]
                \draw[<-] (-0.2,0.2) -- (-0.2,0) arc (180:360:0.2) -- (0.2,0.2);
                \diamdec{0,-0.2}{north}{1}{k-1};
                \multdot{-0.2,0}{east}{-1};
            \end{tikzpicture}
            & \text{if } k > 0,
            \\
            t\
            \begin{tikzpicture}[anchorbase]
                \draw[<-] (-0.2,0.4) to[out=-45,in=up] (0.2,0) arc(360:180:0.2);
                \draw[wipe] (-0.2,0) to[out=up,in=225] (0.2,0.4);
                \draw (-0.2,0) to[out=up,in=225] (0.2,0.4);
            \end{tikzpicture}
            & \text{if } k = 0,
            \\
            t^{-1}\
            \begin{tikzpicture}[anchorbase]
                \draw (-0.2,0) to[out=up,in=225] (0.2,0.4);
                \draw[wipe] (-0.2,0.4) to[out=-45,in=up] (0.2,0) arc(360:180:0.2);
                \draw[<-] (-0.2,0.4) to[out=-45,in=up] (0.2,0) arc(360:180:0.2);
                \multdot{0.2,0}{west}{-k};
            \end{tikzpicture}
            & \text{if } k < 0,
        \end{cases}
        \qquad
        \begin{tikzpicture}[anchorbase]
            \draw[<-] (-0.2,-0.2) -- (-0.2,0) arc (180:0:0.2) -- (0.2,-0.2);
        \end{tikzpicture}
        :=
        \begin{cases}
            t\
            \begin{tikzpicture}[anchorbase]
                \draw (-0.2,0) to[out=down,in=135] (0.2,-0.4);
                \draw[wipe] (-0.2,-0.4) to[out=45,in=down] (0.2,0) arc(0:180:0.2);
                \draw[<-] (-0.2,-0.4) to[out=45,in=down] (0.2,0) arc(0:180:0.2);
                \multdot{0.2,0}{west}{k};
            \end{tikzpicture}
            & \text{if } k \ge 0, \\
            - \displaystyle\frac{t^{-1}}{z}
            \begin{tikzpicture}[anchorbase]
                \draw[<-] (-0.2,-0.2) -- (-0.2,0) arc (180:0:0.2) -- (0.2,-0.2);
                \heartdec{0,0.2}{south}{1}{0};
            \end{tikzpicture}
            & \text{if } k < 0.
        \end{cases}
    \end{equation}
    We then impose one additional relation:
    \begin{equation} \label{impose}
        \begin{tikzpicture}[centerzero]
            \draw[<-] (0,0.2) arc(90:450:0.2);
            \token{0.2,0}{west}{a};
        \end{tikzpicture}
        \ = \frac{t}{z} \tr(a) 1_\one \text{ if } k > 0,
        \quad
        \begin{tikzpicture}[centerzero]
            \draw[<-] (0,0.2) arc(90:450:0.2);
            \token{0.2,0}{west}{a};
        \end{tikzpicture}
        \ = \frac{t-t^{-1}}{z} \tr(a) 1_\one \text{ if } k = 0,
        \quad
        \ccbubble{a}{-k}
        \ = \frac{t}{z} \tr(a) 1_\one \text{ if } k < 0.
    \end{equation}
    This concludes the definition of $\Heis_k(A;z,t)$.
\end{defin}

One left crossing has been defined as the first entry in the matrix appearing in \cref{invrel2}.  We define the other left crossing so that the Frobenius skein relation holds:
\begin{equation} \label{lskein}
    \begin{tikzpicture}[anchorbase]
        \draw[<-] (-0.3,-0.4) -- (0.3,0.4);
        \draw[wipe] (0.3,-0.4) -- (-0.3,0.4);
        \draw[->] (0.3,-0.4) -- (-0.3,0.4);
    \end{tikzpicture}
    \ -\
    \begin{tikzpicture}[anchorbase]
        \draw[->] (0.3,-0.4) -- (-0.3,0.4);
        \draw[wipe] (-0.3,-0.4) -- (0.3,0.4);
        \draw[<-] (-0.3,-0.4) -- (0.3,0.4);
    \end{tikzpicture}
    = \
    \begin{tikzpicture}[anchorbase]
        \draw[<-] (0.1,0.4) -- (0.1,0.2) arc(180:360:0.2) -- (0.5,0.4);
        \draw[<-] (-0.5,-0.4) -- (-0.5,-0.2) arc(180:0:0.2) -- (-0.1,-0.4);
        \teleport{-0.1,-0.2}{0.1,0.2};
    \end{tikzpicture}
    \ .
\end{equation}
For $a \in A$, we also set
\begin{gather} \label{nakano1}
    \begin{tikzpicture}[>=To,baseline={([yshift=1ex]current bounding box.center)}]
        \draw[<-] (-0.2,0.2) -- (-0.2,0) arc (180:360:0.2) -- (0.2,0.2);
        \heartdec{0,-0.2}{north}{a}{0};
    \end{tikzpicture}
    :=
    \begin{tikzpicture}[>=To,baseline={([yshift=1ex]current bounding box.center)}]
        \draw[<-] (-0.2,0.2) -- (-0.2,0) arc (180:360:0.2) -- (0.2,0.2);
        \diamdec{0,-0.2}{north}{a}{0};
    \end{tikzpicture}
    + z\
    \begin{tikzpicture}[anchorbase]
        \draw[<-] (-0.2,0.4) to[out=-45,in=up] (0.2,0) arc(360:180:0.2);
        \draw[wipe] (-0.2,0) to[out=up,in=225] (0.2,0.4);
        \draw (-0.2,0) to[out=up,in=225] (0.2,0.4);
        \token{0.2,0}{west}{a};
    \end{tikzpicture}
    \ \text{if } k > 0,
    \qquad
    \begin{tikzpicture}[>=To,baseline={([yshift=1ex]current bounding box.center)}]
        \draw[<-] (-0.2,0.2) -- (-0.2,0) arc (180:360:0.2) -- (0.2,0.2);
        \heartdec{0,-0.2}{north}{a}{n};
    \end{tikzpicture}
    \ :=\
    \begin{tikzpicture}[>=To,baseline={([yshift=1ex]current bounding box.center)}]
        \draw[<-] (-0.2,0.2) -- (-0.2,0) arc (180:360:0.2) -- (0.2,0.2);
        \diamdec{0,-0.2}{north}{a}{n};
    \end{tikzpicture}
    \ \text{if } 0 < n < k,
    \\ \label{nakano3}
    \begin{tikzpicture}[>=To,baseline={([yshift=-2ex]current bounding box.center)}]
        \draw[<-] (-0.2,-0.2) -- (-0.2,0) arc (180:0:0.2) -- (0.2,-0.2);
        \diamdec{0,0.2}{south}{a}{0};
    \end{tikzpicture}
    \ :=\
    \begin{tikzpicture}[>=To,baseline={([yshift=-2ex]current bounding box.center)}]
        \draw[<-] (-0.2,-0.2) -- (-0.2,0) arc (180:0:0.2) -- (0.2,-0.2);
        \heartdec{0,0.2}{south}{a}{0};
    \end{tikzpicture}
    + (-1)^{\bar a} z\
    \begin{tikzpicture}[anchorbase]
        \draw[<-] (-0.2,-0.4) to[out=45,in=down] (0.2,0) arc(0:180:0.2);
        \draw[wipe] (-0.2,0) to[out=down,in=135] (0.2,-0.4);
        \draw (-0.2,0) to[out=down,in=135] (0.2,-0.4);
        \token{0.2,0}{west}{a};
    \end{tikzpicture}
    \ \text{if } k < 0,
    \qquad
    \begin{tikzpicture}[>=To,baseline={([yshift=-2ex]current bounding box.center)}]
        \draw[<-] (-0.2,-0.2) -- (-0.2,0) arc (180:0:0.2) -- (0.2,-0.2);
        \diamdec{0,0.2}{south}{a}{n};
    \end{tikzpicture}
    \ :=\
    \begin{tikzpicture}[>=To,baseline={([yshift=-2ex]current bounding box.center)}]
        \draw[<-] (-0.2,-0.2) -- (-0.2,0) arc (180:0:0.2) -- (0.2,-0.2);
        \heartdec{0,0.2}{south}{a}{n};
    \end{tikzpicture}
    \ \text{if } 0 < n < -k.
\end{gather}

For $n \le 0$ and $a \in A$, we define the \emph{$(+)$-bubbles}:
\begin{equation} \label{fake0}
    \plusccbubble{a}{n}
    :=
    \begin{cases}
        -\displaystyle{\frac{t}{z}}
        \begin{tikzpicture}[anchorbase]
            \draw[<-] (0,0.2) arc(90:450:0.2);
            \diamdec{0,-0.2}{north}{a}{-n};
            \multdot{-0.2,0}{east}{k};
        \end{tikzpicture}
        & \text{if } n > -k,
        \\
        \displaystyle{\frac{t}{z} \tr(a) 1_\one} & \text{if } n = -k,
        \\
        0 & \text{if } n < -k,
    \end{cases}
    \qquad
    \pluscbubble{a}{n}
    :=
    \begin{cases}
        \displaystyle{\frac{t^{-1}}{z}}
        \begin{tikzpicture}[centerzero]
            \draw[->] (0,-0.2) arc(-90:270:0.2);
            \diamdec{0,0.2}{south}{a}{-n};
            \multdot{0.2,0}{west}{-k};
        \end{tikzpicture}
        & \text{if } n > k,
        \\
        \displaystyle{-\frac{t^{-1}}{z} \tr(a) 1_\one} & \text{if } n = k,
        \\
        0 & \text{if } n < k.
    \end{cases}
\end{equation}
For $a \in A$, we define
\begin{align} \label{fake2}
    \pluscbubble{a}{n} &:= \cbubble{a}{n}
    \quad \text{if } n > 0,
    &
    \plusccbubble{a}{n} &:= \ccbubble{a}{n}
    \quad \text{if } n > 0.
\end{align}
Then we define the \emph{$(-)$-bubbles} so that, for all $a \in A$, $n \in \Z$,
\begin{align} \label{fake1}
    \cbubble{a}{n}
    &= \pluscbubble{a}{n} + \minuscbubble{a}{n}
    \ ,&
    \ccbubble{a}{n}
    &= \plusccbubble{a}{n} + \minusccbubble{a}{n}
    \quad \text{for all } n \in \Z.
\end{align}
It follows from the definitions that the $(\pm)$-bubbles are linear in the label $a \in A$.  Note that although our notation for the $(\pm)$-bubbles involves tokens and dots, these bubbles are not built from cups, caps, tokens, and dots in general.  Furthermore, we allow ourselves the freedom to draw the tokens and dots on $(\pm)$-bubbles in any position, and we use the relations \cref{tokrel,antimatter,rcslide} to interpret $(\pm)$-bubbles with multiple tokens.  For example,
\begin{equation}
    \begin{tikzpicture}[centerzero]
        \draw[->-=0.1] (0,0.2) arc(90:-270:0.2);
        \node at (0,0) {\dotlabel{+}};
        \multdot{0.2,0}{west}{n};
        \token{-0.141,0.141}{east}{a};
        \token{-0.141,-0.141}{east}{b};
    \end{tikzpicture}
    = \pluscbubble{ab}{n},
    \qquad
    \begin{tikzpicture}[centerzero]
        \draw[->] (0,-0.3) -- (0,0.3);
        \minusleftblank{0.5,0};
        \teleport{0,0}{0.3,0};
        \multdot{0.7,0}{west}{n};
    \end{tikzpicture}
    =
    \begin{tikzpicture}[centerzero]
        \draw[->] (0,-0.3) -- (0,0.3);
        \minusleftblank{0.9,0};
        \token{0,0}{east}{b};
        \token{0.7,0}{east}{b^\vee};
        \multdot{1.1,0}{west}{n};
    \end{tikzpicture}
    \ .
\end{equation}

When $k = 0$, the assertion that \cref{invrel,invrel2} are two-sided inverses means that
\begin{equation} \label{lunch}
    \begin{tikzpicture}[anchorbase]
        \draw[->] (0.2,-0.5) \braidup (-0.2,0) \braidup (0.2,0.5);
        \draw[wipe] (-0.2,-0.5) \braidup (0.2,0) \braidup (-0.2,0.5);
        \draw[<-] (-0.2,-0.5) \braidup (0.2,0) \braidup (-0.2,0.5);
    \end{tikzpicture}
    \ =\
    \begin{tikzpicture}[anchorbase]
        \draw[<-] (-0.2,-0.5) -- (-0.2,0.5);
        \draw[->] (0.2,-0.5) -- (0.2,0.5);
    \end{tikzpicture}
    \quad \text{if } k = 0,
    \qquad \qquad
    \begin{tikzpicture}[anchorbase]
        \draw[->] (-0.2,-0.5) \braidup (0.2,0) \braidup (-0.2,0.5);
        \draw[wipe] (0.2,-0.5) \braidup (-0.2,0) \braidup (0.2,0.5);
        \draw[<-] (0.2,-0.5) \braidup (-0.2,0) \braidup (0.2,0.5);
    \end{tikzpicture}
    \ =\
    \begin{tikzpicture}[anchorbase]
        \draw[->] (-0.2,-0.5) -- (-0.2,0.5);
        \draw[<-] (0.2,-0.5) -- (0.2,0.5);
    \end{tikzpicture}
    \quad \text{if } k = 0.
\end{equation}

When $k > 0$, the inversion relation implies that, for all $a \in A$,
\begin{gather} \label{tea1}
    \begin{tikzpicture}[anchorbase]
        \draw[->] (0.2,-0.5) \braidup (-0.2,0) \braidup (0.2,0.5);
        \draw[wipe] (-0.2,-0.5) \braidup (0.2,0) \braidup (-0.2,0.5);
        \draw[<-] (-0.2,-0.5) \braidup (0.2,0) \braidup (-0.2,0.5);
    \end{tikzpicture}
    \ =\
    \begin{tikzpicture}[anchorbase]
        \draw[<-] (-0.2,-0.5) -- (-0.2,0.5);
        \draw[->] (0.2,-0.5) -- (0.2,0.5);
    \end{tikzpicture}
    \quad \text{if } k > 0,
    \qquad \qquad
    \begin{tikzpicture}[anchorbase]
        \draw[->] (-0.2,-0.5) \braidup (0.2,0) \braidup (-0.2,0.5);
        \draw[wipe] (0.2,-0.5) \braidup (-0.2,0) \braidup (0.2,0.5);
        \draw[<-] (0.2,-0.5) \braidup (-0.2,0) \braidup (0.2,0.5);
    \end{tikzpicture}
    \ =\
    \begin{tikzpicture}[anchorbase]
        \draw[->] (-0.2,-0.5) -- (-0.2,0.5);
        \draw[<-] (0.2,-0.5) -- (0.2,0.5);
    \end{tikzpicture}
    - \sum_{r=0}^{k-1}
    \begin{tikzpicture}[anchorbase]
        \draw[->] (-0.2,-0.8) -- (-0.2,-0.4) arc(180:0:0.2) -- (0.2,-0.8);
        \draw[<-] (-0.2,0.6) -- (-0.2,0.4) arc(180:360:0.2) -- (0.2,0.6);
        \multdot{-0.2,-0.4}{west}{r};
        \token{-0.2,-0.6}{east}{b^\vee};
        \diamdec{0,0.2}{north}{b}{r};
    \end{tikzpicture}
    \ \text{if } k > 0,
    \\ \label{tea0}
    \begin{tikzpicture}[anchorbase]
        \draw (-0.2,0.4) to[out=-45,in=up] (0.2,0) arc(360:180:0.2);
        \draw[wipe] (-0.2,0) to[out=up,in=225] (0.2,0.4);
        \draw[->] (-0.2,0) to[out=up,in=225] (0.2,0.4);
        \token{-0.2,0}{east}{a};
    \end{tikzpicture}
    \ =\
    \begin{tikzpicture}[anchorbase]
        \draw (-0.2,0.4) to[out=-45,in=up] (0.2,0) arc(360:180:0.2);
        \draw[wipe] (-0.2,0) to[out=up,in=225] (0.2,0.4);
        \draw[->] (-0.2,0) to[out=up,in=225] (0.2,0.4);
        \token{0.2,0}{west}{a};
    \end{tikzpicture}
    =
    \begin{tikzpicture}[anchorbase]
        \draw (0.2,-0.4) to[out=135,in=down] (-0.2,0) -- (-0.2,0.2) arc(180:0:0.2) -- (0.2,0);
        \draw[wipe] (0.2,0) to[out=down,in=45] (-0.2,-0.4);
        \draw[->] (0.2,0) to[out=down,in=45] (-0.2,-0.4);
        \token{-0.2,0}{east}{a};
        \multdot{-0.2,0.2}{east}{r};
    \end{tikzpicture}
    \ =0,\ 0 \le r < k,
    \qquad
    \cbubble{a}{r}
    = -\delta_{r,k} \frac{t^{-1}}{z} \tr(a) 1_\one,\ 0 < r \le k.
\end{gather}
Recall, in \cref{tea1} and in analogous expressions below, our convention that there is an implicit sum over $b \in \BA$.

Using \cref{teapos,tea0,rskein,adecomp,impose} we also have
\begin{equation} \label{exercise}
    \begin{tikzpicture}[anchorbase]
        \draw[->] (-0.2,0) to[out=up,in=225] (0.2,0.4);
        \draw[wipe] (-0.2,0.4) to[out=-45,in=up] (0.2,0) arc(360:180:0.2);
        \draw (-0.2,0.4) to[out=-45,in=up] (0.2,0) arc(360:180:0.2);
        \multdot{-0.2,0}{east}{n};
    \end{tikzpicture}
    \ = \delta_{n,0} t\
    \begin{tikzpicture}[anchorbase]
        \draw[->] (-0.2,0.2) -- (-0.2,0) arc (180:360:0.2) -- (0.2,0.2);
    \end{tikzpicture}
    \quad \text{and} \quad
    \begin{tikzpicture}[anchorbase]
        \draw (-0.2,0.4) to[out=-45,in=up] (0.2,0) arc(360:180:0.2);
        \draw[wipe] (-0.2,0) to[out=up,in=225] (0.2,0.4);
        \draw[->] (-0.2,0) to[out=up,in=225] (0.2,0.4);
        \multdot{-0.2,0}{east}{n};
    \end{tikzpicture}
    \ = \delta_{n,k} t^{-1}\
    \begin{tikzpicture}[anchorbase]
        \draw[->] (-0.2,0.2) -- (-0.2,0) arc (180:360:0.2) -- (0.2,0.2);
    \end{tikzpicture}
    \quad \text{for } 0 \le n \le k.
\end{equation}
\details{
    Consider the first relation in \cref{exercise}.  If $k = 0$, we obtain this relation by attaching a positive right crossing to the top of the left cup, as defined by \cref{leftwards}, and then use \cref{lunch}.  Now suppose $k > 0$.  Then, for $0 < n \le k$, we have
    \[
        \begin{tikzpicture}[anchorbase]
            \draw[->] (-0.2,0) to[out=up,in=225] (0.2,0.4);
            \draw[wipe] (-0.2,0.4) to[out=-45,in=up] (0.2,0) arc(360:180:0.2);
            \draw (-0.2,0.4) to[out=-45,in=up] (0.2,0) arc(360:180:0.2);
            \multdot{-0.2,0}{east}{n};
        \end{tikzpicture}
        \stackrel{\cref{teapos}}{=}
        \begin{tikzpicture}[anchorbase]
            \draw (-0.2,0.6) -- (-0.2,0.4) to[out=down,in=up] (0.2,0) arc(360:180:0.2);
            \draw[wipe] (-0.2,0) to[out=up,in=down] (0.2,0.4) -- (0.2,0.6);
            \draw[->] (-0.2,0) to[out=up,in=down] (0.2,0.4) -- (0.2,0.6);
            \multdot{0.2,0.4}{west}{n};
        \end{tikzpicture}
        - z \sum_{\substack{r + s = n \\ r,s > 0}}
        \begin{tikzpicture}[anchorbase]
            \draw[->] (-0.2,0.4) to (-0.2,0.2) arc(180:360:0.2) to (0.2,0.4);
            \token{-0.2,0.2}{east}{b};
            \multdot{0.2,0.2}{west}{r};
            \bubright{0,-0.4}{b^\vee}{s};
        \end{tikzpicture}
        \stackrel{\cref{tea0}}{=} 0.
  \]
  If $n=0$, we have
  \[
        \begin{tikzpicture}[anchorbase]
            \draw[->] (-0.2,0) to[out=up,in=225] (0.2,0.4);
            \draw[wipe] (-0.2,0.4) to[out=-45,in=up] (0.2,0) arc(360:180:0.2);
            \draw (-0.2,0.4) to[out=-45,in=up] (0.2,0) arc(360:180:0.2);
        \end{tikzpicture}
        \stackrel{\cref{rskein}}{=}
        \begin{tikzpicture}[anchorbase]
            \draw (-0.2,0.4) to[out=-45,in=up] (0.2,0) arc(360:180:0.2);
            \draw[wipe] (-0.2,0) to[out=up,in=225] (0.2,0.4);
            \draw[->] (-0.2,0) to[out=up,in=225] (0.2,0.4);
        \end{tikzpicture}
        + z
        \begin{tikzpicture}[anchorbase]
            \draw[->] (-0.2,0.4) to (-0.2,0.2) arc(180:360:0.2) to (0.2,0.4);
            \token{0.2,0.2}{west}{b};
            \draw[->] (0,-0.2) arc(90:-270:0.2);
            \filldraw (-0.2,-0.4) circle (1.5pt) node[anchor=east] {\dotlabel{b^\vee}};
        \end{tikzpicture}
        \stackrel[\cref{tea0}]{\substack{\cref{impose} \\ \cref{adecomp}}}{=}
        t\
        \begin{tikzpicture}[anchorbase]
            \draw[->] (-0.2,0.2) -- (-0.2,0) arc (180:360:0.2) -- (0.2,0.2);
        \end{tikzpicture}
        \ .
    \]

    Now consider the second relation in \cref{exercise}.  When $k=0$, we use \cref{rskein} to flip the crossing, then use the first relation in \cref{exercise}, together with \cref{impose,adecomp}.  Now suppose $k>0$.  For $n=0$, this relation follows immediately from \cref{tea0}.  If $0 < n \le k$, we have
    \[
        \begin{tikzpicture}[anchorbase]
            \draw (-0.2,0.4) to[out=-45,in=up] (0.2,0) arc(360:180:0.2);
            \draw[wipe] (-0.2,0) to[out=up,in=225] (0.2,0.4);
            \draw[->] (-0.2,0) to[out=up,in=225] (0.2,0.4);
            \multdot{-0.2,0}{east}{n};
        \end{tikzpicture}
        \stackrel{\cref{rskein}}{=}
        \begin{tikzpicture}[anchorbase]
            \draw[->] (-0.2,0) to[out=up,in=225] (0.2,0.4);
            \draw[wipe] (-0.2,0.4) to[out=-45,in=up] (0.2,0) arc(360:180:0.2);
            \draw (-0.2,0.4) to[out=-45,in=up] (0.2,0) arc(360:180:0.2);
            \multdot{-0.2,0}{east}{n};
        \end{tikzpicture}
        - z
        \begin{tikzpicture}[anchorbase]
            \draw[->] (-0.2,0.4) to (-0.2,0.2) arc(180:360:0.2) to (0.2,0.4);
            \token{-0.2,0.2}{east}{b};
            \bubright{0,-0.4}{b^\vee}{n};
        \end{tikzpicture}
        \stackrel[\cref{adecomp}]{\cref{tea0}}{=}
        t^{-1} \delta_{n,k}\
        \begin{tikzpicture}[anchorbase]
            \draw[->] (-0.2,0.2) -- (-0.2,0) arc (180:360:0.2) -- (0.2,0.2);
        \end{tikzpicture}
        \ .
    \]
}

When $k < 0$, we have
\begin{gather} \label{tea3}
    \begin{tikzpicture}[anchorbase]
        \draw[<-] (-0.2,-0.5) to[out=up,in=down] (0.2,0);
        \draw[->] (-0.2,0) to[out=up,in=down] (0.2,0.5);
        \draw[wipe] (0.2,0) to[out=up,in=down] (-0.2,0.5);
        \draw (0.2,0) to[out=up,in=down] (-0.2,0.5);
        \draw[wipe] (0.2,-0.5) to[out=up,in=down] (-0.2,0);
        \draw (0.2,-0.5) to[out=up,in=down] (-0.2,0);
    \end{tikzpicture}
    \ =\
    \begin{tikzpicture}[anchorbase]
        \draw[<-] (-0.2,-0.5) -- (-0.2,0.5);
        \draw[->] (0.2,-0.5) -- (0.2,0.5);
    \end{tikzpicture}
    - \sum_{r=0}^{-k-1}
    \begin{tikzpicture}[anchorbase]
        \draw[<-] (-0.2,-0.6) -- (-0.2,-0.4) arc(180:0:0.2) -- (0.2,-0.6);
        \draw[->] (-0.2,0.8) -- (-0.2,0.4) arc(180:360:0.2) -- (0.2,0.8);
        \multdot{0.2,0.6}{east}{r};
        \token{0.2,0.4}{west}{b^\vee};
        \heartdec{0,-0.2}{south}{b}{r};
    \end{tikzpicture}
    \quad \text{if } k < 0,
    \qquad \qquad
    \begin{tikzpicture}[anchorbase]
        \draw (-0.2,-0.5) to[out=up,in=down] (0.2,0);
        \draw (-0.2,0) to[out=up,in=down] (0.2,0.5);
        \draw[wipe] (0.2,0) to[out=up,in=down] (-0.2,0.5);
        \draw[->] (0.2,0) to[out=up,in=down] (-0.2,0.5);
        \draw[wipe] (0.2,-0.5) to[out=up,in=down] (-0.2,0);
        \draw[<-] (0.2,-0.5) to[out=up,in=down] (-0.2,0);
    \end{tikzpicture}
    \ =\
    \begin{tikzpicture}[anchorbase]
        \draw[->] (-0.2,-0.5) -- (-0.2,0.5);
        \draw[<-] (0.2,-0.5) -- (0.2,0.5);
    \end{tikzpicture}
    \quad \text{if } k < 0,
    \\ \label{tea2}
    \begin{tikzpicture}[anchorbase]
        \draw (-0.2,-0.4) to[out=45,in=down] (0.2,0) arc(0:180:0.2);
        \draw[wipe] (-0.2,0) to[out=down,in=135] (0.2,-0.4);
        \draw[->] (-0.2,0) to[out=down,in=135] (0.2,-0.4);
    \end{tikzpicture}
    \ = 0 \text{ if } k < 0,
    \quad
    \begin{tikzpicture}[anchorbase]
        \draw (-0.2,0) to[out=up,in=225] (0.2,0.4);
        \draw[wipe] (-0.2,0.4) to[out=-45,in=up] (0.2,0) arc(360:180:0.2);
        \draw[<-] (-0.2,0.4) to[out=-45,in=up] (0.2,0) arc(360:180:0.2);
        \multdot{0.2,0}{west}{r};
        \token{-0.2,0}{east}{a};
    \end{tikzpicture}
    \ = 0 \text{ if } 0 \le r < -k,
    \quad
    \ccbubble{a}{r}
    = -\delta_{r,0} \frac{t^{-1}}{z} \tr(a) 1_\one \text{ if } 0 \le r < -k,
\end{gather}
for all $a \in A$.

\begin{lem}
    The following relations hold for all $a \in A$:
    \begin{equation} \label{ltokcross}
        \begin{tikzpicture}[centerzero]
            \draw[->] (0.3,-0.4) -- (-0.3,0.4);
            \draw[wipe] (-0.3,-0.4) -- (0.3,0.4);
            \draw[<-] (-0.3,-0.4) -- (0.3,0.4);
            \token{0.15,-0.2}{west}{a};
        \end{tikzpicture}
        \ =\
        \begin{tikzpicture}[centerzero]
            \draw[->] (0.3,-0.4) -- (-0.3,0.4);
            \draw[wipe] (-0.3,-0.4) -- (0.3,0.4);
            \draw[<-] (-0.3,-0.4) -- (0.3,0.4);
            \token{-0.15,0.2}{east}{a};
        \end{tikzpicture}
        \ ,\quad
        \begin{tikzpicture}[centerzero]
            \draw[<-] (-0.3,-0.4) -- (0.3,0.4);
            \draw[wipe] (0.3,-0.4) -- (-0.3,0.4);
            \draw[->] (0.3,-0.4) -- (-0.3,0.4);
            \token{0.15,0.2}{west}{a};
        \end{tikzpicture}
        \ =\
        \begin{tikzpicture}[centerzero]
            \draw[<-] (-0.3,-0.4) -- (0.3,0.4);
            \draw[wipe] (0.3,-0.4) -- (-0.3,0.4);
            \draw[->] (0.3,-0.4) -- (-0.3,0.4);
            \token{-0.15,-0.2}{east}{a};
        \end{tikzpicture}
        \ ,\quad
        \begin{tikzpicture}[anchorbase]
            \draw[<-] (-0.3,-0.4) -- (0.3,0.4);
            \draw[wipe] (0.3,-0.4) -- (-0.3,0.4);
            \draw[->] (0.3,-0.4) -- (-0.3,0.4);
            \token{0.15,-0.2}{west}{a};
        \end{tikzpicture}
        \ =\
        \begin{tikzpicture}[anchorbase]
            \draw[<-] (-0.3,-0.4) -- (0.3,0.4);
            \draw[wipe] (0.3,-0.4) -- (-0.3,0.4);
            \draw[->] (0.3,-0.4) -- (-0.3,0.4);
            \token{-0.15,0.2}{east}{a};
        \end{tikzpicture}
        \ ,\quad
        \begin{tikzpicture}[anchorbase]
            \draw[->] (0.3,-0.4) -- (-0.3,0.4);
            \draw[wipe] (-0.3,-0.4) -- (0.3,0.4);
            \draw[<-] (-0.3,-0.4) -- (0.3,0.4);
            \token{0.15,0.2}{west}{a};
        \end{tikzpicture}
        \ =\
        \begin{tikzpicture}[anchorbase]
            \draw[->] (0.3,-0.4) -- (-0.3,0.4);
            \draw[wipe] (-0.3,-0.4) -- (0.3,0.4);
            \draw[<-] (-0.3,-0.4) -- (0.3,0.4);
            \token{-0.15,-0.2}{east}{a};
        \end{tikzpicture}
        \ .
  \end{equation}
\end{lem}

\begin{proof}
  First suppose $k \ge 0$.  The first relation in \cref{ltokcross} follows from composing the second relation in \cref{rtokcross} on the top and bottom with the negative leftwards crossing, then using \cref{lunch,tea1,tea0}.  The fourth relation in \cref{ltokcross} follows similarly from composing the third relation in \cref{rtokcross} on the top and bottom with the negative leftwards crossing, then using \cref{lunch,tea1,tea0}.  Then the second and third relations in \cref{ltokcross} follow from the first and second using \cref{lskein,teleport}.  The case $k < 0$ is similar.
\end{proof}

Now consider the analogue of the morphism \cref{invrel} defined using the negative instead of the positive rightward crossing:
\begin{equation} \label{invrel1}
    \begin{aligned}
        \left[
            \begin{tikzpicture}[anchorbase]
                \draw[<-] (0.6,0) -- (0,0.6);
                \draw[wipe] (0,0) -- (0.6,0.6);
                \draw[->] (0,0) -- (0.6,0.6);
            \end{tikzpicture}
            \
            \begin{tikzpicture}[anchorbase]
                \draw[->] (-0.2,-0.4) -- (-0.2,0) arc(180:0:0.2) -- (0.2,-0.4);
                \multdot{-0.2,0}{west}{r};
                \token{-0.2,-0.2}{east}{b^\vee};
            \end{tikzpicture}\ ,\
            0 \le r \le k-1,\ b \in \BA
        \right]^T
        \ &\colon \uparrow \downarrow \to \downarrow \uparrow \oplus \one^{\oplus k \dim A} \quad \text{if } k \ge 0,
        \\
        \left[
            \begin{tikzpicture}[anchorbase]
                \draw[<-] (0.6,0) -- (0,0.6);
                \draw[wipe] (0,0) -- (0.6,0.6);
                \draw[->] (0,0) -- (0.6,0.6);
            \end{tikzpicture}
            \quad
            \begin{tikzpicture}[anchorbase]
                \draw[->] (-0.2,0.4) -- (-0.2,0) arc(180:360:0.2) -- (0.2,0.4);
                \multdot{0.2,0.2}{east}{r};
                \token{0.2,0}{west}{b^\vee};
            \end{tikzpicture}\ ,\
            0 \le r \le -k-1,\ b \in \BA
        \right]
        \ &\colon \uparrow \downarrow \oplus \one^{\oplus (-k \dim A)} \to \downarrow \uparrow \quad \text{if } k < 0.
    \end{aligned}
\end{equation}

\begin{lem} \label{irrelevant}
    The morphism \cref{invrel1} is invertible with two-sided inverse
    \begin{equation} \label{invrel3}
        \begin{aligned}
            \left[
                \begin{tikzpicture}[anchorbase]
                    \draw[->] (0.6,0) -- (0,0.6);
                    \draw[wipe] (0,0) -- (0.6,0.6);
                    \draw[<-] (0,0) -- (0.6,0.6);
                \end{tikzpicture}
                \
                \begin{tikzpicture}[centerzero={(0,-0.1)}]
                    \draw[<-] (-0.2,0.2) -- (-0.2,0) arc (180:360:0.2) -- (0.2,0.2);
                    \heartdec{0,-0.2}{north}{b}{r};
                \end{tikzpicture},\
                0 \le r \le k-1,\ b \in \BA
            \right]
            \ &\colon \downarrow \uparrow \oplus \one^{\oplus k \dim A} \to \uparrow \downarrow
            \text{ if } k > 0,
            \\
            \left[
                \begin{tikzpicture}[anchorbase]
                    \draw[<-] (0,0) -- (0.6,0.6);
                    \draw[wipe] (0.6,0) -- (0,0.6);
                    \draw[->] (0.6,0) -- (0,0.6);
                \end{tikzpicture}
                \
                \begin{tikzpicture}[centerzero={(0,0.1)}]
                    \draw[<-] (-0.2,-0.2) -- (-0.2,0) arc (180:0:0.2) -- (0.2,-0.2);
                    \diamdec{0,0.2}{south}{b}{r};
                \end{tikzpicture},\
                0 \le r \le -k-1,\ b \in \BA
            \right]^T
            \ &\colon \downarrow \uparrow \to \uparrow \downarrow \oplus \one^{\oplus (-k \dim A)}
            \quad \text{if } k \le 0.
        \end{aligned}
    \end{equation}
    Moreover, we have that
    \begin{gather} \label{colder}
        \cbubble{a}{k}
        \ = -\frac{t^{-1}}{z} \tr(a) 1_\one \text{ if } k > 0,
        \quad
        \begin{tikzpicture}[anchorbase]
            \draw[->] (0,0.2) arc(90:450:0.2);
            \token{0.2,0}{west}{a};
        \end{tikzpicture}
        \ = \frac{t-t^{-1}}{z} \tr(a) 1_\one \text{ if } k = 0,
        \quad
        \begin{tikzpicture}[anchorbase]
            \draw[->] (0,0.2) arc(90:450:0.2);
            \token{0.2,0}{west}{a};
        \end{tikzpicture}
        \ = -\frac{t^{-1}}{z} \tr(a) 1_\one \text{ if } k < 0,
        \\ \label{coldest}
        \begin{tikzpicture}[anchorbase]
            \draw[<-] (-0.2,0.2) -- (-0.2,0) arc (180:360:0.2) -- (0.2,0.2);
        \end{tikzpicture}
        =
        \begin{cases}
            \displaystyle{\frac{t}{z}}\
            \begin{tikzpicture}[centerzero={(0,-0.2)}]
                \draw[<-] (-0.2,0.2) -- (-0.2,0) arc (180:360:0.2) -- (0.2,0.2);
                \heartdec{0,-0.2}{north}{1}{0};
            \end{tikzpicture}
            & \text{if } k > 0, \\
            t^{-1}\
            \begin{tikzpicture}[anchorbase]
                \draw (-0.2,0) to[out=up,in=225] (0.2,0.4);
                \draw[wipe] (-0.2,0.4) to[out=-45,in=up] (0.2,0) arc(360:180:0.2);
                \draw[<-] (-0.2,0.4) to[out=-45,in=up] (0.2,0) arc(360:180:0.2);
                \multdot{0.2,0}{west}{-k};
            \end{tikzpicture}
            & \text{if } k \le 0,
        \end{cases}
        \qquad \qquad
        \begin{tikzpicture}[anchorbase]
            \draw[<-] (-0.2,-0.2) -- (-0.2,0) arc (180:0:0.2) -- (0.2,-0.2);
        \end{tikzpicture}
        =
        \begin{cases}
            t\
            \begin{tikzpicture}[anchorbase]
                \draw (-0.2,0) to[out=down,in=135] (0.2,-0.4);
                \draw[wipe] (-0.2,-0.4) to[out=45,in=down] (0.2,0) arc(0:180:0.2);
                \draw[<-] (-0.2,-0.4) to[out=45,in=down] (0.2,0) arc(0:180:0.2);
                \multdot{-0.2,0}{east}{k};
            \end{tikzpicture}
            & \text{if } k > 0,
            \\
            t^{-1}\
            \begin{tikzpicture}[anchorbase]
                \draw[<-] (-0.2,-0.4) to[out=45,in=down] (0.2,0) arc(0:180:0.2);
                \draw[wipe] (-0.2,0) to[out=down,in=135] (0.2,-0.4);
                \draw (-0.2,0) to[out=down,in=135] (0.2,-0.4);
            \end{tikzpicture}
            & \text{if } k=0,
            \\
            \displaystyle\frac{t}{z}
            \begin{tikzpicture}[centerzero={(0,0.2)}]
                \draw[<-] (-0.2,-0.2) -- (-0.2,0) arc (180:0:0.2) -- (0.2,-0.2);
                \multdot{-0.2,0}{east}{-1};
                \diamdec{0,0.2}{south}{1}{-k-1};
            \end{tikzpicture}
            & \text{if } k < 0.
        \end{cases}
    \end{gather}
\end{lem}

\begin{proof}
    First consider the case $k=0$.  Note that
    \[
        \begin{tikzpicture}[anchorbase]
            \draw[->] (0,0.2) arc(90:450:0.2);
            \token{0.2,0}{west}{a};
        \end{tikzpicture}
        \stackrel{\cref{leftwards}}{=} t\
        \begin{tikzpicture}[anchorbase]
            \draw (-0.2,0.2) to[out=down,in=up] (0.2,-0.2) to[out=down,in=down,looseness=1.5] (-0.2,-0.2);
            \draw[wipe] (-0.2,-0.2) to[out=up,in=down] (0.2,0.2) to[out=up,in=up,looseness=1.5] (-0.2,0.2);
            \draw[<-] (-0.2,-0.2) to[out=up,in=down] (0.2,0.2) to[out=up,in=up,looseness=1.5] (-0.2,0.2);
            \token{0.2,-0.2}{west}{a};
        \end{tikzpicture}
        \stackrel[\cref{rcslide}]{\cref{ltokcross}}{=} t\
        \begin{tikzpicture}[anchorbase]
           \draw (-0.2,0.2) to[out=down,in=up] (0.2,-0.2) to[out=down,in=down,looseness=1.5] (-0.2,-0.2);
            \draw[wipe] (-0.2,-0.2) to[out=up,in=down] (0.2,0.2) to[out=up,in=up,looseness=1.5] (-0.2,0.2);
            \draw[<-] (-0.2,-0.2) to[out=up,in=down] (0.2,0.2) to[out=up,in=up,looseness=1.5] (-0.2,0.2);
            \token{0.2,0.2}{west}{a};
        \end{tikzpicture}
        \stackrel{\cref{leftwards}}{=}
        \begin{tikzpicture}[anchorbase]
            \draw[<-] (0,0.2) arc(90:450:0.2);
            \token{0.2,0}{west}{a};
        \end{tikzpicture}
        \stackrel{\cref{impose}}{=}
        \frac{t-t^{-1}}{z} \tr(f) 1_\one,
    \]
    proving \cref{colder}.  Also,
    \[
        \begin{tikzpicture}[anchorbase]
            \draw[<-] (-0.2,0.2) -- (-0.2,0) arc (180:360:0.2) -- (0.2,0.2);
        \end{tikzpicture}
        \stackrel{\cref{leftwards}}{=}
        t\
        \begin{tikzpicture}[anchorbase]
            \draw[<-] (-0.2,0.4) to[out=-45,in=up] (0.2,0) arc(360:180:0.2);
            \draw[wipe] (-0.2,0) to[out=up,in=225] (0.2,0.4);
            \draw (-0.2,0) to[out=up,in=225] (0.2,0.4);
        \end{tikzpicture}
        \stackrel[\substack{\cref{colder} \\ \cref{adecomp}}]{\cref{lskein}}{=} t\
        \begin{tikzpicture}[anchorbase]
            \draw (-0.2,0) to[out=up,in=225] (0.2,0.4);
            \draw[wipe] (-0.2,0.4) to[out=-45,in=up] (0.2,0) arc(360:180:0.2);
            \draw[<-] (-0.2,0.4) to[out=-45,in=up] (0.2,0) arc(360:180:0.2);
        \end{tikzpicture}
        - (t^2-1)
        \begin{tikzpicture}[anchorbase]
            \draw[<-] (-0.2,0.2) -- (-0.2,0) arc (180:360:0.2) -- (0.2,0.2);
        \end{tikzpicture},
    \]
    from which the left-hand relation in \cref{coldest} follows.  The proof of the right-hand relation in \cref{coldest} is analogous.
    \details{
        We have
        \[
            \begin{tikzpicture}[anchorbase]
                \draw[<-] (-0.2,-0.2) -- (-0.2,0) arc (180:0:0.2) -- (0.2,-0.2);
            \end{tikzpicture}
            \stackrel{\cref{leftwards}}{=}
            t\
            \begin{tikzpicture}[anchorbase]
                \draw (-0.2,0) to[out=down,in=135] (0.2,-0.4);
                \draw[wipe] (-0.2,-0.4) to[out=45,in=down] (0.2,0) arc(0:180:0.2);
                \draw[<-] (-0.2,-0.4) to[out=45,in=down] (0.2,0) arc(0:180:0.2);
            \end{tikzpicture}
            \stackrel{\cref{lskein}}{=}
            t\
            \begin{tikzpicture}[anchorbase]
                \draw[<-] (-0.2,-0.4) to[out=45,in=down] (0.2,0) arc(0:180:0.2);
                \draw[wipe] (-0.2,0) to[out=down,in=135] (0.2,-0.4);
                \draw (-0.2,0) to[out=down,in=135] (0.2,-0.4);
            \end{tikzpicture}
            - (t^2-1)
            \begin{tikzpicture}[anchorbase]
                \draw[<-] (-0.2,-0.2) -- (-0.2,0) arc (180:0:0.2) -- (0.2,-0.2);
            \end{tikzpicture}
            \ .
        \]
        Solving for the left cap then gives the right relation in \cref{coldest}.
    }

    Now suppose $k > 0$.  By \cref{rskein}, we have
    \[
        \left[
            \begin{tikzpicture}[anchorbase]
                \draw[<-] (0.6,0) -- (0,0.6);
                \draw[wipe] (0,0) -- (0.6,0.6);
                \draw[->] (0,0) -- (0.6,0.6);
            \end{tikzpicture}
            \
            \begin{tikzpicture}[anchorbase]
                \draw[->] (0,0) -- (0,0.4) arc (180:0:0.2) -- (0.4,0);
                \multdot{0,0.4}{east}{r};
                \token{0,0.15}{east}{b^\vee};
            \end{tikzpicture}\ ,\
            0 \le r \le k-1,\ b \in \BA
        \right]^T
        =
        M
        \left[
            \begin{tikzpicture}[anchorbase]
                \draw[->] (0,0) -- (0.6,0.6);
                \draw[wipe] (0.6,0) -- (0,0.6);
                \draw[<-] (0.6,0) -- (0,0.6);
            \end{tikzpicture}
            \
            \begin{tikzpicture}[anchorbase]
                \draw[->] (0,0) -- (0,0.4) arc (180:0:0.2) -- (0.4,0);
                \multdot{0,0.4}{east}{r};
                \token{0,0.15}{east}{b^\vee};
            \end{tikzpicture}\ ,\
            0 \le r \le k-1,\ b \in \BA
        \right]^T,
    \]
    where $M$ is the $J \times J$ matrix whose $(\star,\star)$ and $(\star,(0,b))$ entries are
    \[
        \begin{tikzpicture}[anchorbase]
            \draw[->] (0,0.2) to (0,-0.2);
            \draw[<-] (0.3,0.2) to (0.3,-0.2);
        \end{tikzpicture}
        \qquad \text{and} \qquad
        \ -z\
        \begin{tikzpicture}[anchorbase]
            \draw[->] (-0.2,0.2) -- (-0.2,0) arc (180:360:0.2) -- (0.2,0.2);
            \token{0.2,0}{west}{b};
        \end{tikzpicture}
        \ ,
    \]
    whose $((r,b),(s,a))$ entry is $\delta_{r,s} \delta_{a,b}$ for all $0 \le r,s \le k-1$, $a,b \in \BA$, and whose other entries are zero.  The matrix $M$ is upper unitriangular, and hence invertible.  Its inverse is the matrix whose $(\star,\star)$ and $(\star,(0,b))$ entries are
    \[
        \begin{tikzpicture}[anchorbase]
            \draw[->] (0,0.2) to (0,-0.2);
            \draw[<-] (0.3,0.2) to (0.3,-0.2);
        \end{tikzpicture}
        \qquad \text{and} \qquad
        \ z\
        \begin{tikzpicture}[anchorbase]
            \draw[->] (-0.2,0.2) -- (-0.2,0) arc (180:360:0.2) -- (0.2,0.2);
            \token{0.2,0}{west}{b};
        \end{tikzpicture}
        \ ,
    \]
    whose $((r,b),(s,a))$ entry is $\delta_{r,s} \delta_{a,b}$ for all $0 \le r,s \le k-1$, $a,b \in \BA$, and whose other entries are zero.  Thus,
    \begin{align*}
        \left( \left[
            \begin{tikzpicture}[anchorbase]
                \draw[<-] (0.6,0) -- (0,0.6);
                \draw[wipe] (0,0) -- (0.6,0.6);
                \draw[->] (0,0) -- (0.6,0.6);
            \end{tikzpicture}
            \
            \begin{tikzpicture}[anchorbase]
                \draw[->] (0,0) -- (0,0.4) arc (180:0:0.2) -- (0.4,0);
                \multdot{0,0.4}{east}{r};
                \token{0,0.15}{east}{b^\vee};
            \end{tikzpicture}\ ,\
            0 \le r \le k-1,\ b \in \BA
        \right]^T \right)^{-1}
        &=
        \left[
            \begin{tikzpicture}[anchorbase]
                \draw[->] (0.6,0) -- (0,0.6);
                \draw[wipe] (0,0) -- (0.6,0.6);
                \draw[<-] (0,0) -- (0.6,0.6);
            \end{tikzpicture}
            \
            \begin{tikzpicture}[centerzero={(0,-0.1)}]
                \draw[<-] (-0.2,0.2) -- (-0.2,0) arc (180:360:0.2) -- (0.2,0.2);
                \diamdec{0,-0.2}{north}{b}{r};
            \end{tikzpicture},\
            0 \le r \le k-1,\ b \in \BA
        \right]
        M^{-1}
        \\
        &\stackrel{\mathclap{\cref{nakano1}}}{=}\ \
        \left[
            \begin{tikzpicture}[anchorbase]
                \draw[->] (0.6,0) -- (0,0.6);
                \draw[wipe] (0,0) -- (0.6,0.6);
                \draw[<-] (0,0) -- (0.6,0.6);
            \end{tikzpicture}
            \quad
            \begin{tikzpicture}[centerzero={(0,-0.1)}]
                \draw[<-] (-0.2,0.2) -- (-0.2,0) arc (180:360:0.2) -- (0.2,0.2);
                \heartdec{0,-0.2}{north}{b}{r};
            \end{tikzpicture},\
            0 \le r \le k-1,\ b \in \BA
        \right].
    \end{align*}

    The relation \cref{colder} follows from \cref{tea0}, and the right-hand relation in \cref{coldest} is the right-hand relation in \cref{leftwards}.  To prove the left-hand relation in \cref{coldest}, we compose both sides on the top with the isomorphism \cref{invrel1}, to see that it suffices to show that
    \[
        \begin{tikzpicture}[anchorbase]
            \draw (-0.2,0.4) to[out=-45,in=up] (0.2,0) arc(360:180:0.2);
            \draw[wipe] (-0.2,0) to[out=up,in=225] (0.2,0.4);
            \draw[->] (-0.2,0) to[out=up,in=225] (0.2,0.4);
        \end{tikzpicture}
        \ =0,
        \qquad
        \cbubble{b^\vee}{n}
        = \delta_{n,0} \frac{t}{z} \tr(b^\vee) 1_\one,\quad 0 \le n < k,\ b \in \BA.
    \]
    These relations follow immediately from \cref{tea0,impose}.

    The case $k < 0$ is similar to the case $k > 0$ and is left to the reader.
    \details{
        The relation \cref{colder} follows from \cref{tea2}, and the left-hand relation in \cref{coldest} is the left-hand relation in \cref{leftwards}.  To prove the right-hand relation in \cref{coldest}, we compose both sides on the bottom with the isomorphism \cref{invrel}, to see that it suffices to show that
        \[
            \begin{tikzpicture}[anchorbase]
                \draw (-0.2,-0.4) to[out=45,in=down] (0.2,0) arc(0:180:0.2);
                \draw[wipe] (-0.2,0) to[out=down,in=135] (0.2,-0.4);
                \draw[->] (-0.2,0) to[out=down,in=135] (0.2,-0.4);
            \end{tikzpicture}
            \ = \frac{t}{z}\
            \begin{tikzpicture}[centerzero={(0,-0.1)}]
                \draw (-0.2,-0.4) to[out=45,in=down] (0.2,0) arc(0:180:0.2);
                \draw[wipe] (-0.2,0) to[out=down,in=135] (0.2,-0.4);
                \draw[->] (-0.2,0) to[out=down,in=135] (0.2,-0.4);
                \multdot{-0.2,0}{east}{-1};
                \diamdec{0,0.2}{south}{1}{-k-1};
            \end{tikzpicture}
             \ ,
            \qquad
            \ccbubble{b^\vee}{n}
            = \delta_{n,-k} \frac{t}{z} \tr(b^\vee) 1_\one,\quad 0 < n \le -k,\ b \in \BA.
        \]
        To see the first relation, we compute
        \[
            \begin{tikzpicture}[centerzero={(0,-0.1)}]
                \draw (-0.2,-0.4) to[out=45,in=down] (0.2,0) arc(0:180:0.2);
                \draw[wipe] (-0.2,0) to[out=down,in=135] (0.2,-0.4);
                \draw[->] (-0.2,0) to[out=down,in=135] (0.2,-0.4);
                \multdot{-0.2,0}{east}{-1};
                \diamdec{0,0.2}{south}{1}{-k-1};
            \end{tikzpicture}
            \ =\
            \begin{tikzpicture}[centerzero={(0,-0.2)}]
                \draw[->] (-0.2,0) to[out=down,in=up] (0.2,-0.4) -- (0.2,-0.6);
                \draw[wipe] (-0.2,-0.6) -- (-0.2,-0.4) to[out=up,in=down] (0.2,0) arc(0:180:0.2);
                \draw (-0.2,-0.6) -- (-0.2,-0.4) to[out=up,in=down] (0.2,0) arc(0:180:0.2);
                \multdot{0.2,-0.4}{west}{-1};
                \diamdec{0,0.2}{south}{1}{-k-1};
            \end{tikzpicture}
            \ = 0
            \stackrel{\cref{tea2}}{=}
            \begin{tikzpicture}[anchorbase]
                \draw (-0.2,-0.4) to[out=45,in=down] (0.2,0) arc(0:180:0.2);
                \draw[wipe] (-0.2,0) to[out=down,in=135] (0.2,-0.4);
                \draw[->] (-0.2,0) to[out=down,in=135] (0.2,-0.4);
            \end{tikzpicture}
            \ .
        \]
        The other relations follow immediately from \cref{impose,tea2}.
    }
\end{proof}

\section{Second approach\label{sec:second}}

We now give a second definition of $\Heis_k(A;z,t)$.  Intuitively, this differs from \cref{def1} by replacing the positive crossing in the inversion relation by the negative crossing.


\begin{defin} \label{def2}
    The \emph{quantum Frobenius Heisenberg category} $\Heis_k(A;z,t)$ is the strict $\kk$-linear monoidal supercategory obtained from $\QAW(A;z,t)$ by adjoining a right dual $\downarrow$ to $\uparrow$, together with matrix entries of \cref{invrel3}, which we declare to be a two-sided inverse to \cref{invrel1}.  In addition, we impose the relation \cref{colder} for the leftwards cups and caps, which are defined in this approach by \cref{coldest}.
\end{defin}

Introduce the other leftward crossing not appearing in \cref{invrel3} so that \cref{lskein} holds, and set
\begin{gather}
    \begin{tikzpicture}[centerzero]
        \draw[<-] (-0.2,0.2) -- (-0.2,0) arc (180:360:0.2) -- (0.2,0.2);
        \diamdec{0,-0.2}{north}{a}{0};
    \end{tikzpicture}
    \ :=\
    \begin{tikzpicture}[>=To,baseline={([yshift=1ex]current bounding box.center)}]
        \draw[<-] (-0.2,0.2) -- (-0.2,0) arc (180:360:0.2) -- (0.2,0.2);
        \heartdec{0,-0.2}{north}{a}{0};
    \end{tikzpicture}
    \ - z\
    \begin{tikzpicture}[anchorbase]
        \draw[<-] (-0.2,0.4) to[out=-45,in=up] (0.2,0) arc(360:180:0.2);
        \draw[wipe] (-0.2,0) to[out=up,in=225] (0.2,0.4);
        \draw (-0.2,0) to[out=up,in=225] (0.2,0.4);
        \token{0.2,0}{west}{a};
    \end{tikzpicture}
    \ \text{if } k > 0,
    \qquad
    \begin{tikzpicture}[centerzero]
        \draw[<-] (-0.2,0.2) -- (-0.2,0) arc (180:360:0.2) -- (0.2,0.2);
        \diamdec{0,-0.2}{north}{a}{n};
    \end{tikzpicture}
    \ :=\
    \begin{tikzpicture}[>=To,baseline={([yshift=1ex]current bounding box.center)}]
        \draw[<-] (-0.2,0.2) -- (-0.2,0) arc (180:360:0.2) -- (0.2,0.2);
        \heartdec{0,-0.2}{north}{a}{n};
    \end{tikzpicture}
    \ \text{if } 0 < n < k,
    \\ \label{nakano2}
    \begin{tikzpicture}[>=To,baseline={([yshift=-2ex]current bounding box.center)}]
        \draw[<-] (-0.2,-0.2) -- (-0.2,0) arc (180:0:0.2) -- (0.2,-0.2);
        \heartdec{0,0.2}{south}{a}{0};
    \end{tikzpicture}
    \ :=\
    \begin{tikzpicture}[>=To,baseline={([yshift=-2ex]current bounding box.center)}]
        \draw[<-] (-0.2,-0.2) -- (-0.2,0) arc (180:0:0.2) -- (0.2,-0.2);
        \diamdec{0,0.2}{south}{a}{0};
    \end{tikzpicture}
    \ - (-1)^{\bar a} z\
    \begin{tikzpicture}[anchorbase]
        \draw[<-] (-0.2,-0.4) to[out=45,in=down] (0.2,0) arc(0:180:0.2);
        \draw[wipe] (-0.2,0) to[out=down,in=135] (0.2,-0.4);
        \draw (-0.2,0) to[out=down,in=135] (0.2,-0.4);
        \token{0.2,0}{west}{a};
    \end{tikzpicture}
    \ \text{if } k < 0,
    \qquad
    \begin{tikzpicture}[>=To,baseline={([yshift=-2ex]current bounding box.center)}]
        \draw[<-] (-0.2,-0.2) -- (-0.2,0) arc (180:0:0.2) -- (0.2,-0.2);
        \heartdec{0,0.2}{south}{a}{n};
    \end{tikzpicture}
    \ :=\
    \begin{tikzpicture}[>=To,baseline={([yshift=-2ex]current bounding box.center)}]
        \draw[<-] (-0.2,-0.2) -- (-0.2,0) arc (180:0:0.2) -- (0.2,-0.2);
        \diamdec{0,0.2}{south}{a}{n};
    \end{tikzpicture}
    \ \text{if } 0 < n < -k.
\end{gather}
Finally, we define the fake bubbles from \cref{fake0,fake1,fake2} as before.

\begin{theo}
    \Cref{def1,def2} give two different presentations for the same monoidal supercategory, and all of the named morphisms introduced in the two definitions are the same.  Moreover, there is a unique isomorphism of $\kk$-linear monoidal supercategories
    \begin{equation} \label{om}
        \Omega_k \colon \Heis_k(A;z,t) \to \Heis_{-k}(A;z,t^{-1})^\op,
    \end{equation}
    determined by
    \begin{gather*}
        \begin{tikzpicture}[anchorbase]
            \draw[->] (0,-0.2) -- (0,0.2);
            \token{0,0}{east}{a};
        \end{tikzpicture}
        \mapsto
        \begin{tikzpicture}[anchorbase]
            \draw[<-] (0,-0.2) -- (0,0.2);
            \token{0,0}{west}{a};
        \end{tikzpicture}
        \ ,\ a \in A,\qquad
        \begin{tikzpicture}[anchorbase]
            \draw[->] (0,-0.2) -- (0,0.2);
            \singdot{0,0};
        \end{tikzpicture}
        \mapsto
        \begin{tikzpicture}[anchorbase]
            \draw[<-] (0,-0.2) -- (0,0.2);
            \singdot{0,0};
        \end{tikzpicture}
        \ ,\qquad
        \begin{tikzpicture}[anchorbase]
            \draw[->] (0.2,-0.2) -- (-0.2,0.2);
            \draw[wipe] (-0.2,-0.2) -- (0.2,0.2);
            \draw[->] (-0.2,-0.2) -- (0.2,0.2);
        \end{tikzpicture}
        \mapsto -
        \begin{tikzpicture}[anchorbase]
            \draw[<-] (-0.2,-0.2) -- (0.2,0.2);
            \draw[wipe] (0.2,-0.2) -- (-0.2,0.2);
            \draw[<-] (0.2,-0.2) -- (-0.2,0.2);
        \end{tikzpicture}
        \ ,\qquad
        \begin{tikzpicture}[anchorbase]
            \draw[->] (-0.2,-0.2) -- (-0.2,0) arc (180:0:0.2) -- (0.2,-0.2);
        \end{tikzpicture}
        \mapsto
        \begin{tikzpicture}[anchorbase]
            \draw[->] (-0.2,0.2) -- (-0.2,0) arc (180:360:0.2) -- (0.2,0.2);
        \end{tikzpicture}
        \ ,\qquad
        \begin{tikzpicture}[anchorbase]
            \draw[->] (-0.2,0.2) -- (-0.2,0) arc (180:360:0.2) -- (0.2,0.2);
        \end{tikzpicture}
        \mapsto
        \begin{tikzpicture}[anchorbase]
            \draw[->] (-0.2,-0.2) -- (-0.2,0) arc (180:0:0.2) -- (0.2,-0.2);
        \end{tikzpicture}
        \ .
    \end{gather*}
    The isomorphism $\Omega_k$ acts on the other morphisms as follows:
    \begin{gather*}
        \begin{tikzpicture}[anchorbase]
            \draw[<-] (0,-0.2) -- (0,0.2);
            \token{0,0}{east}{a};
        \end{tikzpicture}
        \mapsto
        \begin{tikzpicture}[anchorbase]
            \draw[->] (0,-0.2) -- (0,0.2);
            \token{0,0}{west}{a};
        \end{tikzpicture}
        \ ,\ a \in A,\qquad
        \begin{tikzpicture}[anchorbase]
            \draw[<-] (0,-0.2) -- (0,0.2);
            \singdot{0,0};
        \end{tikzpicture}
        \mapsto
        \begin{tikzpicture}[anchorbase]
            \draw[->] (0,-0.2) -- (0,0.2);
            \singdot{0,0};
        \end{tikzpicture}
        \ ,\qquad
        \begin{tikzpicture}[anchorbase]
            \draw[->] (-0.2,-0.2) -- (0.2,0.2);
            \draw[wipe] (0.2,-0.2) -- (-0.2,0.2);
            \draw[->] (0.2,-0.2) -- (-0.2,0.2);
        \end{tikzpicture}
        \mapsto -
        \begin{tikzpicture}[anchorbase]
            \draw[<-] (0.2,-0.2) -- (-0.2,0.2);
            \draw[wipe] (-0.2,-0.2) -- (0.2,0.2);
            \draw[<-] (-0.2,-0.2) -- (0.2,0.2);
        \end{tikzpicture}
        \ ,\qquad
        \begin{tikzpicture}[anchorbase]
            \draw[<-] (0.2,-0.2) -- (-0.2,0.2);
            \draw[wipe] (-0.2,-0.2) -- (0.2,0.2);
            \draw[->] (-0.2,-0.2) -- (0.2,0.2);
        \end{tikzpicture}
        \mapsto -
        \begin{tikzpicture}[anchorbase]
            \draw[->] (-0.2,-0.2) -- (0.2,0.2);
            \draw[wipe] (0.2,-0.2) -- (-0.2,0.2);
            \draw[<-] (0.2,-0.2) -- (-0.2,0.2);
        \end{tikzpicture}
        \ ,\qquad
        \begin{tikzpicture}[anchorbase]
            \draw[->] (-0.2,-0.2) -- (0.2,0.2);
            \draw[wipe] (0.2,-0.2) -- (-0.2,0.2);
            \draw[<-] (0.2,-0.2) -- (-0.2,0.2);
        \end{tikzpicture}
        \mapsto -
        \begin{tikzpicture}[anchorbase]
            \draw[<-] (0.2,-0.2) -- (-0.2,0.2);
            \draw[wipe] (-0.2,-0.2) -- (0.2,0.2);
            \draw[->] (-0.2,-0.2) -- (0.2,0.2);
        \end{tikzpicture}
        \ ,
        \\
        \begin{tikzpicture}[anchorbase]
            \draw[<-] (0.2,-0.2) -- (-0.2,0.2);
            \draw[wipe] (-0.2,-0.2) -- (0.2,0.2);
            \draw[<-] (-0.2,-0.2) -- (0.2,0.2);
        \end{tikzpicture}
        \mapsto -
        \begin{tikzpicture}[anchorbase]
            \draw[->] (-0.2,-0.2) -- (0.2,0.2);
            \draw[wipe] (0.2,-0.2) -- (-0.2,0.2);
            \draw[->] (0.2,-0.2) -- (-0.2,0.2);
        \end{tikzpicture}
        \ ,\qquad
        \begin{tikzpicture}[anchorbase]
            \draw[<-] (-0.2,-0.2) -- (0.2,0.2);
            \draw[wipe] (0.2,-0.2) -- (-0.2,0.2);
            \draw[<-] (0.2,-0.2) -- (-0.2,0.2);
        \end{tikzpicture}
        \mapsto -
        \begin{tikzpicture}[anchorbase]
            \draw[->] (0.2,-0.2) -- (-0.2,0.2);
            \draw[wipe] (-0.2,-0.2) -- (0.2,0.2);
            \draw[->] (-0.2,-0.2) -- (0.2,0.2);
        \end{tikzpicture}
        \ ,\qquad
        \begin{tikzpicture}[anchorbase]
            \draw[->] (0.2,-0.2) -- (-0.2,0.2);
            \draw[wipe] (-0.2,-0.2) -- (0.2,0.2);
            \draw[<-] (-0.2,-0.2) -- (0.2,0.2);
        \end{tikzpicture}
        \mapsto -
        \begin{tikzpicture}[anchorbase]
            \draw[<-] (-0.2,-0.2) -- (0.2,0.2);
            \draw[wipe] (0.2,-0.2) -- (-0.2,0.2);
            \draw[->] (0.2,-0.2) -- (-0.2,0.2);
        \end{tikzpicture}
        \ ,\qquad
        \begin{tikzpicture}[anchorbase]
            \draw[<-] (-0.2,-0.2) -- (0.2,0.2);
            \draw[wipe] (0.2,-0.2) -- (-0.2,0.2);
            \draw[->] (0.2,-0.2) -- (-0.2,0.2);
        \end{tikzpicture}
        \mapsto -
        \begin{tikzpicture}[anchorbase]
            \draw[->] (0.2,-0.2) -- (-0.2,0.2);
            \draw[wipe] (-0.2,-0.2) -- (0.2,0.2);
            \draw[<-] (-0.2,-0.2) -- (0.2,0.2);
        \end{tikzpicture}
        \ ,
        \\
        \begin{tikzpicture}[centerzero]
            \draw[<-] (-0.2,0.2) -- (-0.2,0) arc (180:360:0.2) -- (0.2,0.2);
            \diamdec{0,-0.2}{north}{a}{n};
        \end{tikzpicture}
        \mapsto (-1)^{\bar a}\
        \begin{tikzpicture}[centerzero]
            \draw[<-] (-0.2,-0.2) -- (-0.2,0) arc (180:0:0.2) -- (0.2,-0.2);
            \diamdec{0,0.2}{south}{a}{n};
        \end{tikzpicture}
        \ ,\quad
        \begin{tikzpicture}[centerzero]
            \draw[<-] (-0.2,-0.2) -- (-0.2,0) arc (180:0:0.2) -- (0.2,-0.2);
            \diamdec{0,0.2}{south}{a}{n};
        \end{tikzpicture}
        \mapsto (-1)^{\bar a}\
        \begin{tikzpicture}[centerzero]
            \draw[<-] (-0.2,0.2) -- (-0.2,0) arc (180:360:0.2) -- (0.2,0.2);
            \diamdec{0,-0.2}{north}{a}{n};
        \end{tikzpicture}
        \ ,\quad
        \begin{tikzpicture}[centerzero]
            \draw[<-] (-0.2,0.2) -- (-0.2,0) arc (180:360:0.2) -- (0.2,0.2);
            \heartdec{0,-0.2}{north}{a}{n};
        \end{tikzpicture}
        \mapsto (-1)^{\bar a}\
        \begin{tikzpicture}[centerzero]
            \draw[<-] (-0.2,-0.2) -- (-0.2,0) arc (180:0:0.2) -- (0.2,-0.2);
            \heartdec{0,0.2}{south}{a}{n};
        \end{tikzpicture}
        \ ,\quad
        \begin{tikzpicture}[centerzero]
            \draw[<-] (-0.2,-0.2) -- (-0.2,0) arc (180:0:0.2) -- (0.2,-0.2);
            \heartdec{0,0.2}{south}{a}{n};
        \end{tikzpicture}
        \mapsto (-1)^{\bar a}\
        \begin{tikzpicture}[centerzero]
            \draw[<-] (-0.2,0.2) -- (-0.2,0) arc (180:360:0.2) -- (0.2,0.2);
            \heartdec{0,-0.2}{north}{a}{n};
        \end{tikzpicture}
        \ ,
        \\
        \begin{tikzpicture}[anchorbase]
            \draw[<-] (-0.2,-0.2) -- (-0.2,0) arc (180:0:0.2) -- (0.2,-0.2);
        \end{tikzpicture}
        \mapsto -
        \begin{tikzpicture}[anchorbase]
            \draw[<-] (-0.2,0.2) -- (-0.2,0) arc (180:360:0.2) -- (0.2,0.2);
        \end{tikzpicture}
        \ ,\qquad
        \begin{tikzpicture}[anchorbase]
            \draw[<-] (-0.2,0.2) -- (-0.2,0) arc (180:360:0.2) -- (0.2,0.2);
        \end{tikzpicture}
        \mapsto -
        \begin{tikzpicture}[anchorbase]
            \draw[<-] (-0.2,-0.2) -- (-0.2,0) arc (180:0:0.2) -- (0.2,-0.2);
        \end{tikzpicture}
        \ ,\qquad
        \begin{tikzpicture}[anchorbase]
            \draw[->] (0,0.2) arc(90:450:0.2);
            \node at (0,0) {\dotlabel{\pm}};
            \token{-0.2,0}{east}{a};
            \multdot{0.2,0}{west}{n};
        \end{tikzpicture}
        \mapsto -
        \begin{tikzpicture}[anchorbase]
            \draw[->] (0,0.2) arc(90:-270:0.2);
            \node at (0,0) {\dotlabel{\pm}};
            \token{0.2,0}{west}{a};
            \multdot{-0.2,0}{east}{n};
        \end{tikzpicture}
        \ .
    \end{gather*}
    In particular, $\Omega_k^2$ is the identity.
\end{theo}

Note that the isomorphism $\Omega_k$ is given by reflecting diagrams in a horizontal plane and multiplying by $(-1)^{c+d+\binom{y}{2}}$, where $c$ is the number of crossings, $d$ is the number of left cups and caps (including left cups and caps in fake bubbles, but not ones labelled by $\color{violet}\scriptstyle{\diamondsuit}$ or $\color{violet}\scriptstyle{\heartsuit}$), and $y$ is the number of odd tokens.

\begin{proof}
  To avoid confusion, we denote the category $\Heis_k(A;z,t)$ from \cref{def1} by $\Heis_k^\old(A;z,t)$ and the one from \cref{def2} by $\Heis_k^\new(A;z,t)$.  The relations and other definitions for the category $\Heis_k^\new(A;z,t)$ in \cref{def2} and the ones for $\Heis_k^\old(A;z,t^{-1})$ from \cref{def1} are related by reflecting all diagrams in a horizontal plane and multiplying by $(-1)^{c + d + \binom{y}{2}}$, where $c$ is the number of crossings and $d$ is the number of left cups and caps (including left cups and caps in fake bubbles, but not ones labelled by $\color{violet}\scriptstyle{\diamondsuit}$ or $\color{violet}\scriptstyle{\heartsuit}$), and $y$ is the number of odd tokens.  It follows that there are mutually inverse isomorphisms
  \[
    \Heis_{-k}^\old(A;z,t^{-1}) \stackrel[\Omega_+]{\Omega_-}{\rightleftarrows} \Heis_k^\new(A;z,t)^\op
  \]
  both defined in the same way as the functor $\Omega_k$ in the statement of the theorem.  Now, by \cref{irrelevant,def2}, there exists a strict $\kk$-linear monoidal functor
  \[
    \Theta_k \colon \Heis_k^\new(A;z,t) \to \Heis_k^\old(A;z,t)
  \]
  that is the identity on diagrams.  This functor is an isomorphism because it has a two-sided inverse, namely $\Omega_- \circ \Theta_{-k} \circ \Omega_-$.  Thus, using $\Theta_k$, we may identify $\Heis_k^\new(A;z,t)$ and $\Heis_k^\old(A;z,t)$.  Finally, $\Omega_k := \Omega_+$ gives the required symmetry.
\end{proof}


We now prove some important relations.  In particular, we will show how we can remove the decorated left cups and caps from our string diagrams.  Then we prove the important infinite Grassmannian relations.

\begin{lem}
    For all $a \in A$, we have
    \begin{equation} \label{diamondcup}
        \begin{tikzpicture}[centerzero]
            \draw[<-] (-0.2,0.2) -- (-0.2,0) arc (180:360:0.2) -- (0.2,0.2);
            \diamdec{0,-0.2}{north}{a}{k-1};
        \end{tikzpicture}
        \ = - zt\
        \begin{tikzpicture}[centerzero={(0,0.1)}]
            \draw[<-] (-0.2,0.4) -- (-0.2,0) arc (180:360:0.2) -- (0.2,0.4);
            \singdot{-0.2,0.2};
            \token{-0.2,0}{east}{a};
        \end{tikzpicture}
        \quad \text{if } k > 0,
        \qquad
        \begin{tikzpicture}[centerzero]
            \draw[<-] (-0.2,-0.2) -- (-0.2,0) arc (180:0:0.2) -- (0.2,-0.2);
            \diamdec{0,0.2}{south}{a}{-k-1};
        \end{tikzpicture}
        = (-1)^{\bar a}zt^{-1}\
        \begin{tikzpicture}[centerzero={(0,-0.2)}]
            \draw[<-] (-0.2,-0.4) -- (-0.2,0) arc (180:0:0.2) -- (0.2,-0.4);
            \singdot{-0.2,-0.2};
            \token{-0.2,0}{east}{a};
        \end{tikzpicture}
        \quad \text{if } k < 0.
    \end{equation}
\end{lem}

\begin{proof}
    To prove the first relation in \cref{diamondcup}, compose the second relation in \cref{tea1} on the bottom with
    \[
        \begin{tikzpicture}[anchorbase]
            \draw[<-] (-0.2,0.4) -- (-0.2,0) arc (180:360:0.2) -- (0.2,0.4);
            \token{-0.2,0}{east}{a};
            \singdot{-0.2,0.2};
        \end{tikzpicture}
    \]
    and use \cref{rdotcross,tea0}.  The second relation in \cref{diamondcup} then follows by applying $\Omega_k$.
\end{proof}

\begin{lem}
    The following relations hold for all $f \in F$:
    \begin{gather} \label{ldotcross}
        \begin{tikzpicture}[anchorbase]
            \draw[->] (0.3,-0.3) -- (-0.3,0.3);
            \draw[wipe] (-0.3,-0.3) -- (0.3,0.3);
            \draw[<-] (-0.3,-0.3) -- (0.3,0.3);
            \singdot{0.17,-0.17};
        \end{tikzpicture}
        \ =\
        \begin{tikzpicture}[anchorbase]
            \draw[<-] (-0.3,-0.3) -- (0.3,0.3);
            \draw[wipe] (0.3,-0.3) -- (-0.3,0.3);
            \draw[->] (0.3,-0.3) -- (-0.3,0.3);
            \singdot{-0.15,0.15};
        \end{tikzpicture}
        \ ,\quad
        \begin{tikzpicture}[anchorbase]
            \draw[<-] (-0.3,-0.3) -- (0.3,0.3);
            \draw[wipe] (0.3,-0.3) -- (-0.3,0.3);
            \draw[->] (0.3,-0.3) -- (-0.3,0.3);
            \singdot{0.17,0.17};
        \end{tikzpicture}
        \ =\
        \begin{tikzpicture}[anchorbase]
            \draw[->] (0.3,-0.3) -- (-0.3,0.3);
            \draw[wipe] (-0.3,-0.3) -- (0.3,0.3);
            \draw[<-] (-0.3,-0.3) -- (0.3,0.3);
            \singdot{-0.15,-0.15};
        \end{tikzpicture}
        \ ,
        \\ \label{piv}
        \begin{tikzpicture}[anchorbase]
            \draw[<-] (-0.2,-0.2) -- (-0.2,0) arc (180:0:0.2) -- (0.2,-0.2);
            \token{-0.2,0}{east}{a};
        \end{tikzpicture}
        \ =\
        \begin{tikzpicture}[anchorbase]
            \draw[<-] (-0.2,-0.2) -- (-0.2,0) arc (180:0:0.2) -- (0.2,-0.2);
            \token{0.2,0}{west}{a};
        \end{tikzpicture}
        \ ,\qquad
        \begin{tikzpicture}[anchorbase]
            \draw[<-] (-0.2,0.2) -- (-0.2,0) arc (180:360:0.2) -- (0.2,0.2);
            \token{-0.2,0}{east}{a};
        \end{tikzpicture}
        \ =\
        \begin{tikzpicture}[anchorbase]
            \draw[<-] (-0.2,0.2) -- (-0.2,0) arc (180:360:0.2) -- (0.2,0.2);
            \token{0.2,0}{west}{a};
        \end{tikzpicture}
        \ ,\qquad
        \begin{tikzpicture}[anchorbase]
            \draw[<-] (-0.2,-0.2) -- (-0.2,0) arc (180:0:0.2) -- (0.2,-0.2);
            \singdot{-0.2,0};
        \end{tikzpicture}
        \ =\
        \begin{tikzpicture}[anchorbase]
            \draw[<-] (-0.2,-0.2) -- (-0.2,0) arc (180:0:0.2) -- (0.2,-0.2);
            \singdot{0.2,0};
        \end{tikzpicture}
        \ ,\qquad
        \begin{tikzpicture}[anchorbase]
            \draw[<-] (-0.2,0.2) -- (-0.2,0) arc (180:360:0.2) -- (0.2,0.2);
            \singdot{-0.2,0};
        \end{tikzpicture}
        \ =\
        \begin{tikzpicture}[anchorbase]
            \draw[<-] (-0.2,0.2) -- (-0.2,0) arc (180:360:0.2) -- (0.2,0.2);
            \singdot{0.2,0};
        \end{tikzpicture}
        \ .
    \end{gather}
\end{lem}

\begin{proof}
    Using $\Omega_k$, it suffices to consider the case $k \ge 0$.  First suppose $k=0$.  Composing the first relation in \cref{rdotcross} on the top with the leftward negative crossing and on the bottom with the leftward positive crossing, we have
    \[
        \begin{tikzpicture}[anchorbase]
            \draw[<-] (-0.2,0) to[out=up,in=down] (0.2,0.5);
            \draw (-0.2,0.5) to[out=up,in=down] (0.2,1);
            \draw[->] (0.2,1) to[out=up,in=down] (-0.2,1.5);
            \draw[wipe] (0.2,0) to[out=up,in=down] (-0.2,0.5);
            \draw (0.2,0) to[out=up,in=down] (-0.2,0.5);
            \draw[wipe] (0.2,0.5) to[out=up,in=down] (-0.2,1);
            \draw (0.2,0.5) to[out=up,in=down] (-0.2,1);
            \draw[wipe] (-0.2,1) to[out=up,in=down] (0.2,1.5);
            \draw (-0.2,1) to[out=up,in=down] (0.2,1.5);
            \singdot{-0.2,0.5};
        \end{tikzpicture}
        \ =\
        \begin{tikzpicture}[anchorbase]
            \draw[<-] (-0.2,0) to[out=up,in=down] (0.2,0.5);
            \draw (0.2,0.5) to[out=up,in=down] (-0.2,1);
            \draw[->] (0.2,1) to[out=up,in=down] (-0.2,1.5);
            \draw[wipe] (0.2,0) to[out=up,in=down] (-0.2,0.5);
            \draw (0.2,0) to[out=up,in=down] (-0.2,0.5);
            \draw[wipe] (-0.2,0.5) to[out=up,in=down] (0.2,1);
            \draw (-0.2,0.5) to[out=up,in=down] (0.2,1);
            \draw[wipe] (-0.2,1) to[out=up,in=down] (0.2,1.5);
            \draw (-0.2,1) to[out=up,in=down] (0.2,1.5);
            \singdot{0.2,1};
        \end{tikzpicture}
        \stackrel{\cref{lunch}}{\implies}
        \begin{tikzpicture}[anchorbase]
            \draw[<-] (-0.3,-0.3) -- (0.3,0.3);
            \draw[wipe] (0.3,-0.3) -- (-0.3,0.3);
            \draw[->] (0.3,-0.3) -- (-0.3,0.3);
            \singdot{-0.15,0.15};
        \end{tikzpicture}
        \ =\
        \begin{tikzpicture}[anchorbase]
            \draw[->] (0.3,-0.3) -- (-0.3,0.3);
            \draw[wipe] (-0.3,-0.3) -- (0.3,0.3);
            \draw[<-] (-0.3,-0.3) -- (0.3,0.3);
            \singdot{0.17,-0.17};
        \end{tikzpicture}
        \ .
    \]
    The proof of the second relation in \cref{ldotcross} is obtained similarly from composing the second relation in \cref{rdotcross} on the top with the leftward negative crossing and on the bottom with the leftward positive crossing.  The relations in \cref{piv} then follow from the definition \cref{leftwards}, together with \cref{ltokcross,ldotcross,rcslide}.

    Now suppose $k > 0$.  Composing the first relation in \cref{rdotcross} on the top and bottom with the leftward negative crossing, we have
    \begin{multline*}
        0 =\
        \begin{tikzpicture}[anchorbase]
            \draw (0.2,0) to[out=up,in=down] (-0.2,0.5);
            \draw (-0.2,0.5) to[out=up,in=down] (0.2,1);
            \draw[->] (0.2,1) to[out=up,in=down] (-0.2,1.5);
            \draw[wipe] (-0.2,0) to[out=up,in=down] (0.2,0.5);
            \draw[<-] (-0.2,0) to[out=up,in=down] (0.2,0.5);
            \draw[wipe] (0.2,0.5) to[out=up,in=down] (-0.2,1);
            \draw (0.2,0.5) to[out=up,in=down] (-0.2,1);
            \draw[wipe] (-0.2,1) to[out=up,in=down] (0.2,1.5);
            \draw (-0.2,1) to[out=up,in=down] (0.2,1.5);
            \singdot{-0.2,0.5};
        \end{tikzpicture}
        \ -\
        \begin{tikzpicture}[anchorbase]
            \draw (0.2,0) to[out=up,in=down] (-0.2,0.5);
            \draw (0.2,0.5) to[out=up,in=down] (-0.2,1);
            \draw[->] (0.2,1) to[out=up,in=down] (-0.2,1.5);
            \draw[wipe] (-0.2,0) to[out=up,in=down] (0.2,0.5);
            \draw[<-] (-0.2,0) to[out=up,in=down] (0.2,0.5);
            \draw[wipe] (-0.2,0.5) to[out=up,in=down] (0.2,1);
            \draw (-0.2,0.5) to[out=up,in=down] (0.2,1);
            \draw[wipe] (-0.2,1) to[out=up,in=down] (0.2,1.5);
            \draw (-0.2,1) to[out=up,in=down] (0.2,1.5);
            \singdot{0.2,1};
        \end{tikzpicture}
        \ \stackrel[\cref{rskein}]{\cref{tea1}}{=}\
        \begin{tikzpicture}[anchorbase]
            \draw[->] (0.3,-0.3) -- (-0.3,0.3);
            \draw[wipe] (-0.3,-0.3) -- (0.3,0.3);
            \draw[<-] (-0.3,-0.3) -- (0.3,0.3);
            \singdot{-0.15,0.15};
        \end{tikzpicture}
        \ - \sum_{r=0}^{k-1}
        \begin{tikzpicture}[anchorbase]
            \draw (-0.2,0) to (-0.2,0.2) to[out=up,in=up,looseness=1.5] (0.2,0.2) to (0.2,0);
            \draw (-0.2,0) to[out=down,in=up] (0.2,-0.6);
            \draw[wipe] (0.2,0) to[out=down,in=up] (-0.2,-0.6);
            \draw[->] (0.2,0) to[out=down,in=up] (-0.2,-0.6);
            \draw[->] (0.2,1.2) -- (0.2,1.1) arc(360:180:0.2) -- (-0.2,1.2);
            \multdot{-0.2,0.2}{east}{r+1};
            \token{-0.2,0}{north east}{b^\vee};
            \diamdec{0,0.9}{north}{b}{r};
        \end{tikzpicture}
        \ -\
        \begin{tikzpicture}[anchorbase]
            \draw (0.2,0) to[out=up,in=down] (-0.2,0.5);
            \draw (-0.2,0.5) to[out=up,in=down] (0.2,1);
            \draw[->] (0.2,1) to[out=up,in=down] (-0.2,1.5);
            \draw[wipe] (-0.2,0) to[out=up,in=down] (0.2,0.5);
            \draw[<-] (-0.2,0) to[out=up,in=down] (0.2,0.5);
            \draw[wipe] (0.2,0.5) to[out=up,in=down] (-0.2,1);
            \draw (0.2,0.5) to[out=up,in=down] (-0.2,1);
            \draw[wipe] (-0.2,1) to[out=up,in=down] (0.2,1.5);
            \draw (-0.2,1) to[out=up,in=down] (0.2,1.5);
            \singdot{0.2,1};
        \end{tikzpicture}
        \ - z\
        \begin{tikzpicture}[anchorbase]
            \draw (0.2,-0.4) to[out=up,in=up,looseness=1.5] (-0.2,-0.4);
            \draw (-0.2,-0.4) to[out=down,in=up] (0.2,-1);
            \draw[wipe] (-0.2,-1) to[out=up,in=down] (0.2,-0.4);
            \draw[<-] (-0.2,-1) to[out=up,in=down] (0.2,-0.4);
            \draw[<-] (-0.2,1) to[out=down,in=up] (0.2,0.4) to (0.2,0.2);
            \draw (0.2,0.2) to[out=down,in=down,looseness=1.5] (-0.2,0.2);
            \draw[wipe] (-0.2,0.2) to (-0.2,0.4) to[out=up,in=down] (0.2,1);
            \draw (-0.2,0.2) to (-0.2,0.4) to[out=up,in=down] (0.2,1);
            \singdot{0.2,0.4};
            \token{0.2,0.2}{west}{b};
            \token{-0.2,-0.4}{east}{b^\vee};
        \end{tikzpicture}
    \\
        \stackrel[\cref{tea1}]{\cref{tea0}}{=}\
        \begin{tikzpicture}[anchorbase]
            \draw[->] (0.3,-0.3) -- (-0.3,0.3);
            \draw[wipe] (-0.3,-0.3) -- (0.3,0.3);
            \draw[<-] (-0.3,-0.3) -- (0.3,0.3);
            \singdot{-0.15,0.15};
        \end{tikzpicture}
        -
        \begin{tikzpicture}[anchorbase]
            \draw (-0.2,0) to (-0.2,0.2) to[out=up,in=up,looseness=1.5] (0.2,0.2) to (0.2,0);
            \draw (-0.2,0) to[out=down,in=up] (0.2,-0.6);
            \draw[wipe] (0.2,0) to[out=down,in=up] (-0.2,-0.6);
            \draw[->] (0.2,0) to[out=down,in=up] (-0.2,-0.6);
            \draw[->] (0.2,1.2) -- (0.2,1.1) arc(360:180:0.2) -- (-0.2,1.2);
            \multdot{-0.2,0.2}{east}{k};
            \token{-0.2,0}{north east}{b^\vee};
            \diamdec{0,0.9}{north}{b}{k-1};
        \end{tikzpicture}
        \ -\
        \begin{tikzpicture}[anchorbase]
            \draw[->] (0.3,-0.3) -- (-0.3,0.3);
            \draw[wipe] (-0.3,-0.3) -- (0.3,0.3);
            \draw[<-] (-0.3,-0.3) -- (0.3,0.3);
            \singdot{0.15,-0.15};
        \end{tikzpicture}
        \ \stackrel{\cref{lskein}}{=}\
        \begin{tikzpicture}[anchorbase]
            \draw[<-] (-0.3,-0.3) -- (0.3,0.3);
            \draw[wipe] (0.3,-0.3) -- (-0.3,0.3);
            \draw[->] (0.3,-0.3) -- (-0.3,0.3);
            \singdot{-0.15,0.15};
        \end{tikzpicture}
        - z\
        \begin{tikzpicture}[anchorbase]
            \draw[<-] (-0.2,1.1) -- (-0.2,0.7) arc(180:360:0.2) -- (0.2,1.1);
            \draw[<-] (-0.2,-0.1) -- (-0.2,0.1) arc(180:0:0.2) -- (0.2,-0.1);
            \token{-0.2,0.7}{east}{b};
            \token{0.2,0.1}{west}{b^\vee};
            \singdot{-0.2,0.9};
        \end{tikzpicture}
        -
        \begin{tikzpicture}[anchorbase]
            \draw (-0.2,0) to (-0.2,0.2) to[out=up,in=up,looseness=1.5] (0.2,0.2) to (0.2,0);
            \draw (-0.2,0) to[out=down,in=up] (0.2,-0.6);
            \draw[wipe] (0.2,0) to[out=down,in=up] (-0.2,-0.6);
            \draw[->] (0.2,0) to[out=down,in=up] (-0.2,-0.6);
            \draw[->] (0.2,1.2) -- (0.2,1.1) arc(360:180:0.2) -- (-0.2,1.2);
            \multdot{-0.2,0.2}{east}{k};
            \token{-0.2,0}{north east}{b^\vee};
            \diamdec{0,0.9}{north}{b}{k-1};
        \end{tikzpicture}
        \ -\
        \begin{tikzpicture}[anchorbase]
            \draw[->] (0.3,-0.3) -- (-0.3,0.3);
            \draw[wipe] (-0.3,-0.3) -- (0.3,0.3);
            \draw[<-] (-0.3,-0.3) -- (0.3,0.3);
            \singdot{0.15,-0.15};
        \end{tikzpicture}
        \stackrel[\cref{diamondcup}]{\substack{\cref{ltokcross} \\ \cref{leftwards}}}{=}
        \begin{tikzpicture}[anchorbase]
            \draw[<-] (-0.3,-0.3) -- (0.3,0.3);
            \draw[wipe] (0.3,-0.3) -- (-0.3,0.3);
            \draw[->] (0.3,-0.3) -- (-0.3,0.3);
            \singdot{-0.15,0.15};
        \end{tikzpicture}
        \ -\
        \begin{tikzpicture}[anchorbase]
            \draw[->] (0.3,-0.3) -- (-0.3,0.3);
            \draw[wipe] (-0.3,-0.3) -- (0.3,0.3);
            \draw[<-] (-0.3,-0.3) -- (0.3,0.3);
            \singdot{0.15,-0.15};
        \end{tikzpicture}
        \ .
    \end{multline*}
    The proof of the second relation in \cref{ldotcross} is similar, starting with the second relation in \cref{rdotcross}.

    The first relation in \cref{piv} follows immediately from the definition \cref{leftwards}, together with \cref{ltokcross,QAWC}.  To prove the second relation in \cref{piv}, compose the second relation in \cref{tea1} on the bottom with
    $
        \begin{tikzpicture}[anchorbase]
            \draw[<-] (-0.2,0.2) -- (-0.2,0) arc (180:360:0.2) -- (0.2,0.2);
            \token{0.2,0}{west}{a};
            \singdot{-0.2,0};
        \end{tikzpicture}
    $
    and use \cref{rdotcross,tea0} to obtain
    \[
        \begin{tikzpicture}[centerzero]
            \draw[<-] (-0.2,0.2) -- (-0.2,0) arc (180:360:0.2) -- (0.2,0.2);
            \token{0.2,0}{west}{a};
            \singdot{-0.2,0};
        \end{tikzpicture}
        = -z^{-1}t^{-1}
        \begin{tikzpicture}[centerzero]
            \draw[<-] (-0.2,0.2) -- (-0.2,0) arc (180:360:0.2) -- (0.2,0.2);
            \diamdec{0,-0.2}{north}{a}{k-1};
        \end{tikzpicture}
        \stackrel{\cref{diamondcup}}{=}
        \begin{tikzpicture}[centerzero={(0,0.2)}]
            \draw[<-] (-0.2,0.4) -- (-0.2,0) arc (180:360:0.2) -- (0.2,0.4);
            \token{-0.2,0}{east}{a};
            \singdot{-0.2,0.2};
        \end{tikzpicture}
        \ .
    \]
    Then we add an inverse dot to the top of the left strand.

    To prove the third relation in \cref{piv}, we note that
    \[
        \begin{tikzpicture}[anchorbase]
            \draw[<-] (-0.2,-0.2) -- (-0.2,0) arc (180:0:0.2) -- (0.2,-0.2);
            \singdot{-0.2,0};
        \end{tikzpicture}
        \stackrel{\cref{leftwards}}{=}
        t\
        \begin{tikzpicture}[centerzero={(0,-0.2)}]
            \draw (-0.2,0) to[out=down,in=up] (0.2,-0.4) -- (0.2,-0.6);
            \draw[wipe] (-0.2,-0.6) -- (-0.2,-0.4) to[out=up,in=down] (0.2,0) arc(0:180:0.2);
            \draw[<-] (-0.2,-0.6) -- (-0.2,-0.4) to[out=up,in=down] (0.2,0) arc(0:180:0.2);
            \singdot{-0.2,-0.4};
            \multdot{-0.2,0}{east}{k};
        \end{tikzpicture}
        \stackrel[\cref{rcslide}]{\cref{ldotcross}}{=}
        t\
        \begin{tikzpicture}[centerzero={(0,-0.2)}]
            \draw[<-] (-0.2,-0.6) -- (-0.2,-0.4) to[out=up,in=down] (0.2,0) arc(0:180:0.2);
            \draw[wipe] (-0.2,0) to[out=down,in=up] (0.2,-0.4) -- (0.2,-0.6);
            \draw (-0.2,0) to[out=down,in=up] (0.2,-0.4) -- (0.2,-0.6);
            \multdot{-0.2,0}{east}{k+1};
        \end{tikzpicture}
        \stackrel{\cref{ldotcross}}{=}
        t\
        \begin{tikzpicture}[centerzero={(0,-0.2)}]
            \draw (-0.2,0) to[out=down,in=up] (0.2,-0.4) -- (0.2,-0.6);
            \draw[wipe] (-0.2,-0.6) -- (-0.2,-0.4) to[out=up,in=down] (0.2,0) arc(0:180:0.2);
            \draw[<-] (-0.2,-0.6) -- (-0.2,-0.4) to[out=up,in=down] (0.2,0) arc(0:180:0.2);
            \singdot{0.2,-0.4};
            \multdot{-0.2,0}{east}{k};
        \end{tikzpicture}
        \stackrel{\cref{leftwards}}{=}
        \begin{tikzpicture}[anchorbase]
            \draw[<-] (-0.2,-0.2) -- (-0.2,0) arc (180:0:0.2) -- (0.2,-0.2);
            \singdot{0.2,0};
        \end{tikzpicture}
        \ .
    \]
    Finally, to prove the fourth relation in \cref{piv}, we compose the second relation in \cref{tea1} on the bottom with
    $
        \begin{tikzpicture}[anchorbase]
            \draw[<-] (-0.2,0.2) -- (-0.2,0) arc (180:360:0.2) -- (0.2,0.2);
            \singdot{0.2,0};
        \end{tikzpicture}
    $
    and use \cref{rdotcross,tea0} to get
    \[
        \begin{tikzpicture}[centerzero]
            \draw[<-] (-0.2,0.2) -- (-0.2,0) arc (180:360:0.2) -- (0.2,0.2);
            \singdot{0.2,0};
        \end{tikzpicture}
        = - z^{-1}t^{-1}
        \begin{tikzpicture}[centerzero]
            \draw[<-] (-0.2,0.2) -- (-0.2,0) arc (180:360:0.2) -- (0.2,0.2);
            \diamdec{0,-0.2}{north}{1}{k-1};
        \end{tikzpicture}
        \stackrel{\cref{diamondcup}}{=}
        \begin{tikzpicture}[centerzero]
            \draw[<-] (-0.2,0.2) -- (-0.2,0) arc (180:360:0.2) -- (0.2,0.2);
            \singdot{-0.2,0};
        \end{tikzpicture}
        \ . \qedhere
    \]
\end{proof}

From now on, we will slide tokens and dots over cups and caps, and tokens through crossings, without citing the relevant relations.

\begin{theo}[Infinite Grassmannian relations]
    For any $n \in \Z$, we have
    \begin{equation} \label{infgrass}
        \sum_{r + s = n}
        \begin{tikzpicture}[centerzero]
            \plusleftblank{-0.4,0};
            \plusrightblank{0.4,0};
            \token{-0.6,0}{east}{a};
            \token{0.6,0}{west}{b};
            \multdot{-0.4,-0.2}{north}{r};
            \multdot{0.4,-0.2}{north}{s};
            \teleport{-0.2,0}{0.2,0};
        \end{tikzpicture}
        \ = \sum_{r + s = n}
        \begin{tikzpicture}[centerzero]
            \minusleftblank{-0.4,0};
            \minusrightblank{0.4,0};
            \token{-0.6,0}{east}{a};
            \token{0.6,0}{west}{b};
            \multdot{-0.4,-0.2}{north}{r};
            \multdot{0.4,-0.2}{north}{s};
            \teleport{-0.2,0}{0.2,0};
        \end{tikzpicture}
        \ = - \delta_{n,0} \frac{1}{z} \tr(ab) 1_\one.
    \end{equation}
    Furthermore,
    \begin{align} \label{tanks}
        \plusccbubble{a}{n}
        \ &= \delta_{n,-k} \frac{t}{z} \tr(a) 1_\one \quad \text{if } n \le -k,
        &
        \pluscbubble{a}{n}
        \ &= -\delta_{n,k} \frac{t^{-1}}{z} \tr(a) 1_\one \quad \text{if } n \le k,
        \\ \label{tanks2}
        \minuscbubble{a}{n}
        \ &= \delta_{n,0} \frac{t}{z} \tr(a) 1_\one \quad \text{if } n \ge 0,
        &
        \minusccbubble{a}{n}
        \ &= -\delta_{n,0} \frac{t^{-1}}{z} \tr(a) 1_\one \quad \text{if } n \ge 0.
    \end{align}
\end{theo}

\begin{proof}
    Using $\Omega_k$, it suffices to consider the cases $k \ge 0$.  In addition, using \cref{teleport} and the fact that the $(\pm)$-bubbles are linear in $a$, it suffices to consider \cref{infgrass} in the case where $a=1$.

    First suppose $k = 0$.  Then \cref{tanks} follows immediately from \cref{fake0}, while \cref{tanks2} follows from \cref{fake1,fake2,impose,colder}.  If $n \le 0$, it follows from \cref{fake0,adecomp} that the left-hand side of \cref{infgrass} is equal to the right-hand side.  If $n > 0$, we have
    \[
        t^{-1}\
        \begin{tikzpicture}[centerzero]
            \draw[<-] (0,0.2) arc(90:-270:0.2);
            \multdot{-0.2,0}{east}{n};
            \token{0.2,0}{west}{b};
        \end{tikzpicture}
        \stackrel{\cref{leftwards}}{=}
        \begin{tikzpicture}[centerzero]
            \draw (-0.2,-0.2) to[out=up,in=down] (0.2,0.2) to[out=up,in=up,looseness=1.5] (-0.2,0.2);
            \draw[<-] (-0.2,0.2) to[out=down,in=up] (0.2,-0.2) to[out=down,in=down,looseness=1.5] (-0.2,-0.2);
            \draw[wipe] (-0.2,-0.2) to[out=up,in=down] (0.2,0.2) to[out=up,in=up,looseness=1.5] (-0.2,0.2);
            \draw (-0.2,-0.2) to[out=up,in=down] (0.2,0.2) to[out=up,in=up,looseness=1.5] (-0.2,0.2);
            \multdot{-0.2,-0.2}{east}{n};
            \token{0.2,-0.2}{west}{b};
        \end{tikzpicture}
        \stackrel{\cref{teapos}}{=}
        \begin{tikzpicture}[centerzero]
            \draw (-0.2,-0.2) to[out=up,in=down] (0.2,0.2) to[out=up,in=up,looseness=1.5] (-0.2,0.2);
            \draw (-0.2,-0.2) to[out=up,in=down] (0.2,0.2) to[out=up,in=up,looseness=1.5] (-0.2,0.2);
            \draw[wipe] (-0.2,0.2) to[out=down,in=up] (0.2,-0.2) to[out=down,in=down,looseness=1.5] (-0.2,-0.2);
            \draw[<-] (-0.2,0.2) to[out=down,in=up] (0.2,-0.2) to[out=down,in=down,looseness=1.5] (-0.2,-0.2);
            \multdot{0.2,0.2}{west}{n};
            \token{0.2,-0.2}{west}{b};
        \end{tikzpicture}
        +z \sum_{\substack{r+s=n \\ r,s > 0}}
        \begin{tikzpicture}[centerzero]
            \draw[->] (0,0.5) arc(90:450:0.2);
            \draw[<-] (0,-0.1) arc(90:450:0.2);
            \multdot{0.2,0.3}{west}{r};
            \token{-0.2,0.3}{east}{c};
            \multdot{-0.2,-0.3}{east}{s};
            \token{0.2,-0.3}{west}{c^\vee b};
        \end{tikzpicture}
        = t\, \cbubble{b}{n}
        \ +z \sum_{\substack{r+s=n \\ r,s > 0}}\
        \begin{tikzpicture}[centerzero]
            \draw (-0.2,-0.2) to[out=up,in=down] (0.2,0.2) to[out=up,in=up,looseness=1.5] (-0.2,0.2);
            \draw (-0.2,-0.2) to[out=up,in=down] (0.2,0.2) to[out=up,in=up,looseness=1.5] (-0.2,0.2);
            \draw[wipe] (-0.2,0.2) to[out=down,in=up] (0.2,-0.2) to[out=down,in=down,looseness=1.5] (-0.2,-0.2);
            \draw[<-] (-0.2,0.2) to[out=down,in=up] (0.2,-0.2) to[out=down,in=down,looseness=1.5] (-0.2,-0.2);
            \multdot{0.2,0.2}{west}{n};
            \token{0.2,-0.2}{west}{b};
        \end{tikzpicture},
    \]
    where the final equality follows from the fact that
    \[
        \begin{tikzpicture}[centerzero]
            \draw (-0.2,-0.2) to[out=up,in=down] (0.2,0.2) to[out=up,in=up,looseness=1.5] (-0.2,0.2);
            \draw (-0.2,-0.2) to[out=up,in=down] (0.2,0.2) to[out=up,in=up,looseness=1.5] (-0.2,0.2);
            \draw[wipe] (-0.2,0.2) to[out=down,in=up] (0.2,-0.2) to[out=down,in=down,looseness=1.5] (-0.2,-0.2);
            \draw[<-] (-0.2,0.2) to[out=down,in=up] (0.2,-0.2) to[out=down,in=down,looseness=1.5] (-0.2,-0.2);
            \multdot{0.2,0.2}{west}{n};
            \token{0.2,-0.2}{west}{b};
        \end{tikzpicture}
        \ \stackrel{\cref{lskein}}{=}\
        \begin{tikzpicture}[centerzero]
            \draw (-0.2,-0.2) to[out=up,in=down] (0.2,0.2) to[out=up,in=up,looseness=1.5] (-0.2,0.2);
            \draw[<-] (-0.2,0.2) to[out=down,in=up] (0.2,-0.2) to[out=down,in=down,looseness=1.5] (-0.2,-0.2);
            \draw[wipe] (-0.2,-0.2) to[out=up,in=down] (0.2,0.2) to[out=up,in=up,looseness=1.5] (-0.2,0.2);
            \draw (-0.2,-0.2) to[out=up,in=down] (0.2,0.2) to[out=up,in=up,looseness=1.5] (-0.2,0.2);
            \multdot{0.2,0.2}{west}{n};
            \token{0.2,-0.2}{west}{b};
        \end{tikzpicture}
        \ + z\
        \begin{tikzpicture}[centerzero]
            \draw[<-] (0,0.5) arc(90:450:0.2);
            \draw[->] (0,-0.1) arc(90:450:0.2);
            \multdot{-0.2,0.3}{east}{n};
            \token{0.2,0.3}{west}{c};
            \token{0.2,-0.3}{west}{c^\vee b};
        \end{tikzpicture}
        \ \stackrel[\cref{colder}]{\cref{leftwards}}{=} t^{-1} \cbubble{b}{n} +(t-t^{-1}) \cbubble{b}{n}
        = t\, \cbubble{b}{n}.
    \]
    Thus we have
    \[
        0
        = -\frac{t^{-1}}{z}
        \begin{tikzpicture}[centerzero]
            \draw[<-] (0,0.2) arc(90:-270:0.2);
            \multdot{-0.2,0}{east}{n};
            \token{0.2,0}{west}{b};
        \end{tikzpicture}
        \ + \frac{t}{z}
        \begin{tikzpicture}[centerzero]
            \draw[->] (0,0.2) arc(90:-270:0.2);
            \multdot{-0.2,0}{east}{n};
            \token{0.2,0}{west}{b};
        \end{tikzpicture}
        + \sum_{\substack{r+s=n \\ r,s > 0}}
        \begin{tikzpicture}[centerzero]
            \draw[->] (0,0.5) arc(90:450:0.2);
            \draw[<-] (0,-0.1) arc(90:450:0.2);
            \multdot{0.2,0.3}{west}{r};
            \token{-0.2,0.3}{east}{c};
            \multdot{-0.2,-0.3}{east}{s};
            \token{0.2,-0.3}{west}{c^\vee b};
        \end{tikzpicture}
        \stackrel[\cref{fake2}]{\cref{fake0}}{=} \sum_{\substack{r,s \ge 0 \\ r + s = n}}
        \begin{tikzpicture}[centerzero]
            \plusleft{0,0.3}{c}{r};
            \plusright{0,-0.3}{c^\vee b}{s};
        \end{tikzpicture}
        \ .
    \]

    If $n \ge 0$, then the $(-)$-bubble sum in \cref{infgrass} is equal to $-\delta_{n,0} z^{-2} \tr(ab)$ by \cref{tanks2}.  Now suppose $n < 0$.  Then, as above, we have
    \[
        t^{-1}
        \begin{tikzpicture}[centerzero]
            \draw[<-] (0,0.2) arc(90:-270:0.2);
            \multdot{-0.2,0}{east}{n};
            \token{0.2,0}{west}{b};
        \end{tikzpicture}
        \stackrel{\cref{leftwards}}{=}
        \begin{tikzpicture}[centerzero]
            \draw (-0.2,-0.2) to[out=up,in=down] (0.2,0.2) to[out=up,in=up,looseness=1.5] (-0.2,0.2);
            \draw[<-] (-0.2,0.2) to[out=down,in=up] (0.2,-0.2) to[out=down,in=down,looseness=1.5] (-0.2,-0.2);
            \draw[wipe] (-0.2,-0.2) to[out=up,in=down] (0.2,0.2) to[out=up,in=up,looseness=1.5] (-0.2,0.2);
            \draw (-0.2,-0.2) to[out=up,in=down] (0.2,0.2) to[out=up,in=up,looseness=1.5] (-0.2,0.2);
            \multdot{-0.2,-0.2}{east}{n};
            \token{0.2,-0.2}{west}{b};
        \end{tikzpicture}
        \ \stackrel{\cref{teapos}}{=}\
        \begin{tikzpicture}[centerzero]
            \draw (-0.2,-0.2) to[out=up,in=down] (0.2,0.2) to[out=up,in=up,looseness=1.5] (-0.2,0.2);
            \draw (-0.2,-0.2) to[out=up,in=down] (0.2,0.2) to[out=up,in=up,looseness=1.5] (-0.2,0.2);
            \draw[wipe] (-0.2,0.2) to[out=down,in=up] (0.2,-0.2) to[out=down,in=down,looseness=1.5] (-0.2,-0.2);
            \draw[<-] (-0.2,0.2) to[out=down,in=up] (0.2,-0.2) to[out=down,in=down,looseness=1.5] (-0.2,-0.2);
            \multdot{0.2,0.2}{west}{n};
            \token{0.2,-0.2}{west}{b};
        \end{tikzpicture}
        - z \sum_{\substack{r+s=n \\ r,s \le 0}}
        \begin{tikzpicture}[centerzero={(0,0.1)}]
            \draw[->] (0,0.5) arc(90:450:0.2);
            \draw[<-] (0,-0.1) arc(90:450:0.2);
            \multdot{0.2,0.3}{west}{r};
            \token{-0.2,0.3}{east}{c};
            \multdot{-0.2,-0.3}{east}{s};
            \token{0.2,-0.3}{west}{c^\vee b};
        \end{tikzpicture}
        = t\, \cbubble{b}{n}
        \ - z \sum_{\substack{r+s=n \\ r,s \le 0}}
        \begin{tikzpicture}[centerzero={(0,0.1)}]
            \bubleft{0,0.3}{c}{r};
            \bubright{0,-0.3}{c^\vee b}{s};
        \end{tikzpicture}
        \ .
    \]
    Thus, by \cref{impose,colder,fake1,tanks,tanks2}, we have
    \[
        \sum_{\substack{r+s=n \\ r,s < 0}}
        \begin{tikzpicture}[centerzero]
            \plusleft{0,0.3}{c}{r};
            \plusright{0,-0.3}{c^\vee b}{s};
        \end{tikzpicture}
        = 0.
    \]

    Now suppose $k > 0$.  Note that
    \begin{equation} \label{slide8}
        z^{-1}
        \begin{tikzpicture}[centerzero]
            \draw[->] (-0.2,-0.2) to[out=up,in=down] (0.2,0.2) to[out=up,in=up,looseness=1.5] (-0.2,0.2);
            \draw[wipe] (-0.2,0.2) to[out=down,in=up] (0.2,-0.2) to[out=down,in=down,looseness=1.5] (-0.2,-0.2);
            \draw (-0.2,0.2) to[out=down,in=up] (0.2,-0.2) to[out=down,in=down,looseness=1.5] (-0.2,-0.2);
            \multdot{-0.2,-0.2}{east}{n};
            \token{0.2,-0.2}{west}{b};
        \end{tikzpicture}
        \stackrel{\cref{leftwards}}{=} z^{-1} t
        \begin{tikzpicture}[centerzero]
            \draw[->] (-0.2,-0.5) to[out=up,in=down] (0.2,0) to[out=up,in=down] (-0.2,0.5) to[out=up,in=up,looseness=1.5] (0.2,0.5);
            \draw[wipe] (0.2,0.5) to[out=down,in=up] (-0.2,0) to[out=down,in=up] (0.2,-0.5) to[out=down,in=down,looseness=1.5] (-0.2,-0.5);
            \draw (0.2,0.5) to[out=down,in=up] (-0.2,0) to[out=down,in=up] (0.2,-0.5) to[out=down,in=down,looseness=1.5] (-0.2,-0.5);
            \multdot{-0.2,-0.5}{east}{n};
            \multdot{-0.2,0.5}{east}{k};
            \token{0.2,-0.5}{west}{b};
        \end{tikzpicture}
        \stackrel{\cref{tea1}}{=} z^{-1}t\ \cbubble{b}{k+n}
        -z^{-1}t \sum_{r=0}^{k-1}
        \begin{tikzpicture}[centerzero]
            \draw[<-] (0,0.6) arc(90:450:0.2);
            \multdot{-0.2,0.4}{east}{k};
            \diamdec{0,0.2}{north east}{c}{r};
            \bubright{0,-0.4}{c^\vee b}{r+n};
        \end{tikzpicture}
        \stackrel{\cref{fake0}}{=}
        \sum_{\substack{r+s=n \\ r \le 0}}
        \begin{tikzpicture}[centerzero]
            \plusleft{0,0.3}{c}{r};
            \bubright{0,-0.3}{c^\vee b}{s};
        \end{tikzpicture}
        .
    \end{equation}
    Now consider the second relation in \cref{tanks}.  If $n \le 0$, this relation holds immediately by the definition \cref{fake0}.  On the other hand, if $0 < n \le k$, then it holds by \cref{fake2,tea0}.  The first relation in \cref{tanks} follows immediately from \cref{fake0}.

    Now we prove that the left-hand side of \cref{infgrass} is equal to the right-hand side.  If $n < 0$, then, by \cref{fake0}, we have
    \[
        \pluscbubble{c^\vee b}{s}
        \ = 0 \text{ when } s < k
        \qquad \text{and} \qquad
        \plusccbubble{ac}{r}
        \ = 0 \text{ when } r < -k,
    \]
    and so the sum on the left-hand side of \cref{infgrass} zero.  Similarly, if $n=0$, then only the $s=k$, $r=-k$ terms survives, and so the sum is equal to
    \[
        - \frac{1}{z} \tr(ac) \tr(c^\vee b) 1_\one
        \stackrel{\cref{adecomp}}{=} -\frac{1}{z} \tr(ab) 1_\one.
    \]
    Now assume $n>0$.  Using \cref{teapos}, the left-hand side of \cref{slide8} becomes
    \[
        z^{-1}
        \begin{tikzpicture}[centerzero]
            \draw (-0.2,0.2) to[out=down,in=up] (0.2,-0.2) to[out=down,in=down,looseness=1.5] (-0.2,-0.2);
            \draw[wipe] (-0.2,-0.2) to[out=up,in=down] (0.2,0.2) to[out=up,in=up,looseness=1.5] (-0.2,0.2);
            \draw[->] (-0.2,-0.2) to[out=up,in=down] (0.2,0.2) to[out=up,in=up,looseness=1.5] (-0.2,0.2);
            \multdot{0.2,0.2}{west}{n};
            \token{0.2,-0.2}{west}{b};
        \end{tikzpicture}
        - \sum_{\substack{r+s=n \\ r,s > 0}}
        \begin{tikzpicture}[centerzero]
            \bubleft{0,0.3}{c}{r};
            \bubright{0,-0.3}{c^\vee b}{s};
        \end{tikzpicture}
        \stackrel[\cref{fake2}]{\cref{tea0}}{=}
        - \sum_{\substack{r+s=n \\ r,s > 0}}
        \begin{tikzpicture}[centerzero]
            \plusleft{0,0.3}{c}{r};
            \bubright{0,-0.3}{c^\vee b}{s};
        \end{tikzpicture}
        .
    \]
    Then it follows from \cref{fake2,tanks} that
    \[
        \sum_{r + s = n}
        \begin{tikzpicture}[centerzero]
            \plusleft{0,0.3}{c}{r};
            \plusright{0,-0.3}{c^\vee b}{s};
        \end{tikzpicture}
        = 0,
    \]
    as desired.

    Next we prove that the $(-)$-bubble sum in \cref{infgrass} is equal to the right-hand side.   By \cref{fake1,fake2}, we have
    \[
        \minuscbubble{a}{r} = 0
        \qquad \text{and} \qquad
        \plusccbubble{a}{r} = 0
        \qquad \text{when } r > 0.
    \]
    Thus the identity holds when $n > 0$.  Now suppose $n \le 0$.  Then, using \cref{teapos} on the left-hand side, the relation \cref{slide8} becomes
    \[
        \sum_{\substack{r+s=n \\ r \le 0}}
        \begin{tikzpicture}[centerzero]
            \plusleft{0,0.3}{c}{r};
            \bubright{0,-0.3}{c^\vee b}{s};
        \end{tikzpicture}
        \ =\
        z^{-1}
        \begin{tikzpicture}[centerzero]
            \draw (-0.2,0.2) to[out=down,in=up] (0.2,-0.2) to[out=down,in=down,looseness=1.5] (-0.2,-0.2);
            \draw[wipe] (-0.2,-0.2) to[out=up,in=down] (0.2,0.2) to[out=up,in=up,looseness=1.5] (-0.2,0.2);
            \draw[->] (-0.2,-0.2) to[out=up,in=down] (0.2,0.2) to[out=up,in=up,looseness=1.5] (-0.2,0.2);
            \multdot{0.2,0.2}{west}{n};
        \end{tikzpicture}
        \ + \sum_{\substack{r+s=n \\ r,s \le 0}}
        \begin{tikzpicture}[centerzero]
            \bubleft{0,0.3}{c}{r};
            \bubright{0,-0.3}{c^\vee b}{s};
        \end{tikzpicture}
        \ \stackrel{\cref{tea0}}{=} \sum_{\substack{r+s=n \\ r,s \le 0}}
        \begin{tikzpicture}[centerzero]
            \bubleft{0,0.3}{c}{r};
            \bubright{0,-0.3}{c^\vee b}{s};
        \end{tikzpicture}
        \ .
    \]
    Using \cref{fake1}, we then have
    \[
        \sum_{\substack{r+s=n \\ r \le 0}}
        \begin{tikzpicture}[centerzero]
            \plusleft{0,0.3}{c}{r};
            \bubright{0,-0.3}{c^\vee b}{s};
        \end{tikzpicture}
        \ =\
        \sum_{\substack{r+s=n \\ r,s \le 0}}
        \left(
            \begin{tikzpicture}[centerzero]
                \plusleft{0,0.3}{c}{r};
                \bubright{0,-0.3}{c^\vee b}{s};
            \end{tikzpicture}
            \ +\
            \begin{tikzpicture}[centerzero]
                \minusleft{0,0.3}{c}{r};
                \plusright{0,-0.3}{c^\vee b}{s};
            \end{tikzpicture}
            \ +\
            \begin{tikzpicture}[centerzero]
                \minusleft{0,0.3}{c}{r};
                \minusright{0,-0.3}{c^\vee b}{s};
            \end{tikzpicture}
        \right).
    \]
    Thus, by \cref{tanks}, we have
    \[
        \sum_{r+s=n}
        \begin{tikzpicture}[centerzero]
            \minusleft{0,0.3}{c}{r};
            \minusright{0,-0.3}{c^\vee b}{s};
        \end{tikzpicture}
        = \sum_{\substack{r+s=n \\ r \le 0,\ s > 0}}
        \begin{tikzpicture}[centerzero]
            \plusleft{0,0.3}{c}{r};
            \bubright{0,-0.3}{c^\vee b}{s};
        \end{tikzpicture}
        \stackrel[\cref{tea0}]{\cref{tanks}}{=} - \delta_{n,0} \frac{1}{z^2} \tr(c) \tr(c^\vee b) 1_\one
        \stackrel{\cref{adecomp}}{=} - \delta_{n,0} \frac{1}{z^2} \tr(b) 1_\one.
    \]

    It remains to prove \cref{tanks2}.  When $n > 0$, these relations follows immediately from \cref{fake1,fake2}.  Now,
    \[
        \frac{t}{z} \tr(a) 1_\one
        \stackrel{\cref{impose}}{=}
        \begin{tikzpicture}[centerzero]
            \draw[<-] (0,0.2) arc(90:450:0.2);
            \token{0.2,0}{west}{a};
        \end{tikzpicture}
        \stackrel{\cref{fake1}}{=}
        \pluscbubble{a}{0} + \minuscbubble{a}{0}
        \stackrel{\cref{tanks}}{=} \minuscbubble{a}{0}.
    \]
    Then, by \cref{infgrass} with $n=0$, we have
    \[
        -\frac{t^{-1}}{z} \tr(a) 1_\one
        = z t^{-1}
        \begin{tikzpicture}[centerzero]
            \minusleft{0,0.3}{ac}{0};
            \minusright{0,-0.3}{c^\vee}{0};
        \end{tikzpicture}
        = \tr(c^\vee) \minusccbubble{ac}{0}
        \stackrel{\cref{adecomp}}{=} \minusccbubble{a}{0}. \qedhere
    \]
\end{proof}

It will be useful to expression some of our relations in terms of generating functions.  For $a \in A$ and an indeterminate $w$, let
\begin{align}
    \\ \label{bubgen1}
    \begin{tikzpicture}[centerzero]
        \plusgenleft{0,0};
        \token{-0.2,0}{east}{a};
    \end{tikzpicture}
    &:= t^{-1}z \sum_{r \in \Z} \plusccbubble{a}{r} w^{-r} \in \tr(a) 1_\one + w^{k-1} \End_{\Heis_k(A;z,t)}(\one) \llbracket w^{-1} \rrbracket,
    \\ \label{bubgen2}
    \begin{tikzpicture}[centerzero]
        \plusgenright{0,0};
        \token{0.2,0}{west}{a};
    \end{tikzpicture}
    &:= -tz \sum_{r \in \Z} \pluscbubble{a}{r} w^{-r} \in \tr(a) u^{-k} 1_\one + w^{-k-1} \End_{\Heis_k(A;z,t)}(\one) \llbracket w^{-1} \rrbracket,
    \\ \label{bubgen3}
    \begin{tikzpicture}[centerzero]
        \minusgenleft{0,0};
        \token{-0.2,0}{east}{a};
    \end{tikzpicture}
    &:=  -tz \sum_{r \in \Z} \minusccbubble{a}{r} w^{-r} \in \tr(a) w^k 1_\one + w \End_{\Heis_k(A;z,t)}(\one) \llbracket w \rrbracket,
    \\ \label{bubgen4}
    \begin{tikzpicture}[centerzero]
        \minusgenright{0,0};
        \token{0.2,0}{west}{a};
    \end{tikzpicture}
    &:= t^{-1} z \sum_{r \in \Z} \minuscbubble{a}{r} w^{-r} \in \tr(a) 1_\one + w \End_{\Heis_k(A;z,t)}(\one) \llbracket w \rrbracket.
\end{align}

\begin{cor}
    We have
    \begin{equation} \label{infweeds}
        \begin{tikzpicture}[centerzero]
            \plusgenleft{-0.4,0};
            \plusgenright{0.4,0};
            \token{-0.6,0}{east}{a};
            \token{0.6,0}{west}{b};
            \teleport{-0.2,0}{0.2,0};
        \end{tikzpicture}
        =
        \begin{tikzpicture}[centerzero]
            \minusgenleft{-0.4,0};
            \minusgenright{0.4,0};
            \token{-0.6,0}{east}{a};
            \token{0.6,0}{west}{b};
            \teleport{-0.2,0}{0.2,0};
        \end{tikzpicture}
        \ = z \tr(ab) 1_\one.
    \end{equation}
\end{cor}

We adopt the convention that determinants of matrices whose entries lie in a superalgebra are to be computed using the usual Laplace expansions, keeping the non-commuting variables in each monomial ordered in the same way as the columns from which they are taken (see \cite[(17)]{Sav19}).

\begin{lem} \label{dlem}
    For all $a \in A$, we have
    \begin{align} \label{d1}
        \plusccbubble{a}{r-k}
        &= z^{r-1} t^{r+1} \sum_{b_1,\dotsc,b_{r-1} \in \BA} \det
        \left( \cbubble{b^\vee_{j-1}b_j}{i-j+k+1} \right)_{i,j=1}^r,
        \quad r \le k,
        \\ \label{d2}
        \pluscbubble{a}{r+k}
        &= - z^{r-1} t^{-r-1} \sum_{b_1,\dotsc,b_{r-1} \in \BA} \det
        \left( - \ccbubble{b^\vee_{j-1}b_j}{i-j-k+1} \right)_{i,j=1}^r,
        \quad r \le -k,
    \end{align}
    adopting the convention that $b^\vee_0 := a$ and $b_r:=1$.  We interpret these determinants as $\tr(a)$ if $r=0$ or as $0$ if $r<0$.
\end{lem}

\begin{proof}
    It suffices to prove \cref{d1}, since \cref{d2} then follows by applying $\Omega_k$.  If $k < 0$, then the right-hand side in \cref{d1} is zero by our convention, and so the result follows by \cref{tanks}.  Now suppose $k \ge 0$.  The cases $r \le 0$ follow immediately from \cref{fake0}.  Now assume that $r > 0$ and that the result holds for $r-1$.  For $b_1,\dotsc,b_{r-1} \in \BA$, define the matrix
    \[
        A = \left( \cbubble{b^\vee_{j-1}b_j}{i-j+k+1} \right)_{i,j=1}^r.
    \]
    (We leave the dependence on $b_1,\dotsc,b_{r-1} \in \BA$ implicit to simplify the notation.)  We have
    \[
        \det A = \sum_{m=1}^r (-1)^{m+1} \cbubble{a b_1}{m+k} \det A_{m,1},
    \]
    where $A_{m,1}$ is the $(r-1) \times (r-1)$ matrix obtained from $A$ by removing the $m$-th row and first column.  If we consider $A_{m,1}$ as a block matrix with upper-left block of size $(m-1) \times (m-1)$ and lower-right block of size $(r-m) \times (r-m)$, we see that it is block lower triangular.  By \cref{tea0}, the upper-left block is lower triangular with diagonal entries
    \[
        -t^{-1}z^{-1} \tr(b^\vee_1b_2),\
        -t^{-1}z^{-1} \tr(b^\vee_2b_3),\
        \dotsc\ ,\
        -t^{-1}z^{-1} \tr(b^\vee_{m-1}b_m).
    \]
    On the other hand, the lower-right block is the matrix
    \[
        \left( \cbubble{b^\vee_{m+j-1} b_{m+j}}{i-j+k+1} \right)_{i,j=1}^{r-m}.
    \]
    Thus, using the induction hypothesis and \cref{chkdef}, we have
    \[
        z^{r-1} t^{r+1} \sum_{b_1,\dotsc,b_{r-1} \in \BA} \det A
        = zt \sum_{m=1}^r
        \begin{tikzpicture}[centerzero]
            \bubright{0,0.3}{ab_1}{m+k};
            \plusleft{0,-0.3}{b^\vee_1}{r-k-m};
        \end{tikzpicture}
        \stackrel{\cref{fake2}}{=} zt \sum_{m=1}^r
        \begin{tikzpicture}[centerzero]
            \plusright{0,0.3}{ab_1}{m+k};
            \plusleft{0,-0.3}{b^\vee_1}{r-k-m};
        \end{tikzpicture}
        \overset{\substack{\cref{infgrass} \\ \cref{tanks}}}{\underset{\cref{adecomp}}{=}}
        \plusccbubble{a}{r-k}\ . \qedhere
    \]
\end{proof}

\begin{lem}\label{boomerang}
    For all $r \in \Z$ and $a,b \in A$, we have
    \begin{align} \label{boomerang+}
        \plusccbubble{ab}{r}
        &= (-1)^{\bar{a} \bar{b}} \plusccbubble{ba}{r}
        \ ,&
        \pluscbubble{ab}{r}
        &= (-1)^{\bar{a} \bar{b}} \pluscbubble{ba}{r}
        \ ,
        \\ \label{boomerang-}
        \minusccbubble{ab}{r}
        &= (-1)^{\bar{a} \bar{b}} \minusccbubble{ba}{r}
        \ ,&
        \minuscbubble{ab}{r}
        &= (-1)^{\bar{a} \bar{b}} \minuscbubble{ba}{r}
        \ .
    \end{align}
\end{lem}

\begin{proof}
    Writing out the columns of the matrix in \cref{d1}, we have
    \begin{align*}
        &z^{-r+1} t^{-r-1} \plusccbubble{ab}{r-k}
        \\
        &\quad \overset{\mathclap{\cref{d1}}}{=} \sum_{b_1,\dotsc,b_{n+k-1} \in \BA} \det
        \begin{pmatrix}
            \cbubble{abb_1}{i+k} &
            \cbubble{b_1^\vee b_2}{i-1+k} &
            \cdots &
            \cbubble{b_{r-1}^\vee}{i-r+k+1}
        \end{pmatrix}_{i=1}^r
        \\
        &\quad \overset{\mathclap{\cref{beam}}}{=} \sum_{b_1,\dotsc,b_{r-1} \in \BA} \det
        \begin{pmatrix}
            \cbubble{ab_1}{i+k} &
            \cbubble{b_1^\vee b_2}{i-1+k} &
            \cdots &
            \cbubble{b_{r-1}^\vee b}{i-r+k+1}
        \end{pmatrix}_{i=1}^r
        \\
        &\quad = \sum_{b_1,\dotsc,b_{n+k-1} \in \BA}  (-1)^{\bar{b} \bar{b}_{r-1}} \det
        \begin{pmatrix}
            \cbubble{ab_1}{i+k} &
            \cbubble{b_1^\vee b_2}{i-1+k} &
            \cdots &
            \cbubble{b b_{r-1}^\vee}{i-r+k+1}
        \end{pmatrix}_{i=1}^r
        \\
        &\quad \stackrel{\mathclap{\cref{beam}}}{=} \sum_{b_1,\dotsc,b_{r-1} \in \BA}  (-1)^{\bar{b} \bar{b}_1} \det
        \begin{pmatrix}
            \cbubble{ab_1b}{i+k} &
            \cbubble{b_1^\vee b_2}{i-1+k} &
            \cdots &
            \cbubble{b_{r-1}^\vee}{i-r+k+1}
        \end{pmatrix}_{i=1}^r
        \\
        &\quad = \sum_{b_1,\dotsc,b_{n+k-1} \in \BA}  (-1)^{\bar{a} \bar{b}} \det
        \begin{pmatrix}
            \cbubble{bab_1}{i+k} &
            \cbubble{b_1^\vee b_2}{i-1+k} &
            \cdots &
            \cbubble{b_{r-1}^\vee}{i-r+k+1}
        \end{pmatrix}_{i=1}^r
        \\
        &\quad \stackrel{\mathclap{\cref{d1}}}{=}\ (-1)^{\bar{a} \bar{b}} z^{-r+1} t^{-r-1} \plusccbubble{ba}{r-k}.
    \end{align*}
    This completes the proof of the first relation in \cref{boomerang+}.  The second follows by applying $\Omega_k$.  Finally, \cref{boomerang-} follows from \cref{boomerang+,fake1}.
\end{proof}

It follows from \cref{boomerang} that we can label the tokens in $(\pm)$-bubbles by elements of the cocenter $C(A)$ of $A$.  The following lemma, together with \cref{nakano1,nakano2}, allows us to eliminate the leftwards caps decorated by $\color{violet}\scriptstyle{\diamondsuit}$ and $\color{violet}\scriptstyle{\heartsuit}$ from any diagram.

\begin{lem} \label{redundant}
    The following relations hold:
    \begin{align} \label{nakano4}
        \begin{tikzpicture}[centerzero]
            \draw[<-] (-0.2,0.2) -- (-0.2,0) arc (180:360:0.2) -- (0.2,0.2);
            \diamdec{0,-0.2}{north}{a}{n};
        \end{tikzpicture}
        \ &= -z^2 \sum_{r \ge 1}
        \begin{tikzpicture}[centerzero]
            \draw[<-] (-0.2,0.5) -- (-0.2,0.3) arc (180:360:0.2) -- (0.2,0.5);
            \multdot{-0.2,0.3}{east}{r};
            \token{0.2,0.3}{west}{b};
            \plusleft{0,-0.3}{b^\vee a}{-r-n};
        \end{tikzpicture}
        \quad \text{if } 0 \le n < k,
        \\ \label{nakano5}
        \begin{tikzpicture}[centerzero]
            \draw[<-] (-0.2,-0.2) -- (-0.2,0) arc (180:0:0.2) -- (0.2,-0.2);
            \diamdec{0,0.2}{south}{a}{n};
        \end{tikzpicture}
        \ &= - (-1)^{\bar{a}} z^2 \sum_{r \ge 1}
        \begin{tikzpicture}[centerzero]
            \draw[<-] (-0.2,-0.5) -- (-0.2,-0.3) arc(180:0:0.2) -- (0.2,-0.5);
            \multdot{0.2,-0.3}{west}{r};
            \token{-0.2,-0.3}{east}{b^\vee};
            \plusright{0,0.3}{ab}{-r-n};
        \end{tikzpicture}
        \quad \text{if } 0 \le n < -k.
    \end{align}
\end{lem}

\begin{proof}
    We first prove \cref{nakano4}, and so we assume $0 \le n < k$.  Composing the second relation in \cref{tea1} with the right side of \cref{nakano4} gives
    \[
        -z^2 \sum_{r \ge 1}\
        \begin{tikzpicture}[centerzero]
            \draw[<-] (-0.2,0.5) -- (-0.2,0.3) arc (180:360:0.2) -- (0.2,0.5);
            \multdot{-0.2,0.3}{east}{r};
            \token{0.2,0.3}{west}{b};
            \plusleft{0,-0.3}{b^\vee a}{-r-n};
        \end{tikzpicture}
        = -z^2 \sum_{r \ge 1}
        \begin{tikzpicture}[anchorbase]
            \draw[->] (-0.2,0) to[out=up,in=down] (0.2,0.5) to[out=up,in=down] (-0.2,1);
            \draw[wipe] (0.2,0) to[out=up,in=down] (-0.2,0.5) to[out=up,in=down] (0.2,1);
            \draw (0.2,0) to[out=up,in=down] (-0.2,0.5) to[out=up,in=down] (0.2,1);
            \draw (0.2,0) arc(0:-180:0.2);
            \token{0.2,0}{west}{b};
            \multdot{-0.2,0}{east}{r};
            \plusleft{0,-0.5}{b^\vee a}{-r-n};
        \end{tikzpicture}
        - z^2 \sum_{r \ge 1} \sum_{s=0}^{k-1}
        \begin{tikzpicture}[centerzero]
            \draw[<-] (-0.2,0.8) -- (-0.2,0.6) arc(180:360:0.2) -- (0.2,0.8);
            \diamdec{0,0.4}{west}{c}{s};
            \bubright{0,-0.05}{c^\vee b}{r+s};
            \plusleft{0,-0.6}{b^\vee a}{-r-n};
        \end{tikzpicture}
        .
    \]
    By \cref{tanks}, the $(+)$-bubble in the first sum on the right side above is zero unless $r \le k$.  But for these values of $r$ this term is zero by \cref{exercise}.  Thus the first sum on the right side above vanishes.  Now, in the second (double) sum, \cref{fake2} allows us to replace the clockwise bubble by $\pluscbubble{c^\vee b}{r+s}$.  We can also replace the sum over $r \ge 1$ with a sum over $r \in \Z$, since the additional terms are zero by \cref{tanks}.  Then \cref{nakano4} follows by \cref{infgrass,adecomp}.

    Now, applying to $\Omega_k$ to \cref{nakano4}, we have, for $0 \le n < -k$,
    \[
        \begin{tikzpicture}[centerzero]
            \draw[<-] (-0.2,-0.2) -- (-0.2,0) arc (180:0:0.2) -- (0.2,-0.2);
            \diamdec{0,0.2}{south}{a}{n};
        \end{tikzpicture}
        \ = - (-1)^{\bar{b} + \bar{b}\bar{a} + \bar{a}} z^2 \sum_{r \ge 1}
        \begin{tikzpicture}[centerzero]
            \draw[<-] (-0.2,-0.5) -- (-0.2,-0.3) arc(180:0:0.2) -- (0.2,-0.5);
            \multdot{0.2,-0.3}{west}{r};
            \token{-0.2,-0.3}{east}{b};
            \plusright{0,0.3}{b^\vee a}{-r-n};
        \end{tikzpicture}
        \ \stackrel[\cref{doubledual}]{\cref{boomerang+}}{=} - (-1)^{\bar{a}} z^2 \sum_{r \ge 1}
        \begin{tikzpicture}[centerzero]
            \draw[<-] (-0.2,-0.5) -- (-0.2,-0.3) arc(180:0:0.2) -- (0.2,-0.5);
            \multdot{0.2,-0.3}{west}{r};
            \token{-0.2,-0.3}{east}{b^\vee};
            \plusright{0,0.3}{ab}{-r-n};
        \end{tikzpicture}
        \ .\qedhere
    \]
\end{proof}

\begin{cor}
    The following relations hold for all $a \in A$:
    \begin{align} \label{nakano6}
        \begin{tikzpicture}[>=To,baseline={([yshift=1ex]current bounding box.center)}]
            \draw[<-] (-0.2,0.2) -- (-0.2,0) arc (180:360:0.2) -- (0.2,0.2);
            \diamdec{0,-0.2}{north}{a}{n};
        \end{tikzpicture}
        &=
        \begin{tikzpicture}[>=To,baseline={([yshift=1ex]current bounding box.center)}]
            \draw[<-] (-0.2,0.2) -- (-0.2,0) arc (180:360:0.2) -- (0.2,0.2);
            \diamdec{0,-0.2}{north}{1}{n};
            \token{-0.2,0}{east}{a};
        \end{tikzpicture}
        \ ,&
        \begin{tikzpicture}[>=To,baseline={([yshift=1ex]current bounding box.center)}]
            \draw[<-] (-0.2,0.2) -- (-0.2,0) arc (180:360:0.2) -- (0.2,0.2);
            \heartdec{0,-0.2}{north}{a}{n};
        \end{tikzpicture}
        &=
        \begin{tikzpicture}[>=To,baseline={([yshift=1ex]current bounding box.center)}]
            \draw[<-] (-0.2,0.2) -- (-0.2,0) arc (180:360:0.2) -- (0.2,0.2);
            \heartdec{0,-0.2}{north}{1}{n};
            \token{-0.2,0}{east}{a};
        \end{tikzpicture}
        \ ,\quad 0 \le n < k,
        \\ \label{nakano7}
        \begin{tikzpicture}[>=To,baseline={([yshift=-2ex]current bounding box.center)}]
            \draw[<-] (-0.2,-0.2) -- (-0.2,0) arc (180:0:0.2) -- (0.2,-0.2);
            \diamdec{0,0.2}{south}{a}{n};
        \end{tikzpicture}
        &= (-1)^{\bar a}
        \begin{tikzpicture}[>=To,baseline={([yshift=-2ex]current bounding box.center)}]
            \draw[<-] (-0.2,-0.2) -- (-0.2,0) arc (180:0:0.2) -- (0.2,-0.2);
            \diamdec{0,0.2}{south}{1}{n};
            \token{-0.2,0}{east}{a};
        \end{tikzpicture}
        \ ,&
        \begin{tikzpicture}[>=To,baseline={([yshift=-2ex]current bounding box.center)}]
            \draw[<-] (-0.2,-0.2) -- (-0.2,0) arc (180:0:0.2) -- (0.2,-0.2);
            \heartdec{0,0.2}{south}{a}{n};
        \end{tikzpicture}
        &= (-1)^{\bar a}
        \begin{tikzpicture}[>=To,baseline={([yshift=-2ex]current bounding box.center)}]
            \draw[<-] (-0.2,-0.2) -- (-0.2,0) arc (180:0:0.2) -- (0.2,-0.2);
            \heartdec{0,0.2}{south}{1}{n};
            \token{-0.2,0}{east}{a};
        \end{tikzpicture}
        \ ,\quad 0 \le n < -k.
    \end{align}
\end{cor}

\begin{proof}
    The relations involving the $\color{violet}\scriptstyle{\diamondsuit}$-decorated left cup and cap follow immediately from \cref{nakano4,nakano5,beam}.  Then the relations involving the $\color{violet}\scriptstyle{\heartsuit}$-decorated left cup and cap follow from \cref{nakano1,nakano3}.
\end{proof}


The final goal of this section is to endow $\Heis_k(A;z,t)$ with the structure of a strict pivotal category.

\begin{lem}
    We have
    \begin{equation} \label{water}

        \qquad \text{if } k \ge 0.
    \end{equation}
\end{lem}

\begin{proof}
    Suppose $k \ge 0$.  We first claim that

    \noindent\begin{minipage}{0.5\linewidth}
        \begin{equation} \label{melon1}
            %
            \ .
        \end{equation}
    \end{minipage}\par\vspace{\belowdisplayskip}

    \noindent To prove \cref{melon1}, we note that, composing on the top of the two leftmost strands with the invertible map \cref{invrel}, it suffices to prove that

    \noindent\begin{minipage}{0.5\linewidth}
        \begin{equation} \label{melon1a}
            %
        \ .
    \]
    (We note that the above step is the reason we need to impose the relation \cref{impose}.)  On the other hand, if $0 < n < k$, we have (recall the notation $a^\dagger$ from \cref{adag})
    \begin{multline*}
        %
        \ ,
    \end{multline*}
    where, in the third and fourth equalities, we used that
    \begin{equation} \label{curlpop}
        %
        \stackrel{\cref{exercise}}{=} 0
        \quad \text{for } 0 \le s < k.
    \end{equation}

    To prove \cref{melon2}, we note that, composing on the top of the two rightmost strands with the invertible map \cref{invrel}, it suffices to prove that

    \noindent\begin{minipage}{0.5\linewidth}
        \begin{equation} \label{melon2a}
            %
            \ ,\quad 0 \le n < k,\ a \in A.
        \end{equation}
    \end{minipage}\par\vspace{\belowdisplayskip}

    \noindent The proofs of \cref{melon2a,melon2b} are similar to those of \cref{melon1a,melon1b} and are left to the reader.
    \details{
        To see \cref{melon2a}, we compute
        \begin{multline*}
            %
            \ .
        \]
        On the other hand, if $0 < n < k$, we have
        \begin{multline*}
            %
            \ ,
        \end{multline*}
        where, in the third and fourth equalities, we used the fact that
        \[
            %
            \stackrel{\cref{exercise}}{=} 0
            \quad \text{for } 0 \le s < k.
        \]
    }

    Now, to prove the left-hand relation in \cref{water}, we compose \cref{melon1} on the top with
    %
        \ .
    \]
    Then the left relation in \cref{water} follows from \cref{curlpop}.  The right relation in \cref{water} follows similarly from composing \cref{melon2} on the top with
    %
            \ .
        \]
        The right relation in \cref{water} then follows from the fact that
        \[
            %
            \stackrel{\cref{exercise}}{=} 0
            \quad \text{for } 0 \le s < k.
        \]
    }
\end{proof}

\begin{lem}
    The following relations hold:
    \begin{equation} \label{adjfinal}
        %
        \ .
    \end{equation}
\end{lem}

\begin{proof}
    It suffices to consider the case $k \ge 0$, since the case $k < 0$ then follows after applying $\Omega_k$.  For the first relation, we have
    \[
        %
        \ .
    \]
    For the second relation, we have
    \[
        %
        \ . \qedhere
    \]
\end{proof}

\begin{lem}
    The following relations hold:
    \begin{align} \label{pitchl1}
        %
        \ .
    \end{align}
\end{lem}

\begin{proof}
    It suffices to consider the case $k \ge 0$, since the case $k < 0$ then follows by applying $\Omega_k$.  The second and third relations in \cref{pitchl1} are \cref{water}.  To see the first relation in \cref{pitchl1}, we compute
    \[
        %
        \ .
    \]
    The proof of the fourth relation in \cref{pitchl1} is analogous.

    The relations \cref{pitchl2} now follow by attaching left caps to the relations in \cref{pitchl1} and using \cref{adjfinal}.  For example, attaching left caps to the top left and top right strands of the first relation in \cref{pitchl1} gives
    \[
        %
        \ .
    \]
    Then the first relation in \cref{pitchl2} follows from \cref{adjfinal}.
\end{proof}

It follows from \cref{pitchl1,pitchl2,piv} that the category $\Heis_k(A;z,t)$ is strictly pivotal, with duality function
\begin{equation} \label{pivot}
    * \colon \Heis_k(A;z,t) \xrightarrow{\cong} \left( \left( \Heis_k(A;z,t) \right)^\op \right)^\rev
\end{equation}
defined by rotating diagrams through $180\degree$ and multiplying by $(-1)^{\binom{y}{2}}$, where $y$ is the number of odd tokens in the diagram.  Intuitively, this means that morphisms are invariant under isotopy fixing the endpoints (multiplying by the appropriate sign when odd elements change height).  For this reason, we will allow ourselves to draw tokens and dots at the critical points of cups and caps, and this does not give rise to any ambiguity.

\section{Third approach\label{sec:third}}

In this section we give our third definition of $\Heis_k(A;z,t)$.  In this approach, the inversion relation is replaced by explicit relations amongst generators.  To do this, we must add a left cup and cap as generating morphisms.


\begin{defin} \label{def3}
    The \emph{quantum Frobenius Heisenberg category} $\Heis_k(A;z,t)$ is the strict $\kk$-linear monoidal supercategory obtained from $\QAW(A;z,t)$ by adjoining a right dual $\downarrow$ to $\uparrow$, plus two more generating morphisms
    \begin{tikzpicture}
        \draw[<-] (-0.15,-0.1) -- (-0.15,0) arc(180:0:0.15) -- (0.15,-0.1);
    \end{tikzpicture}
    and
    \begin{tikzpicture}[anchorbase]
        \draw[<-] (-0.15,0.1) -- (-0.15,0) arc(180:360:0.15) -- (0.15,0.1);
    \end{tikzpicture},
    subject to the following additional relations:
    \begin{gather} \label{pos}
        \begin{tikzpicture}[centerzero]
            \draw[->] (0.2,0) \braidup (-0.2,0.5);
            \draw[<-] (0.2,-0.5) \braidup (-0.2,0);
            \draw[wipe] (-0.2,-0.5) \braidup (0.2,0);
            \draw (-0.2,-0.5) \braidup (0.2,0);
            \draw[wipe] (-0.2,0) \braidup (0.2,0.5);
            \draw (-0.2,0) \braidup (0.2,0.5);
        \end{tikzpicture}
        \ =\
        \begin{tikzpicture}[centerzero]
            \draw[->] (-0.2,-0.5) to (-0.2,0.5);
            \draw[<-] (0.2,-0.5) to (0.2,0.5);
        \end{tikzpicture}
        \ - t^{-1}\
        \begin{tikzpicture}[centerzero]
            \draw[<-] (-0.2,0.5) -- (-0.2,0.3) arc(180:360:0.2) -- (0.2,0.5);
            \draw[->] (-0.2,-0.5) -- (-0.2,-0.3) arc(180:0:0.2) -- (0.2,-0.5);
            \teleport{-0.2,-0.3}{-0.2,0.3};
        \end{tikzpicture}
        \ + \sum_{r,s > 0}\
        \begin{tikzpicture}[anchorbase]
            \draw[<-] (-0.2,0.6) to (-0.2,0.35) arc(180:360:0.2) to (0.2,0.6);
            \draw[->] (-0.2,-0.6) to (-0.2,-0.35) arc(180:0:0.2) to (0.2,-0.6);
            \draw[->] (-0.8,0) arc(-180:180:0.2);
            \node at (-0.6,0) {\dotlabel{+}};
            \multdot{-0.4,0}{west}{-r-s};
            \multdot{-0.2,0.42}{east}{r};
            \multdot{-0.2,-0.42}{east}{s};
            \teleport{-0.18,0.25}{-0.459,0.141};
            \teleport{-0.18,-0.25}{-0.459,-0.141};
        \end{tikzpicture}
        \ ,
        \\ \label{neg}
        \begin{tikzpicture}[centerzero]
            \draw[->] (-0.2,0) \braidup (0.2,0.5);
            \draw[wipe] (-0.2,-0.5) \braidup (0.2,0) \braidup (-0.2,0.5);
            \draw[<-] (-0.2,-0.5) \braidup (0.2,0) \braidup (-0.2,0.5);
            \draw[wipe] (0.2,-0.5) \braidup (-0.2,0);
            \draw (0.2,-0.5) \braidup (-0.2,0);
        \end{tikzpicture}
        \ =\
        \begin{tikzpicture}[centerzero]
            \draw[<-] (-0.2,-0.5) to (-0.2,0.5);
            \draw[->] (0.2,-0.5) to (0.2,0.5);
        \end{tikzpicture}
        \ + t\
        \begin{tikzpicture}[centerzero]
            \draw[->] (-0.2,0.5) -- (-0.2,0.3) arc(180:360:0.2) -- (0.2,0.5);
            \draw[<-] (-0.2,-0.5) -- (-0.2,-0.3) arc(180:0:0.2) -- (0.2,-0.5);
            \teleport{0.2,-0.3}{0.2,0.3};
        \end{tikzpicture}
        \ + \sum_{r,s > 0}\
        \begin{tikzpicture}[anchorbase]
            \draw[->] (-0.2,0.6) to (-0.2,0.35) arc(180:360:0.2) to (0.2,0.6);
            \draw[<-] (-0.2,-0.6) to (-0.2,-0.35) arc(180:0:0.2) to (0.2,-0.6);
            \draw[->] (0.8,0) arc(360:0:0.2);
            \node at (0.6,0) {\dotlabel{+}};
            \multdot{0.4,0}{east}{-r-s};
            \multdot{0.2,0.42}{west}{r};
            \multdot{0.2,-0.42}{west}{s};
            \teleport{0.18,0.25}{0.459,0.141};
            \teleport{0.18,-0.25}{0.459,-0.141};
        \end{tikzpicture}
        \ ,
        \\ \label{curls}
        \begin{tikzpicture}[centerzero]
            \draw (-0.2,-0.5) -- (-0.2,-0.35) to[out=up,in=west] (0.05,0.2) to[out=right,in=up] (0.2,0);
            \draw[wipe] (0.2,0) to[out=down,in=east] (0.05,-0.2) to[out=left,in=down] (-0.2,0.35) -- (-0.2,0.5);
            \draw[->] (0.2,0) to[out=down,in=east] (0.05,-0.2) to[out=left,in=down] (-0.2,0.35) -- (-0.2,0.5);
        \end{tikzpicture}
        \ = \delta_{k,0} t^{-1}\
        \begin{tikzpicture}[centerzero]
            \draw[->] (0,-0.5) to (0,0.5);
        \end{tikzpicture}
        \quad \text{if } k \ge 0,
        \qquad \qquad
        \begin{tikzpicture}[centerzero]
            \draw[->] (0.2,0) arc(360:0:0.2);
            \token{-0.141,0.141}{east}{a};
            \multdot{-0.141,-0.141}{east}{n};
        \end{tikzpicture}
        = \frac{\delta_{n,0} t - \delta_{n,k} t^{-1}}{z} \tr(a) 1_\one \quad \text{if } 0 \le n \le k,
        \\ \label{morecurls}
        \begin{tikzpicture}[centerzero]
            \draw (0.2,-0.5) -- (0.2,-0.35) to[out=up,in=east] (-0.05,0.2) to[out=left,in=up] (-0.2,0);
            \draw[wipe] (-0.2,0) to[out=down,in=west] (-0.05,-0.2) to[out=right,in=down] (0.2,0.35) -- (0.2,0.5);
            \draw[->] (-0.2,0) to[out=down,in=west] (-0.05,-0.2) to[out=right,in=down] (0.2,0.35) -- (0.2,0.5);
        \end{tikzpicture}
        \ = \delta_{k,0} t\
        \begin{tikzpicture}[centerzero]
            \draw[->] (0,-0.5) to (0,0.5);
        \end{tikzpicture}
        \quad \text{if } k \le 0,
        \qquad \qquad
        \begin{tikzpicture}[centerzero]
            \draw[->] (-0.2,0) arc(-180:180:0.2);
            \token{0.141,0.141}{west}{a};
            \multdot{0.141,-0.141}{west}{n};
        \end{tikzpicture}
        = \frac{\delta_{n,-k} t - \delta_{n,0} t^{-1}}{z} \tr(a) 1_\one \quad \text{if } 0 \le n \le -k.
    \end{gather}
    Here we have used the leftward crossings, which are defined in this approach by
    \begin{equation} \label{cold}
        \begin{tikzpicture}[anchorbase]
            \draw[<-] (0,0) -- (0.6,0.6);
            \draw[wipe] (0.6,0) -- (0,0.6);
            \draw[->] (0.6,0) -- (0,0.6);
        \end{tikzpicture}
        \ :=\
        \begin{tikzpicture}[anchorbase,scale=0.6]
            \draw[<-] (-0.75,-1) -- (-0.75,-0.5) .. controls (-0.75,-0.2) and (-0.5,0) .. (0,-0.5) .. controls (0.5,-1) and (0.75,-0.8) .. (0.75,-0.5) -- (0.75,0);
            \draw[wipe] (0.3,0) -- (-0.3,-1);
            \draw[<-] (0.3,0) -- (-0.3,-1);
        \end{tikzpicture}
        \qquad , \qquad
        \begin{tikzpicture}[anchorbase]
            \draw[->] (0.6,0) -- (0,0.6);
            \draw[wipe] (0,0) -- (0.6,0.6);
            \draw[<-] (0,0) -- (0.6,0.6);
        \end{tikzpicture}
        \ :=\
        \begin{tikzpicture}[anchorbase,scale=0.6]
            \draw[<-] (0.3,0) -- (-0.3,-1);
            \draw[wipe] (-0.75,-1) -- (-0.75,-0.5) .. controls (-0.75,-0.2) and (-0.5,0) .. (0,-0.5) .. controls (0.5,-1) and (0.75,-0.8) .. (0.75,-0.5) -- (0.75,0);
            \draw[<-] (-0.75,-1) -- (-0.75,-0.5) .. controls (-0.75,-0.2) and (-0.5,0) .. (0,-0.5) .. controls (0.5,-1) and (0.75,-0.8) .. (0.75,-0.5) -- (0.75,0);
        \end{tikzpicture}
        \ ,
    \end{equation}
    and the $(+)$-bubbles, which are defined in this approach by \cref{d1,d2} when they carry a token labelled $a \in A$ and a dot labeled $n \le 0$.  Finally, we define the $(+)$-bubbles with $n > 0$ dots to be the usual bubbles with $n$ dots as in \cref{fake2}, and then define the $(-)$-bubbles, with an arbitrary number of dots, so that \cref{fake1} holds.
\end{defin}

We will prove that \cref{def3} is equivalent to \cref{def1,def2}.  But before doing so, we make some remarks about the relations \cref{neg,pos,curls,morecurls,cold,d1,d2}.  If $k \le 1$, it follows immediately from \cref{d1} that relation \cref{pos} is equivalent to
\begin{equation} \label{posalt}
    \begin{tikzpicture}[centerzero]
        \draw[->] (0.2,0) \braidup (-0.2,0.5);
        \draw[<-] (0.2,-0.5) \braidup (-0.2,0);
        \draw[wipe] (-0.2,-0.5) \braidup (0.2,0);
        \draw (-0.2,-0.5) \braidup (0.2,0);
        \draw[wipe] (-0.2,0) \braidup (0.2,0.5);
        \draw (-0.2,0) \braidup (0.2,0.5);
    \end{tikzpicture}
    \ =\
    \begin{tikzpicture}[centerzero]
        \draw[->] (-0.2,-0.5) to (-0.2,0.5);
        \draw[<-] (0.2,-0.5) to (0.2,0.5);
    \end{tikzpicture}
    \ - t^{-1}\
    \begin{tikzpicture}[centerzero]
        \draw[<-] (-0.2,0.5) -- (-0.2,0.3) arc(180:360:0.2) -- (0.2,0.5);
        \draw[->] (-0.2,-0.5) -- (-0.2,-0.3) arc(180:0:0.2) -- (0.2,-0.5);
        \teleport{-0.2,-0.3}{-0.2,0.3};
    \end{tikzpicture}
    \quad \text{if } k \le 1.
\end{equation}
Similarly, when $k \ge -1$, relation \cref{neg} is equivalent to
\begin{equation} \label{negalt}
    \begin{tikzpicture}[centerzero]
        \draw[->] (-0.2,0) \braidup (0.2,0.5);
        \draw[wipe] (-0.2,-0.5) \braidup (0.2,0) \braidup (-0.2,0.5);
        \draw[<-] (-0.2,-0.5) \braidup (0.2,0) \braidup (-0.2,0.5);
        \draw[wipe] (0.2,-0.5) \braidup (-0.2,0);
        \draw (0.2,-0.5) \braidup (-0.2,0);
    \end{tikzpicture}
    \ =\
    \begin{tikzpicture}[centerzero]
        \draw[<-] (-0.2,-0.5) to (-0.2,0.5);
        \draw[->] (0.2,-0.5) to (0.2,0.5);
    \end{tikzpicture}
    \ + t\
    \begin{tikzpicture}[centerzero]
        \draw[->] (-0.2,0.5) -- (-0.2,0.3) arc(180:360:0.2) -- (0.2,0.5);
        \draw[<-] (-0.2,-0.5) -- (-0.2,-0.3) arc(180:0:0.2) -- (0.2,-0.5);
        \teleport{0.2,-0.3}{0.2,0.3};
    \end{tikzpicture}
    \quad \text{if } k \ge -1.
\end{equation}
In addition, using \cref{rcross,rskein}, we obtain the following relations from \cref{curls,morecurls}:
\begin{align} \label{swirl+}
    \begin{tikzpicture}[centerzero={(0,0.1)}]
        \draw (-0.2,0.4) to[out=-45,in=up] (0.2,0) arc(360:180:0.2);
        \draw[wipe] (-0.2,0) to[out=up,in=225] (0.2,0.4);
        \draw[->] (-0.2,0) to[out=up,in=225] (0.2,0.4);
    \end{tikzpicture}
    &= \delta_{k,0} t^{-1}\
    \begin{tikzpicture}[centerzero]
        \draw[->] (-0.2,0.2) -- (-0.2,0) arc (180:360:0.2) -- (0.2,0.2);
    \end{tikzpicture}
    \ ,&
    \begin{tikzpicture}[centerzero={(0,0.1)}]
        \draw[->] (-0.2,0) to[out=up,in=225] (0.2,0.4);
        \draw[wipe] (-0.2,0.4) to[out=-45,in=up] (0.2,0) arc(360:180:0.2);
        \draw (-0.2,0.4) to[out=-45,in=up] (0.2,0) arc(360:180:0.2);
    \end{tikzpicture}
    &= t\
    \begin{tikzpicture}[centerzero]
        \draw[->] (-0.2,0.2) -- (-0.2,0) arc (180:360:0.2) -- (0.2,0.2);
    \end{tikzpicture}
    \ ,&
    \begin{tikzpicture}[centerzero={(0,-0.1)}]
        \draw (-0.2,0) to[out=down,in=135] (0.2,-0.4);
        \draw[wipe] (-0.2,-0.4) to[out=45,in=down] (0.2,0) arc(0:180:0.2);
        \draw[<-] (-0.2,-0.4) to[out=45,in=down] (0.2,0) arc(0:180:0.2);
    \end{tikzpicture}
    &= \delta_{k,0} t^{-1}
    \begin{tikzpicture}[centerzero]
        \draw[<-] (-0.2,-0.2) -- (-0.2,0) arc (180:0:0.2) -- (0.2,-0.2);
    \end{tikzpicture}
    \ ,&
    \begin{tikzpicture}[centerzero={(0,-0.1)}]
        \draw[<-] (-0.2,-0.4) to[out=45,in=down] (0.2,0) arc(0:180:0.2);
        \draw[wipe] (-0.2,0) to[out=down,in=135] (0.2,-0.4);
        \draw (-0.2,0) to[out=down,in=135] (0.2,-0.4);
    \end{tikzpicture}
    &= t\
    \begin{tikzpicture}[centerzero]
        \draw[<-] (-0.2,-0.2) -- (-0.2,0) arc (180:0:0.2) -- (0.2,-0.2);
    \end{tikzpicture}
    \, &\text{if } k \ge 0,
    \\ \label{swirl-}
    \begin{tikzpicture}[centerzero={(0,-0.1)}]
        \draw (-0.2,-0.4) to[out=45,in=down] (0.2,0) arc(0:180:0.2);
        \draw[wipe] (-0.2,0) to[out=down,in=135] (0.2,-0.4);
        \draw[->] (-0.2,0) to[out=down,in=135] (0.2,-0.4);
    \end{tikzpicture}
    &= \delta_{k,0} t\
    \begin{tikzpicture}[centerzero]
        \draw[->] (-0.2,-0.2) -- (-0.2,0) arc (180:0:0.2) -- (0.2,-0.2);
    \end{tikzpicture}
    \ ,&
    \begin{tikzpicture}[centerzero={(0,-0.1)}]
        \draw[->] (-0.2,0) to[out=down,in=135] (0.2,-0.4);
        \draw[wipe] (-0.2,-0.4) to[out=45,in=down] (0.2,0) arc(0:180:0.2);
        \draw (-0.2,-0.4) to[out=45,in=down] (0.2,0) arc(0:180:0.2);
    \end{tikzpicture}
    &= t^{-1}\
    \begin{tikzpicture}[centerzero]
        \draw[->] (-0.2,-0.2) -- (-0.2,0) arc (180:0:0.2) -- (0.2,-0.2);
    \end{tikzpicture}
    \ ,&
    \begin{tikzpicture}[centerzero={(0,0.1)}]
        \draw (-0.2,0) to[out=up,in=225] (0.2,0.4);
        \draw[wipe] (-0.2,0.4) to[out=-45,in=up] (0.2,0) arc(360:180:0.2);
        \draw[<-] (-0.2,0.4) to[out=-45,in=up] (0.2,0) arc(360:180:0.2);
    \end{tikzpicture}
    &= \delta_{k,0} t\
    \begin{tikzpicture}[centerzero]
        \draw[<-] (-0.2,0.2) -- (-0.2,0) arc (180:360:0.2) -- (0.2,0.2);
    \end{tikzpicture}
    \ ,&
    \begin{tikzpicture}[centerzero={(0,0.1)}]
        \draw[<-] (-0.2,0.4) to[out=-45,in=up] (0.2,0) arc(360:180:0.2);
        \draw[wipe] (-0.2,0) to[out=up,in=225] (0.2,0.4);
        \draw (-0.2,0) to[out=up,in=225] (0.2,0.4);
    \end{tikzpicture}
    &= t^{-1}
    \begin{tikzpicture}[centerzero]
        \draw[<-] (-0.2,0.2) -- (-0.2,0) arc (180:360:0.2) -- (0.2,0.2);
    \end{tikzpicture}
    \ ,& \text{if } k \le 0.
\end{align}
Then, using \cref{swirl-} and \cref{lskein} to convert the negative crossings in \cref{posalt} to positive ones, we see that, when $k < 0$, \cref{posalt} is equivalent to
\begin{equation} \label{posalter}
    \begin{tikzpicture}[centerzero]
        \draw (-0.2,0) \braidup (0.2,0.5);
        \draw[wipe] (-0.2,-0.5) \braidup (0.2,0) \braidup (-0.2,0.5);
        \draw[->] (-0.2,-0.5) \braidup (0.2,0) \braidup (-0.2,0.5);
        \draw[wipe] (0.2,-0.5) \braidup (-0.2,0);
        \draw[<-] (0.2,-0.5) \braidup (-0.2,0);
    \end{tikzpicture}
    \ =\
    \begin{tikzpicture}[anchorbase]
        \draw[->] (-0.2,-0.5) to (-0.2,0.5);
        \draw[<-] (0.2,-0.5) to (0.2,0.5);
    \end{tikzpicture}
    \quad \text{if } k < 0.
\end{equation}
Similarly, when $k > 0$, \cref{negalt} is equivalent to
\begin{equation} \label{negalter}
    \begin{tikzpicture}[centerzero]
        \draw (0.2,0) \braidup (-0.2,0.5);
        \draw (0.2,-0.5) \braidup (-0.2,0);
        \draw[wipe] (-0.2,-0.5) \braidup (0.2,0);
        \draw[<-] (-0.2,-0.5) \braidup (0.2,0);
        \draw[wipe] (-0.2,0) \braidup (0.2,0.5);
        \draw[->] (-0.2,0) \braidup (0.2,0.5);
    \end{tikzpicture}
    \ =\
    \begin{tikzpicture}[anchorbase]
        \draw[<-] (-0.2,-0.5) to (-0.2,0.5);
        \draw[->] (0.2,-0.5) to (0.2,0.5);
    \end{tikzpicture}
    \quad \text{if } k > 0.
\end{equation}
Finally, when $k = 0$, \cref{posalt,negalt} together are equivalent to the single assertion
\begin{equation} \label{bothalt}
    \begin{tikzpicture}[anchorbase]
        \draw[->] (0.2,-0.2) -- (-0.2,0.2);
        \draw[wipe] (-0.2,-0.2) -- (0.2,0.2);
        \draw[<-] (-0.2,-0.2) -- (0.2,0.2);
    \end{tikzpicture}
    \ =\
    \left(
        \begin{tikzpicture}[anchorbase]
            \draw[->] (-0.2,-0.2) -- (0.2,0.2);
            \draw[wipe] (0.2,-0.2) -- (-0.2,0.2);
            \draw[<-] (0.2,-0.2) -- (-0.2,0.2);
        \end{tikzpicture}
    \right)^{-1},
\end{equation}
i.e.\ both of the relations \cref{lunch}.

\begin{theo}
    The category $\Heis_k(A;z,t)$ of \cref{def3} is the same as the one of \cref{def1,def2}, with all morphisms introduced in \cref{def3} being the same as the ones of \cref{def1,def2}.
\end{theo}

\begin{proof}
    To avoid confusion in the proof, we denote the category from the equivalent \cref{def1,def2} by $\Heis_k^\old(A;z,t)$, and the one from \cref{def3} by $\Heis_k^\new(A;z,t)$.  It is clear from the symmetry in the relations \cref{pos,neg,curls,morecurls,cold,d1,d2} that there is an isomorphism
    \[
        \Heis_k^\new(A;z,t) \to \Heis_{-k}^\new(A;z,t^{-1})^\op
    \]
    that reflects diagrams in a horizontal plane and multiplies by $(-1)^{c+l+\binom{y}{2}}$, where $c$ is the number of crossings, $d$ is the number of left cups and caps, and $y$ is the number of odd tokens. Thus, by \cref{om}, it suffices to prove the theorem in the case $k \le 0$.

    We first check that the relations \cref{pos,neg,curls,morecurls,cold} are satisfied in $\Heis_k^\old(A;z,t)$, so that we have a natural strict $\kk$-linear monoidal functor
    \[
        \Theta \colon \Heis_k^\new(A;z,t) \to \Heis_k^\old(A;z,t),
    \]
    which is the identity on diagrams.
    \begin{itemize}[wide]
        \item \cref{pos,neg}:  When $k=0$, \cref{pos,neg} are equivalent to \cref{bothalt}, which holds by \cref{lunch}.  Now suppose $k < 0$.  Then \cref{pos} is equivalent to \cref{posalter}, which holds by the second relation in \cref{tea3}.  To check \cref{neg}, we start with the first relation in \cref{tea3} and expand the $\color{violet}\scriptstyle{\heartsuit}$-decorated left caps using \cref{nakano7,leftwards} when $r=0$, or \cref{nakano3,nakano5} when $r>0$.

        \item \cref{curls,morecurls}:  These relations follow easily from \cref{tea2,impose,leftwards,colder,coldest}.

        \item \cref{cold}: This holds in $\Heis_k^\old(A;z,t)$ since we have shown this category is strictly pivotal.
    \end{itemize}

    Now we want to show that $\Theta$ is an isomorphism.  We do this by using the presentation from \cref{def1} to construct a two-sided inverse
    \[
        \Phi \colon \Heis_k^\old(A;z,t) \to \Heis_k^\new(A;z,t),
    \]
    still assuming $k \le 0$.  We define $\Phi$ on morphisms be declaring that it takes the rightwards cup, the rightwards cap, and all the tokens, dots, and crossings (with any orientation) to the corresponding morphisms in $\Heis_k^\new(A;z,t)$, and also (see \cref{nakano5,nakano2,leftwards,nakano7})
    \begin{gather*}
        \Phi
        \left(
            \begin{tikzpicture}[anchorbase]
                \draw[<-] (-0.2,-0.2) -- (-0.2,0) arc (180:0:0.2) -- (0.2,-0.2);
                \heartdec{0,0.2}{south}{a}{0};
            \end{tikzpicture}
        \right)
        = -(-1)^{\bar a} tz\
        \begin{tikzpicture}[centerzero]
            \draw[<-] (-0.2,-0.2) -- (-0.2,0) arc (180:0:0.2) -- (0.2,-0.2);
            \token{0.2,0}{west}{a};
        \end{tikzpicture}
        \quad \text{if } k < 0,
        \\
        \Phi
        \left(
            \begin{tikzpicture}[anchorbase]
                \draw[<-] (-0.2,-0.2) -- (-0.2,0) arc (180:0:0.2) -- (0.2,-0.2);
                \heartdec{0,0.2}{south}{a}{n};
            \end{tikzpicture}
        \right)
        = - (-1)^{\bar{a}} z^2 \sum_{r \ge 1}
        \begin{tikzpicture}[centerzero]
            \draw[<-] (-0.2,-0.45) -- (-0.2,-0.25) arc(180:0:0.2) -- (0.2,-0.45);
            \token{-0.2,-0.25}{east}{b^\vee};
            \multdot{0.2,-0.25}{west}{r};
            \plusright{0,0.25}{ab}{-r-n};
        \end{tikzpicture}
        \quad \text{if } 0 < n < -k.
    \end{gather*}
    To see that $\Phi$ is well defined, we must verify the relations from \cref{def1}.  For \cref{impose}, this amounts to checking that, in $\Heis_k^\new(A;z,t)$, we have
    \begin{align*}
        t\
        \begin{tikzpicture}[anchorbase]
           \draw (-0.2,0.2) to[out=down,in=up] (0.2,-0.2) to[out=down,in=down,looseness=1.5] (-0.2,-0.2);
            \draw[wipe] (-0.2,-0.2) to[out=up,in=down] (0.2,0.2) to[out=up,in=up,looseness=1.5] (-0.2,0.2);
            \draw[<-] (-0.2,-0.2) to[out=up,in=down] (0.2,0.2) to[out=up,in=up,looseness=1.5] (-0.2,0.2);
            \token{0.2,0.2}{west}{a};
        \end{tikzpicture}
        &= \frac{t-t^{-1}}{z} \tr(a) 1_\one
        \text{ if } k = 0,
        &
        \ccbubble{a}{-k}
        &= \frac{t}{z} \tr(a) 1_\one \text{ if } k < 0.
    \end{align*}
    The first relation follows from \cref{swirl-,curls}, while the second follows from \cref{morecurls}.  Thus, it remains to show that the images under $\Phi$ of the morphisms \cref{invrel,invrel2} are two-sided inverses in $\Heis_k^\new(A;z,t)$.  When $k = 0$, this is immediate from \cref{bothalt}, so suppose $k < 0$.  The images under $\Phi$ of the two equations in \cref{tea3} are precisely \cref{neg} and \cref{posalter}.  We are left with checking that the images under $\Phi$ of the relations
    \begin{align} \label{invcheck}
        \begin{tikzpicture}[centerzero={(0,0.1)}]
            \draw (-0.2,0) to[out=up,in=225] (0.2,0.4);
            \draw[wipe] (-0.2,0.4) to[out=-45,in=up] (0.2,0) arc(360:180:0.2);
            \draw[<-] (-0.2,0.4) to[out=-45,in=up] (0.2,0) arc(360:180:0.2);
            \token{-0.2,0}{east}{a};
            \multdot{0.2,0}{west}{m};
        \end{tikzpicture}
        &=0,&
        \begin{tikzpicture}[centerzero={(0,-0.1)}]
            \draw (-0.2,-0.4) to[out=45,in=down] (0.2,0) arc(0:180:0.2);
            \draw[wipe] (-0.2,0) to[out=down,in=135] (0.2,-0.4);
            \draw[->] (-0.2,0) to[out=down,in=135] (0.2,-0.4);
            \heartdec{0,0.2}{south}{a}{n};
        \end{tikzpicture}
        &=0,&
        \begin{tikzpicture}[centerzero]
            \draw[->] (0,-0.2) arc(-90:270:0.2);
            \token{-0.2,0}{east}{a};
            \multdot{0.2,0}{west}{m};
            \heartdec{0,0.2}{south}{b}{n};
        \end{tikzpicture}
        &= \delta_{m,n} \tr(ab) 1_\one,
     \end{align}
    hold in $\Heis_k^\new(A;z,t)$ for all $a,b \in A$ and $0 \le m,n < -k$.  For the first relation, we can slide the token through the crossing, and so it suffices to prove it without the token.  When $m=0$, it then follows by \cref{swirl-}.  When $0 < m < -k$, we have
    \[
        \begin{tikzpicture}[centerzero={(0,0.1)}]
            \draw (-0.2,0) to[out=up,in=225] (0.2,0.4);
            \draw[wipe] (-0.2,0.4) to[out=-45,in=up] (0.2,0) arc(360:180:0.2);
            \draw[<-] (-0.2,0.4) to[out=-45,in=up] (0.2,0) arc(360:180:0.2);
            \multdot{0.2,0}{west}{m};
        \end{tikzpicture}
        \stackrel[\cref{morecurls}]{\cref{lskein}}{=}
        \begin{tikzpicture}[centerzero={(0,0.1)}]
            \draw[<-] (-0.2,0.4) to[out=-45,in=up] (0.2,0) arc(360:180:0.2);
            \draw[wipe] (-0.2,0) to[out=up,in=225] (0.2,0.4);
            \draw (-0.2,0) to[out=up,in=225] (0.2,0.4);
            \multdot{0.2,0}{west}{m};
        \end{tikzpicture}
        \overset{\substack{\cref{teaneg} \\ \cref{cold}}}{\underset{\cref{morecurls}}{=}}
        \begin{tikzpicture}[anchorbase]
            \draw (-0.2,-0.2) \braidup (0.2,0.2) -- (0.2,0.4);
            \draw[wipe] (-0.2,0.4) -- (-0.2,0.2) \braiddown (0.2,-0.2) arc(360:180:0.2);
            \draw[<-] (-0.2,0.4) -- (-0.2,0.2) \braiddown (0.2,-0.2) arc(360:180:0.2);
            \multdot{-0.2,0.2}{east}{m};
        \end{tikzpicture}
        \ \stackrel{\cref{swirl-}}{=} 0.
    \]
    When $n=0$, the second and third relations in \cref{invcheck} follow from \cref{swirl-,morecurls}.  To prove them when $0 < n < -k$, we must show that
    \begin{align} \label{invcheck2}
        \sum_{r \ge 1}
        \begin{tikzpicture}[anchorbase]
            \draw (-0.2,-0.4) to[out=45,in=down] (0.2,0) arc(0:180:0.2);
            \draw[wipe] (-0.2,0) to[out=down,in=135] (0.2,-0.4);
            \draw[->] (-0.2,0) to[out=down,in=135] (0.2,-0.4);
            \token{-0.2,0}{east}{b^\vee};
            \multdot{0.2,0}{west}{r};
            \plusright{0,0.5}{ab}{-r-n};
        \end{tikzpicture}
        &= 0,&
        \sum_{r \ge 1}
        \begin{tikzpicture}[centerzero]
            \plusright{0,0.3}{bc}{-r-n};
            \bubleft{0,-0.3}{c^\vee a}{r+m};
        \end{tikzpicture}
        &= - (-1)^{\bar b} \delta_{m,n} \frac{1}{z^2} \delta_{m,n} \tr(ab) 1_\one
    \end{align}
    in $\Heis_k^\new(A;z,t)$.  For the first identity, note that the terms with $r \ge -k$ vanish since the $(+)$-bubble is zero by \cref{d1}.  The terms with $0 < r < -k$ also vanish, as can be seen by using the skein relation to flip the crossing, sliding the dots past the crossing, and then using \cref{morecurls,swirl-}.  To prove the second identity in \cref{invcheck}, we first note that, by \cref{beam}, it suffices to consider the case where $b=1$.  When $m \le n$, the identity follows from \cref{d2,morecurls}.  When $m > n$, we compute:
    \[
        \pluscbubble{a}{m-n+k}
        \ \stackrel{\mathclap{\cref{d2}}}{=}\
        - z^{m-n-1} t^{n-m-1} \sum_{b_1,\dotsc,b_{m-n-1} \in \BA} \det
        \left( - \ccbubble{b^\vee_{j-1}b_j}{i-j-k+1} \right)_{i,j=1}^{m-n}
        = -zt^{-1} \sum_{s=1}^{m-n}
        \begin{tikzpicture}[centerzero]
            \bubleft{0,0.3}{ab^\vee}{s-k};
            \plusright{0,-0.3}{b}{m-n+k};
        \end{tikzpicture}
        \ ,
    \]
    where, in the final equality, we expanded the determinant along the first column as in proof of \cref{dlem}.  By \cref{d2}, it follows that
    \[
        0
        = \sum_{s=0}^{m-n}
        \begin{tikzpicture}[centerzero]
            \bubleft{0,0.3}{ab}{s-k};
            \plusright{0,-0.3}{b^\vee}{m-n+k};
        \end{tikzpicture}
        = \sum_{s=0}^{m-n} (-1)^{\bar b}
        \begin{tikzpicture}[centerzero]
            \bubleft{0,0.3}{ba}{s-k};
            \plusright{0,-0.3}{b^\vee}{m-n+k};
        \end{tikzpicture}
        = \sum_{s=0}^{m-n}
        \begin{tikzpicture}[centerzero]
            \plusright{0,0.3}{b^\vee}{m-n+k};
            \bubleft{0,-0.3}{b^\vee a}{s-k};
        \end{tikzpicture}
        \ .
    \]
    Then the second identity in \cref{invcheck2} follows by letting $r=s-m-k$.

    To complete the proof, we must show that $\Theta$ and $\Phi$ are indeed two-sided inverses.  To verify that $\Theta \circ \Phi = \id$, the only nontrivial thing to check is that
    \[
        \Theta
        \left(
            \Phi
            \left(
                \begin{tikzpicture}[centerzero={(0,0.15)}]
                    \draw[<-] (-0.2,-0.2) -- (-0.2,0) arc (180:0:0.2) -- (0.2,-0.2);
                    \heartdec{0,0.2}{south}{a}{n};
                \end{tikzpicture}
            \right)
        \right)
        =
        \begin{tikzpicture}[centerzero={(0,0.15)}]
            \draw[<-] (-0.2,-0.2) -- (-0.2,0) arc (180:0:0.2) -- (0.2,-0.2);
            \heartdec{0,0.2}{south}{a}{n};
        \end{tikzpicture}
        \quad \text{for } a \in A,\ 0 \le n < -k.
    \]
    When $n=0$, this follows immediately from \cref{leftwards,nakano7}.  On the other hand, if $0 < n < -k$, it follows from \cref{nakano3,nakano5}.  To verify that $\Phi \circ \Theta = \id$, the only nontrivial thing to check is that
    \[
        \Phi
        \left(
            \begin{tikzpicture}[anchorbase]
                \draw[<-] (-0.2,-0.2) -- (-0.2,0) arc (180:0:0.2) -- (0.2,-0.2);
            \end{tikzpicture}\,
        \right)
        \ =\
        \begin{tikzpicture}[anchorbase]
            \draw[<-] (-0.2,-0.2) -- (-0.2,0) arc (180:0:0.2) -- (0.2,-0.2);
        \end{tikzpicture}
        \quad \text{and} \quad
        \Phi
        \left(
            \begin{tikzpicture}[anchorbase]
                \draw[<-] (-0.2,0.2) -- (-0.2,0) arc (180:360:0.2) -- (0.2,0.2);
            \end{tikzpicture}\,
        \right)
        \ =\
        \begin{tikzpicture}[anchorbase]
            \draw[<-] (-0.2,0.2) -- (-0.2,0) arc (180:360:0.2) -- (0.2,0.2);
        \end{tikzpicture}
        \ .
    \]
    The first equality follows from \cref{leftwards,swirl+}.  The second follows from \cref{leftwards} when $k=0$.  When $k < 0$, we have
    \[
        \Phi
        \left(
            \begin{tikzpicture}[anchorbase]
                \draw[<-] (-0.2,0.2) -- (-0.2,0) arc (180:360:0.2) -- (0.2,0.2);
            \end{tikzpicture}\,
        \right)
        \ \stackrel{\cref{leftwards}}{=} t^{-1}\
        \begin{tikzpicture}[centerzero={(0,0.1)}]
            \draw (-0.2,0) to[out=up,in=225] (0.2,0.4);
            \draw[wipe] (-0.2,0.4) to[out=-45,in=up] (0.2,0) arc(360:180:0.2);
            \draw[<-] (-0.2,0.4) to[out=-45,in=up] (0.2,0) arc(360:180:0.2);
            \multdot{0.2,0}{west}{-k};
        \end{tikzpicture}
        \ \stackrel[\substack{\cref{morecurls} \\ \cref{adecomp}}]{\cref{lskein}}{=}\
        \begin{tikzpicture}[centerzero={(0,0.1)}]
            \draw[<-] (-0.2,0.4) to[out=-45,in=up] (0.2,0) arc(360:180:0.2);
            \draw[wipe] (-0.2,0) to[out=up,in=225] (0.2,0.4);
            \draw (-0.2,0) to[out=up,in=225] (0.2,0.4);
            \multdot{0.2,0}{west}{-k};
        \end{tikzpicture}
        \ +\
        \begin{tikzpicture}[anchorbase]
            \draw[<-] (-0.2,0.2) -- (-0.2,0) arc (180:360:0.2) -- (0.2,0.2);
        \end{tikzpicture}        \ \stackrel[\substack{\cref{morecurls} \\ \cref{swirl-}}]{\substack{\cref{teaneg} \\ \cref{cold}}}{=}\
        \begin{tikzpicture}[anchorbase]
            \draw[<-] (-0.2,0.2) -- (-0.2,0) arc (180:360:0.2) -- (0.2,0.2);
        \end{tikzpicture}
        \ . \qedhere
    \]
\end{proof}

\begin{rem}
    The Frobeneius Heisenberg categories studied in \cite{Sav19,BSW-foundations} depend only, up to isomorphism, on the underlying algebra $A$, and not on the trace map; see \cite[Lem.~5.3]{BSW-foundations}.  However, in the quantum setting of the current paper, there do not seem to be obvious isomorphisms between quantum Frobenius Heisenberg categories corresponding to the same algebra but with different trace maps.  This is true even for the quantum affine wreath product algebras; see \cite[Rem.~2.2]{RS20}.
\end{rem}

\begin{lem} \label{cherries}
    Suppose that $\cC$ is a strict $\kk$-linear monoidal category containing objects $\uparrow$ and $\downarrow$ and morphisms $\tokup$, $\dotup$, $\poscross$, $\negcross$, $\rightcup$, and $\rightcap$ satisfying \cref{tokrel,braid,skein,QAWC}.  Then $\cC$ contains at most one pair of morphisms $\leftcup$ and $\leftcap$ satisfying \cref{pos,neg,curls,morecurls} (for sideways crossings and the $(+)$-bubbles defined via \cref{rcross,cold,d1,d2}).
\end{lem}

\begin{proof}
    The proof is analogous to that of \cite[Lem.~4.3]{BSW-qheis}.
\end{proof}


We now prove some further useful relations that hold in $\Heis_k(A;z,t)$.  We will state some of these relations with and without the language of generating functions.  To state the relations in terms of generating functions, we need to unify our notation for dots and tokens, adopting the notation
\[
    \begin{tikzpicture}[anchorbase]
        \draw[->] (0,-0.3) to (0,0.3);
        \token{0,0}{east}{a x^n};
    \end{tikzpicture}
    :=
    \begin{tikzpicture}[anchorbase]
        \draw[->] (0,-0.3) to (0,0.3);
        \token{0,0.1}{east}{a};
        \multdot{0,-0.1}{east}{n};
    \end{tikzpicture}
    \qquad \text{and} \qquad
    \begin{tikzpicture}[anchorbase]
        \draw[<-] (0,-0.3) to (0,0.3);
        \token{0,0}{east}{a x^n};
    \end{tikzpicture}
    :=
    \begin{tikzpicture}[anchorbase]
        \draw[<-] (0,-0.3) to (0,0.3);
        \token{0,0.1}{east}{a};
        \multdot{0,-0.1}{east}{n};
    \end{tikzpicture}
\]
for $a \in A$, $n \in \Z$.  Then we can also label tokens by polynomials $a_n x^n + \dotsb + a_1 x + a_0 \in A[x]$, meaning the sum of morphisms defined by the tokens labelled $a_n x^n, a_1 x,\dotsc, x_0$, or even by Laurent series in $A[x]\Laurent{w^{-1}}$ or $A[x]\Laurent{w}$.  For example, expanding in $A[x]\Laurent{w^{-1}}$, we have
\[
    \begin{tikzpicture}[anchorbase]
        \draw[->] (0,-0.3) to (0,0.3);
        \token{0,0}{east}{wx(w-x)^{-2}};
    \end{tikzpicture}
    = w^{-1}\,
    \begin{tikzpicture}[anchorbase]
        \draw[->] (0,-0.3) to (0,0.3);
        \singdot{0,0};
    \end{tikzpicture}
    + 2w^{-2}\,
    \begin{tikzpicture}[anchorbase]
        \draw[->] (0,-0.3) to (0,0.3);
        \multdot{0,0}{west}{2};
    \end{tikzpicture}
    + 3w^{-3}\,
    \begin{tikzpicture}[anchorbase]
        \draw[->] (0,-0.3) to (0,0.3);
        \multdot{0,0}{west}{3};
    \end{tikzpicture}
    + 4w^{-4}\,
    \begin{tikzpicture}[anchorbase]
        \draw[->] (0,-0.3) to (0,0.3);
        \multdot{0,0}{west}{4};
    \end{tikzpicture}
    + \dotsb.
\]
The caveat here is that, whereas tokens labelled by elements of $A\Laurent{w^{-1}}$ slide through crossings and can be teleported, the more general tokens labelled by elements of $A[x]\Laurent{w^{-1}}$ no longer have these properties in general.  For a Laurent series $p(w)$, we let $[p(w)]_{w^r}$ denote its $w^r$-coefficient, and we write $[p(w)]_{w^{<0}}$ for $\sum_{n<0} [p(w)]_{w^n} u^n$.  We adopt the convention that, in any equation involving the generating functions \cref{bubgen4,bubgen3}, we expand all rational functions as Laurent series in $A[x]\Laurent{w}$.  In all other equations, we expand rational functions as Laurent series in $A[x]\Laurent{w^{-1}}$.

\begin{lem}[Curl relations]
    The following relations hold for all $n \in \Z$:
    \begin{align} \label{dog1}
        \begin{tikzpicture}[centerzero]
            \draw (0.2,-0.5) -- (0.2,-0.35) to[out=up,in=east] (-0.05,0.2) to[out=left,in=up] (-0.2,0);
            \draw[wipe] (-0.2,0) to[out=down,in=west] (-0.05,-0.2) to[out=right,in=down] (0.2,0.35) -- (0.2,0.5);
            \draw[->] (-0.2,0) to[out=down,in=west] (-0.05,-0.2) to[out=right,in=down] (0.2,0.35) -- (0.2,0.5);
            \multdot{-0.2,0}{east}{n};
        \end{tikzpicture}
        &= \sum_{r \ge 0}
        \begin{tikzpicture}[centerzero]
            \draw[->] (0,-0.5) -- (0,0.5);
            \multdot{0,0.25}{west}{r};
            \plusleftblank{-0.5,0};
            \multdot{-0.7,0}{east}{n-r};
            \teleport{-0.3,0}{0,0};
        \end{tikzpicture}
        - \sum_{r > 0}
        \begin{tikzpicture}[centerzero]
            \draw[->] (0,-0.5) -- (0,0.5);
            \multdot{0,0.25}{west}{-r};
            \minusleftblank{-0.5,0};
            \multdot{-0.7,0}{east}{n+r};
            \teleport{-0.3,0}{0,0};
        \end{tikzpicture}
        \ ,&
        \begin{tikzpicture}[centerzero]
            \draw (0.2,-0.5) -- (0.2,-0.35) to[out=up,in=east] (-0.05,0.2) to[out=left,in=up] (-0.2,0);
            \draw[wipe] (-0.2,0) to[out=down,in=west] (-0.05,-0.2) to[out=right,in=down] (0.2,0.35) -- (0.2,0.5);
            \draw[->] (-0.2,0) to[out=down,in=west] (-0.05,-0.2) to[out=right,in=down] (0.2,0.35) -- (0.2,0.5);
            \multdot{-0.2,0}{east}{(w-x)^{-1}};
        \end{tikzpicture}
        &= t
        \left[
            \begin{tikzpicture}[centerzero]
                \draw[->] (0,-0.5) -- (0,0.5);
                \multdot{0,0.25}{west}{(w-x)^{-1}};
                \plusgenleft{-0.5,0};
                \teleport{-0.3,0}{0,0};
            \end{tikzpicture}
        \right]_{w^{<0}},
        \\ \label{dog2}
        \begin{tikzpicture}[centerzero]
            \draw[->] (-0.2,0) to[out=down,in=west] (-0.05,-0.2) to[out=right,in=down] (0.2,0.35) -- (0.2,0.5);
            \draw[wipe] (0.2,-0.5) -- (0.2,-0.35) to[out=up,in=east] (-0.05,0.2) to[out=left,in=up] (-0.2,0);
            \draw (0.2,-0.5) -- (0.2,-0.35) to[out=up,in=east] (-0.05,0.2) to[out=left,in=up] (-0.2,0);
            \multdot{-0.2,0}{east}{n};
        \end{tikzpicture}
        &= \sum_{r > 0}
        \begin{tikzpicture}[centerzero]
            \draw[->] (0,-0.5) -- (0,0.5);
            \multdot{0,0.25}{west}{r};
            \plusleftblank{-0.5,0};
            \multdot{-0.7,0}{east}{n-r};
            \teleport{-0.3,0}{0,0};
        \end{tikzpicture}
        - \sum_{r \ge 0}
        \begin{tikzpicture}[centerzero]
            \draw[->] (0,-0.5) -- (0,0.5);
            \multdot{0,0.25}{west}{-r};
            \minusleftblank{-0.5,0};
            \multdot{-0.7,0}{east}{n+r};
            \teleport{-0.3,0}{0,0};
        \end{tikzpicture}
        \ ,
        \\ \label{dog3}
        \begin{tikzpicture}[centerzero]
            \draw[->] (0.2,0) to[out=down,in=east] (0.05,-0.2) to[out=left,in=down] (-0.2,0.35) -- (-0.2,0.5);
            \draw[wipe] (-0.2,-0.5) -- (-0.2,-0.35) to[out=up,in=west] (0.05,0.2) to[out=right,in=up] (0.2,0);
            \draw (-0.2,-0.5) -- (-0.2,-0.35) to[out=up,in=west] (0.05,0.2) to[out=right,in=up] (0.2,0);
            \multdot{0.2,0}{west}{n};
        \end{tikzpicture}
        &= \sum_{r \ge 0}
        \begin{tikzpicture}[centerzero]
            \draw[->] (0,-0.5) -- (0,0.5);
            \multdot{0,0.25}{east}{-r};
            \teleport{0,0}{0.3,0};
            \minusrightblank{0.5,0};
            \multdot{0.7,0}{west}{n+r};
        \end{tikzpicture}
        - \sum_{r > 0}
        \begin{tikzpicture}[centerzero]
            \draw[->] (0,-0.5) -- (0,0.5);
            \multdot{0,0.25}{east}{r};
            \teleport{0,0}{0.3,0};
            \plusrightblank{0.5,0};
            \multdot{0.7,0}{west}{n-r};
        \end{tikzpicture}
        \ ,
        \\ \label{dog4}
        \begin{tikzpicture}[centerzero]
            \draw (-0.2,-0.5) -- (-0.2,-0.35) to[out=up,in=west] (0.05,0.2) to[out=right,in=up] (0.2,0);
            \draw[wipe] (0.2,0) to[out=down,in=east] (0.05,-0.2) to[out=left,in=down] (-0.2,0.35) -- (-0.2,0.5);
            \draw[->] (0.2,0) to[out=down,in=east] (0.05,-0.2) to[out=left,in=down] (-0.2,0.35) -- (-0.2,0.5);
            \multdot{0.2,0}{west}{n};
        \end{tikzpicture}
        &= \sum_{r > 0}
        \begin{tikzpicture}[centerzero]
            \draw[->] (0,-0.5) -- (0,0.5);
            \multdot{0,0.25}{east}{-r};
            \teleport{0,0}{0.3,0};
            \minusrightblank{0.5,0};
            \multdot{0.7,0}{west}{n+r};
        \end{tikzpicture}
        - \sum_{r \ge 0}
        \begin{tikzpicture}[centerzero]
            \draw[->] (0,-0.5) -- (0,0.5);
            \multdot{0,0.25}{east}{r};
            \teleport{0,0}{0.3,0};
            \plusrightblank{0.5,0};
            \multdot{0.7,0}{west}{n-r};
        \end{tikzpicture}
        \ ,&
        \begin{tikzpicture}[centerzero]
            \draw (-0.2,-0.5) -- (-0.2,-0.35) to[out=up,in=west] (0.05,0.2) to[out=right,in=up] (0.2,0);
            \draw[wipe] (0.2,0) to[out=down,in=east] (0.05,-0.2) to[out=left,in=down] (-0.2,0.35) -- (-0.2,0.5);
            \draw[->] (0.2,0) to[out=down,in=east] (0.05,-0.2) to[out=left,in=down] (-0.2,0.35) -- (-0.2,0.5);
            \multdot{0.2,0}{west}{(w-x)^{-1}};
        \end{tikzpicture}
       &= t^{-1}
        \left[
            \begin{tikzpicture}[centerzero]
                \draw[->] (0,-0.5) -- (0,0.5);
                \multdot{0,0.25}{east}{(w-x)^{-1}};
                \teleport{0,0}{0.3,0};
                \plusgenright{0.5,0};
            \end{tikzpicture}
        \right]_{w^{<0}}.
    \end{align}
\end{lem}

\begin{proof}
    The right-hand relations in \cref{dog1,dog4} are simply reformulations of the $n \ge 0$ cases of the left-hand relations in terms of generating functions.  Next note that \cref{dog3} follows easily from \cref{dog4,skein,fake1}.  Similarly, \cref{dog2} follows from \cref{dog1,skein,fake1}.  Thus, it suffices to prove \cref{dog1,dog4}.  Now, given one of \cref{dog1} or \cref{dog4} for $k \ge 0$, we can rotate through $180^\circ$ (using the strictly pivotal structure), and then apply $\Omega_k$ to obtain the other relation for $k \le 0$.  Thus, we are reduced to proving \cref{dog4} when $k \ge 0$ and \cref{dog1} when $k > 0$.

    Suppose $k \ge 0$.  When $n \ge 0$, we have
    \begin{multline*}

        \ .
    \end{multline*}
    This gives \cref{dog4}, since the $(-)$-bubbles there are zero when $n \ge 0$, by \cref{fake1,fake2}.  When $n < 0$, we have
    \[
        %
        \ .
    \]
    This gives \cref{dog4}, since the $(+)$-bubbles there are zero when $n<0$, by \cref{tanks,fake1}.

    Now suppose $k > 0$.  When $n=0$, we have
    \[
        %
        \ .
    \]
    Using the strictly pivotal structure, this gives \cref{dog1} since the $(-)$-bubbles appearing there are zero by \cref{fake1,fake2}.

    When $n > 0$, we have
    \begin{multline*}
        %
        \ ,
    \end{multline*}
    where, in the second-to-last equality we used the $n=0$ case of \cref{dog1}.  Thus we have proved \cref{dog1} when $n > 0$ since the $(-)$-bubbles appearing there are zero by \cref{fake1,fake2}.

    Finally, suppose $n<0$.  Then
    \begin{multline*}
        %
        \ ,
    \end{multline*}
    where, in the second equality, we used the $n=0$ case of \cref{dog2}.
\end{proof}

Recall the definition \cref{adag} of the notation $a^\dagger$.

\begin{lem}[Bubble slides]
    The following relations hold for all $n \in \Z$ and $a \in A$:
    \begin{align} \label{bs1}
        %
        \ .
    \end{align}
\end{lem}

\begin{proof}
    Note that, in each of \cref{bs1,bs2,bs3,bs4}, the right-hand relation is simply a restatement of the left-hand relation in terms of generating functions.  It suffices to prove \cref{bs1,bs2,bs3,bs4} for $k \ge 0$, since the identities for $k < 0$ then follow by rotating the pictures $180^\circ$ and applying $\Omega_k$.  So we assume $k \ge 0$.

    The identity \cref{bs1} is trivial when $n \le k$, by \cref{tanks}.  Thus, we suppose $n > k$.  Then we have
    \begin{multline*}
        %
        \ ,
    \end{multline*}
    where, in the final equality, we have used the fact that the $(+)$-bubble there is zero when $t > n$ by \cref{tanks}, together with our assumptions that $k \ge 0$.

    Next we prove \cref{bs2}.  Take the second relation in \cref{bs1} with $a$ replaced by $b^\vee c$, tensor on the left with
    %
    and on the right with
    %
,
    multiply by $z^2$, and then sum over $b,c \in \BA$ to give
    \[
        %
        \ .
    \]
    Then the second relation in \cref{bs2} follows after applying \cref{infweeds}.

    The identity \cref{bs3} is trivial when $n \ge 0$ by \cref{tanks2}.  When $n < 0$, the proof is almost identical to the proof of \cref{bs1} given above.
    \details{
        We have
        \begin{multline*}
            %
            \ ,
        \end{multline*}
        where, in the final equality, we have used the fact that the $(-)$-bubble there is zero when $t > -n$ by \cref{fake1,fake2}.
    }
    We then obtain the second relation in \cref{bs4} from the second relation in \cref{bs3} just as we did for the $(+)$-bubbles.
\end{proof}

\begin{lem}[Alternating braid relation]
    The following relation holds:
    \begin{align} \label{altbraid1}
        %
        \qquad \text{if } k \le 0.
    \end{align}
\end{lem}

\begin{proof}
    We prove \cref{altbraid1}, since \cref{altbraid2} then follows by rotating through $180^\circ$ and applying $\Omega_k$.  So we suppose $k \ge 0$.  Consider the relation
    \[
        %
    \]
    from \cref{spit1} and compose on the top and bottom with
    %
    respectively.  Then we have
    \begin{gather*}
        %
        \ ,
    \end{align*}
    and \cref{altbraid1} follows.
\end{proof}

\section{Action on modules over quantum cyclotomic wreath product algebras\label{sec:action}}

In this section we describe a natural categorical action of $\Heis_k(A;z,t)$ on categories of modules over cyclotomic quotients of quantum affine wreath product algebras.


Recall from \cref{center} the center $Z(A)$ of $A$, with even part $Z(A)_{\bar{0}}$. We fix a polynomial
\begin{equation} \label{parrot}
    f(w) = f_0 w^l + f_1 w^{l-1} + \dotsb + f_l \in Z(A)_{\bar{0}}[w],
\end{equation}
of degree $l \ge 0$ such that $f_0=1$ and $f_l = t^2$.  When $l \ge 1$, the \emph{quantum cyclotomic wreath product algebra} $\QWA_n^f = \QWA_n^f(A;z)$ of level $l$ asssociated to the polynomial $f(w)$  is the quotient of $\QAWA_n(A;z)$ by the two-sided ideal generated by $f(x_1)$.  We also allow $n=0$, with the convention that $\QWA_0^f = \kk$.  When $f(w)=1$, we have that $\QWA_0^f = \kk$ and $\QWA_n^f = 0$ for $n  > 0$.

The basis theorem \cite[Th.~4.10]{RS20} states that
\begin{equation} \label{CQ-basis}
    \left\{x_1^{r_1} \dotsb x_n^{r_n} \ba \sigma_g : 0 \le r_1,\dotsc,r_n < l,\ \ba \in \BA^{\otimes n},\ g \in \fS_n\right\}
\end{equation}
 is a basis of $\QWA_n^f$,
where $\sigma_g = \sigma_{i_1} \dotsm \sigma_{i_m}$ for some reduced expression $g = s_{i_1} \dotsm s_{i_m}$ for $g$ in the symmetric group $\fS_n$.  (The element $\sigma_g$ is independent of the choice of reduced expression.)  It follows that we have an injective homomorphism $\QWA_n^f \to \QWA_{n+1}^f$, sending the generators $x_i$, $\sigma_j$ to the elements of $\QWA_{n+1}^f$ with the same names and sending $\ba \in A^{\otimes n}$ to $1 \otimes \ba \in A^{\otimes (n+1)}$.  In this way, we identify $\QWA_n^f$ with a subalgebra of $\QWA_{n+1}^f$.  We then have induction and restriction superfunctors
\begin{align}
    \ind_n^{n+1} := - \otimes_{\QWA_n^f} \QWA_{n+1}^f &\colon \smod\QWA_n^f \to \smod\QWA_{n+1}^f,
    \\
    \res_n^{n+1} &: \smod\QWA_{n+1}^f \to \smod\QWA_n^f.
\end{align}

Our goal is to endow the abelian category $\bigoplus_{n \ge 0} \smod\QWA_n^f$ with the structure of a left $\Heis_{-l}(A;z,t)$-module category, with $\uparrow$ and $\downarrow$ acting as induction and restriction, respectively.  The key algebraic result that allows us to do this is the cyclotomic Mackey theorem \cite[Prop.~4.13]{RS20}, which states that we have an isomorphism of $(\QWA_n^f,\QWA_n^f)$-bimodules
\begin{equation} \label{marbleface}
        \QWA_n^f \otimes_{\QWA_{n-1}^f} \QWA_n^f \oplus \bigoplus_{\substack{0 \le r < l \\ b \in \BA}} \QWA_n^f \to \QWA_{n+1}^f,
        \quad
        \left( u \otimes v, (w_{r,b}) \right) \mapsto u \sigma_n v + \sum_{\substack{0 \le r < l \\ b \in \BA}} x_{n+1}^r b^{(n+1)} w_r,
\end{equation}
where we recall that $b^{(n+1)} = b \otimes 1^{\otimes n}$.  As explained in \cite[\S4.5]{RS20}, there is a unique homomorphism of $(\QWA_n^f,\QWA_n^f)$-bimodules
\[
    \tr_{n+1}^f \colon \QWA_{n+1}^f \to \QWA_n^f
\]
such that $\tr_{n+1}^f(\sigma_n) = 0$, $\tr_{n+1}^f( a^{(n+1)} x_{n+1}^r ) = \delta_{r,0} \tr(a)$ for $a \in A$, $0 \le r < l$.  (Recall that $\tr$ is the trace map of the Frobenius algebra $A$.)

\begin{lem} \label{fort}
    For any $n \ge 1$, we have $\tr_n^f(f(x_n)) = 0$.
\end{lem}

\begin{proof}
    The proof is almost identical to that of \cite[Lem.~6.1]{BSW-qheis}.
    \details{
        For $u,v \in \QWA_{n+1}^f$, write $u \equiv_n v$ as shorthand for $u=v$ in case $n=0$ or $u-v \in \QWA_n^f \sigma_n \QWA_n^f$ in case $n>0$.  We first show by induction on $n=0,1,\dotsc$ that
        \[ \tag{$\maltese$}
            \sigma_n \dotsm \sigma_1 x_1^m \sigma_1 \dotsm \sigma_n
            \equiv_n
            \begin{cases}
                \displaystyle \sum_{\substack{r + p_1 + \dotsb + p_n = m \\ r>0,\, p_1,\dotsc,p_n \ge 0}} \left( \prod_{i \text{ with } p_i \ne 0} (-z^2 \tau_i^2) \right) x_{n+1}^r x_n^{p_n} \dotsm x_1^{p_1} & \text{if } m > 0,
                \\
                \displaystyle \sum_{\substack{r + p_1 + \dotsb + p_n = m \\ r,p_1,\dotsc,p_n \le 0}} \left( \prod_{i \text{ with } p_i \ne 0} (z^2 \tau_i^2) \right) x_{n+1}^r x_n^{p_n} \dotsm x_1^{p_1} & \text{if } m \le 0.
            \end{cases}
        \]
        We give the details for $m>0$, since the case $m \le 0$ is similar.  (In fact, we will not need the $m<0$ cases.)  The base case is trivial.  For the induction step, we have
        \begin{multline*}
            \sigma_n x_n^m \sigma_n
            \overset{\cref{teaneg}}{=} \sigma_n \sigma_n^{-1} x_{n+1}^m - z \tau_n \sum_{\substack{r+s=m \\ r,s > 0}} \sigma_n x_{n+1}^r x_n^s
            \overset{\cref{skein}}{=} x_{n+1}^m - z \tau_n \sum_{\substack{r+s=m \\ r,s > 0}} \sigma_n^{-1} x_{n+1}^r x_n^s - z^2 \tau_n^2 \sum_{\substack{r+s=m \\ r,s > 0}} x_{n+1}^r x_n^s
            \\
            \overset{\cref{teaneg}}{\equiv_n} x_{n+1}^m - z^2 \tau_n^2 \sum_{\substack{r+s+u=m \\ r,s,u >0}} x_{n+1}^r x_n^{s+u} - z^2 \tau_n^2 \sum_{\substack{r+s=m \\ s,p > 0}} x_{n+1}^r x_n^s
            = x_{n+1}^m - z^2 \tau_n^2 \sum_{\substack{r+s=m \\ r,s > 0}} s x_{n+1}^r x_n^s.
        \end{multline*}
        Now take the expression for $\sigma_{n-1} \dotsm \sigma_1 x_1^m \sigma_1 \dotsm \sigma_{n-1}$ given by the induction hypothesis, multiply on the left and right by $\sigma_n$, and use the above identity plus the observation
        \begin{multline*}
            \sigma_n \left( \QWA_{n-1}^f \sigma_{n-1} \QWA_{n-1}^f \right) \sigma_n
            = \QWA_{n-1}^f \sigma_n \sigma_{n-1} \sigma_n \QWA_{n-1}^f
            \\
            = \QWA_{n-1}^r \sigma_{n-1} \sigma_n \sigma_{n-1} \QWA_{n-1}^f
            \subseteq \QWA_n^f \sigma_n \QWA_n^f.
        \end{multline*}

        Finally, to deduce the lemma, we multiply ($\maltese$) by $f_{l-m}$ and sum over $m=0,1,\dotsc,l$ to show
        \[
            \sigma_n \dotsm \sigma_1 f(x_1) \sigma_1 \dotsm \sigma_n
            \equiv_n f_l + \sum_{m=1}^l f_{l-a} \sum_{\substack{r + p_1 + \dotsb + p_n = m \\ r>0,\, p_1,\dotsc,p_n \ge 0}} \left( \prod_{i \text{ with } p_i \ne 0} (-z^2 \tau_i^2) \right) x_{n+1}^r x_n^{p_n} \dotsm x_1^{p_1}.
        \]
        The left-hand side is zero by the cyclotomic relation in $\QWA_{n+1}^f$.  The right-hand side is equal to $f(x_{n+1})$ plus terms in the kernel of $\tr_{n+1}^f$.
    }
\end{proof}


\begin{theo} \label{dingo}
    There is a unique strict $\kk$-linear monoidal superfunctor
    \[
        \Psi_f \colon \Heis_{-l}(A;z,t)
        \to \SEnd_\kk \left( \bigoplus_{n \ge 0} \smod\QWA_n^f \right)
    \]
    sending the generating object $\uparrow$ (resp.\ $\downarrow$) to the additive endosuperfunctor that takes a $\QWA_n^f$-module $M$ to $\ind_n^{n+1} M$ (resp.\ $\res^n_{n-1}M$), and the generating morphisms to the supernatural transformations defined on the $\QWA_n^f$-module $M$ as follows:
    \begin{itemize}
        \item $\Psi_f(\tokup)_M \colon M \otimes_{\QWA_n^f} \QWA_{n+1}^f \to M \otimes_{\QWA_n^f} \QWA_{n+1}^f$, $u \otimes v \mapsto (-1)^{\bar{a} \bar{u}} u \otimes a^{(n+1)} v$;

        \item $\Psi_f(\dotup)_M \colon M \otimes_{\QWA_n^f} \QWA_{n+1}^f \to M \otimes_{\QWA_n^f} \QWA_{n+1}^f$, $u \otimes v \mapsto u \otimes x_{n+1} v$;

        \item $\Psi_f \left( \poscross \right)_M \colon M \otimes_{\QWA_n^f} \QWA_{n+2}^f \to M \otimes_{\QWA_n^f} \QWA_{n+2}^f$, $u \otimes v \mapsto u \otimes \sigma_{n+1} v$ (where we have identified $\ind_{n+1}^{n+2} \circ \ind_n^{n+1}$ with $\ind_n^{n+2}$ in the obvious way);

        \item $\Psi_f \left( \rightcup \right)_M \colon M \to \res_n^{n+1} \left( M \otimes_{\QWA_n^f} \QWA_{n+1}^f \right)$, $v \mapsto v \otimes 1$, i.e. it is the unit of the canonical adjunction making $(\ind_n^{n+1}, \res_n^{n+1})$ into an adjoint pair of superfunctors;

        \item $\Psi_f \left( \rightcap \right)_M \colon (\res_{n-1}^n M) \otimes_{\QWA_{n-1}^f} \QWA_n^f \to M$, $u \otimes v \mapsto uv$, i.e.\ it is the counit of the canonical adjunction making $(\ind_{n-1}^n, \res_{n-1}^n)$ into an adjoint pair of superfunctors.
    \end{itemize}
\end{theo}

\begin{proof}
    The proof is similar to the proof of \cite[Th.~6.2]{BSW-qheis}, using the presentation of $\Heis_{-l}(A;z,t)$ from \cref{def1}.  The inversion relation follows from \cref{marbleface}, while the relation \cref{impose} follows from \cref{fort}.
    \details{
        We first consider the case $l=0$.  In this case, the polynomial $f(w)$ from \cref{parrot} is $1$ and $t^2=1$.  The category $\bigoplus_{n \ge 0} \smod\QWA_n^f$ is simply the category of right $\kk$-supermodules, and all of the induction and restriction superfunctors are zero.  Thus, almost all of the definition relations of $\Heis_{-1}(A;z,t)$ trivially hold.  The only nontrivial relation is \cref{impose}, which holds since $t^2=1$ implies $t=t^{-1}$.

        Now assume $l>0$.  Then $\Heis_{-1}(A;z,t)$ is generated by the objects $\uparrow$ and $\downarrow$ and morphisms $\tokup$ ($a \in A$), $\dotup$, $\poscross$, $\rightcup$, and $\rightcap$, subject to the relations \cref{tokrel,braid,skein,QAWC,rightadj}, plus two more relations:
        \begin{enumerate}
            \item \label{olaf1} the inversion relation, stating that the second morphism in \cref{invrel} in invertible, where the morphism
              $
                \begin{tikzpicture}[centerzero]
                    \draw[->] (-0.2,-0.2) -- (0.2,0.2);
                    \draw[wipe] (0.2,-0.2) -- (-0.2,0.2);
                    \draw[<-] (0.2,-0.2) -- (-0.2,0.2);
                \end{tikzpicture}
              $
              is defined by \cref{rcross};

            \item \label{olaf2} The relation \cref{impose}, stating that $\ccbubble{a}{l} = tz^{-1} 1_\one$, where $\leftcap$ is defined by \cref{leftwards}, i.e.\ it is $-t^{-1}z^{-1}$ times the $(r,b)=(0,1)$ entry of the inverse \cref{invrel2} of the matrix in \ref{olaf1}.  (Without loss of generality, we can choose the basis $B$ to contain the identity element of $A$.)
        \end{enumerate}

        The relations \cref{tokrel,braid,skein,QAWC,rightadj} are straightforward to verify.  On $\QWA_n^f$-modules,
        $
            \Psi_f
            \left(
                \begin{tikzpicture}[centerzero]
                    \draw[->] (-0.2,-0.2) -- (0.2,0.2);
                    \draw[wipe] (0.2,-0.2) -- (-0.2,0.2);
                    \draw[<-] (0.2,-0.2) -- (-0.2,0.2);
                \end{tikzpicture}
            \right)
        $
        corresponds to the $(\QWA_n^f,\QWA_n^f)$-bimodule homomorphism $\QWA_n^f \otimes_{\QWA_{n-1}^f} \QWA_n^f \to \QWA_{n+1}^f$, $u \otimes v \mapsto u \sigma_n v$.  Then the inversion relation \ref{olaf1} follows from \cref{marbleface}.  Furthermore, we see from \cref{marbleface} that $\Psi_f \left( \leftcap \right)$ corresponds to the $(\QWA_n^f,\QWA_n^f)$-bimodule homomorphisms $-t^{-1}z^{-1} \tr_{n+1}^f \colon \QAW_{n+1}^f \to \QAW_n^f$ for all $n \ge 0$.  So for \cref{olaf2}, we must show that $-t^{-1} z^{-1} \tr_{n+1}^f \left( a x_{n+1}^l \right) = tz^{-1} \tr(a)$.  This follows from \cref{fort}, the definition of $\tr_n^f$, and the assumption that $t^2 = f_l$.
    }
\end{proof}

We can reformulate \cref{dingo} in terms of Heisenberg categories of positive central charge by switching the roles of induction and restriction.  In fact, it is somewhat more natural to replace the induction superfunctor $\ind_n^{n+1}$, which is the canonical left adjoint to restriction, with the \emph{coinduction superfunctor}
\begin{equation}
    \coind_n^{n+1} := \Hom_{\QWA_n^f} \left( \QWA_{n+1}^f, - \right) \colon \smod\QWA_n^f \to \smod\QWA_{n+1}^f,
\end{equation}
which is its canonical right adjoint.

\begin{theo} \label{baby}
    There is a unique strict $\kk$-linear monoidal superfunctor
    \[
        \Psi_f^\vee \colon \Heis_l(A^\op;z,t^{-1})
        \to \SEnd_\kk \left( \bigoplus_{n \ge 0} \smod\QWA_n^f \right)
    \]
    sending the generating object $\uparrow$ (resp.\ $\downarrow$) to the additive endosuperfunctor that takes a $\QWA_n^f$-module $M$ to $\res_n^{n+1} M$ (resp.\ $\coind^n_{n-1}M$), and the generating morphisms to the natural transformations defined on the $\QWA_n^f$-module $M$ as follows:
    \begin{itemize}
        \item $\Psi_f^\vee(\tokup)_M \colon \res_{n-1}^n M \to \res_{n-1}^n M$, $v \mapsto (-1)^{\bar{a} \bar{v}} v a^{(n)}$;

        \item $\Psi_f^\vee(\dotup)_M \colon \res_{n-1}^n M \to \res_{n-1}^n M$, $v \mapsto v x_n$;

        \item $\Psi_f^\vee \left( \poscross \right)_M \colon \res_{n-2}^n M \to \res_{n-2}^n M$, $v \mapsto -v \sigma_{n-1}^{-1}$;

        \item $\Psi_f^\vee \left( \rightcup \right)_M \colon M \to \Hom_{\QWA_{n-1}^f} ( \QWA_n^f, \res_{n-1}^n M)$, $v \mapsto (u \mapsto vu)$,
          i.e.\ it is the unit of the canonical adjunction making $(\res_n^{n+1}, \coind_n^{n+1})$ into an adjoint pair of superfunctors;

        \item $\Psi_f^\vee \left( \rightcap \right)_M \colon \res_n^{n+1} \left( \Hom_{\QWA_n^f} ( \QWA_{n+1}^f, M) \right) \to M$, $\theta \mapsto \theta(1)$, i.e.\ it is the counit of the canonical adjunction making $(\res_{n-1}^n, \coind_{n-1}^n)$ into an adjoint pair of superfunctors.
    \end{itemize}
\end{theo}

\begin{proof}
    The proof is similar to the proof of \cref{dingo}, using instead the presentation for $\Heis_l(A^\op;z,t^{-1})$ from \cref{def2}; see also \cite[Th.~6.3]{BSW-qheis}.
\end{proof}

In fact, we have that $\ind_n^{n+1} \cong \coind_n^{n+1}$.  This follows from the fact that $\QWA_{n+1}^f$ is a Frobenius extension of $\QWA_n^f$; see \cite[Prop.~4.18]{RS20}.  It also follows from the results of the current paper and the uniqueness of adjoints, since \cref{adjfinal} and \cref{dingo} (resp.\ \cref{baby}) imply that $\ind_n^{n+1}$ is right adjoint to restriction as well as being left adjoint (resp.\ $\coind_n^{n+1}$ is left adjoint to restriction as well as being right adjoint).  As a consequence, all three functors (induction, coinduction, and restriction) send finitely-generated projective modules to finitely-generated projective modules.  This gives the following result.

\begin{lem} \label{sunset}
    The restrictions of the functors $\Psi_f$ and $\Psi_f^\vee$ from \cref{dingo,baby} give strict $\kk$-linear monoidal superfunctors
    \[
        \Psi_f \colon \Heis_{-l}(A;z,t) \to \SEnd_\kk \left( \bigoplus_{n \ge 0} \psmod\QWA_n^f \right),\quad
        \Psi_f^\vee \colon \Heis_l(A^\op;z,t^{-1}) \to \SEnd_\kk \left( \bigoplus_{n \ge 0} \psmod\QWA_n^f \right),
    \]
    recalling that $\psmod\QWA_n^f$ denotes the category of finite-generated projective right $\QWA_n^f$-supermodules.
\end{lem}

\section{Categorical comultiplication\label{sec:comult}}

In this section, we construct a quantum analog of the categorical comultiplication of \cite[Th.~5.12]{BSW-foundations}; see also \cite[Th.~8.9]{BSW-qheis} and \cite[Th.~5.4]{BSW-K0}.


Given a diagram $D$ representing a morphism in $\Heis_k(A;z,t)$ and two generic points $P$ and $Q$ on this diagram, we will denote the morphism represented by
\[
    D - (D \text{ with a negative dot at } P \text{ and a positive dot at } Q)
\]
by labelling the points with dots joined by a dotted line oriented from $P$ to $Q$.  We call this oriented dotted line a \emph{spear}.  For example,
\begin{equation} \label{spear}
    \begin{tikzpicture}[centerzero]
        \draw[->] (-0.2,-0.3) -- (-0.2,0.3);
        \draw[->] (0.2,-0.3) -- (0.2,0.3);
        \spear{black}{black}{-0.2,0}{0.2,0};
    \end{tikzpicture}
    :=
    \begin{tikzpicture}[centerzero]
        \draw[->] (-0.2,-0.3) -- (-0.2,0.3);
        \draw[->] (0.2,-0.3) -- (0.2,0.3);
    \end{tikzpicture}
    -
    \begin{tikzpicture}[centerzero]
        \draw[->] (-0.2,-0.3) -- (-0.2,0.3);
        \draw[->] (0.2,-0.3) -- (0.2,0.3);
        \multdot{-0.2,0}{east}{-1};
        \singdot{0.2,0};
    \end{tikzpicture}
    \ .
\end{equation}
It follows that dots and tokens pass through the endpoints of spears.  For example,
\begin{equation}
    \begin{tikzpicture}[centerzero]
        \draw[->] (-0.2,-0.4) -- (-0.2,0.4);
        \draw[->] (0.2,-0.4) -- (0.2,0.4);
        \spear{black}{black}{-0.2,0}{0.2,0};
        \token{-0.2,0.2}{east}{a};
    \end{tikzpicture}
    =
    \begin{tikzpicture}[centerzero]
        \draw[->] (-0.2,-0.4) -- (-0.2,0.4);
        \draw[->] (0.2,-0.4) -- (0.2,0.4);
        \spear{black}{black}{-0.2,0}{0.2,0};
        \token{-0.2,-0.2}{east}{a};
    \end{tikzpicture}
    \ ,\qquad
    \begin{tikzpicture}[centerzero]
        \draw[->] (-0.2,-0.4) -- (-0.2,0.4);
        \draw[->] (0.2,-0.4) -- (0.2,0.4);
        \spear{black}{black}{-0.2,0}{0.2,0};
        \token{0.2,0.2}{west}{a};
    \end{tikzpicture}
    =
    \begin{tikzpicture}[centerzero]
        \draw[->] (-0.2,-0.4) -- (-0.2,0.4);
        \draw[->] (0.2,-0.4) -- (0.2,0.4);
        \spear{black}{black}{-0.2,0}{0.2,0};
        \token{0.2,-0.2}{west}{a};
    \end{tikzpicture}
    \ ,\qquad
    \begin{tikzpicture}[centerzero]
        \draw[->] (-0.2,-0.4) -- (-0.2,0.4);
        \draw[->] (0.2,-0.4) -- (0.2,0.4);
        \spear{black}{black}{-0.2,0}{0.2,0};
        \singdot{-0.2,0.2};
    \end{tikzpicture}
    =
    \begin{tikzpicture}[centerzero]
        \draw[->] (-0.2,-0.4) -- (-0.2,0.4);
        \draw[->] (0.2,-0.4) -- (0.2,0.4);
        \spear{black}{black}{-0.2,0}{0.2,0};
        \singdot{-0.2,-0.2};
    \end{tikzpicture}
    \ ,\qquad
    \begin{tikzpicture}[centerzero]
        \draw[->] (-0.2,-0.4) -- (-0.2,0.4);
        \draw[->] (0.2,-0.4) -- (0.2,0.4);
        \spear{black}{black}{-0.2,0}{0.2,0};
        \singdot{0.2,0.2};
    \end{tikzpicture}
    =
    \begin{tikzpicture}[centerzero]
        \draw[->] (-0.2,-0.4) -- (-0.2,0.4);
        \draw[->] (0.2,-0.4) -- (0.2,0.4);
        \spear{black}{black}{-0.2,0}{0.2,0};
        \singdot{0.2,-0.2};
    \end{tikzpicture}
    \ .
\end{equation}
We also have
\begin{equation} \label{yellow}
    \begin{tikzpicture}[centerzero]
        \draw[->] (-0.2,-0.3) -- (-0.2,0.3);
        \draw[->] (0.2,-0.3) -- (0.2,0.3);
        \spear{black}{black}{-0.2,-0.1}{0.2,-0.1};
        \singdot{-0.2,0.1};
        \multdot{0.2,0.1}{west}{-1};
    \end{tikzpicture}
    = -\
    \begin{tikzpicture}[centerzero]
        \draw[->] (-0.2,-0.3) -- (-0.2,0.3);
        \draw[->] (0.2,-0.3) -- (0.2,0.3);
        \spear{black}{black}{0.2,-0.1}{-0.2,-0.1};
    \end{tikzpicture}
    \ ,\qquad
    \begin{tikzpicture}[centerzero]
        \draw[->] (-0.2,-0.3) -- (-0.2,0.3);
        \draw[->] (0.2,-0.3) -- (0.2,0.3);
        \spear{black}{black}{-0.2,0.1}{0.2,0.1};
        \spear{black}{black}{0.2,-0.1}{-0.2,-0.1};
    \end{tikzpicture}
    =
    \begin{tikzpicture}[centerzero]
        \draw[->] (-0.2,-0.3) -- (-0.2,0.3);
        \draw[->] (0.2,-0.3) -- (0.2,0.3);
        \spear{black}{black}{-0.2,0}{0.2,0};
    \end{tikzpicture}
    +
    \begin{tikzpicture}[centerzero]
        \draw[->] (-0.2,-0.3) -- (-0.2,0.3);
        \draw[->] (0.2,-0.3) -- (0.2,0.3);
        \spear{black}{black}{0.2,0}{-0.2,0};
    \end{tikzpicture}
    \ .
\end{equation}

Now define the \emph{neutral crossing}
\begin{equation} \label{intertwine}
    \begin{tikzpicture}[centerzero]
        \draw[->] (-0.2,-0.4) -- (0.2,0.4);
        \draw[->] (0.2,-0.4) -- (-0.2,0.4);
    \end{tikzpicture}
    :=
    \begin{tikzpicture}[centerzero]
        \draw[->] (-0.2,-0.4) -- (0.2,0.4);
        \draw[wipe] (0.2,-0.4) -- (-0.2,0.4);
        \draw[->] (0.2,-0.4) -- (-0.2,0.4);
    \end{tikzpicture}
    -
    \begin{tikzpicture}[centerzero]
        \draw[->] (0.2,-0.4) -- (0.2,-0.2) \braidup (-0.2,0.4);
        \draw[wipe] (-0.2,-0.4) -- (-0.2,-0.2) \braidup (0.2,0.4);
        \draw[->] (-0.2,-0.4) -- (-0.2,-0.2) \braidup (0.2,0.4);
        \multdot{-0.2,-0.2}{east}{-1};
        \singdot{0.2,-0.2};
    \end{tikzpicture}
    =
    \begin{tikzpicture}[centerzero]
        \draw[->] (0.2,-0.4) -- (-0.2,0.4);
        \draw[wipe] (-0.2,-0.4) -- (0.2,0.4);
        \draw[->] (-0.2,-0.4) -- (0.2,0.4);
    \end{tikzpicture}
    -
    \begin{tikzpicture}[centerzero]
        \draw[->] (-0.2,-0.4) \braidup (0.2,0.2) -- (0.2,0.4);
        \draw[wipe] (0.2,-0.4) \braidup (-0.2,0.2) -- (-0.2,0.4);
        \draw[->] (0.2,-0.4) \braidup (-0.2,0.2) -- (-0.2,0.4);
        \multdot{0.2,0.2}{west}{-1};
        \singdot{-0.2,0.2};
    \end{tikzpicture}
    =
    \begin{tikzpicture}[centerzero]
        \draw[->] (0.2,-0.4) -- (0.2,-0.2) \braidup (-0.2,0.4);
        \draw[wipe] (-0.2,-0.4) -- (-0.2,-0.2) \braidup (0.2,0.4);
        \draw[->] (-0.2,-0.4) -- (-0.2,-0.2) \braidup (0.2,0.4);
        \spear{black}{black}{-0.2,-0.2}{0.2,-0.2};
    \end{tikzpicture}
    -
    \begin{tikzpicture}[centerzero]
        \draw[->] (-0.2,-0.4) -- (-0.2,0.4);
        \draw[->] (0.2,-0.4) -- (0.2,0.4);
        \teleport{-0.2,0}{0.2,0};
    \end{tikzpicture}
    =
    \begin{tikzpicture}[centerzero]
        \draw[->] (-0.2,-0.4) \braidup (0.2,0.2) -- (0.2,0.4);
        \draw[wipe] (0.2,-0.4) \braidup (-0.2,0.2) -- (-0.2,0.4);
        \draw[->] (0.2,-0.4) \braidup (-0.2,0.2) -- (-0.2,0.4);
        \spear{black}{black}{0.2,0.2}{-0.2,0.2};
    \end{tikzpicture}
    +
    \begin{tikzpicture}[centerzero]
        \draw[->] (-0.2,-0.4) -- (-0.2,0.4);
        \draw[->] (0.2,-0.4) -- (0.2,0.4);
        \teleport{-0.2,0}{0.2,0};
    \end{tikzpicture}
    \ .
\end{equation}
Direct computation shows that
\begin{equation} \label{interprop}
    \begin{tikzpicture}[centerzero]
        \draw[->] (-0.2,-0.4) -- (0.2,0.4);
        \draw[->] (0.2,-0.4) -- (-0.2,0.4);
        \token{-0.1,-0.2}{east}{a};
    \end{tikzpicture}
    =
    \begin{tikzpicture}[centerzero]
        \draw[->] (-0.2,-0.4) -- (0.2,0.4);
        \draw[->] (0.2,-0.4) -- (-0.2,0.4);
        \token{0.1,0.2}{west}{a};
    \end{tikzpicture}
    \ ,\quad
    \begin{tikzpicture}[centerzero]
        \draw[->] (-0.2,-0.4) -- (0.2,0.4);
        \draw[->] (0.2,-0.4) -- (-0.2,0.4);
        \token{0.1,-0.2}{west}{a};
    \end{tikzpicture}
    =
    \begin{tikzpicture}[centerzero]
        \draw[->] (-0.2,-0.4) -- (0.2,0.4);
        \draw[->] (0.2,-0.4) -- (-0.2,0.4);
        \token{-0.1,0.2}{east}{a};
    \end{tikzpicture}
    \ ,\quad
    \begin{tikzpicture}[centerzero]
        \draw[->] (-0.2,-0.4) -- (0.2,0.4);
        \draw[->] (0.2,-0.4) -- (-0.2,0.4);
        \singdot{-0.1,-0.2};
    \end{tikzpicture}
    =
    \begin{tikzpicture}[centerzero]
        \draw[->] (-0.2,-0.4) -- (0.2,0.4);
        \draw[->] (0.2,-0.4) -- (-0.2,0.4);
        \singdot{0.1,0.2};
    \end{tikzpicture}
    \ ,\quad
    \begin{tikzpicture}[centerzero]
        \draw[->] (-0.2,-0.4) -- (0.2,0.4);
        \draw[->] (0.2,-0.4) -- (-0.2,0.4);
        \singdot{0.1,-0.2};
    \end{tikzpicture}
    =
    \begin{tikzpicture}[centerzero]
        \draw[->] (-0.2,-0.4) -- (0.2,0.4);
        \draw[->] (0.2,-0.4) -- (-0.2,0.4);
        \singdot{-0.1,0.2};
    \end{tikzpicture}
    \ ,\quad
    \begin{tikzpicture}[anchorbase]
        \draw[->] (0.4,-0.4) -- (-0.4,0.4);
        \draw[->] (0,-0.4) \braidup (-0.3,0) \braidup (0,0.4);
        \draw[->] (-0.4,-0.4) -- (0.4,0.4);
    \end{tikzpicture}
    \ =\
    \begin{tikzpicture}[anchorbase]
        \draw[->] (0.4,-0.4) -- (-0.4,0.4);
        \draw[->] (0,-0.4) \braidup (0.3,0) \braidup (0,0.4);
        \draw[->] (-0.4,-0.4) -- (0.4,0.4);
    \end{tikzpicture}
    \ ,\quad
    \begin{tikzpicture}[anchorbase]
        \draw[->] (0.2,-0.4) \braidup (-0.2,0)\braidup (0.2,0.4);
        \draw[->] (-0.2,-0.4) \braidup (0.2,0) \braidup (-0.2,0.4);
    \end{tikzpicture}
    \ =\
    \begin{tikzpicture}[anchorbase]
        \draw[->] (-0.2,-0.4) -- (-0.2,0.4);
        \draw[->] (0.2,-0.4) -- (0.2,0.4);
        \spear{black}{black}{-0.2,0.15}{0.2,0.15};
        \spear{black}{black}{0.2,-0.15}{-0.2,-0.15};
    \end{tikzpicture}
    +
    \begin{tikzpicture}[anchorbase]
        \draw[->] (-0.2,-0.4) -- (-0.2,0.4);
        \draw[->] (0.2,-0.4) -- (0.2,0.4);
        \teleport{-0.2,0.15}{0.2,0.15};
        \teleport{0.2,-0.15}{-0.2,-0.15};
    \end{tikzpicture}
    \ .
\end{equation}
(The proof of the fifth relation is lengthy but straightforward.)  The relations \cref{interprop} show that the neutral crossing is an analogue of the intertwining operators that play an important role in the study of (degenerate) affine Hecke algebras.  In addition, the neutral crossing is a quantum analogue of the intertwining operators for affine wreath product algebras introduced in \cite[\S4.7]{Sav20}.

We would like use the neutral crossing as a generator instead of the positive and negative crossings.  However, in order to write the positive and negative crossings in terms of the neutral crossings, we need to invert the spear.  Let $\invQAW(A;z)$ denote the linear monoidal supercategory obtained from $\QAW(A;z)$ by adjoining a two-sided inverse to the spear:
\begin{equation} \label{dart}
    \begin{tikzpicture}[centerzero]
        \draw[->] (-0.2,-0.3) -- (-0.2,0.3);
        \draw[->] (0.2,-0.3) -- (0.2,0.3);
        \dart{black}{black}{-0.2,0}{0.2,0};
    \end{tikzpicture}
    :=
    \left(
        \begin{tikzpicture}[centerzero]
            \draw[->] (-0.2,-0.3) -- (-0.2,0.3);
            \draw[->] (0.2,-0.3) -- (0.2,0.3);
            \spear{black}{black}{-0.2,0}{0.2,0};
        \end{tikzpicture}
    \right)^{-1}.
\end{equation}
We call this inverse a \emph{dart}.  Using \cref{yellow}, we then automatically have the other dart
\begin{tikzpicture}[centerzero]
    \draw[->] (-0.2,-0.2) -- (-0.2,0.2);
    \draw[->] (0.2,-0.2) -- (0.2,0.2);
    \dart{black}{black}{0.2,0}{-0.2,0};
\end{tikzpicture}
:=
$\left(
    \begin{tikzpicture}[centerzero]
        \draw[->] (-0.2,-0.2) -- (-0.2,0.2);
        \draw[->] (0.2,-0.2) -- (0.2,0.2);
        \spear{black}{black}{0.2,0}{-0.2,0};
    \end{tikzpicture}
\right)^{-1}$.
It also follows from \cref{yellow} that
\begin{gather} \label{green}
    \begin{tikzpicture}[centerzero]
        \draw[->] (-0.2,-0.3) -- (-0.2,0.3);
        \draw[->] (0.2,-0.3) -- (0.2,0.3);
        \dart{black}{black}{0.2,-0.1}{-0.2,-0.1};
        \singdot{-0.2,0.1};
        \multdot{0.2,0.1}{west}{-1};
    \end{tikzpicture}
    = -\
    \begin{tikzpicture}[centerzero]
        \draw[->] (-0.2,-0.3) -- (-0.2,0.3);
        \draw[->] (0.2,-0.3) -- (0.2,0.3);
        \dart{black}{black}{-0.2,-0.1}{0.2,-0.1};
    \end{tikzpicture}
    \ ,\qquad
    \begin{tikzpicture}[centerzero]
        \draw[->] (-0.2,-0.3) -- (-0.2,0.3);
        \draw[->] (0.2,-0.3) -- (0.2,0.3);
    \end{tikzpicture}
    =
    \begin{tikzpicture}[centerzero]
        \draw[->] (-0.2,-0.3) -- (-0.2,0.3);
        \draw[->] (0.2,-0.3) -- (0.2,0.3);
        \dart{black}{black}{-0.2,0}{0.2,0};
    \end{tikzpicture}
    +
    \begin{tikzpicture}[centerzero]
        \draw[->] (-0.2,-0.3) -- (-0.2,0.3);
        \draw[->] (0.2,-0.3) -- (0.2,0.3);
        \dart{black}{black}{0.2,0}{-0.2,0};
    \end{tikzpicture}
    \ ,\qquad
    \begin{tikzpicture}[centerzero]
        \draw[->] (-0.2,-0.3) -- (-0.2,0.3);
        \draw[->] (0.2,-0.3) -- (0.2,0.3);
        \dart{black}{black}{-0.2,-0.1}{0.2,-0.1};
        \multdot{-0.2,0.1}{east}{n};
    \end{tikzpicture}
    =
    \begin{tikzpicture}[centerzero]
        \draw[->] (-0.2,-0.3) -- (-0.2,0.3);
        \draw[->] (0.2,-0.3) -- (0.2,0.3);
        \dart{black}{black}{-0.2,-0.1}{0.2,-0.1};
        \multdot{0.2,0.1}{west}{n};
    \end{tikzpicture}
    + \sum_{\substack{r+s=n \\ r>0,\, s \ge 0}}
    \begin{tikzpicture}[centerzero]
        \draw[->] (-0.2,-0.3) -- (-0.2,0.3);
        \draw[->] (0.2,-0.3) -- (0.2,0.3);
        \multdot{-0.2,0}{east}{r};
        \multdot{0.2,0}{west}{s};
    \end{tikzpicture}
    - \sum_{\substack{r+s=n \\ r \le 0,\, s < 0}}
    \begin{tikzpicture}[centerzero]
        \draw[->] (-0.2,-0.3) -- (-0.2,0.3);
        \draw[->] (0.2,-0.3) -- (0.2,0.3);
        \multdot{-0.2,0}{east}{r};
        \multdot{0.2,0}{west}{s};
    \end{tikzpicture}
    \ ,
    \\ \label{greener}
    \begin{tikzpicture}[centerzero]
        \draw[->] (-0.2,-0.3) -- (-0.2,0.3);
        \draw[->] (0.2,-0.3) -- (0.2,0.3);
        \dart{black}{black}{-0.2,-0.1}{0.2,-0.1};
        \multdot{-0.2,0.1}{east}{n};
    \end{tikzpicture}
    +
    \begin{tikzpicture}[centerzero]
        \draw[->] (-0.2,-0.3) -- (-0.2,0.3);
        \draw[->] (0.2,-0.3) -- (0.2,0.3);
        \dart{black}{black}{0.2,-0.1}{-0.2,-0.1};
        \multdot{0.2,0.1}{west}{n};
    \end{tikzpicture}
    = \sum_{\substack{r+s=n \\ r, s \ge 0}}
    \begin{tikzpicture}[centerzero]
        \draw[->] (-0.2,-0.3) -- (-0.2,0.3);
        \draw[->] (0.2,-0.3) -- (0.2,0.3);
        \multdot{-0.2,0}{east}{r};
        \multdot{0.2,0}{west}{s};
    \end{tikzpicture}
    - \sum_{\substack{r+s=n \\ r,s < 0}}
    \begin{tikzpicture}[centerzero]
        \draw[->] (-0.2,-0.3) -- (-0.2,0.3);
        \draw[->] (0.2,-0.3) -- (0.2,0.3);
        \multdot{-0.2,0}{east}{r};
        \multdot{0.2,0}{west}{s};
    \end{tikzpicture}
    \ ,\quad
    \begin{tikzpicture}[centerzero]
        \draw[->] (-0.2,-0.3) -- (-0.2,0.3);
        \draw[->] (0.2,-0.3) -- (0.2,0.3);
        \dart{black}{black}{0.2,-0.1}{-0.2,-0.1};
        \multdot{-0.2,0.1}{east}{n};
    \end{tikzpicture}
    +
    \begin{tikzpicture}[centerzero]
        \draw[->] (-0.2,-0.3) -- (-0.2,0.3);
        \draw[->] (0.2,-0.3) -- (0.2,0.3);
        \dart{black}{black}{-0.2,-0.1}{0.2,-0.1};
        \multdot{0.2,0.1}{west}{n};
    \end{tikzpicture}
    = - \sum_{\substack{r+s=n \\ r, s > 0}}
    \begin{tikzpicture}[centerzero]
        \draw[->] (-0.2,-0.3) -- (-0.2,0.3);
        \draw[->] (0.2,-0.3) -- (0.2,0.3);
        \multdot{-0.2,0}{east}{r};
        \multdot{0.2,0}{west}{s};
    \end{tikzpicture}
    + \sum_{\substack{r+s=n \\ r,s \le 0}}
    \begin{tikzpicture}[centerzero]
        \draw[->] (-0.2,-0.3) -- (-0.2,0.3);
        \draw[->] (0.2,-0.3) -- (0.2,0.3);
        \multdot{-0.2,0}{east}{r};
        \multdot{0.2,0}{west}{s};
    \end{tikzpicture}
    \ ,
\end{gather}
for all $n \in \Z$.

\begin{theo} \label{invQAWpres}
    The strict $\kk$-linear monoidal supercategory $\invQAW(A;z)$ is generated by one object $\uparrow$ and morphisms
    \begin{equation}
        \begin{tikzpicture}[centerzero]
            \draw[->] (0.2,-0.2) -- (-0.2,0.2);
            \draw[->] (-0.2,-0.2) -- (0.2,0.2);
        \end{tikzpicture}
        \ ,\
        \begin{tikzpicture}[centerzero]
            \draw[->] (-0.2,-0.2) -- (-0.2,0.2);
            \draw[->] (0.2,-0.2) -- (0.2,0.2);
            \dart{black}{black}{-0.2,0}{0.2,0};
        \end{tikzpicture}
        \ \colon
        \uparrow \otimes \uparrow \to \uparrow \otimes \uparrow,
        \qquad
        \begin{tikzpicture}[anchorbase]
            \draw[->] (0,-0.2) -- (0,0.2);
            \token{0,0}{west}{a};
        \end{tikzpicture}
        \colon \uparrow \to \uparrow,\quad a \in A,
    \end{equation}
    where the crossings are even and the parity of the morphism
    \begin{tikzpicture}[anchorbase]
        \draw[->] (0,-0.2) -- (0,0.2);
        \token{0,0}{west}{a};
    \end{tikzpicture}
    is the same as the parity of $A$.  A complete set of relations is given by \cref{tokrel,interprop,dart} (where in \cref{dart} we use the definition \cref{spear} of the spear).
\end{theo}

\begin{proof}
    Let $\cC$ be the category whose presentation is given in the statement of the theorem.  As noted above, we have a strict monoidal superfunctor $\cC \to \invQAW(A;z)$ sending each generating morphism to the morphism of $\invQAW(A;z)$ represented by the same string diagram.  Its inverse is given by mapping tokens and dots to the same tokens and dots, and
    \begin{align}
        \begin{tikzpicture}[centerzero]
            \draw[->] (0.2,-0.2) -- (-0.2,0.2);
            \draw[wipe] (-0.2,-0.2) -- (0.2,0.2);
            \draw[->] (-0.2,-0.2) -- (0.2,0.2);
        \end{tikzpicture}
        &\mapsto
        \begin{tikzpicture}[centerzero]
            \draw[->] (0.2,-0.4) -- (0.2,-0.2) \braidup (-0.2,0.4);
            \draw[->] (-0.2,-0.4) -- (-0.2,-0.2) \braidup (0.2,0.4);
            \dart{black}{black}{-0.2,-0.2}{0.2,-0.2};
        \end{tikzpicture}
        +
        \begin{tikzpicture}[centerzero]
            \draw[->] (-0.2,-0.4) -- (-0.2,0.4);
            \draw[->] (0.2,-0.4) -- (0.2,0.4);
            \teleport{-0.2,0.15}{0.2,0.15};
            \dart{black}{black}{-0.2,-0.15}{0.2,-0.15};
        \end{tikzpicture}
        \ ,&
        \begin{tikzpicture}[centerzero]
            \draw[->] (-0.2,-0.2) -- (0.2,0.2);
            \draw[wipe] (0.2,-0.2) -- (-0.2,0.2);
            \draw[->] (0.2,-0.2) -- (-0.2,0.2);
        \end{tikzpicture}
        &\mapsto
        \begin{tikzpicture}[centerzero]
            \draw[->] (-0.2,-0.4) \braidup (0.2,0.2) -- (0.2,0.4);
            \draw[->] (0.2,-0.4) \braidup (-0.2,0.2) -- (-0.2,0.4);
            \dart{black}{black}{0.2,0.2}{-0.2,0.2};
        \end{tikzpicture}
        -
        \begin{tikzpicture}[centerzero]
            \draw[->] (-0.2,-0.4) -- (-0.2,0.4);
            \draw[->] (0.2,-0.4) -- (0.2,0.4);
            \teleport{-0.2,0.15}{0.2,0.15};
            \dart{black}{black}{0.2,-0.15}{-0.2,-0.15};
        \end{tikzpicture}
        \ ,
    \end{align}
    where the left-pointing spear is defined using the first relation in \cref{green}.  Direct computation shows that this respects the relations \cref{braid,skein}, hence is well defined.
\end{proof}

Define the \emph{box dumbbell}
\begin{equation} \label{ups}
    \begin{tikzpicture}[centerzero]
        \draw[->] (-0.2,-0.4) -- (-0.2,0.4);
        \draw[->] (0.2,-0.4) -- (0.2,0.4);
        \boxbell{black}{black}{-0.2,0}{0.2,0};
    \end{tikzpicture}
    \ :=\
    \begin{tikzpicture}[centerzero]
        \draw[->] (-0.2,-0.4) -- (-0.2,0.4);
        \draw[->] (0.2,-0.4) -- (0.2,0.4);
        \spear{black}{black}{-0.2,0.15}{0.2,0.15};
        \spear{black}{black}{0.2,-0.15}{-0.2,-0.15};
    \end{tikzpicture}
    +
    \begin{tikzpicture}[centerzero]
        \draw[->] (-0.2,-0.4) -- (-0.2,0.4);
        \draw[->] (0.2,-0.4) -- (0.2,0.4);
        \teleport{-0.2,0.15}{0.2,0.15};
        \teleport{0.2,-0.15}{-0.2,-0.15};
    \end{tikzpicture}
    \ ,\qquad \text{so that}\quad
    \begin{tikzpicture}[centerzero]
        \draw[->] (0.2,-0.4) \braidup (-0.2,0)\braidup (0.2,0.4);
        \draw[->] (-0.2,-0.4) \braidup (0.2,0) \braidup (-0.2,0.4);
    \end{tikzpicture}
    \ =\
    \begin{tikzpicture}[centerzero]
        \draw[->] (-0.2,-0.4) -- (-0.2,0.4);
        \draw[->] (0.2,-0.4) -- (0.2,0.4);
        \boxbell{black}{black}{-0.2,0}{0.2,0};
    \end{tikzpicture}
    \ .
\end{equation}
Straightforward computation shows that the box dumbbell is central in $\End_{\QAW(A;z)}(\uparrow^{\otimes 2})$:
\begin{equation}
    \begin{tikzpicture}[centerzero]
        \draw[->] (-0.2,-0.4) -- (-0.2,0.4);
        \draw[->] (0.2,-0.4) -- (0.2,0.4);
        \boxbell{black}{black}{-0.2,0}{0.2,0};
        \token{-0.2,-0.2}{west}{a};
    \end{tikzpicture}
    \ =\
    \begin{tikzpicture}[centerzero]
        \draw[->] (-0.2,-0.4) -- (-0.2,0.4);
        \draw[->] (0.2,-0.4) -- (0.2,0.4);
        \boxbell{black}{black}{-0.2,0}{0.2,0};
        \token{-0.2,0.2}{west}{a};
    \end{tikzpicture}
    \ ,\ \
    \begin{tikzpicture}[centerzero]
        \draw[->] (-0.2,-0.4) -- (-0.2,0.4);
        \draw[->] (0.2,-0.4) -- (0.2,0.4);
        \boxbell{black}{black}{-0.2,0}{0.2,0};
        \token{0.2,-0.2}{east}{a};
    \end{tikzpicture}
    \ =\
    \begin{tikzpicture}[centerzero]
        \draw[->] (-0.2,-0.4) -- (-0.2,0.4);
        \draw[->] (0.2,-0.4) -- (0.2,0.4);
        \boxbell{black}{black}{-0.2,0}{0.2,0};
        \token{0.2,0.2}{east}{a};
    \end{tikzpicture}
    \ ,\ \
    \begin{tikzpicture}[centerzero]
        \draw[->] (-0.2,-0.4) -- (-0.2,0.4);
        \draw[->] (0.2,-0.4) -- (0.2,0.4);
        \boxbell{black}{black}{-0.2,0}{0.2,0};
        \singdot{-0.2,-0.2};
    \end{tikzpicture}
    \ =\
    \begin{tikzpicture}[centerzero]
        \draw[->] (-0.2,-0.4) -- (-0.2,0.4);
        \draw[->] (0.2,-0.4) -- (0.2,0.4);
        \boxbell{black}{black}{-0.2,0}{0.2,0};
        \singdot{-0.2,0.2};
    \end{tikzpicture}
    \ ,\ \
    \begin{tikzpicture}[centerzero]
        \draw[->] (-0.2,-0.4) -- (-0.2,0.4);
        \draw[->] (0.2,-0.4) -- (0.2,0.4);
        \boxbell{black}{black}{-0.2,0}{0.2,0};
        \singdot{0.2,-0.2};
    \end{tikzpicture}
    \ =\
    \begin{tikzpicture}[centerzero]
        \draw[->] (-0.2,-0.4) -- (-0.2,0.4);
        \draw[->] (0.2,-0.4) -- (0.2,0.4);
        \boxbell{black}{black}{-0.2,0}{0.2,0};
        \singdot{0.2,0.2};
    \end{tikzpicture}
    \ ,\ \
    \begin{tikzpicture}[centerzero]
        \draw[->] (0.2,-0.4) -- (0.2,-0.2) \braidup (-0.2,0.4);
        \draw[wipe] (-0.2,-0.4) -- (-0.2,-0.2) \braidup (0.2,0.4);
        \draw[->] (-0.2,-0.4) -- (-0.2,-0.2) \braidup (0.2,0.4);
        \boxbell{black}{black}{-0.2,-0.2}{0.2,-0.2};
    \end{tikzpicture}
    =
    \begin{tikzpicture}[centerzero]
        \draw[->] (0.2,-0.4) \braidup (-0.2,0.2) -- (-0.2,0.4);
        \draw[wipe] (-0.2,-0.4) \braidup (0.2,0.2) -- (0.2,0.4);
        \draw[->] (-0.2,-0.4) \braidup (0.2,0.2) -- (0.2,0.4);
        \boxbell{black}{black}{-0.2,0.2}{0.2,0.2};
    \end{tikzpicture}
    \ ,\ \
    \begin{tikzpicture}[centerzero]
        \draw[->] (-0.2,-0.4) -- (-0.2,-0.2) \braidup (0.2,0.4);
        \draw[wipe] (0.2,-0.4) -- (0.2,-0.2) \braidup (-0.2,0.4);
        \draw[->] (0.2,-0.4) -- (0.2,-0.2) \braidup (-0.2,0.4);
        \boxbell{black}{black}{-0.2,-0.2}{0.2,-0.2};
    \end{tikzpicture}
    =
    \begin{tikzpicture}[centerzero]
        \draw[->] (-0.2,-0.4) \braidup (0.2,0.2) -- (0.2,0.4);
        \draw[wipe] (0.2,-0.4) \braidup (-0.2,0.2) -- (-0.2,0.4);
        \draw[->] (0.2,-0.4) \braidup (-0.2,0.2) -- (-0.2,0.4);
        \boxbell{black}{black}{-0.2,0.2}{0.2,0.2};
    \end{tikzpicture}
    \ ,\ \
    \begin{tikzpicture}[centerzero]
        \draw[->] (-0.2,-0.4) -- (-0.2,-0.2) \braidup (0.2,0.4);
        \draw[->] (0.2,-0.4) -- (0.2,-0.2) \braidup (-0.2,0.4);
        \boxbell{black}{black}{-0.2,-0.2}{0.2,-0.2};
    \end{tikzpicture}
    =
    \begin{tikzpicture}[centerzero]
        \draw[->] (-0.2,-0.4) \braidup (0.2,0.2) -- (0.2,0.4);
        \draw[->] (0.2,-0.4) \braidup (-0.2,0.2) -- (-0.2,0.4);
        \boxbell{black}{black}{-0.2,0.2}{0.2,0.2};
    \end{tikzpicture}
    \ .
\end{equation}
We also introduce the \emph{box darts}
\begin{equation} \label{boxdart}
    \begin{tikzpicture}[centerzero]
        \draw[->] (-0.2,-0.4) -- (-0.2,0.4);
        \draw[->] (0.2,-0.4) -- (0.2,0.4);
        \boxdart{black}{black}{-0.2,0}{0.2,0};
    \end{tikzpicture}
    :=
    \begin{tikzpicture}[centerzero]
        \draw[->] (-0.2,-0.4) -- (-0.2,0.4);
        \draw[->] (0.2,-0.4) -- (0.2,0.4);
        \boxbell{black}{black}{-0.2,0.15}{0.2,0.15};
        \dart{black}{black}{-0.2,-0.15}{0.2,-0.15};
    \end{tikzpicture}
    =
    \begin{tikzpicture}[centerzero]
        \draw[->] (-0.2,-0.4) -- (-0.2,0.4);
        \draw[->] (0.2,-0.4) -- (0.2,0.4);
        \spear{black}{black}{0.2,0}{-0.2,0};
    \end{tikzpicture}
    +
    \begin{tikzpicture}[centerzero]
        \draw[->] (-0.2,-0.4) -- (-0.2,0.4);
        \draw[->] (0.2,-0.4) -- (0.2,0.4);
        \dart{black}{black}{-0.2,-0.2}{0.2,-0.2};
        \teleport{-0.2,0.2}{0.2,0.2};
        \teleport{-0.2,0}{0.2,0};
    \end{tikzpicture}
    \ ,
\end{equation}
and similarly for the other orientation of the box darts or, more generally, a box dart joining any two generic points in a string diagram.

The categorical comultiplication to be defined shortly will go from $\Heis_k(A;z,t)$ to a certain monoidal supercategory built from $\blue{\Heis_l(A;z,u)}$ and $\red{\Heis_m(A;z,v)}$ for $\blue{l}, \red{m} \in \Z$ and $\blue{u}, \red{v} \in \kk^\times$ chosen so that
\begin{equation}
    k = \blue{l}+\red{m}, \qquad
    t =\blue{u}\red{v}.
\end{equation}
To avoid confusion between these different categories, the reader will want to view the subsequent material in this section in color.

Let $\blue{\Heis_l(A;z,u)}\odot \red{\Heis_m(A;z,v)}$ be the symmetric product of $\blue{\Heis_l(A;z,u)}$ and $\red{\Heis_m(A;z,v)}$ as defined in \cite[$\S$3]{BSW-K0}.  This is the strict $\kk$-linear monoidal category defined by first taking the free product of $\blue{\Heis_l(A;z,u)}$ and $\red{\Heis_m(A;z,v)}$, i.e.\
the strict $\kk$-linear monoidal category defined by the disjoint union of the given generators and relations of $\blue{\Heis_l(A;z,u)}$ and of $\red{\Heis_m(A;z,v)}$, then adjoining even isomorphisms $\sigma_{X,Y} \colon X \otimes Y \xrightarrow{\cong} Y \otimes X$ for each pair of objects $X \in \blue{\Heis_l(A;z,u)}$ and $Y \in \red{\Heis_m(A;z,v)}$ subject to the relations
\begin{align*}
    \sigma_{X_1 \otimes X_2, Y}
    &= (\sigma_{X_1,Y} \otimes 1_{X_2}) \circ (1_{X_1} \otimes \sigma_{X_2,Y}),
    &
    \sigma_{X_2,Y} \circ (f \otimes 1_Y)
    &= (1_Y \otimes f) \circ \sigma_{X_1,Y},
    \\
    \sigma_{X, Y_1 \otimes Y_2}
    &= (1_{Y_1} \otimes \sigma_{X,Y_2}) \circ (\sigma_{X, Y_1} \otimes 1_{Y_2}),
    &
    \sigma_{X,Y_2} \circ (1_X \otimes g) &= (g \otimes 1_X)\circ \sigma_{X,Y_1},
\end{align*}
for all $X, X_1,X_2 \in\blue{\Heis_l(A;z,u)}$, $Y, Y_1,Y_2 \in \red{\Heis_m(A;z,v)}$ and $f \colon X_1 \rightarrow X_2$, $g \colon Y_1\rightarrow Y_2$.  Morphisms in $\blue{\Heis_l(A;z,u)} \odot \red{\Heis_m(A;z,v)}$ are linear combinations of diagrams colored both blue and red. In these diagrams,
as well as the generating morphisms of $\blue{\Heis_l(A;z,u)}$ and $\red{\Heis_m(A;z,v)}$, we have the additional two-color crossings
\[

\]
which represent the isomorphisms $\sigma_{X,Y}$ for $X \in \{\upblue, \downblue\}$ and $Y \in \{\upred, \downred\}$, and their inverses
\[
    %
    \ .
\]
We also have the symmetric product $\blue{\QAW(A;z)} \odot \red{\QAW(A;z)}$, defined in an analogous manner, where we now only have the upward two-colored crossings.

Let $\blue{\Heis_l(A;z,u)} \barodot \red{\Heis_m(A;z,v)}$ be the strict $\kk$-linear monoidal category obtained from $\blue{\Heis_l(A;z,u)} \odot \red{\Heis_m(A;z,v)}$ by localizing at the spear
%
.
This means that we adjoin a two-sided inverse to this morphism, which we denote again by a dart:
\begin{equation} \label{chromedart}
    %
    \right)^{-1}.
\end{equation}
As with the dumbbell defined in \cite[\S4,\S5]{BSW-K0}, all morphisms whose string diagram is that of an identity morphism with a spear joining two points of different colors are also automatically invertible in the localized category.  We also denote the inverses of such morphisms by using a dart in place of the spear. For instance:
\[
    %
    \right)^{-1}.
\]
We will use two-colored box dumbbells and box darts, defined in the obvious way.

We have the following relations to move darts past one-color crossings:
\begin{align} \label{dartslide1}
    %
    \ .
\end{align}
The same relations hold if we interchange the colors of the strands.

\begin{theo} \label{marble}
    There is a strict $\kk$-linear monoidal superfunctor
    \[
        \Delta \colon \invQAW(A;z) \to \Add \left( \blue{\QAW(A;z)} \barodot \red{\QAW(A;z)} \right)
    \]
    such that $\uparrow \mapsto \upblue \oplus \upred$, $\downarrow \mapsto \downblue \oplus \downred$, and on morphisms
    \begin{equation} \label{snowball}
        %
        \ .
    \end{equation}
\end{theo}

\begin{proof}
    It is straightforward to verify that the defining relations from \cref{invQAWpres} are preserved.  For example
    \[
        \Delta
        \left(
            %
        \right)
        \ . \qedhere
    \]
\end{proof}

\begin{cor} \label{racing}
    There is a strict $\kk$-linear monoidal superfunctor
    \[
        \Delta \colon \QAW(A;z) \to \Add \left( \blue{\QAW(A;z)} \barodot \red{\QAW(A;z)} \right)
    \]
    such that $\uparrow \mapsto \upblue \oplus \upred$, $\downarrow \mapsto \downblue \oplus \downred$, and on morphisms
    \begin{gather} \label{pitlane}
        %
        \ .
    \end{gather}
\end{cor}

\begin{proof}
    We compose the canonical superfunctor $\QAW(A;z) \to \invQAW(A;z)$ with the superfunctor of \cref{marble} and compute
    \begin{multline*}
        \Delta
        \left(
            %
        \ .
    \end{multline*}
    The computation of the image of the negative crossing is similar.
\end{proof}

\begin{rem}
    The categorical comultiplications of \cref{marble,racing} are not unique.  For example, we could place the box dumbbell on the other two-colored crossing in \cref{snowball}.  In addition, when $A = \kk$ and $z = q - q^{-1}$ for some $q \in \kk^\times$, we can factor the box dumbbell as
    \[
        %
        \right).
    \]
    Placing the first factor at the bottom of the first two-colored crossing in \cref{snowball} and the second factor at the top of the second two-colored crossing in \cref{snowball} (and removing the box dumbbell) yields the categorical comultiplication of \cite[Th.~8.9]{BSW-qheis}.
\end{rem}

Our goal is to extend the categorical comultiplication of \cref{racing} to Frobenius Heisenberg categories.  We will need the following morphisms, which we refer to as \emph{internal bubbles}:
\begin{align} \label{blueint}
    %
    \ .
\end{align}

The category $\blue{\Heis_l(A;z,u)} \barodot \red{\Heis_m(A;z,v)}$ has several useful symmetries.  Coming from \cref{om}, we have the strict $\kk$-linear monoidal isomorphism
\begin{equation}
    \Omega_{\blue{l}|\red{m}} \colon \blue{\Heis_l(A;z,u)} \barodot \red{\Heis_m(A;z,v)} \xrightarrow{\cong} \left( \blue{\Heis_{-l}(A;z,u^{-1})} \barodot \red{\Heis_{-m}(A;z,v^{-1})} \right)^\op,
\end{equation}
which reflects a diagram in a horizontal plane and multiplies by $(-1)^{c+d+\binom{y}{2}}$, where $c$ is the number of one-colored crossings, $d$ is the number of left cups and caps (including ones in $(+)$-, $(-)$-, and internal bubbles), and $y$ is the number of odd tokens.  We also have
\begin{equation}
    \flip \colon \blue{\Heis_l(A;z,u)} \barodot \red{\Heis_m(A;z,v)} \xrightarrow{\cong} \blue{\Heis_l(A;z,u)} \barodot \red{\Heis_m(A;z,v)}
\end{equation}
defined on diagrams by switching the colors blue and red.  Finally, the category $\blue{\Heis_l(A;z,u)} \barodot \red{\Heis_m(A;z,v)}$ is strictly pivotal, with duality superfunctor
\[
    * \colon \blue{\Heis_l(A;z,u)} \barodot \red{\Heis_m(A;z,v)} \xrightarrow{\cong} \left( \left( \blue{\Heis_l(A;z,u)} \barodot \red{\Heis_m(A;z,v)} \right)^\op \right)^\rev
\]
defined, as in \cref{pivot}, by rotating diagrams through $180\degree$ and multiplying by $(-1)^{\binom{y}{2}}$, where $y$ is the number of odd tokens in the diagram.

We denote the duals of the internal bubbles \cref{blueint,redint} by
\[
    \begin{tikzpicture}[centerzero]
        \draw[red,<-] (0,-0.4) -- (0,0.4);
        \intleft{blue}{0,0};
    \end{tikzpicture}
    \ ,\qquad
    \begin{tikzpicture}[centerzero]
        \draw[red,<-] (0,-0.4) -- (0,0.4);
        \intright{blue}{0,0};
    \end{tikzpicture}
    \ ,\qquad
    \begin{tikzpicture}[centerzero]
        \draw[blue,<-] (0,-0.4) -- (0,0.4);
        \intleft{red}{0,0};
    \end{tikzpicture}
    \ ,\qquad
    \begin{tikzpicture}[centerzero]
        \draw[blue,<-] (0,-0.4) -- (0,0.4);
        \intright{red}{0,0};
    \end{tikzpicture}
    \ .
\]
It follows that internal bubbles slide past all cups and caps.

We now state the main result of the section, which is an extension of the categorical comultiplication of \cref{racing} to the quantum Frobenius Heisenberg category.  The proof of \cref{comult} will be based on a series of lemmas.

\begin{theo} \label{comult}
    For $k = \blue{l} + \red{m}$ and $t = \blue{u} \red{v}$, there is a unique strict $\kk$-linear monoidal superfunctor
    \[
        \cDelt \colon \Heis_k(A;z,t) \to \Add \left( \blue{\Heis_l(A;z,u)} \barodot \red{\Heis_m(A;z,v)} \right)
    \]
    such that $\uparrow \mapsto \upblue \oplus \upred$, $\downarrow \mapsto \downblue \oplus \downred$, and on morphisms by \cref{pitlane}, \cref{galaxy} and
    \begin{align} \label{hazers}

        \ .
    \end{align}
    Moreover, we have that
    \begin{align} \label{wisps}
        \cDelt
        \left(
            %
        \ .
    \end{align}
    Also, the following hold for all $n \in \Z$ and $a \in A$:
    \begin{align} \label{ducks}
        \cDelt
        \left(
            \plusccbubble{a}{n}
        \right)
        &= \sum_{r \in \Z}
        %
        \ .
    \end{align}
    Equivalently, in terms of the generating functions \cref{bubgen2,bubgen1,bubgen4,bubgen3} and their analogues in $\blue{\Heis_l(A;z,u)}$ and $\red{\Heis_m(A;z,v)}$:
    \begin{align} \label{catseye}
        \cDelt
        \left(
            %
        \right)^{-1}.
    $
\end{lem}

\begin{proof}
    We first compute
    \[
        %
        \overset{\cref{infgrass}}{=} -1.
    \end{multline*}
\end{proof}

\begin{lem} \label{slalom}
    For any $n \in \Z$ and $a \in A$, we have\footnote{There is an error in \cite[Lem.~8.3]{BSW-qheis}.  Both sums there should be over all $b \in \Z$, so that the statement matches the $A = \kk$ case of the current lemma.}
    \[
        %
        \ .
    \end{multline*}
\end{proof}

\begin{lem} \label{slippy}
    The following relations hold:
    \[
        %
        \ ,
    \]
    which is equal to the right-hand side of the first relation.  The second relation is proved similarly.
\end{proof}

\begin{lem} \label{funnel}
    We have
    \[
        %
        \ .
    \]
\end{lem}

\begin{proof}
    First note that, using \cref{funnel}, we have
    \[
        %
        \ ,
    \end{multline*}
    where, in the last equality, we used \cref{slalom}, \cref{tanks2}, and the fact that $t=uv$.
\end{proof}

\begin{lem} \label{curling}
    We have
    \[
        %
        \ .
    \]
\end{lem}

\begin{proof}
    First, a lengthy but straightforward computation using \cref{pos,neg,dartslide1,dartslide2,green} shows that
    \begin{equation} \label{cochon}
        %
        \ .
    \end{align*}
    Then, adding a counter-clockwise interior bubble to bottom of the blue strand completes the proof of the lemma.
\end{proof}

\begin{proof}[Proof of \cref{comult}]
    In light of the uniqueness from \cref{cherries}, we can take \cref{pitlane,galaxy,hazers,wisps} as the definition of $\cDelt$ on generating morphisms, and must check that the images of the relations \cref{tokrel,braid,skein,QAWC,rightadj,pos,neg,curls,morecurls} are all satisfied in $\Add \left( \blue{\Heis_l(A;z,u)} \barodot \red{\Heis_m(A;z,v)} \right)$.  We must also check \cref{ducks,jumpers}.  Relations \cref{tokrel,braid,skein,QAWC} follow from \cref{racing} and \cref{rightadj} is straightforward; so it remains to check the others.

    We first check \cref{ducks,jumpers}.  First assume $k \ge 0$.  Consider the bubble \pluscbubble{a}{n}.  When $n < 0$, both sides of the second relation in \cref{ducks} are zero by \cref{tanks} and our assumption that $k \ge 0$.  When $n = 0$, we have
    \[
        \cDelt \left( \pluscbubble{a}{n} \right)
        \overset{\cref{tanks}}{=} -\delta_{0,k} \frac{t^{-1}}{z} \tr(a) 1_\one
        = -\delta_{0,k} \frac{u^{-1} v^{-1}}{z} \tr(a) 1_\one
        \overset{\cref{tanks}}{\underset{\cref{adecomp}}{=}} -
        \begin{tikzpicture}[centerzero]
            \plusrightblank[blue]{-0.4,0};
            \plusrightblank[red]{0.4,0};
            \token[blue]{-0.6,0}{east}{a};
            \multdot[blue]{-0.4,-0.2}{north}{n};
            \multdot[red]{0.4,-0.2}{north}{n-r};
            \telecolor{blue}{red}{-0.2,0}{0.2,0};
        \end{tikzpicture}
        \ .
    \]
    Now assume $n>0$, so that $\pluscbubble{a}{n} = \cbubble{a}{n}$.  Then, using \cref{slalom}, we have
    \[
        \cDelt \left( \pluscbubble{a}{n} \right)
        = \cDelt \left( \cbubble{a}{n} \right)
        = -
        \begin{tikzpicture}[centerzero]
            \draw[blue,->] (0,0.3) arc(90:-270:0.3);
            \intright{red}{-0.212,-0.212};
            \multdot[blue]{-0.212,0.212}{east}{n};
            \token[blue]{0.3,0}{west}{a};
        \end{tikzpicture}
        -
        \begin{tikzpicture}[centerzero]
            \draw[red,->] (0,0.3) arc(90:-270:0.3);
            \intright{blue}{-0.212,-0.212};
            \multdot[red]{-0.212,0.212}{east}{n};
            \token[red]{0.3,0}{west}{a};
        \end{tikzpicture}
        =
        - \sum_{r \in \Z}
        \begin{tikzpicture}[centerzero]
            \plusrightblank[blue]{-0.4,0};
            \plusrightblank[red]{0.4,0};
            \token[blue]{-0.6,0}{east}{a};
            \multdot[blue]{-0.4,-0.2}{north}{r};
            \multdot[red]{0.4,-0.2}{north}{n-r};
            \telecolor{blue}{red}{-0.2,0}{0.2,0};
        \end{tikzpicture}
        \ .
    \]
    This establishes the right-hand identity in \cref{ducks}, hence the right-hand identity in \cref{catseye}.   Next we compute
    \begin{gather*}
        \cDelt
        \left(
            \begin{tikzpicture}[centerzero]
                \plusgenleft{0,0};
                \token{-0.2,0}{east}{ab};
            \end{tikzpicture}
        \right)
        \circ
        \left(
            \begin{tikzpicture}[centerzero]
                \plusgenright[blue]{-0.4,0};
                \plusgenright[red]{0.4,0};
                \telecolor{blue}{red}{-0.2,0}{0.2,0};
                \token[blue]{-0.6,0}{east}{b^\vee};
            \end{tikzpicture}
        \right)
        =
        z \cDelt
        \left(
            \begin{tikzpicture}[centerzero]
                \plusgenleft{0,0};
                \token{-0.2,0}{east}{ab};
            \end{tikzpicture}
        \right)
        \circ \cDelt
        \left(
            \begin{tikzpicture}[centerzero]
                \plusgenright{0,0};
                \token{-0.2,0}{east}{b^\vee};
            \end{tikzpicture}
        \right)
        =
        \cDelt
        \left(
            \begin{tikzpicture}[centerzero]
                \plusgenright{-0.4,0};
                \plusgenright{0.4,0};
                \teleport{-0.2,0}{0.2,0};
                \token{-0.6,0}{east}{a};
            \end{tikzpicture}
        \right)
        \overset{\cref{infweeds}}{=} z \tr(a) 1_\one,
        \\
        z^{-1}
        \begin{tikzpicture}[centerzero]
            \plusgenleft[blue]{-0.4,0};
            \plusgenleft[red]{0.4,0};
            \telecolor{blue}{red}{-0.2,0}{0.2,0};
            \token[blue]{-0.6,0}{east}{ab};
        \end{tikzpicture}
        \begin{tikzpicture}[centerzero]
            \plusgenright[blue]{-0.4,0};
            \plusgenright[red]{0.4,0};
            \telecolor{blue}{red}{-0.2,0}{0.2,0};
            \token[blue]{-0.6,0}{east}{b^\vee};
        \end{tikzpicture}
        \overset{\cref{laser}}{=} z^{-1}
        \begin{tikzpicture}[centerzero]
            \plusgenleft[blue]{-0.4,0.3};
            \plusgenright[blue]{0.4,0.3};
            \plusgenleft[red]{-0.4,-0.3};
            \plusgenright[red]{0.4,-0.3};
            \telecolor{blue}{blue}{-0.2,0.3}{0.2,0.3};
            \telecolor{red}{red}{-0.2,-0.3}{0.2,-0.3};
            \token[blue]{-0.6,0.3}{east}{ab};
            \token[red]{-0.6,-0.3}{east}{b^\vee};
        \end{tikzpicture}
        \overset{\cref{infweeds}}{\underset{\cref{adecomp}}{=}}
        z \tr(a) 1_\one.
    \end{gather*}
    Now we observe that there is a unique morphism in $f(a;w) \in \End_{\blue{\Heis_l(A;z,u)} \barodot \red{\Heis_m(A;z,v)}}(\one) \llbracket w^{-1} \rrbracket$ such that
    $
        f(ab;w)
        \left(
            \begin{tikzpicture}[centerzero]
                \plusgenright[blue]{-0.4,0};
                \plusgenright[red]{0.4,0};
                \telecolor{blue}{red}{-0.2,0}{0.2,0};
                \token[blue]{-0.6,0}{east}{b^\vee};
            \end{tikzpicture}
        \right)
        = z \tr(a) 1_\one;
    $
    see the proof of \cite[Lem.~7.1]{BSW-foundations} for a similar situation.  It follows that
    \[
        f(a;w) =
        \cDelt
        \left(
            \begin{tikzpicture}[centerzero]
                \plusgenleft{0,0};
                \token{-0.2,0}{east}{a};
            \end{tikzpicture}
        \right)
        = z^{-1}
        \begin{tikzpicture}[centerzero]
            \plusgenleft[blue]{-0.4,0};
            \plusgenleft[red]{0.4,0};
            \telecolor{blue}{red}{-0.2,0}{0.2,0};
            \token[blue]{-0.6,0}{east}{a};
        \end{tikzpicture}
        \ .
    \]
    This establishes the left-hand identity in \cref{catseye}, hence the left-identity in \cref{ducks}.  Next consider the bubble \plusccbubble{a}{n}.  When $n>0$, both sides of the second relation in \cref{jumpers} are zero by \cref{tanks2}, and this relation is straightforward to verify when $n=0$, again using \cref{tanks2}.  Now assume $n<0$, so that $\minuscbubble{a}{n} = \cbubble{a}{n}$.  Then, using \cref{slalom}, we have
    \[
        \cDelt \left( \minuscbubble{a}{n} \right)
        = \cDelt \left( \cbubble{a}{n} \right)
        = -
        \begin{tikzpicture}[centerzero]
            \draw[blue,->] (0,0.3) arc(90:-270:0.3);
            \intright{red}{-0.212,-0.212};
            \multdot[blue]{-0.212,0.212}{east}{n};
            \token[blue]{0.3,0}{west}{a};
        \end{tikzpicture}
        -
        \begin{tikzpicture}[centerzero]
            \draw[red,->] (0,0.3) arc(90:-270:0.3);
            \intright{blue}{-0.212,-0.212};
            \multdot[red]{-0.212,0.212}{east}{n};
            \token[red]{0.3,0}{west}{a};
        \end{tikzpicture}
        =
        \sum_{r \in \Z}
        \begin{tikzpicture}[centerzero]
            \minusrightblank[blue]{-0.4,0};
            \minusrightblank[red]{0.4,0};
            \token[blue]{-0.6,0}{east}{a};
            \multdot[blue]{-0.4,-0.2}{north}{r};
            \multdot[red]{0.4,-0.2}{north}{n-r};
            \telecolor{blue}{red}{-0.2,0}{0.2,0};
        \end{tikzpicture}
        \ .
    \]
    Then we complete the proof of \cref{jumpers} as above, using \cref{infweeds,mellow}.  It remains to treat the case $k < 0$.  This follows by similar arguments; one starts by considering the counterclockwise $(+)$- and $(-)$-bubbles using the identities obtained by applying $\Omega_{\blue{l}|\red{m}}$ to \cref{slalom}, then gets the clockwise ones using \cref{infweeds}.

    Now consider \cref{curls,morecurls}.  The relations involving bubbles follow easily from \cref{ducks,jumpers}.
    \details{
        If $0 \le n \le k$, we have
        \begin{multline*}
            \cDelt \left( \cbubble{a}{n} \right)
            = \cDelt \left( \pluscbubble{a}{n} \right) + \cDelt \left( \minuscbubble{a}{n} \right)
            \overset{\cref{ducks}}{\underset{\cref{jumpers}}{=}}=  \sum_{r \in \Z}
            \begin{tikzpicture}[centerzero]
                \minusrightblank[blue]{-0.4,0};
                \minusrightblank[red]{0.4,0};
                \token[blue]{-0.6,0}{east}{a};
                \multdot[blue]{-0.4,-0.2}{north}{r};
                \multdot[red]{0.4,-0.2}{north}{n-r};
                \telecolor{blue}{red}{-0.2,0}{0.2,0};
            \end{tikzpicture}
            - \sum_{r \in \Z}
            \begin{tikzpicture}[centerzero]
                \plusrightblank[blue]{-0.4,0};
                \plusrightblank[red]{0.4,0};
                \token[blue]{-0.6,0}{east}{a};
                \multdot[blue]{-0.4,-0.2}{north}{r};
                \multdot[red]{0.4,-0.2}{north}{n-r};
                \telecolor{blue}{red}{-0.2,0}{0.2,0};
            \end{tikzpicture}
            \\
            = \sum_{r \in Z} \delta_{r,0} \delta_{n-r,0} \frac{uv}{z} \tr(a) 1_\one - \sum_{r \in \Z} \delta_{r,l} \delta_{n-r,m} \frac{u^{-1}v^{-1}}{z} \tr(a) 1_\one
            = \delta_{n,0} \frac{t}{z} \tr(a) 1_\one - \delta_{n,k} \frac{t^{-1}}{z} \tr(a) 1_\one.
        \end{multline*}
    }
    Next consider the right curl relation in \cref{curls}, so $k \ge 0$.  Considering the image of this relation under $\cDelt$, we must show that
    \begin{equation} \label{momay}
        -
        \begin{tikzpicture}[centerzero]
            \draw[blue] (-0.2,-0.5) -- (-0.2,-0.35) to[out=up,in=west] (0.05,0.2) to[out=right,in=up] (0.2,0);
            \draw[wipe] (0.2,0) to[out=down,in=east] (0.05,-0.2) to[out=left,in=down] (-0.2,0.35) -- (-0.2,0.5);
            \draw[blue,->] (0.2,0) to[out=down,in=east] (0.05,-0.2) to[out=left,in=down] (-0.2,0.35) -- (-0.2,0.5);
            \intrightsm{red}{0.2,0};
        \end{tikzpicture}
        -
        \begin{tikzpicture}[centerzero]
            \draw[red] (-0.2,-0.5) -- (-0.2,-0.35) to[out=up,in=west] (0.05,0.2) to[out=right,in=up] (0.2,0);
            \draw[wipe] (0.2,0) to[out=down,in=east] (0.05,-0.2) to[out=left,in=down] (-0.2,0.35) -- (-0.2,0.5);
            \draw[red,->] (0.2,0) to[out=down,in=east] (0.05,-0.2) to[out=left,in=down] (-0.2,0.35) -- (-0.2,0.5);
            \intrightsm{blue}{0.2,0};
        \end{tikzpicture}
        +
        \begin{tikzpicture}[centerzero]
            \draw[blue,->] (0,-0.5) -- (0,0.5);
            \draw[red,->-=0.1] (0.5,0.2) arc(90:-270:0.2);
            \telecolor{blue}{red}{0,0.141}{0.359,0.141};
            \dart{red}{blue}{0.359,-0.141}{0,-0.141};
            \intrightsm{blue}{0.7,0};
        \end{tikzpicture}
        +
        \begin{tikzpicture}[centerzero]
            \draw[red,->] (0,-0.5) -- (0,0.5);
            \draw[blue,->-=0.1] (0.5,0.2) arc(90:-270:0.2);
            \telecolor{red}{blue}{0,0.141}{0.359,0.141};
            \dart{blue}{red}{0.359,-0.141}{0,-0.141};
            \intrightsm{red}{0.7,0};
        \end{tikzpicture}
        =
        \delta_{k,0} t^{-1}
        \begin{tikzpicture}[centerzero]
            \draw[blue,->] (0,-0.5) -- (0,0.5);
        \end{tikzpicture}
        + \delta_{k,0} t^{-1}
        \begin{tikzpicture}[centerzero]
            \draw[red,->] (0,-0.5) -- (0,0.5);
        \end{tikzpicture}
        \ .
    \end{equation}
    This follows from \cref{funnel}, together with its image under $\flip$, noting that the only nonzero term in the summation on the right-hand side of that identity is the one with $r=s=0$, due to \cref{tanks} and the assumption that $k \ge 0$.   For the left curl relation in \cref{morecurls}, we see, after applying $\cDelt$, that we must show, for $k \le 0$, that
    \[
        \begin{tikzpicture}[centerzero]
            \draw[blue] (0.2,-0.5) -- (0.2,-0.35) to[out=up,in=east] (-0.05,0.2) to[out=left,in=up] (-0.2,0);
            \draw[wipe] (-0.2,0) to[out=down,in=west] (-0.05,-0.2) to[out=right,in=down] (0.2,0.35) -- (0.2,0.5);
            \draw[blue,->] (-0.2,0) to[out=down,in=west] (-0.05,-0.2) to[out=right,in=down] (0.2,0.35) -- (0.2,0.5);
            \intleftsm{red}{-0.2,0};
        \end{tikzpicture}
        +
        \begin{tikzpicture}[centerzero]
            \draw[red] (0.2,-0.5) -- (0.2,-0.35) to[out=up,in=east] (-0.05,0.2) to[out=left,in=up] (-0.2,0);
            \draw[wipe] (-0.2,0) to[out=down,in=west] (-0.05,-0.2) to[out=right,in=down] (0.2,0.35) -- (0.2,0.5);
            \draw[red,->] (-0.2,0) to[out=down,in=west] (-0.05,-0.2) to[out=right,in=down] (0.2,0.35) -- (0.2,0.5);
            \intleftsm{blue}{-0.2,0};
        \end{tikzpicture}
        +
        \begin{tikzpicture}[centerzero]
            \draw[blue,->] (0,-0.4) -- (0,0.4);
            \draw[red,->-=0.1] (-0.5,0.2) arc(90:450:0.2);
            \telecolor{red}{blue}{-0.359,0.141}{0,0.141};
            \dart{red}{blue}{-0.359,-0.141}{0,-0.141};
            \intleftsm{blue}{-0.7,0};
        \end{tikzpicture}
        +
        \begin{tikzpicture}[centerzero]
            \draw[red,->] (0,-0.4) -- (0,0.4);
            \draw[blue,->-=0.1] (-0.5,0.2) arc(90:450:0.2);
            \telecolor{blue}{red}{-0.359,0.141}{0,0.141};
            \dart{blue}{red}{-0.359,-0.141}{0,-0.141};
            \intleftsm{red}{-0.7,0};
        \end{tikzpicture}
        = \delta_{k,0} t\
        \begin{tikzpicture}[centerzero]
            \draw[blue,->] (0,-0.5) -- (0,0.5);
        \end{tikzpicture}
        + \delta_{k,0} t\
        \begin{tikzpicture}[centerzero]
            \draw[red,->] (0,-0.5) -- (0,0.5);
        \end{tikzpicture}
        \ .
    \]
    This follows from \cref{momay} for $\blue{\Heis_{-l}(A;z,u^{-1})} \barodot \red{\Heis_{-m}(A;z,v^{-1})}$ after applying $* \circ \Omega_{\blue{(-l)} | \red{(-m)}}$.

    Now consider \cref{neg}.  By definitions \cref{rcross,cold}, we see that
    \begin{align*}
        \cDelt
        \left(
            \begin{tikzpicture}[centerzero]
                \draw[->] (-0.3,-0.3) -- (0.3,0.3);
                \draw[wipe] (0.3,-0.3) -- (-0.3,0.3);
                \draw[<-] (0.3,-0.3) -- (-0.3,0.3);
            \end{tikzpicture}
        \right)
        &=
        \begin{tikzpicture}[centerzero]
            \draw[blue,->] (-0.3,-0.3) -- (0.3,0.3);
            \draw[wipe] (0.3,-0.3) -- (-0.3,0.3);
            \draw[blue,<-] (0.3,-0.3) -- (-0.3,0.3);
        \end{tikzpicture}
        +
        \begin{tikzpicture}[centerzero]
            \draw[red,->] (-0.3,-0.3) -- (0.3,0.3);
            \draw[wipe] (0.3,-0.3) -- (-0.3,0.3);
            \draw[red,<-] (0.3,-0.3) -- (-0.3,0.3);
        \end{tikzpicture}
        +
        \begin{tikzpicture}[centerzero]
            \draw[red,->] (-0.3,-0.5) -- (0.3,0.5);
            \draw[blue,<-] (0.3,-0.5) -- (-0.3,0.5);
            \dart{blue}{red}{0.15,-0.25}{-0.15,-0.25};
        \end{tikzpicture}
        +
        \begin{tikzpicture}[centerzero]
            \draw[blue,->] (-0.3,-0.5) -- (0.3,0.5);
            \draw[red,<-] (0.3,-0.5) -- (-0.3,0.5);
            \boxdart{red}{blue}{0.15,-0.25}{-0.15,-0.25};
        \end{tikzpicture}
        +
        \begin{tikzpicture}[centerzero]
            \draw[blue,->] (-0.2,0.55) -- (-0.2,0.35) arc(180:360:0.2) -- (0.2,0.55);
            \draw[red,->] (-0.2,-0.55) -- (-0.2,-0.35) arc(180:0:0.2) -- (0.2,-0.55);
            \telecolor{blue}{red}{0.141,0.209}{0.141,-0.209};
            \dart{blue}{red}{-0.141,0.209}{-0.141,-0.209};
        \end{tikzpicture}
        +
        \begin{tikzpicture}[centerzero]
            \draw[red,->] (-0.2,0.55) -- (-0.2,0.35) arc(180:360:0.2) -- (0.2,0.55);
            \draw[blue,->] (-0.2,-0.55) -- (-0.2,-0.35) arc(180:0:0.2) -- (0.2,-0.55);
            \telecolor{red}{blue}{0.141,0.209}{0.141,-0.209};
            \dart{red}{blue}{-0.141,0.209}{-0.141,-0.209};
        \end{tikzpicture}
        \ ,
        \\
        \cDelt
        \left(
            \begin{tikzpicture}[centerzero]
                \draw[<-] (-0.3,-0.3) -- (0.3,0.3);
                \draw[wipe] (0.3,-0.3) -- (-0.3,0.3);
                \draw[->] (0.3,-0.3) -- (-0.3,0.3);
            \end{tikzpicture}
        \right)
        &= -
        \begin{tikzpicture}[centerzero]
            \draw[blue,<-] (-0.3,-0.5) -- (0.3,0.5);
            \draw[wipe] (0.3,-0.5) -- (-0.3,0.5);
            \draw[blue,->] (0.3,-0.5) -- (-0.3,0.5);
            \intrightsm{red}{0.15,0.25};
            \intleftsm{red}{-0.15,-0.25};
        \end{tikzpicture}
        -
        \begin{tikzpicture}[centerzero]
            \draw[red,<-] (-0.3,-0.5) -- (0.3,0.5);
            \draw[wipe] (0.3,-0.5) -- (-0.3,0.5);
            \draw[red,->] (0.3,-0.5) -- (-0.3,0.5);
            \intrightsm{blue}{0.15,0.25};
            \intleftsm{blue}{-0.15,-0.25};
        \end{tikzpicture}
        -
        \begin{tikzpicture}[centerzero]
            \draw[red,<-] (-0.3,-0.6) -- (0.3,0.6);
            \draw[blue,->] (0.3,-0.6) -- (-0.3,0.6);
            \intrightsm{blue}{0.15,0.3};
            \intleftsm{blue}{-0.1,-0.2};
            \dart{blue}{red}{0.2,-0.4}{-0.2,-0.4};
        \end{tikzpicture}
        -
        \begin{tikzpicture}[centerzero]
            \draw[blue,<-] (-0.3,-0.6) -- (0.3,0.6);
            \draw[red,->] (0.3,-0.6) -- (-0.3,0.6);
            \intrightsm{red}{0.15,0.3};
            \intleftsm{red}{-0.1,-0.2};
            \boxdart{red}{blue}{0.2,-0.4}{-0.2,-0.4};
        \end{tikzpicture}
        -
        \begin{tikzpicture}[centerzero]
            \draw[red,<-] (-0.2,0.6) -- (-0.2,0.35) arc(180:360:0.2) -- (0.2,0.6);
            \draw[blue,<-] (-0.2,-0.6) -- (-0.2,-0.35) arc(180:0:0.2) -- (0.2,-0.6);
            \telecolor{red}{blue}{0.141,0.209}{0.141,-0.209};
            \dart{blue}{red}{-0.141,-0.209}{-0.141,0.209};
            \intrightsm{blue}{0.2,0.42};
            \intleftsm{red}{0.2,-0.42};
        \end{tikzpicture}
        -
        \begin{tikzpicture}[centerzero]
            \draw[blue,<-] (-0.2,0.6) -- (-0.2,0.35) arc(180:360:0.2) -- (0.2,0.6);
            \draw[red,<-] (-0.2,-0.6) -- (-0.2,-0.35) arc(180:0:0.2) -- (0.2,-0.6);
            \telecolor{blue}{red}{0.141,0.209}{0.141,-0.209};
            \dart{red}{blue}{-0.141,-0.209}{-0.141,0.209};
            \intrightsm{red}{0.2,0.42};
            \intleftsm{blue}{0.2,-0.42};
        \end{tikzpicture}
        \ .
    \end{align*}
    We can then compute the image under $\cDelt$ of both sides of \cref{neg}.  Comparing the matrix entries of the resulting morphisms, we are reduced to verifying the identities
    \begin{align*}
        -
        \begin{tikzpicture}[anchorbase]
            \draw[blue,->] (-0.2,0) \braidup (0.2,0.5);
            \draw[wipe] (-0.2,-0.5) \braidup (0.2,0) \braidup (-0.2,0.5);
            \draw[blue,<-] (-0.2,-0.8) -- (-0.2,-0.5) \braidup (0.2,0) \braidup (-0.2,0.5);
            \draw[wipe] (0.2,-0.5) \braidup (-0.2,0);
            \draw[blue] (0.2,-0.8) -- (0.2,-0.5) \braidup (-0.2,0);
            \intrightsm{red}{0.2,0};
            \intleftsm{red}{-0.2,-0.5};
        \end{tikzpicture}
        -
        \begin{tikzpicture}[anchorbase]
            \draw[blue,->] (-0.2,0.8) -- (-0.2,0.3) arc(180:360:0.2) -- (0.2,0.8);
            \draw[blue,<-] (-0.2,-0.8) -- (-0.2,-0.3) arc(180:0:0.2) -- (0.2,-0.8);
            \draw[red,->] (0.7,0.7) arc(90:0:0.2) -- (0.9,-0.5) arc(360:180:0.2) -- (0.5,0.5) arc(180:80:0.2);
            \telecolor{blue}{red}{0.2,0.5}{0.5,0.5};
            \dart{blue}{red}{0.2,0.3}{0.5,0.3};
            \telecolor{blue}{red}{0.2,-0.5}{0.5,-0.5};
            \dart{blue}{red}{0.2,-0.3}{0.5,-0.3};
            \intrightsm{blue}{0.9,0};
            \intleftsm{red}{-0.2,-0.4};
        \end{tikzpicture}
        &=
        \begin{tikzpicture}[centerzero]
            \draw[blue,<-] (-0.2,-0.65) -- (-0.2,0.65);
            \draw[blue,->] (0.2,-0.65) -- (0.2,0.65);
        \end{tikzpicture}
        + t
        \begin{tikzpicture}[anchorbase]
            \draw[blue,->] (-0.2,0.5) -- (-0.2,0.35) arc(180:360:0.2) -- (0.2,0.5);
            \draw[blue,<-] (-0.2,-0.7) -- (-0.2,-0.35) arc(180:0:0.2) -- (0.2,-0.7);
            \telecolor{blue}{blue}{0,0.15}{0,-0.15};
            \intleftsm{red}{-0.2,-0.4};
        \end{tikzpicture}
        - \sum_{\substack{r,s > 0 \\ p \in \Z}}
        \begin{tikzpicture}[anchorbase]
            \draw[blue,->] (-0.2,0.6) -- (-0.2,0.3) arc(180:360:0.2) -- (0.2,0.6);
            \draw[blue,<-] (-0.2,-0.7) -- (-0.2,-0.3) arc(180:0:0.2) -- (0.2,-0.7);
            \plusrightblank[blue]{0.7,0.3};
            \draw[red,->] (0.7,-0.6) arc(270:-90:0.2);
            \node at (0.7,-0.4) {{\color{red} \dotlabel{+}}};
            \telecolor{blue}{blue}{0.2,0.3}{0.5,0.3};
            \telecolor{blue}{red}{0.2,-0.4}{0.5,-0.4};
            \telecolor{red}{blue}{0.7,-0.2}{0.7,0.1};
            \intleftsm{red}{-0.2,-0.45};
            \multdot[blue]{-0.2,0.3}{east}{r};
            \multdot[blue]{-0.141,-0.159}{east}{s};
            \multdot[blue]{0.9,0.3}{west}{p};
            \multdot[red]{0.9,-0.4}{west}{-r-s-p};
        \end{tikzpicture}
        ,\\
        -
        \begin{tikzpicture}[anchorbase]
            \draw[red,<-] (-0.2,-1.1) -- (-0.2,-0.8) \braidup (0.2,-0.4) -- (0.2,-0.2) arc(0:180:0.2);
            \draw[wipe] (0.2,-1.1) -- (0.2,-0.8) \braidup (-0.2,-0.4) -- (-0.2,-0.2);
            \draw[red] (0.2,-1.1) -- (0.2,-0.8) \braidup (-0.2,-0.4) -- (-0.2,-0.2);
            \draw[blue,->] (-0.2,0.6) -- (-0.2,0.5) arc(180:360:0.2) -- (0.2,0.6);
            \dart{blue}{red}{-0.141,0.359}{-0.141,-0.059};
            \telecolor{blue}{red}{0.141,0.359}{0.141,-0.059};
            \intrightsm{blue}{0.2,-0.35};
            \intleftsm{blue}{-0.2,-0.8};
        \end{tikzpicture}
        -
        \begin{tikzpicture}[anchorbase]
            \draw[red,<-] (-0.2,-0.8) -- (-0.2,-0.5) arc(180:0:0.2) -- (0.2,-0.8);
            \draw[blue,<-] (0.2,0.9) -- (0.2,0.8) \braiddown (-0.2,0.4) -- (-0.2,0.2) arc(180:360:0.2);
            \draw[wipe] (0.2,0.2) -- (0.2,0.4) \braidup (-0.2,0.8) -- (-0.2,0.9);
            \draw[blue] (0.2,0.2) -- (0.2,0.4) \braidup (-0.2,0.8) -- (-0.2,0.9);
            \intleftsm{blue}{0.2,-0.6};
            \telecolor{red}{blue}{0.141,-0.359}{0.141,0.059};
            \dart{red}{blue}{-0.141,-0.359}{-0.141,0.059};
            \intrightsm{red}{0.2,0.3};
        \end{tikzpicture}
        &= t
        \begin{tikzpicture}[anchorbase]
            \draw[blue,->] (-0.2,0.5) -- (-0.2,0.35) arc(180:360:0.2) -- (0.2,0.5);
            \draw[red,<-] (-0.2,-0.7) -- (-0.2,-0.35) arc(180:0:0.2) -- (0.2,-0.7);
            \telecolor{blue}{red}{0,0.15}{0,-0.15};
            \intleftsm{blue}{0.2,-0.4};
        \end{tikzpicture}
         - \sum_{\substack{r,s > 0 \\ p \in \Z}}
        \begin{tikzpicture}[anchorbase]
            \draw[blue,->] (-0.2,0.6) -- (-0.2,0.3) arc(180:360:0.2) -- (0.2,0.6);
            \draw[red,<-] (-0.2,-0.7) -- (-0.2,-0.3) arc(180:0:0.2) -- (0.2,-0.7);
            \plusrightblank[blue]{0.7,0.3};
            \draw[red,->] (0.7,-0.6) arc(270:-90:0.2);
            \node at (0.7,-0.4) {{\color{red} \dotlabel{+}}};
            \telecolor{blue}{blue}{0.2,0.3}{0.5,0.3};
            \telecolor{blue}{red}{0.2,-0.4}{0.5,-0.4};
            \telecolor{red}{blue}{0.7,-0.2}{0.7,0.1};
            \intleftsm{blue}{-0.2,-0.45};
            \multdot[blue]{-0.2,0.3}{east}{r};
            \multdot[red]{-0.141,-0.159}{east}{s};
            \multdot[blue]{0.9,0.3}{west}{p};
            \multdot[red]{0.9,-0.4}{west}{-r-s-p};
        \end{tikzpicture}
        ,
        \\
        -\
        \begin{tikzpicture}[anchorbase]
            \draw[blue,<-] (-0.2,-0.1) -- (-0.2,0.4) \braidup (0.2,0.8) -- (0.2,1.1) \braidup (-0.2,1.5) -- (-0.2,1.6);
            \draw[red,->] (0.2,-0.1) -- (0.2,0.4) \braidup (-0.2,0.8) -- (-0.2,1.1) \braidup (0.2,1.5) -- (0.2,1.6);
            \boxdart{red}{blue}{0.2,0.1}{-0.2,0.1};
            \intleftsm{red}{-0.2,0.4};
            \intrightsm{red}{0.2,0.8};
            \dart{blue}{red}{0.2,1.1}{-0.2,1.1};
        \end{tikzpicture}
        &=
        \begin{tikzpicture}[anchorbase]
            \draw[blue,<-] (-0.2,-0.1) -- (-0.2,1.6);
            \draw[red,->] (0.2,-0.1) -- (0.2,1.6);
        \end{tikzpicture}
        \ ,\qquad -\
        \begin{tikzpicture}[anchorbase]
            \draw[red,<-] (-0.2,-0.1) -- (-0.2,0.4) \braidup (0.2,0.8) -- (0.2,1.1) \braidup (-0.2,1.5) -- (-0.2,1.6);
            \draw[blue,->] (0.2,-0.1) -- (0.2,0.4) \braidup (-0.2,0.8) -- (-0.2,1.1) \braidup (0.2,1.5) -- (0.2,1.6);
            \dart{blue}{red}{0.2,0.1}{-0.2,0.1};
            \intleftsm{blue}{-0.2,0.4};
            \intrightsm{blue}{0.2,0.8};
            \boxdart{red}{blue}{0.2,1.1}{-0.2,1.1};
        \end{tikzpicture}
        =
        \begin{tikzpicture}[anchorbase]
            \draw[red,<-] (-0.2,-0.1) -- (-0.2,1.6);
            \draw[blue,->] (0.2,-0.1) -- (0.2,1.6);
        \end{tikzpicture}
        \ ,
    \end{align*}
    plus the images of the first two under the symmetry $\flip$.  To prove the first two identities, simplify them by multiplying on the bottom of the left string by a clockwise internal bubble and using \cref{dubs}.  The resulting identities then follow from \cref{bobsled,biathlon}.  For the third identity, composing on the top of both sides of the identity by
    $
        \begin{tikzpicture}[centerzero]
            \draw[blue,<-] (-0.2,-0.2) -- (0.2,0.2);
            \draw[red,->] (0.2,-0.2) -- (-0.2,0.2);
        \end{tikzpicture},
    $
    then multiplying on the top of the blue strands by a counter-clockwise internal bubble, on the bottom of the blue strands by a clockwise internal bubble, and finally using \cref{dubs,boxdart}, we see that we must show that
    \[
        -\
        \begin{tikzpicture}[centerzero]
            \draw[red,->] (0.5,-0.6) -- (-0.5,0.6);
            \draw[blue,<-] (-0.5,-0.6) -- (0.5,0.6);
        \end{tikzpicture}
        -\
        \begin{tikzpicture}[centerzero]
            \draw[red,->] (0.5,-0.6) -- (-0.5,0.6);
            \draw[blue,<-] (-0.5,-0.6) -- (0.5,0.6);
            \telecolor{red}{blue}{-0.333,0.4}{0.333,0.4};
            \telecolor{red}{blue}{-0.167,0.2}{0.167,0.2};
            \dart{blue}{red}{-0.167,-0.2}{0.167,-0.2};
            \dart{red}{blue}{0.333,-0.4}{-0.333,-0.4};
        \end{tikzpicture}
        =
        \begin{tikzpicture}[centerzero]
            \draw[red,->] (0.5,-0.6) -- (-0.5,0.6);
            \draw[blue,<-] (-0.5,-0.6) -- (0.5,0.6);
            \intleft{red}{0.25,0.3};
            \intright{red}{-0.25,-0.3};
        \end{tikzpicture}
        \ .
    \]
    This follows from \cref{curling}. The proof of the fourth identity is similar.
\end{proof}

\section{Generalized cyclotomic quotients\label{sec:GCQ}}

We now construct some important strict module supercategories over $\Heis_k(A;z,t)$ known as generalized cyclotomic quotients.


We fix a supercommutative superalgebra $R$ over our usual ground ring $\kk$.  Defining the $R$-superalgebra $A_R := R \otimes_\kk A$, we can extend the trace map of $A$ to an $R$-linear map $\tr_R = \id \otimes \tr \colon A_R \to R$.  We can then also base change to obtain $\Heis_k(A_R;z,t) = R \otimes_\kk \Heis_k(A;z,t)$, which is a strict $R$-linear monoidal supercategory.  Since scalars in $R$ are not necessarily even, we must take extra care with the potential additional signs.  In diagrams for morphisms in $\Heis_k(A_R;z,t)$, we can label tokens by elements of $A_R$.  Recall also the generating function formalism introduced in \cref{sec:third}.

\begin{lem}
    For a polynomial $p(w) \in A_R[w]$, we have
    \begin{align} \label{road1}
        \begin{tikzpicture}[centerzero]
            \draw[->] (0,-0.4) to (0,0.4);
            \token{0,0}{west}{p(x)};
        \end{tikzpicture}
        &=
        \left[
            \begin{tikzpicture}[centerzero]
                \draw[->] (0,-0.4) to (0,0.4);
                \token{0,0.15}{east}{(w-x)^{-1}};
                \token{0,-0.15}{west}{p(w)};
            \end{tikzpicture}
        \right]_{w^{-1}},
        &
        \begin{tikzpicture}[centerzero]
            \draw[<-] (0,-0.4) to (0,0.4);
            \token{0,0}{west}{p(x)};
        \end{tikzpicture}
        &=
        \left[
            \begin{tikzpicture}[centerzero]
                \draw[<-] (0,-0.4) to (0,0.4);
                \token{0,0.15}{east}{(w-x)^{-1}};
                \token{0,-0.15}{west}{p(w)};
            \end{tikzpicture}
        \right]_{w^{-1}},
        \\ \label{road2}
        \begin{tikzpicture}[centerzero]
            \bubrightblank{0,0};
            \token{0.2,0}{west}{p(x)};
        \end{tikzpicture}
        &= tz^{-1} \tr(p(0)) 1_\one - t^{-1}z^{-1}
        \left[
            \begin{tikzpicture}[centerzero]
                \plusgenright{0,0};
                \token{0.2,0}{west}{p(w)};
            \end{tikzpicture}
        \right]_{w^0},
        &
        \begin{tikzpicture}[centerzero]
            \bubleftblank{0,0};
            \token{0.2,0}{west}{p(x)};
        \end{tikzpicture}
        &= tz^{-1}
        \left[
            \begin{tikzpicture}[centerzero]
                \plusgenleft{0,0};
                \token{0.2,0}{west}{p(w)};
            \end{tikzpicture}
        \right]_{w^0}
        - t^{-1}z^{-1} \tr(p(0)) 1_\one,
        \\ \label{road3}
        \begin{tikzpicture}[centerzero]
            \draw (-0.2,-0.5) -- (-0.2,-0.35) to[out=up,in=west] (0.05,0.2) to[out=right,in=up] (0.2,0);
            \draw[wipe] (0.2,0) to[out=down,in=east] (0.05,-0.2) to[out=left,in=down] (-0.2,0.35) -- (-0.2,0.5);
            \draw[->] (0.2,0) to[out=down,in=east] (0.05,-0.2) to[out=left,in=down] (-0.2,0.35) -- (-0.2,0.5);
            \token{0.2,0}{west}{p(x)};
        \end{tikzpicture}
       &= t^{-1}
        \left[
            \begin{tikzpicture}[centerzero]
                \draw[->] (0,-0.5) -- (0,0.5);
                \multdot{0,0.25}{east}{(w-x)^{-1}};
                \teleport{0,0}{0.3,0};
                \plusgenright{0.5,0};
                \token{0,-0.25}{east}{p(u)};
            \end{tikzpicture}
        \right]_{w^{-1}},
        &
        \begin{tikzpicture}[centerzero]
            \draw (0.2,0) to[out=down,in=east] (0.05,-0.2) to[out=left,in=down] (-0.2,0.35) -- (-0.2,0.5);
            \draw[wipe] (-0.2,-0.5) -- (-0.2,-0.35) to[out=up,in=west] (0.05,0.2) to[out=right,in=up] (0.2,0);
            \draw[<-] (-0.2,-0.5) -- (-0.2,-0.35) to[out=up,in=west] (0.05,0.2) to[out=right,in=up] (0.2,0);
            \token{0.2,0}{west}{p(x)};
        \end{tikzpicture}
        &= t
        \left[
            \begin{tikzpicture}[centerzero]
                \draw[->] (0,-0.5) -- (0,0.5);
                \multdot{0,0.25}{east}{(w-x)^{-1}};
                \teleport{0,0}{0.3,0};
                \plusgenleft{0.5,0};
                \token{0,-0.25}{east}{p(w)};
            \end{tikzpicture}
        \right]_{w^{-1}}.
    \end{align}
\end{lem}

\begin{proof}
    By linearity, it suffices to consider the case $p(x) = ax^n$, $a \in A_R$, $n \ge 0$.  In that case, \cref{road1,road2} follow immediately from computing the $w^{-1}$ coefficient on the right-hand side.  For \cref{road2}, we use that
    \[
        \begin{tikzpicture}[centerzero]
            \bubrightblank{0,0};
            \token{-0.2,0}{east}{a};
        \end{tikzpicture}
        = tz^{-1} \tr(a) 1_\one + \pluscbubble{a}{0}
        \quad \text{and} \quad
        \begin{tikzpicture}[centerzero]
            \bubleftblank{0,0};
            \token{0.2,0}{west}{a};
        \end{tikzpicture}
        = \plusccbubble{a}{0} - t^{-1} z^{-1} \tr(a) 1_\one.
    \]
    by \cref{fake1,tanks2}.  Finally, \cref{road3} follows from \cref{dog1,dog4,road1}
\end{proof}

The even part $Z(A_R)_{\bar{0}}$ of the center $Z(A_R) = R \otimes_\kk Z(A)$ of $A_R$ is a commutative algebra.  Fix a pair of monic polynomials
\begin{align}
    f(w) &= f_0 w^l + f_1 w^{l-1} + \dotsb + f_l \in Z(A_R)_{\bar{0}}[w], \\
    g(w) &= g_0 w^m + g_1 w^{m-1} + \dotsb + g_m \in Z(A_R)_{\bar{0}}[w]
\end{align}
such that $k=m-l$, $f_l,g_m \in \kk^\times$, and $t^2 = f_l/g_m$.  Define
\begin{align} \label{zebra1}
    \OO^+(w) &:= t^{-1}z \sum_{n \in \Z} \OO^+_n w^{-n}
    := g(w)/f(w) = w^k + w^{k-1} Z(A_R)_{\bar{0}} \llbracket w^{-1} \rrbracket,
    \\ \label{zebra2}
    \tOO^+(w) &:= -tz \sum_{n \in \Z} \tOO^+_n w^{-n}
    := f(w)/g(w) = w^{-k} + w^{-k-1} Z(A_R)_{\bar{0}} \llbracket w^{-1} \rrbracket,
    \\ \label{zebra3}
    \OO^-(w) &:= -tz \sum_{n \in \Z} \OO^-_n w^{-n}
    := t^2 g(w)/f(w) = 1 + w Z(A_R)_{\bar{0}} \llbracket w \rrbracket,
    \\ \label{zebra4}
    \tOO^-(w) &:= t^{-1}z \sum_{n \in \Z} \tOO^-_n w^{-n}
    := t^{-2} f(w)/g(w) = 1 + w Z(A_R)_{\bar{0}} \llbracket w \rrbracket,
\end{align}
cf.\ \cref{bubgen1,bubgen2,bubgen3,bubgen4}.

\begin{lem} \label{Detroit}
    The $R$-linear left tensor ideal $\cI_R(f|g)$ of $\Heis_k(A_R;z,t)$ generated by
    \begin{equation} \label{Detroit1}
        \tokup[f(x)] \quad \text{and} \quad \plusccbubble{a}{n} - \tr_R(\OO_n^+a) 1_\one,\quad -k < n < l,\ a \in A_R,
    \end{equation}
    is equal to the $R$-linear left tensor ideal $\cI_R(f|g)$ of $\Heis_k(A_R;z,t)$ generated by
    \begin{equation} \label{Detroit2}
        \tokdown[g(x)] \quad \text{and} \quad \pluscbubble{a}{n} - \tr_R(\tOO_n^+a) 1_\one,\quad k < n < m,\ a \in A_R.
    \end{equation}
    Moreover, this ideal contains
    \begin{equation} \label{Detroit3}
        \plusccbubble{a}{n} - \tr_R(\OO_n^+a) 1_\one,\
        \pluscbubble{a}{n} - \tr_R(\tOO_n^+a) 1_\one,\
        \minusccbubble{a}{n} - \tr_R(\OO_n^-a) 1_\one,\
        \minuscbubble{a}{n} - \tr_R(\tOO_n^-a) 1_\one,
    \end{equation}
    for all $a \in A$, $n \in \Z$.
\end{lem}

\begin{proof}
    The proof is similar to those of \cite[Lem.~9.2]{BSW-qheis} and \cite[Lem.~6.2]{BSW-foundations}.
    \details{
        For morphisms $\theta,\phi \colon X \to Y$, we will write $\theta \equiv \phi$ as shorthand for $\theta - \phi \in \cI_R(f|g)$.  By \cref{tanks,tanks2}, we have that $\plusccbubble{a}{n} \equiv \tr_R(\OO_n^+a) 1_\one$ when $n \le -k$, $\pluscbubble{a}{n} \equiv \tr_R(\tOO_n^+a) 1_\one$ when $n \le k$, $\minusccbubble{a}{n} \equiv \tr_R(\OO_n^-a) 1_\one$ when $n \ge 0$, and $\minuscbubble{a}{n} \equiv \tr_R(\tOO_n^-a) 1_\one$ when $n \ge 0$.

        We now use ascending induction on $n$ to show that $\plusccbubble{a}{n} \equiv \tr_R(\OO_n^+a) 1_\one$ for all $n \in \Z$.  By the above and \cref{Detroit1}, this holds for $n < l$, so assume that $n \ge l$.  The fact that $\tokup[f(x)] \equiv 0$ implies that
        \[
            \sum_{r=0}^l \plusccbubble{af_r}{n-r} + \sum_{r=0}^l \minusccbubble{af_r}{n-r} = \sum_{r=0}^l \ccbubble{af_r}{n-r} \equiv 0.
        \]
        On the left-hand side of the above, the only nonzero $(-)$-bubble occurs when $n=r=l$, hence we have $\sum_{r=0}^l \plusccbubble{af_r}{n-r} \equiv \delta_{l,n} t^{-1} z^{-1} \tr_R(af_l) 1_\one$.  Using the induction hypothesis and the fact that $f_l = g_m t^2$, we see that $\plusccbubble{a}{n} + \sum_{r=1}^l \tr_R(\OO_{n-r}^+ af_r) 1_\one \equiv \delta_{l,n} tz^{-1} \tr_R(ag_m) 1_\one$.  Equating $w^{l-n}$-coefficients in $f(w) \OO^+(w) = g(w)$, we get that $\sum_{r=0}^l f_r \OO_{n-r}^+ = \delta_{l,n} t z^{-1} g_m$.  Thus $\plusccbubble{a}{n} = \tr_R(a\OO_n^+) 1_\one$, as desired.

        Next, we show by descending induction on $n$ that $\minusccbubble{a}{n} \equiv \tr_R(\OO_n^-a) 1_\one$ for all $n \in \Z$.  Since this is clear for $n \ge 0$, assume $n<0$.  Equating $w^{-n}$-coefficients in $f(w) \OO^+(w) = t^{-2} f(w) \OO^-(w)$ gives that $\sum_{r=0}^l f_{l-r} \OO_{r+n}^+ = - \sum_{r=0}^l f_{l-r} \OO_{r+n}^-$.  Using the induction hypothesis and the previous paragraph, we deduce that
        \[
            \sum_{r=0}^l \plusccbubble{a f_{l-r}}{r+n} + \tr_R(\OO_n^- a f_l) 1_\one + \sum_{r=1}^l \minusccbubble{af_{l-r}}{r+n} \equiv 0.
        \]
        Since $\tokup[f(x)] \equiv 0$, we also have that
        \[
            \sum_{r=0}^l \plusccbubble{af_{l-r}}{r+n} + \sum_{r=0}^l \minusccbubble{af_{l-r}}{r+n}
            = \sum_{r=0}^l \ccbubble{af_{l-r}}{r+n}
            \equiv 0.
        \]
        Taking the difference of these two identities and using the fact that $f_l \in \kk^\times$ establishes the induction step.

        We have now shown that
        $
            \begin{tikzpicture}[centerzero]
                \plusgenleft{0,0};
                \token{-0.2,0}{east}{a};
            \end{tikzpicture}
            = \tr_R \left( \OO^+(w)a \right) 1_\one
        $
        and
        $
            \begin{tikzpicture}[centerzero]
                \minusgenleft{0,0};
                \token{-0.2,0}{east}{a};
            \end{tikzpicture}
            = \tr_R \left( \OO^-(w)a \right) 1_\one.
        $
        It then follows from \cref{infweeds} that
        $
            \begin{tikzpicture}[centerzero]
                \plusgenright{0,0};
                \token{0.2,0}{west}{a};
            \end{tikzpicture}
            = \tr_R \left( \tOO^+(w)a \right) 1_\one
        $
        and
        $
            \begin{tikzpicture}[centerzero]
                \minusgenright{0,0};
                \token{0.2,0}{west}{a};
            \end{tikzpicture}
            = \tr_R \left( \tOO^-(w)a \right) 1_\one.
        $
        (See \cite[Lem.~7.1]{BSW-foundations} for a similar argument.)  This completes the proof of the last statement in the lemma.

        Now note that
        \[
            \begin{tikzpicture}[centerzero]
                \draw[<-] (0,-0.4) to (0,0.4);
                \token{0,0}{east}{g(x)};
            \end{tikzpicture}
            \overset{\cref{road1}}{=}
            \left[
                \begin{tikzpicture}[anchorbase]
                    \draw[<-] (0,-0.4) to (0,0.4);
                    \token{0,0.15}{east}{(w-x)^{-1}};
                    \token{0,-0.15}{west}{g(w)};
                \end{tikzpicture}
            \right]_{w^{-1}}
            \overset{\cref{zebra1}}{=}
            \left[
                \begin{tikzpicture}[anchorbase]
                    \draw[<-] (0,-0.4) to (0,0.4);
                    \token{0,0.15}{east}{(w-x)^{-1}};
                    \token{0,-0.15}{west}{\OO^+(w) f(w)};
                \end{tikzpicture}
            \right]_{w^{-1}}
            \overset{\cref{Detroit3}}{\underset{\cref{adecomp}}{\equiv}}
            \left[
                \begin{tikzpicture}[anchorbase]
                    \draw[<-] (0,-0.4) to (0,0.4);
                    \plusgenleft{0.5,0};
                    \token{0,-0.2}{east}{f(w)};
                    \telecolor{black}{black}{0,0}{0.3,0};
                    \token{0,0.2}{east}{(w-x)^{-1}};
                \end{tikzpicture}
            \right]_{w^{-1}}
            \overset{\cref{road3}}{=} t^{-1}
            \begin{tikzpicture}[centerzero]
                \draw (0.2,0) to[out=down,in=east] (0.05,-0.2) to[out=left,in=down] (-0.2,0.35) -- (-0.2,0.5);
                \draw[wipe] (-0.2,-0.5) -- (-0.2,-0.35) to[out=up,in=west] (0.05,0.2) to[out=right,in=up] (0.2,0);
                \draw[<-] (-0.2,-0.5) -- (-0.2,-0.35) to[out=up,in=west] (0.05,0.2) to[out=right,in=up] (0.2,0);
                \token{0.2,0}{west}{f(x)};
            \end{tikzpicture}
            = 0.
        \]
        So we have shown that the $\kk$-linear left tensor ideal generated by \cref{Detroit1} contains the elements \cref{Detroit2}.  A similar argument shows that the $\kk$-linear left tensor ideal generated by \cref{Detroit2} contains the elements \cref{Detroit1}, completing the proof.
    }
\end{proof}

The \emph{generalized cyclotomic quotient} of $\Heis_k(A_R;z,t)$ corresponding to $f,g$ is the $R$-linear category
\begin{equation}
    \cH_R(f|g):= \Heis_k(A_R;z,t)/\cI_R(f|g).
\end{equation}


Assume that we are given a factorization $t = u v^{-1}$ for $u,v \in \kk^\times$ such that $u^2 = \tr(f_l)$ and $t^2 g_m = f_l$, which implies that $v^2 = \tr(g_m)$.  Let
\begin{equation}
    \cV(f) := \bigoplus_{n \ge 0} \psmod \QWA_n^f(A_R;z)
    \quad \text{and} \quad
    \cV(g)^\vee := \bigoplus_{n \ge 0} \psmod \QWA_n^f(A_R^\op;z),
\end{equation}
viewed as supermodule categories over $\blue{\Heis_{-l}(A_R;z,u)}$ and $\red{\Heis_m(A_R;z,v^{-1})}$ via the monoidal superfunctors $\Psi_f$ and $\Psi_g^\vee$ from \cref{sunset}.  Let
\begin{equation}
    \cV(f|g) := \cV(f) \boxtimes_R \cV(g)^\vee
\end{equation}
be their linearized Cartesian product, i.e.\ the $R$-linear supercategory with objects that are pairs $(X,Y)$ for $X \in \cV(f)$, $Y \in \cV(g)^\vee$, and morphisms
\[
    \Hom_{\cV(f|g)}((X,Y), (U,V)) := \Hom_{\cV(f)}(X,Y) \otimes_R \Hom_{\cV(g)^\vee}(Y,V),
\]
with the obvious composition law.  There is an equivalence of categories
\[
    \cV(f|g) \to \bigoplus_{n,m \ge 0} \psmod \left( \QWA_n^f(A_R;z) \otimes_R \QWA_m^g(A_R^\op;z) \right),
\]
hence $\cV(f|g)$ is additive and idempotent complete.  Furthermore, $\cV(f|g)$ is a module category over the symmetric product $\blue{\Heis_{-l}(A_R;z,u)} \odot \red{\Heis_m(A_R;z,v^{-1})}$.

\emph{We assume for the remainder of this section} that the base ring $R$ is a finite-dimensional supercommutative superalgebra over an algebraically closed field $\KK \supseteq \kk$; by ``eigenvalue'' we mean eigenvalue in this field $\KK$.   Recall the definition of $\tau_i$ from \cref{tau}.  The element $z \tau_1$ induces a $\KK$-linear endomorphism of $A_R \otimes_R A_R$ by left multiplication.  Since $A_R \otimes_R A_R$ is finite dimensional over $\KK$, this endomorphism has a minimal polynomial $m(w) \in \KK[w]$.  Let $\Gamma_R$ be the multiplicative subgroup of $\KK$ generated by the finite set
\begin{equation} \label{KD}
    \left\{ \frac{\mu}{\mu-\eta} : \mu,\eta \in \KK,\ m(\eta)=0,\ \mu^2 - \mu \eta - 1 = 0 \right\}.
\end{equation}
For example, if $A=\kk$ and $z=q-q^{-1}$ for some $q \in \kk^\times$ then $\tau_1 = 1 \otimes 1$, hence $m(w)=w-z$, and it follows that $\Gamma_R = \{q^{2n} : n \in \Z\}$.
Schur's lemma implies that $Z(A_R)_{\bar{0}}$ acts on any irreducible $A_R$-supermodule $L$ via a central character $\chi_L \colon Z(A_R)_{\bar{0}} \to \KK$.  For $p(w) \in Z(A_R)_{\bar{0}}[w]$, we have $\chi_L(p(w)) \in \KK[w]$.  Let $\Xi_p$ be the set of all roots of the polynomials $\chi_L(p(w))$ for all irreducible $A_R$-supermodules $L$. 

\begin{lem}
    Let $V$ be a finite-dimensional $\QAWA_2(A_R;z)$-module.  All eigenvalues of $x_2$ on $V$ are of the form $\lambda$ or $\gamma \lambda$ for $\gamma \in \Gamma_R$ and an eigenvalue $\lambda$ of $x_1$ on $V$.
\end{lem}

\begin{proof}
    Let $\lambda_2$ be an eigenvalue of $x_2$.  Since $x_1$, $x_2$, and $z \tau_1$ all commute, we can choose a vector $y \in V$ in the $\lambda_2$-eigenspace of $x_2$ that is also an eigenvector of $x_1$ with some eigenvalue $\lambda_1$ and an eigenvector of $z \tau_1$ with some eigenvalue $\eta$.  First suppose that $y$ is an eigenvector for $\sigma_1$ with eigenvalue $\mu$.  It follows from \cref{skein} that the element $\sigma_1 \in \QAWA_n(A_R;z)$ satisfies the equation $\sigma_1^2 - z \tau_1 \sigma_1 - 1= 0$.  Thus we have $0 = (\sigma_1^2 - z \tau_1 \sigma_1 - 1)y = (\mu^2 - \mu \eta - 1)y$.  So $\mu,\eta$ satisfy the conditions in the definition of the set \cref{KD}.  Furthermore, it follows from \cref{QAWC,skein} that $\sigma_1 x_1 = x_2 \sigma_1 - z x_2 \tau_1$.  Thus $\mu \lambda_1 = \mu \lambda_2 - \eta \lambda_2 = 0$, which implies that $\lambda_2 = \frac{\mu}{\mu-\eta} \lambda_1$.

    On the other hand, if $y$ is not an eigenvalue for $\sigma_1$, then $y$ and $\sigma_1 y$ are linearly independent.  Then the matrix describing the action of $x_1$ on the subspace with basis $\{y,\sigma_1 y\}$ is $\begin{pmatrix} \lambda_1 & -\eta \lambda_2 \\ 0 & \lambda_2 \end{pmatrix}$.  Hence $\lambda_2$ is an eigenvalue for the action of $x_1$ on $V$.
\end{proof}

\begin{cor} \label{hot}
    The eigenvalues of $x_1,\dotsc,x_n$ on any $\QWA_n^f(A_R;z)$-module lie in $\Gamma_R \Xi_f$.
\end{cor}

\begin{lem} \label{yard}
    Suppose that $f(w), g(w) \in Z(A_R)_{\bar{0}}[w]$ 
are generic in the sense that the sets $\Gamma_R \Xi_f$ and $\Gamma_R \Xi_g$ are disjoint. Then the categorical action of the symmetric product $\blue{\Heis_{-l}(A_R;z,u)} \odot \red{\Heis_m(A_R;z,v^{-1})}$ on $\cV(f|g)$ defined above extends to an action of the localization $\blue{\Heis_{-l}(A_R;z,u)} \barodot \red{\Heis_m(A_R;z,v^{-1})}$ from \cref{sec:comult}.
\end{lem}

\begin{proof}
    The proof is analogous to that of \cite[Lem.~9.4]{BSW-qheis}, using \cref{hot}.
\end{proof}

When the genericity assumption of \cref{yard} is satisfied, the lemma
implies that there is a strict $\kk$-linear monoidal superfunctor $\Psi_f \barodot \Psi_g^\vee \colon \blue{\Heis_{-l}(A_R;z,u)} \barodot \red{\Heis_m(A_R;z,v^{-1})} \to \SEnd_R(\cV(f|g))$.  Composing this superfunctor with the superfunctor $\cDelt[-l]$ from \cref{comult} yields a strict $R$-linear monoidal endofunctor
\begin{equation} \label{sticky}
    \Psi_{f|g} := \left( \Psi_f \barodot \Psi_g^\vee \right) \circ \cDelt[-l] \colon \Heis_z(A_R;z,t) \to \SEnd_R(\cV(f|g)).
\end{equation}
In this way, we make $\cV(f|g)$ into a module supercategory over $\Heis_k(A_R;z,t)$.

\begin{lem} \label{dimsum}
    Recalling the $(+)$-bubble generating functions from \cref{bubgen1,bubgen2}, we have
    \begin{align*}
        \Psi_{f|g}
        \left. \left(
            \begin{tikzpicture}[centerzero]
                \plusgenleft{0,0};
                \token{-0.2,0}{east}{a};
            \end{tikzpicture}
        \right) \right|_{(\QWA_0^f(A_R;z),\QWA_0^g(A_R^\op;z))}
        &= \tr_R \left( g(w)f(w)^{-1} a \right) \in w^{m-l} R \llbracket w^{-1} \rrbracket,
        \\
        \Psi_{f|g}
        \left. \left(
            \begin{tikzpicture}[centerzero]
                \plusgenright{0,0};
                \token{0.2,0}{west}{a};
            \end{tikzpicture}
        \right) \right|_{(\QWA_0^f(A_R;z),\QWA_0^g(A_R^\op;z))}
        &= \tr_R \left( f(w)g(w)^{-1} a \right) \in w^{l-m} R \llbracket w^{-1} \rrbracket.
    \end{align*}
\end{lem}

\begin{proof}
    The proof is almost identical to that of \cite[Lem.~6.8]{BSW-foundations}, using \cref{Detroit,catseye}.
    \details{
        Applying \cref{Detroit} with $g(w)=1$, we have
        \[
            \Psi_f
            \left(
                \begin{tikzpicture}[centerzero]
                    \plusgenleft[blue]{0,0};
                    \token[blue]{-0.2,0}{east}{a};
                \end{tikzpicture}
            \right)_{\QWA_0^f(A_R;z)}
            = \tr_R \left( f(w)^{-1} a \right),
            \qquad
            \Psi_f
            \left(
                \begin{tikzpicture}[centerzero]
                    \plusgenright[blue]{0,0};
                    \token[blue]{0.2,0}{west}{a};
                \end{tikzpicture}
            \right)_{\QWA_0^f(A_R;z)}
            = \tr_R \left( f(w) a \right).
        \]
        Similarly, applying it with $f(w)=1$, we have
        \[
            \Psi_g
            \left(
                \begin{tikzpicture}[centerzero]
                    \plusgenleft[red]{0,0};
                    \token[red]{-0.2,0}{east}{a};
                \end{tikzpicture}
            \right)_{\QWA_0^g(A_R^\op;z)}
            = \tr_R \left( g(w) a \right),
            \qquad
            \Psi_g
            \left(
                \begin{tikzpicture}[centerzero]
                    \plusgenright[red]{0,0};
                    \token[red]{0.2,0}{west}{a};
                \end{tikzpicture}
            \right)_{\QWA_0^g(A_R^\op;z)}
            = \tr_R \left( g(w)^{-1} a \right).
        \]
        Thus, by \cref{catseye}, we have
        \begin{align*}
            \Psi_{f|g}
            \left. \left(
                \begin{tikzpicture}[centerzero]
                    \plusgenleft{0,0};
                    \token{-0.2,0}{east}{a};
                \end{tikzpicture}
            \right) \right|_{(\QWA_0^f(A_R;z),\QWA_0^g(A_R^\op;z))}
            &= \tr_R \left( f(w)^{-1}ab \right) \tr_R \left( g(w)f(w)^{-1} \right)
            \overset{\cref{adecomp}}{=} \tr_R \left( g(w)f(w)^{-1} a \right),
            \\
            \Psi_{f|g}
            \left. \left(
                \begin{tikzpicture}[centerzero]
                    \plusgenright{0,0};
                    \token{0.2,0}{west}{a};
                \end{tikzpicture}
            \right) \right|_{(\QWA_0^f(A_R;z),\QWA_0^g(A_R^\op;z))}
            &= \tr_R \left( f(w) b^\vee a \right) \tr_R \left( g(w){-1} b \right)
            \overset{\cref{adecomp}}{=} \tr_R \left( f(w)g(w)^{-1} a \right).
        \end{align*}
    }
\end{proof}

\begin{theo} \label{GCQ}
    Suppose $f(w), g(w)$ are generic as in \cref{yard}. Let $\Ev \colon \SEnd_R(\cV(f|g)) \to \cV(f|g)$ be the $R$-linear superfunctor defined by evaluation on $\left( \QWA_0^f(A_R;z), \QWA_0^g(A_R^\op;z) \right) \in \cV(f|g)$.  Then $\Ev \circ \Psi_{f|g}$ factors through the generalized cyclotomic quotient $\cH_R(f|g)$ to induce an equivalence of $\Heis_k(A_R;z,t)$-module supercategories
    \[
        \psi_{f|g} \colon \Kar(\cH_R(f|g)_\pi) \to \cV(f|g),
    \]
    where $\Kar(\cH_R(f|g)_\pi)$ denotes the Karoubi envelope of the $\Pi$-envelope $\cH_R(f|g)_\pi$ of $\cH_R(f|g)$ (see \cite[Def.~1.10]{BE17}).
\end{theo}

\begin{proof}
    It is clear that $\tokup[f(x)]$ acts as zero on $\left( \QWA_0^f(A_R;z), \QWA_0^g(A_R^\op;z) \right)$.  Together with \cref{dimsum}, this implies that $\Ev \circ \Psi_{f|g}$ factors through the quotient $\cH_R(f|g)$.  Since $\cV(f|g)$ is an idempotent complete $\Pi$-supercategory (in the sense of \cite[Def.~3.1]{BE17}) we obtain the induced $R$-linear superfunctor $\psi_{f|g}$ from the statement of the theorem.  To show that $\psi_{f|g}$ is an equivalence, it remains to prove that it is full, faithful and dense.  This argument is almost identical to the one in the proof of \cite[Th.~9.5]{BSW-qheis}.
\end{proof}

\begin{rem}
    When $g(w)=1$, the genericity assumption is vacuous.  Hence \cref{GCQ} yields an equivalence of categories $\psi_{f|1} \colon \Kar(\cH_R(f|1)_\pi) \to \cV(f)$.  So the generalized cyclotomic quotient $\cH_R(f|1)$ is Morita equivalent to the usual cyclotomic quotient, that is, the locally unital superalgebra $\bigoplus_{n \ge 0} \QWA_n^f(A_R;z)$.
\end{rem}

\section{Basis theorem\label{sec:basis}}

In this final section we prove a basis theorem for the morphism spaces in $\Heis_k(A;z,t)$.  Unsurprisingly, our method of proof is a hybrid of the proofs for the quantum Heisenberg category \cite[Th.~10.1]{BSW-qheis} and the Frobenius Heisenberg category \cite[Th.~7.2]{BSW-foundations}; see also \cite[Th.~6.4]{BSW-K0}.

Recall from \cref{cocenter} the definition of the cocenter $C(A)$ of $A$.  For $a \in A$, we let $\cocenter{a}$ denote its image in $C(A)$.  For a general superalgebra $A$, the cocenter $C(A)$ is merely a vector superspace.  However, under our assumption that $A$ is a Frobenius superalgebra, the cocenter can be endowed with the structure of an associative commutative, but not necessarily unital, superalgebra with product $\diamond$; see \cite[Prop.~2.1]{ReeksS20}.  We will not use this product explicitly in this paper; however we note that $a^\dagger = z a \diamond 1$ (see \cref{adag}).  There is a well-defined nondegenerate pairing
\begin{equation} \label{slipstream}
    Z(A) \times C(A) \to \kk,\quad  (a,\cocenter{b}) \mapsto \tr(ab),\quad a,b \in A.
\end{equation}
(For a proof, see \cite[Lem~4.1]{BSW-foundations}.)

Let $\Sym(A)$ denote the symmetric superalgebra generated by the vector superspace $C(A)[x]$, where $x$ here is an even indeterminate.  For $n \in \Z$ and $a \in A$, let $e_n(a) \in \Sym(A)$ denote
\begin{equation}
    e_n(a) :=
    \begin{cases}
        0 & \text{if $n < 0$}, \\
        \tr(a) & \text{if $n=0$}, \\
        \cocenter{a} x^{n-1} & \text{if $n > 0$}.
    \end{cases}
\end{equation}
This defines a parity-preserving linear map $e_n \colon A \rightarrow \Sym(A)$.

\begin{lem}[{\cite[Lem.~7.1]{BSW-foundations}}] \label{controversy}
    For each $n \in \Z$, there is a unique parity-preserving linear map $h_n \colon A \rightarrow \Sym(A)$ such that
    \begin{equation} \label{banana}
        e(ac;-u) h(c^\vee b; u) = \tr(ab),
        \qquad \text{for all } a,b \in A,
    \end{equation}
    where we are using the generating functions
    \begin{equation} \label{moon}
        e(a;u) := \sum_{n \geq 0} e_n(a) u^{-n}, \quad
        h(a;u) := \sum_{n \geq 0} h_n(a) u^{-n}
        \in \Sym(A)\llbracket u^{-1}\rrbracket.
    \end{equation}
\end{lem}

By \cref{controversy,infweeds}, we have a well-defined homomorphism of superalgebras
\begin{equation} \label{beta}
    \beta \colon \Sym(A) \otimes \Sym(A) \to \End_{\Heis_k(A;z,t)}(\one),
\end{equation}
\begin{align} \label{beta-e}
    e_n(a) \otimes 1 &\mapsto (-1)^{n-1} t^{-1}z \plusccbubble{a}{n-k},
    &
    1 \otimes e_n(a) &\mapsto (-1)^{n-1} tz \minusccbubble{a}{-n},
    \\ \label{beta-h}
    h_n(a) \otimes 1 &\mapsto tz \pluscbubble{a}{n+k},
    &
    1 \otimes h_n(a) &\mapsto t^{-1}z \minuscbubble{a}{-n}.
\end{align}
We will show in \cref{colossus} that $\beta$ is an isomorphism.  We have that, for $X,Y \in \Heis_k(A;z,t)$, the superspace $\Hom_{\Heis_k(A;z,t)}(X,Y)$ is a right $\Sym(A) \otimes \Sym(A)$-module under the action
\[
    \phi \theta := \phi \otimes \beta(\theta),\quad
    \phi \in \Hom_{\Heis{A}{k}}(X,Y),\ \theta \in \Sym(A) \otimes \Sym(A).
\]

Let $X = X_n \otimes \dotsb \otimes X_1$ and $Y = Y_m \otimes \dotsb \otimes Y_1$ be objects of $\Heis_k(A;z,t)$ for $X_i, Y_j \in \{\uparrow, \downarrow\}$.  An \emph{$(X,Y)$-matching} is a bijection between the sets
\[
    \{i : X_i = \uparrow\} \sqcup \{j : Y_j = \downarrow\}
    \quad \text{and} \quad
    \{i : X_i = \downarrow\} \sqcup \{j : Y_j = \uparrow\}.
\]
A \emph{reduced lift} of an $(X,Y)$-matching is a string diagram representing a morphism $X \to Y$ such that
\begin{itemize}
    \item the endpoints of each string are points which correspond under the given matching;
    \item there are no floating bubbles and no dots or tokens on any string;
    \item there are no self-intersections of strings and no two strings cross each other more than once.
\end{itemize}
For each $(X,Y)$, fix a set $B(X,Y)$ consisting of a choice of reduced lift for each $(X,Y)$-matching.  Then let $B_\circ(X,Y)$ denote the set of all morphisms that can be obtained from the elements of $B(X,Y)$ by adding an integer number of dots near the terminus of each string and one element of $\BA$ to each string.

\begin{theo} \label{basis-thm}
Assume that the ground ring $\kk$ is an integral domain and that $z,t \in \kk^\times$ are arbitrary. For objects $X,Y \in \Heis_k(A;z,t)$, the morphism space $\Hom_{\Heis_k(A;z,t)}(X,Y)$ is a free right $\Sym(A) \otimes \Sym(A)$-module with basis $B_\circ(X,Y)$.
\end{theo}

\begin{proof}
    It suffices to consider the case $k \le 0$, since the result for $k \ge 0$ then follows by applying $\Omega_k$ from \cref{om}.  Let $X = X_r \otimes \dotsb \otimes X_1$ and $Y = Y_s \otimes \dotsb \otimes Y_1$ be two objects, with $X_i,Y_i \in \{\uparrow,\downarrow\}$.

     The defining relations and the additional relations proved in \cref{sec:first,sec:second,sec:third} give Reidemeister-type relations modulo terms with fewer crossing, plus a skein relation and bubble, dot, and token sliding relations.  These relations allow diagrams for morphisms in $\Heis_k(A;z,t)$ to be manipulated in a similar way to the way oriented tangles are simplified in skein categories, modulo diagrams with fewer crossings.  Thus, we have a straightening algorithm to rewrite any diagram representing a morphism $X \to Y$ as a linear combination of the ones in $B_\circ(X,Y)$. Hence $B_\circ(X,Y)$ spans $\Hom_{\Heis_k(A;z,t)}$ as a right $\Sym(A) \otimes \Sym(A)$-module. 

    It remains to prove linear independence of $B_\circ(X,Y)$.
    For this, recalling that $\kk$ is an integral domain by assumption, we can replace $\kk$ by the algebraic closure of its field of fractions to assume without loss of generality that $\kk$ is actually an algebraically closed field.
                We say $\phi \in B_\circ(X,Y)$ is \emph{positive} if it only involves nonnegative powers of dots.  It suffices to show that the positive morphisms in $B_\circ(X,Y)$ are linearly independent.  Indeed, given any linear relation $\sum_{i=1}^N \phi_i \otimes \beta(\theta_i)=0$ for morphisms $\phi_i \in B_\circ(X,Y)$ and coefficients $\theta_i \in \Sym(A) \otimes \Sym(A)$, we can multiply the termini of the strings by sufficiently large positive powers of dots to reduce to the positive case.

    We begin with the case $X = Y =\ \uparrow^{\otimes n}$.  Consider a linear relation $\sum_{s=1}^N \phi_s \beta(\theta_s)$ for some positive $\phi_s \in B_\circ(X,Y)$ and $\theta_s \in \Sym(A) \otimes \Sym(A)$.  Fix a homogeneous basis $a_1,\dotsc,a_r$ of $C(A)$ with $a_1,\dotsc,a_{r'}$ even and $a_{r'+1},\dotsc,a_r$ odd.  Choose $m_1, m_2 \ge 1$ so that all elements $\theta_s \in \Sym(A)$ are polynomials in $e_1(a_j) \otimes 1, \dotsc, e_{m_1}(a_j) \otimes 1$ and $1 \otimes e_1(a_j), \dotsc, 1 \otimes e_{m_2}(a_j)$, $1 \le j \le r$.  Then choose $l, m \ge 0$ such that
    \begin{itemize}
        \item $k=m-l$;
        \item the multiplicities of dots in all $\phi_s$ arising in the linear relation are $< l$;
        \item $m_1 + m_2 < m$.
    \end{itemize}

    Let $u_{i,j}$, $1 \le i \le m_1$, $1 \le j \le r$, and $v_{i,j}$, $1 \le i \le m_2$, $1 \le j \le r$, be indeterminates, with $u_{i,j}$ and $v_{i,j}$ even for $1 \le j \le r'$ and odd for $r' < j \le r$. Let $\KK$ be the algebraic closure of $\kk(u_{i,j}, v_{p,j} : 1 \le i \le m_1,\, 1 \le p \le m_2,\ 1 \le j \le r')$, and define $R$ to be the free supercommutive $\KK$-superalgebra generated by $u_{i,j}$, $1 \le i \le m_1$, $r' < j \le r$, and $v_{i,j}$, $1 \le i \le m_2$, $r' < j \le r$.  Since the $u_{i,j}$ and $v_{i,j}$ are odd for $r' < j \le r$, $R$ is finite dimensional over $\KK$.  We will now work with algebras/categories that are linear over $R$, as in \cref{sec:GCQ}.  Let $a_1^\vee,\dotsc,a_r^\vee$ be a basis of $Z(A)$ dual to the basis $a_1,\dotsc,a_r$ of $C(A)$ under the pairing of \cref{slipstream}.  Consider the cyclotomic wreath product algebras $\QWA_n^f(A_R;z)$ and $\QWA_n^g(A_R^\op;z)$ associated to the polynomials
    \begin{multline*}
        f(w) := w^l + t^2, \quad
        g(w) := w^m + (u_{1,1} a_1^\vee + \dotsb + u_{1,r} a_r^\vee) w^{m-1} + \dotsb + (u_{m_1,1} a_1^\vee + \dotsb + w_{m_1,r} a_r^\vee) w^{m-m_1} \\
        + (v_{m_2,1} a_1^\vee + \dotsb + v_{m_2,r} a_r^\vee) w^{m_2} + \dotsb + (v_{1,1} a_1^\vee + \dotsb v_{1,r} a_r^\vee) w + 1
        \quad \in Z(A_R)_{\bar{0}}[w].
    \end{multline*}
    The roots of $f(w)$ are contained in $\kk$. Also, since $A$ is defined over the algebraically closed field $\kk$, the set $\Gamma_R$ from \cref{KD} is actually contained in $\kk^\times$. For any irreducible $A_R$-module $L$, the evaluation of $\chi_L(g(w))$ at any element of $\Gamma_R$ involves at least one of the even $u_{i,j}$ or $v_{i,j}$ with a nonzero coefficient, hence is not contained in $\kk$.  This shows that $f(w)$, $g(w)$ is generic  in the sense of \cref{yard}.  Hence we can use the superfunctor $\Psi_{f|g}$ of \cref{sticky} to make $\cV(f|g)$ into a $\Heis_k(A_R;z,t)$-module supercategory. Using the canonical $\kk$-linear monoidal superfunctor $\Heis_k(A;z,t) \to \Heis_k(A_R;z,t)$, we can view $\cV(f|g)$ as a module supercategory over $\Heis_k(A;z,t)$.  We now evaluate the relation $\sum_{s=1}^N \phi_s \otimes \beta(\theta_s) = 0$ on $(\QWA_0^f(A_R;z), \QWA_0^g(A_R^\op;z)) \in \cV(f|g)$ to obtain a relation in $\QWA_n^f(A_R;z)$.  It follows from \cref{CQ-basis} and the choice of $l$ that the images of $\phi_1,\dotsc,\phi_N$ in $\QWA_n^f(A_R;z)$ are linearly independent over $R$.  Thus the image of $\beta(\theta_s) \in R$ is zero for each $s$.  We wish to deduce that each $\theta_s = 0$.  By our choice of $m$, each $\theta_s$ is a supercommutative polynomial in $e_1(a_j) \otimes 1, \dotsc, e_{m_1}(a_j) \otimes 1$ and $1 \otimes e_1(a_j), \dotsc, 1 \otimes e_{m_2}(a_j)$ for $1 \le j \le r$. So it suffices to show that the images of these elements under $\beta$ generate a free supercommutative superalgebra.  In fact, we claim that these images are the elements $u_{i,j}$ and $v_{i,j}$, respectively, up to a sign. To see this, note that the low degree terms of the series $\OO^\pm(w)$ from \cref{zebra1,zebra3} are
    \begin{align*}
        \OO^+(w)
        &=  w^k + (u_{1,1} a_1^\vee + \dotsb + u_{1,r} a_r^\vee) w^{k-1} + \dotsb + (u_{m_1,1} a_1^\vee + \dotsb + w_{m_1,r} a_r^\vee) w^{k-m_1} + \dotsb \in w^k \KK \llbracket w^{-1} \rrbracket,
        \\
        \OO^-(w)
        &=
        1 + (v_{1,1} a_1^\vee + \dotsb v_{1,r} a_r^\vee) w + \dotsb + (v_{m_2,1} a_1^\vee + \dotsb + v_{m_2,r} a_r^\vee) w^{m_1} + \dotsb \in \KK \llbracket w \rrbracket.
    \end{align*}
    Thus, the claim follows from \cref{beta-e,zebra1,zebra3,Detroit}.

    We have now proved the linear independence when $X = Y = \uparrow^{\otimes n}$.  The linear independence for more general $X$ and $Y$ follows from this as in the proof of \cite[Th.~10.1]{BSW-qheis}.
\end{proof}

\begin{cor} \label{colossus}
    The map $\beta \colon \End_{\Heis_k(A;z,t)}(\one) \cong \Sym(A) \otimes \Sym(A)$ is an isomorphism of superalgebras.
\end{cor}

\bibliographystyle{alphaurl}
\bibliography{qFrobHeis}

\end{document}